\documentclass[a4paper, 11pt]{book}

\usepackage[latin1]{inputenc}
\usepackage[english]{babel}
\usepackage[nottoc]{tocbibind}
\usepackage{emptypage}

\usepackage{amsmath}
\usepackage{amssymb}
\usepackage{amsfonts}
\usepackage{mathrsfs}

\usepackage[final]{showkeys}

\usepackage{color}
\usepackage{graphicx}

\usepackage[all]{xy}

\usepackage{cite}
\usepackage{epigraph}

\usepackage{dsfont}

\usepackage{geometry}
\newcommand{\C}{\mathbb{C}}
\newcommand{\rg}{\mathrm{reg}}
\newcommand{\sg}{\mathrm{sing}}
\newcommand{\supp}{\mathrm{supp}\,}
\newcommand{\loc}{\mathrm{loc}}
\newcommand{\Int}{\mathrm{Int}}
\renewcommand{\H}{\mathcal{H}}
\newcommand{\R}{\mathbb{R}}
\newcommand{\Z}{\mathbb{Z}}
\newcommand{\Ci}{\mathcal{C}}
\newcommand{\E}{\mathcal{E}}
\newcommand{\Di}{\mathcal{D}}
\newcommand{\D}{\mathscr{D}}
\newcommand{\B}{\mathscr{B}}
\newcommand{\N}{\mathbb{N}}
\newcommand{\EL}{\mathfrak{L}}
\newcommand{\CP}{\mathbb{CP}}
\newcommand{\A}{\mathscr{A}}
\newcommand{\Ol}{\mathcal{O}}
\newcommand{\DD}{\mathbb{D}}

\newcommand{\Lip}{\mathrm{Lip}}
\newcommand{\End}{\mathrm{End}}
\newcommand{\Hom}{\mathrm{Hom}}
\newcommand{\dist}{\mathrm{dist}}
\def\q{{\mathrm{q-loc}}}

\newcommand{\mass}{\mathbf{M}}
\renewcommand{\flat}[1]{\mathbf{F}(#1)}
\newcommand{\pair}[1]{\left\langle #1 \right\rangle}

\newcommand{\de}{\partial}
\newcommand{\debar}{\overline{\de}}
\renewcommand{\div}{\partial\mathrm{iv}}
\newcommand{\dbiv}{\debar\mathrm{iv}}
\newcommand{\dbad}{\overline{\partial}^{*}\!\!\mathrm{iv}}

\newcommand{\boh}{good}

\newtheorem{Teo}{Theorem}[section]
\newtheorem{Rem}{Remark}[section]
\newtheorem{Prp}[Teo]{Proposition}
\newtheorem{Cor}[Teo]{Corollary}
\newtheorem{Lmm}[Teo]{Lemma}

\begin{document}
\setlength{\epigraphwidth}{7.7cm}
\newenvironment{myparastyle}{\setlength{\parindent}{0pt}}{}
\renewcommand{\textflush}{myparastyle}
\renewcommand{\epigraphsize}{\footnotesize}

\begin{titlepage}
\begin{center}
\textsc{\Large Scuola Normale Superiore di Pisa\\}
\vspace{0.2cm}
\vspace{0.2cm}

\textsc{\large Classe di Scienze, Disciplina Matematica}\\ \vspace{0.5cm}
\end{center}
\begin{figure}
\centering
\includegraphics[width=2cm]{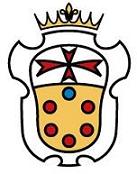}
\end{figure}
\begin{center}

\vspace{4.0cm}

\textsc{\Large Tesi di Perfezionamento}\\

\vspace{1.7cm}

\LARGE{\textbf{Applications of metric currents to complex analysis}}\\ 

\vspace{1.7cm}

%\begin{center}
{\large Dr.}\\ \vspace{-0.2cm} \textbf{{\large Samuele Mongodi}}\\ \vspace{-0.2cm} {\small
\ttfamily 
}
%\end{center}
\vspace{1.7cm}

\begin{minipage}{0.45\linewidth}
\centering {\large Advisor}\\ \vspace{-0.2cm} \textbf{{\large Prof. Giuseppe Tomassini}} \\
\vspace{-0.2cm} {\large Scuola Normale Superiore}\\
\end{minipage}

\vspace{1.5cm}

\end{center}

\cleardoublepage
\thispagestyle{empty} %\phantom{}	%aggiunge una pagina vuota
\newpage
\phantom{k}
\vspace{10cm}
\begin{flushright}
\emph{To my family\hspace{2.1cm}\phantom{l}}\\

\bigskip

\emph{Alla mia famiglia\hspace{1.5cm}\phantom{l}}\\
\end{flushright}
\newpage
\thispagestyle{empty} %\phantom{}	%aggiunge una pagina vuota
\end{titlepage}

\newpage

\vskip5cm
\frontmatter

\tableofcontents

\cleardoublepage
\thispagestyle{empty}
\phantom{k}
\cleardoublepage
\thispagestyle{empty}\epigraphhead[260]{\epigraph{With magic, you can turn a frog into a prince. With science, you can turn a frog into a Ph.D and you still have the frog you started with.}{T. Pratchett, I. Stewart, J. Cohen - \emph{The Science of Discworld}}}
\phantom{k}

\nocite{*}
\setlength{\epigraphwidth}{9.5cm}
\chapter{Introduction}
%\pagestyle{headings}
%\pagenumbering{arabic}
\epigraphhead[45]{\epigraph{Le but de ce travail est de munir son auteur du grade de docteur ès-sciences mathématiques et l'ensemble $H(X)$ des sous-espaces analytiques compacts de $X$ d'une structure d'espace analytique.}{A. Douady - \emph{Le problème des modules pour les sous-espaces analytiques compacts d'un espace analytique donné}}}
%\addcontentsline{toc}{chapter}{Introduction}

The theory of currents, introduced by De Rham in \cite{derham1} to study harmonic forms and developed by Federer, Fleming, Whitney and others \cite{federer2, federer3, federer4, whitney1} in the 50s and 60s, found deep and important applications in complex analysis and geometry, thanks to the work of Harvey, King, Lawson, Shiffman and others.

Just to mention a few of them, we recall the characterization of holomorphic chains by King in \cite{king1} and by Harvey and Shiffman in \cite{harvey2}, the removal of singularities for analytic functions and sets by Shiffman in \cite{shiffman1} and the boundary problem for holomorphic chains by Harvey and Lawson in \cite{harvey1, harvey3}.

\medskip

However, the main purpose of the geometric integration theory was the weak formulation of variational problems of geometric nature, whose main example is, probably, the Plateau problem; currents being the dual space of compactly supported smooth differential forms, their definition can be given on a smooth $n-$manifold, whereas many geometrical problems make sense in much more general settings, such as singular spaces, some classes of metric spaces, infinite dimensional manifolds.

\subsection*{Giants with broad shoulders}

In this direction of greater generality, a significant step was moved in 2000 by Ambrosio and Kirchheim, with the paper \emph{Currents in metric spaces} \cite{ambrosio1}, where they settled the foundations of a theory of currents on metric spaces and used it to pose and solve Plateau problem in a wide class of Banach spaces. They gave a new definition of current which was meaningful on any metric space. A \emph{metric current} is a multilinear functional on $(k+1)-$tuples of Lipschitz functions (with the first one bounded)
$$(f,\pi_1,\ldots, \pi_k)\mapsto T(f,\pi_1,\ldots, \pi_k)$$
satisfying a continuity property, a locality property and a finite mass property. 

\medskip

Locality property requires that $T(f,\pi_1,\ldots, \pi_k)=0$ whenever there exists a $\pi_j$ which is constant on a neighborhood of $\supp f$. Such a condition allows us to interpret $(f,\pi_1,\ldots,\pi_k)$ as the ``differential form'' $fd\pi_1\wedge\ldots\wedge d\pi_k$.

Employing the same idea, we can prescribe the vanishing of such a functional whenever the ``differentials'' belong to some particular class of functions. With the class of constant functions, this gives an analogue of the exterior differentiation; on a metric space with a complex structure, with the class of holomorphic or antiholomorphic functions, we obtain a way to introduce the concept of bidimension and the operators $\de$ and $\debar$.

\medskip

The finite mass condition asks for the existence of a finite Radon measure $\mu$ such that
$$|T(f,\pi_1,\ldots, \pi_k)|\leq \prod_{j=1}^k\Lip(\pi_j)\int |f|d\mu\;,$$
which is more or less the same definition of \emph{mass} given in the classical theory of currents.
On the other hand, the continuity we impose on metric currents is stronger than the usual one: as soon as $(f^n, \pi_1^n,\ldots, \pi_k^n)$ tend pointwise to $(f,\pi_1,\ldots, \pi_k)$, with a uniform bound on Lipschitz constants, we have $T(f^n,\pi_1^n,\ldots, \pi_k^n)\to T(f,\pi_1,\ldots, \pi_k)$. Together with the previous properties, continuity implies a chain rule for differentials and a change of variables formula.

This stronger continuity assumption makes the metric current behave like classical flat currents: we can define their pushforward through any Lipschitz map and the value of $T(f,\pi_1,\ldots, \pi_k)$ is determined by the values of $(f,\pi_1,\ldots, \pi_k)$ on the support of $T$ (and not on a neighborhood of it!).

Such a tight link between the current and its support can be specified in structure theorems, whenever the geometry of the underlying space allows it (e.g. reducible singularities); another consequence is that we can work safely with embeddings, because a control on the support is enough to go back to the original space.

\medskip

Two of these properties are inherently \emph{global} conditions (this is apparent for the mass finiteness, not so evident, but still true, for the continuity) and globality is an essential characteristic, if we want to treat also metric spaces which lack of local compactness, as it is the case for infinite dimensional Banach spaces. However, the theory of metric currents has a local counterpart, developed by U. Lang in \cite{lang1}, in which the ``differential forms''are $(k+1)-$tuples of \emph{locally} Lipschitz functions, with the first one compactly supported (hence Lipschitz and bounded). The mass condition is dropped and continuity is required only on sequences of forms whose ``support'' is contained in a given compact.

Such a variation, obviously, makes sense only in locally compact spaces, but is quite more flexible than the original version and allows us to define of a sheaf of currents. Even relaxed in such a way, the continuity property still gives a strong similarity with the classical locally flat currents and the previous considerations still hold. We mention that, even in $\R^n$, it is still an open problem whether metric currents coincide with classical flat currents or not, although it is known that the former class contains the latter.

Another advantage of the local setting is that we don't have to concern ourselves too much with the distance we are considering: locally equivalent distances will give locally equivalent theories. 

\subsection*{Getting complex}

A finite dimensional complex space can be given a metric space structure in many ways, e.g. by patching the metric given by local embeddings in $\C^n$, by k\"ahlerianity, by Kobayashi hyperbolicity; therefore we can consider on it the space of metric currents or, better, their local version.

As we mentioned before, we can define $(p,q)-$currents by requiring them to vanish on forms with $p+1$ holomorphic differentials or $q+1$ antiholomorphic differentials; this gives, with some caution, a bidimension theory for metric currents, together with a Dolbeault decomposition (unfortunately, not for all the currents), allowing us to define the $\de$ and $\debar$ operators. If we use local metric currents, it is possible to define sheafs of currents and to consider a version of the Dolbeault complex.

The natural questions which arise at this point are the characterization of holomorphic currents, i.e. $\debar-$closed $(p,n)-$currents, and the local (and global) exactness of the Dolbeault complex, i.e. the Cauchy-Riemann equation. 

\medskip

It turns out that holomorphic currents can be characterized by growth conditions around the singular set and this characterization can be made explicit by resolving the singularities; we do the computations in the case of complex curves, where a possible desingularization is given by the normalization, thus simplifying many technical details. A simple, but meaningful observation is that, whenever the space is locally Lipschitz contractible, the Poincar\'e lemma holds for metric holomorphic currents.

\medskip

The Cauchy-Riemann equation on singular complex spaces has been widely studied since the work by Henkin and Polyakov on complete intersections \cite{henk1}; Fornaess and Gavosto tried another approach to the problem, solving it for complex curves in \cite{forn3}. In more recent years, the problem has been tackled by Fornaess, {\O}vrelid, Ruppenthal and Vassiliadou in a series of papers, e.g. see \cite{forn1, forn2, vassil1, ruppenthal1, ruppenthal2, ruppenthal3}.

Recently, the representation formulas approach of Henkin and Polyakov was generalized by Andersson and Samuelsson in \cite{ander1}, leading to a solution of the Cauchy-Riemann equation in a suitable class of forms (or currents), see \cite{ander2}.

\medskip

A first and, maybe, na\"ive approach to such a problem is to try and apply the $L^2$ techniques due to H\"ormander to the regular part of the singular space; we therefore spend some efforts in defining Sobolev spaces on a singular complex space, investigating density theorems for functions and the vector-fields supported away from the singularity. Such vector-fields can obviously be interpreted as metric currents with $L^2$ (or, more generally, $L^p$) coefficients. 

H\"ormander's techniques give also an insight on the regularity of holomorphic currents, in case some density conditions hold; such density requests are equivalent to the coincidence of the minimal and maximal $L^2-$extensions of the $\debar$ operator on smooth forms with compact support in the regular part. Therefore, the failure of such regularity for holomorphic currents permits to identify situations in which the two extensions may not coincide.

\medskip

Again on the Cauchy-Riemann equation, the similarity with locally flat currents stressed before brings, as a consequence, a strong dependence of the structure of currents from the structure of the metric space, a fact, this one, also previously noted. We apply these considerations to locally ``completely'' reducible spaces (which can be locally viewed as union of smooth complex manifolds) and obtain an answer to the previous two questions in this case: Cauchy-Riemann equation is locally solvable on these spaces and the global solvability is controlled by the Dolbeault cohomology of the normalization, thus permitting to characterize holomorphic currents as the pushforward of holomorphic currents on the normalization.

It is worth noting that, for such spaces, the density hypotheses we were speaking before are satisfied.

\medskip

Many evidences can be found of a link between the degeneracy of the metric around the singular set and the behavior of the metric currents in relation to the $\debar-$equation; such a phenomenon is explored in some detail again for complex curves, exploiting the $L^p-$solvability of the Cauchy-Riemann equation in one complex variable, with subharmonic weights, which is proved in \cite{fornsib} by Fornaess and Sibony.

\medskip

The peculiar behavior of the Cauchy-Riemann equation in $L^p$ and the flatness of classical currents with $L^p$ coefficients suggests another interesting case to study: complex hypersurfaces in $\C^n$. Unfortunately, to solve the Cauchy-Riemann equation for metric currents on such spaces, we need $L^p$ solvability in the complement of the hypersurface in a ball.

\subsection*{Back on track}

Metric currents were originally employed to solve the Plateau problem in Banach spaces;  Ambrosio and Kirchheim adapted the classical proof, developing a theory of integral currents, with closure and compactness theorems, which are the main ingredients for any extremal problem. 

Another important ingredient of the proof is the cone construction, which shows that the set where we are looking for the minimum is not empty; the version of the cone construction for Banach spaces given in \cite{ambrosio1} is essentially a form of convolution with the ``primitive'' of $\delta_1-\delta_0$, that is the characteristic function of the interval $[0,1]$, together with a contraction in $0$. 

\medskip

If we are interested in the $\debar$ operator on a Banach space, the ``primitive'' of $\delta_1$ is $(\pi(t-1))^{-1}$. The main difference with the previous convolution kernel is the non-compactness of the support, which influences also the finite-dimensional solutions of the $\debar-$equation: the compactness of the support of the data is not enough to ensure the existence of a compactly supported solution.

The main consequence of this difference is the lack of a global mass estimate for the metric functional so obtained; on the finite mass condition relies the proof of the continuity property, which also breaks down here. However, as $1/t$ is locally integrable on $\C$, we have mass estimates on every bounded set and, consequently, continuity on sequences of forms with supports in a prescribed ball.

\medskip

An infinite dimensional Banach space isn't locally compact, therefore we cannot use local currents, but we can try and define a rough analogue of them, substituting  compactness with boundedness. We called these objects \emph{quasi-local} currents, halfway between local and global metric currents.

We can solve the $\debar$ for a boundedly supported current in a Banach space, if we allow the result to be a quasi-local current, applying the variation of the cone construction explained above.

\medskip

From the geometric viewpoint, in the complex case, Plateau problem is replaced by the boundary problem for holomorphic chains. Its solution, given by Harvey and Lawson, exploits heavily the theory of rectifiable currents and slicing.

Having these ingredients at our disposal, we approach the boundary problem for holomorphic chains in Hilbert spaces, where the inner product permits to recover something similar to the finite dimensional Wirtinger formula.

With the same (few) ideas, we can tackle also the generalizations of some characterization results on holomorphic chains and positive currents.

\subsection*{Per aspera ad aspera}

We organized the contents as follows.

\medskip

In the first chapter, we present the basic notions in complex analysis and geometry and in geometric measure theory, which we will need in the following.

\medskip

The second chapter is devoted to the theory of local metric currents. We introduce the basic concepts and adapt the theory to the complex case, defining the bidimension, the Dolbeault decomposition and related notions. A characterization of holomorphic currents is given.

\medskip

The development of a Sobolev theory on singular space is our main concern in the third chapter. We give a characterization of Sobolev functions in terms of their behavior and growth on the regular part; this leads to a capacity theory for the singular set which allows us to obtain an approximation result, using functions with support in the regular part. The second part of the chapter deals with the $L^2$ theory on singular spaces.

\medskip

In the fourth chapter, some applications of the theory are discussed. We solve Cauchy-Riemann equation on completely reducible singularities, by means of a structure theorem for metric currents; we also treat the equation in $L^p$ on complex curves and outline a possible approach for the same problem on complex spaces which can be embedded as divisors in $\C^n$.

\medskip

The final chapter tries to spread some light on the complex geometry in infinite dimensional spaces. After solving the Cauchy-Riemann equation in Banach spaces, in terms of the quasi-local metric currents, we turn to the study of positive currents, holomorphic chains and their boundaries.

\medskip

Some of the results exposed in this work can be found in \cite{mongodi1}. Many problems are left unanswered, some others are newly raised. We try to summarize this situation in the final remarks.

\subsection*{Acknowledgements}
\addcontentsline{toc}{section}{Acknowledgements}

Anyone who has ever tried to sit down and write an \emph{Acknowledgements} section knows the twofold nature of the inspiration that, sooner or later, comes to the writer.

\medskip

On one, more pragmatic, side, I feel the need to express my gratitude to those who helped me with this mathematical work. That's why I thank my advisor, prof. Giuseppe Tomassini, for his guidance, patience and attentiveness in following and helping me in my efforts at doing these tiny bits of mathematics; almost the same amount of gratitude and thankfulness I owe to prof. Luigi Ambrosio, who helped me through my (many) doubts in geometric measure theory, always willing to debate with me, giving useful advice and precious insights.

My most sincere thanks also go to my former colleague Matteo Scienza, with whom I had many good  chats, both about my thesis' subject and other bits of  mathematics; it is well known, i believe, that one of the best ways to gather one's thoughts is to explain them aloud to someone capable of active criticism.

Moreover, I would like to thank prof. Fabrizio Broglia, prof. Francesca Acquistapace, prof. Valentino Magnani for the aid they gave me in finding references and precise information in fields of mathematics that are even less familiar to me than complex analysis.

I should also thank everyone I pestered by talking about metric currents and related problems, independently of the answers or pieces of advice I got back, because, as I already said, talking is a major tool in proof-checking. Alas, the list would be long and I would forget someone, so I will stop here with this first half of thanks.

\medskip

On the other side, a PhD thesis marks the end of 3 (or more, as this is the case) years of study and life. Somewhere I read that you cannot be a PhD student only from 9 to 17 on weekdays and, to me, that's absolutely true; I won't say that I have been working 24/7 all these years, but it is surely more correct to say that I \emph{lived} as a PhD student, than that I studied or worked as one.

Therefore it's mandatory, for me, to thank all the people that put up with me, lived, played, laughed, cried, drank with me in these 3 years and more. Some of them maybe won't read these words, some of them surely won't, but that's not the point.

I thank Mitch, L\`o, Dade and Cri,  the best flatmates I could hope to meet.

I thank Giancarlo, Patrizia, their pub (\emph{the} Pub) and its waitresses.

I thank Vale, Giangi and Teo (again!).

I thank all the people into Italian Mathematical Olympiads: Bobo, Max, Ludo, Francesco, Michele, Maria, Lab, Pol, Venez, Marco, Simo, Ale, Jack, Cla, Gdt and all the others.

I thank Paolo, Ast, Alberto, Maurizio, Alessio, Vins, Francesco, Al, Stefano, Annalisa, Filippo and all the others, students and former students in Pisa, whom I met at university and with whom I always enjoyed speaking about maths or anything else.

I thank all the guys from the ju jitsu classes at the Mithos.

I thank Gau and Luisa and many other pub-friends.

I thank Teo, Frenzis, Ndrew, Ely, Luca, Barbara for their friendship.

I thank everyone I am forgetting and I apologize for this fault of mine.

\medskip

And obviously, as I am not that good at counting, there is a third half of thanks, really short, but not at all small. As predictable as it is, I have to thank my parents Aldina and Zeffirino, my brother Gionata and my grandmother Iva. Too many reasons there are, to list them all here, many of which have not even to do with this thesis.

\medskip

To all these people, thank you.

\mainmatter 

\setlength{\epigraphwidth}{7.7cm}
\chapter{Basic notions}
%\pagenumbering{arabic}
\epigraphhead[55]{\epigraph{The White Rabbit put on his spectacles. 

``Where shall I begin, please your Majesty?'' he asked.

``Begin at the beginning,'' the King said gravely, ``and go on till you come to the end: then stop.''}%
{Lewis Carroll - \emph{Alice's Adventures in Wonderland}}}

Summarizing the main concepts of complex geometry and geometric measure theory in little more than a dozen pages is somewhat ridiculous, considering the extension, both in time and papers, of the work done in both by mathematicians. However, we feel the need to fix some notations, recall the basic ideas and give some bibliographical indications, which will be forcefully incomplete.

\section{Complex manifolds and complex spaces}

A \emph{complex analytic manifold} of (complex) dimension $n$ is a differentiable manifold $M$, $2n-$dimensional as a real manifold, equipped with an atlas $\{(U_\alpha,\tau_\alpha)\}_{\alpha\in A}$ which is holomorphic with values in $\C^n$; by definition, this means that the transition functions $\tau_{\alpha\beta}$ are holomorphic.

The tangent space $T_xM$ has a natural complex vector space structure, given by the isomorphism $d\tau_\alpha: T_xM\to \C^n$, for $x\in U_\alpha$; if $T^\R_xM$ is the underlying real vector space, the multiplication by $i$ induces $J_x\in\End(T^\R_xM)$ such that $J^2_x=-I_x$ and the distribution $J_x$ is integrable. Such an endomorphism is called \emph{complex structure}.

Given an open set $U\subset M$ and analytic coordinates $z_1,\ldots, z_n$ on $U$, with $z_k=x_k+iy_k$, the real tangent space $T^\R M\vert_U$ admits the basis
$$\left\{\frac{\de}{\de x_1},\ldots, \frac{\de}{\de x_n},\frac{\de}{\de y_1},\ldots, \frac{\de}{\de y_n}\right\}$$
and the complex structure $J$ can be described by
$$J(\de/\de x_i)=\de/\de y_i\qquad J(\de/\de y_i)=-\de/\de x_i\;.$$

Let $T^\C M$ be the complexification of $T^\R M$, that is $T^\C M=T^\R M\otimes_\R\C=T^\R M\oplus iT^\R M$. The endomorphism $J$ extends naturally to $J\otimes 1\in\End(T^\C M)$, which satisfies again $(J\otimes1)^2=-I$; therefore, this endomorphism admits two $n-$dimensional eigenspaces, associated to the eigenvalues $i$ and $-i$. We denote them by $T^{1,0}M$ (\emph{holomorphic vectors}) and $T^{0,1}M$ (\emph{antiholomorphic vectors}) respectively and we   endow $T^{1,0}M\oplus{T^{0,1}M}$ with projections 
$$v\mapsto \frac{1}{2}(v-i(J\otimes1)v)\qquad v\mapsto \frac{1}{2}(v+i(J\otimes1)v)$$
 It is easy to see that $T^{1,0}M$ and $\overline{T^{0,1}M}$ are both canonically isomorphic to $TM$ through these projections.

We define
$$\frac{\de}{\de z_j}=\frac{1}{2}\left(\frac{\de}{\de x_j}-i(J\otimes i)\frac{\de}{\de x_j}\right)=\frac{1}{2}\left(\frac{\de}{\de x_j}-i\frac{\de}{\de y_j}\right)\;,$$
$$\frac{\de}{\de \bar{z}_j}=\frac{1}{2}\left(\frac{\de}{\de x_j}+i(J\otimes i)\frac{\de}{\de x_j}\right)=\frac{1}{2}\left(\frac{\de}{\de x_j}+i\frac{\de}{\de y_j}\right)\;.$$

In this way, we have a canonical decomposition $T^\C M=TM\oplus \overline{TM}$, which carries on to the cotangent bundle:
$$\Hom_\R(T^\R M, \C)\cong\Hom_\C(TM\otimes_\R\C, \C)\cong T^*M\oplus \overline{T^*M}$$
where $T^*M$ is the space of $\C-$linear forms. With the previous notations, $(dx_k, dy_k)$ is a basis for $\Hom_\R(T^\R M, \C)$, $(dz_k)$ a basis for $T^*M$, $(d\bar{z}_k)$ a basis for $\overline{T^*M}$. The differential of a function $f\in\Ci^1(U,\C)$ can be written
$$df=\sum\frac{\de f}{\de x_k}dx_k +\sum\frac{\de f}{\de y_k }dy_k=\sum\frac{\de f}{\de z_k}dz_k+\sum\frac{\de f}{\de \bar{z}_k}d\bar{z}_k$$
and the function $f$ is holomorphic if and only if $df$ is $\C-$linear i.e. if and only if $f$ is a solution of  the \emph{Cauchy-Riemann equations} $\de f/\de \bar{z}_k=0$ for $k=1,\ldots, n$. 

The ring of holomorphic functions on an open set $U\subset M$ will be denoted by $\Ol(U)$.

The above decomposition of $df$ gives
$$\Lambda^k_\R T^\R M\otimes\C=\Lambda^k_\C(TM\otimes\C)=\bigoplus_{p+q=k}\Lambda^{p,q}TM=\bigoplus_{p+q=k}\Lambda^p TM\oplus \Lambda^q\overline{TM}$$
$$\Lambda^k_\R (T^\R M)^*\otimes\C=\Lambda^k_\C(TM\otimes\C)^*=\bigoplus_{p+q=k}\Lambda^{p,q}T^*M=\bigoplus_{p+q=k}\Lambda^p T^*M\oplus \Lambda^q\overline{T^*M}\;.$$
A complex vector field is said to be of \emph{type $p,q$} if its value at every point of $M$ lies in $\Lambda^{p,q}TM$; similarly, a complex differential form is said to be of \emph{bidegree $p,q$} if its value at every point of $M$ lies in $\Lambda^{p,q}T^*M$. A general complex differential form (or vector field) splits into its $(p,q)-$components. The exterior differential on forms (and divergence on vector fields) splits into two components: on the $(p,q)-$forms, we have
$$d:\Lambda^{p,q}T^*M\to\Lambda^{p+1,q}T^*M\oplus \Lambda^{p,q+1}T^*M$$
and we set 
$$\de=\pi_{p+1,q}\circ d:\Lambda^{p,q}T^*M\to\Lambda^{p+1,q}T^*M$$
$$\debar=\pi_{p,q+1}\circ d:\Lambda^{p,q}T^*M\to\Lambda^{p,q+1}T^*M\;.$$
We define these operators on a generic form by linearity.

The identity $d^2=0$ implies $\de^2=\debar^2=\de\debar+\debar\de=0$, so, in particular, for each $p$ the operator $\debar$ is associated to a complex, called \emph{Dolbeault complex}:
$$\Ci^\infty(M,\Lambda^{p,0}T^*M)\xrightarrow{\debar}\cdots\xrightarrow{\debar}\Ci^\infty(M\Lambda^{p,q}T^*M)\xrightarrow{\debar}\Ci^\infty(M\Lambda^{p,q+1}T^*M)\xrightarrow{\debar}\cdots$$
and to the corresponding \emph{Dolbeault cohomology groups}
$$H^{p,q}(M,\C)=\frac{\ker \debar^{p,q}}{\mathrm{Im}\debar^{p,q-1}}$$
with the convention that the image of $\debar$ is $0$ for $q=0$. The groups $H^{p,0}(M,\C)$ correspond to the spaces of \emph{holomorphic $p-$forms} on $M$.

A holomorphic map $F:M_1\to M_2$ induces a homomorphism
 $$
 F_1^*:H^{p,q}(M_2,\C)\to H^{p,q}(M_1,\C).
 $$
\subsection{Analytic sets and complex spaces}

Let $M$ be a complex analytic manifold; a subset $A\subset M$ is called \emph{analytic subset} of $M$ if it is closed and if for every point $x_0\in A$ there exist a neighbourhood $U$ and holomorphic functions $g_1,\ldots, g_m\in \Ol(U)$ such that
$$A\cap U=\{z\in U\ :\ g_1(z)=\ldots=g_m(z)=0\}\;.$$

Unions and intersections of analytic sets are again analytic sets and the analytic continuation principle shows that if $A$ is an analytic subset of a connected manifold $M$, then either $A=M$ or $A$ has no interior points.

\medskip

In what follows we will need the concepts of germ (of a function or of a set) and of sheaf; we refer the reader to \cite{godement1, iversen1} for the basic definitions and properties.

\medskip

Let $(A,x)$ be the germ of the set $A$ at the point $x$; if $A$ is an analytic subset of an open set $U\subset M$, with $x\in U$, let $\mathscr{I}_{A,x}$ be the ideal of germs $f\in\Ol_{M,x}$ which vanish on $(A,x)$. Conversely, if $\mathscr{I}=(g_1,\ldots, g_m)$ is an ideal of $\Ol_{M,x}$, we denote by $(V(\mathscr{I}, x)$ the germ at $x$ of the zero variety $V(\mathscr{I})=\{z\in U\ :\ g_1(z)=\ldots=g_m(z)=0\}$. A local version of Hilbert's Nullstellenstatz holds for analytic sets. 

\begin{Teo}Let $\Ol_n=\Ol_{\C^n,0}$, then for every ideal $\mathscr{I}\subset\Ol_n$,
$$\mathscr{I}_{V(\mathscr(I),x}=\sqrt{\mathscr{I}}\;.$$
\end{Teo}

A germ of analytic set is irreducible if it has no decomposition $(A,x)=(A_1,x)\cup (A_2,x)$ in different analytic sets; $(A,x)$ is irreducible if and only if $\mathscr{I}_{A,x}$ is a prime ideal in $\Ol_{M,x}$. We have the following local parametrization result, due to Ruckert.

\begin{Teo}Let $\mathscr{I}$ be a prime ideal in $\Ol_n$ and $A=V(\mathscr{I})$. Then there exist an integer $d$, a choice of coordinates $(z';z'')=(z_1,\ldots, z_d;z_{d+1},\ldots, z_n)$, polydiscs $\Delta'$, $\Delta''$ in $\C^d$, $\C^{n-d}$ with sufficiently small radii such that the projection $\pi:A\cap(\Delta'\times \Delta'')\to\Delta'$ is a ramified covering with $q$ sheets, whose ramification locus is contained in $S=\{z'\:\ \delta(z')=0\}$ where $\delta\in\Ol_d$.\end{Teo}

We remark that, in the setting of the previous Theorem, $\Ol_n/\mathscr{I}$ is a finite extension of $\Ol_d$ and $q$ is the degree of such extension. Moreover, $\delta(z')$ is the discriminant of the minimal polynomial of a primitive element of such extension.

\medskip

A point $x\in A$ is said to be \emph{regular} if there exists a neighbourhood $\Omega$ of $x$ in $M$ such that $\Omega\cap A$ is an analytic submanifold of $\Omega$; otherwise, $x$ is called \emph{singular point}. The corresponding subsets of $A$ will be denoted by $A_\rg$ and $A_\sg$, respectively. For every point $x\in A$ such that $(A,x)$ is irreducible, there exist a family of neighbourhoods $\Omega$ of $x$ such that $A_\rg\cap\Omega$ is dense and connected in $A\cap \Omega$.

The dimension of an irreducible germ of an analytic set $(A,x)$ is defined to be $\dim (A,x)=\dim (A_\rg, x)$.

\medskip

If $x\in A$, we define the local ring $\Ol_{A,x}$ of germs of functions on $(A,x)$ which can be extended to germs of holomorphic functions on $(M,x)$; $\Ol_{A,x}=\Ol_{M,x}/\mathscr{I}_{A,x}$ and it is called the {\em ring of germs of holomorphic functions on $(A,x)$}. This definition allows us to define the analytic subsets of an analytic set in just the same way as we already did for manifolds. We have the following important result.

\begin{Teo}$A_\sg$ is an analytic subset of $A$.\end{Teo}

An immediate consequence of this theorem is that every analytic set admits a stratification in analytic subsets $A=A^0\supset A^1\supset A^2\ldots\supset A^m$ with $A^j=A^{j-1}_\sg$, therefore $A^j_\rg=A^j\setminus A^{j+1}$.

\bigskip

Given two analytic sets, $A\subset\Omega\subset\C^n$ and $B\subset\Omega'\subset\C^p$, a map $F:A\to B$ is said to be an \emph{analytic morphism} (or \emph{holomorphic map between analytic sets}) if for every $x\in A$ there exists a neighborhood $U$ of $x$ in $\Omega$ and a holomorphic map $\tilde{F}:\Omega\to\C^p$ such that $\tilde{F}\vert_{A\cap\Omega}=F\vert_{A\cap\Omega}$.

Equivalently, if $F$ is continuous and for every $x\in A$ and $g\in\Ol_{B,F(x)}$, we have $g\circ F\in\Ol_{A,x}$. The induced map
$$F^*_x:\Ol_{B,F(x)}\to\Ol_{A,x}$$
is called the \emph{comorphism} of $F$ at point $x$.

\medskip

A \emph{complex space} $X$ is a locally compact Hausdorff space, countable at infinity, together with a sheaf $\Ol_X$ of continuous functions on $X$, such that there exists an open covering $\{U_\lambda\}_{\lambda}$ of $X$ and for each $\lambda$ an homeomorphism $F_\lambda:U_\lambda\to A_\lambda$ onto an analytic set $A_\lambda\subset\Omega_\lambda\subset\C^{n(\lambda)}$ such that the comorphism $F^*_\lambda:\Ol_{A_\lambda}\to\Ol_{X\vert U_\lambda}$ is an isomorphism of sheaves of rings. $\Ol_X$ is called the {\em structure sheaf} of $X$.

\medskip

By definition, a complex space $X$ is locally isomorphic to an analytic set, therefore the concepts of holomorphic and meromorphic functions on $X$, analytic subsets, analytic morhpisms, regular and singular points, etc. are meaningful. For instance, the analogous result for analytic sets implies that $X_\sg$ is an analytic subset of $X$.

\subsection{Normalization and resolution of singularities}

In \cite{hironaka1}, Hironaka showed that every singular algebraic variety $X$ admits a resolution of singularities, i.e. there exist a manifold $Y$ and a proper morphism $\pi:Y\to X$ such that $\pi^{-1}(X_\sg)$ is a normal crossing divisor $E$ in $Y$ (called the \emph{exceptional divisor}) and $\pi$ is an isomorphism outside $E$.

The manifold $Y$ can be obtained by repeatedly blowing up the space $X$ along its singular locus; we refer to \cite{kollar1} for an excellent explanation of the basic (and not so basic) concepts.

\medskip

We will not need the full extent of Hironaka's result; in fact, we will be only concerned with the resolution of singularities for a complex curve. In such a particular case, we can obtain the desired result with significantly less efforts.

A complex space $X$ with structure sheaf $\Ol_X$ is \emph{normal} at $x\in X$ if $\Ol_{X,x}$ is a normal ring, i.e. if $\Ol_{X,x}$ is an integrally closed integral domain. Let $N(X)$ be the set of non normal points; we say that $X$ is normal if $N(X)=\emptyset$. In general, we have that (see \cite{abhyankar1}
\begin{enumerate}
\item $N(X)$ is a closed analytic subspace of $X$ and $N(X)\subset X_\sg$
\item for $x\in X\setminus N(X)$,
$$\dim_{x}X_\sg\leq \dim_x X-2$$
\item if $X$ is a complete intersection at $x$ and if the above inequality holds, then $X$ is normal at $x$.
\end{enumerate}

By the above, we deduce that if a complex curve is normal, then it is smooth and viceversa. 

A \emph{normalization} of a reduced analytic space $X$ is a pair $(X^\nu,\pi)$, where $X^\nu$ is a normal analytic space and $\pi$ is a finite surjective analytic mapping inducing an isomorphism of the open sets
$$X^\nu\setminus\pi^{-1}(N(X))\to X\setminus N(X)\;.$$

Therefore, to resolve the singularities of a complex curve it is enough to find a normalization.

\medskip

The normalization is uniquely determined up to an isomorphism, that is, if $(X^\nu_1,\pi_1)$ and $(X^\nu_2,\pi_2)$ are two normalizations, then there exists a unique analytic isomorphism $\phi:X_1^\nu\to X_2^\nu$ such that $\pi_1=\pi_2\circ\phi$.

The normalization exists and has the following properties. For every point $x\in X$  the set of irreducible components of $X$ at $x$ is in one-to-one correspondence with $\pi^{-1}(x)$. The fibre at $x$ of the direct image $\pi_*\Ol_{X^\nu}$ of the structure sheaf $\Ol_{X^\nu}$ is naturally isomorphic to the integral closure of the ring $\Ol_{X,x}$ in its complete ring of fractions.

The concept of a normal analytic space over $\C$ can be introduced in terms of analytic continuation of holomorphic functions. Namely, a reduced complex space is normal if and only if Riemann's first theorem on the removal of singularities holds for it: if $U\subset X$ is an open subset and $A\subset U$ is a closed analytic subset not containing irreducible components of $U$, then any function that is holomorphic on $U\setminus A$ and locally bounded on $U$ has a unique analytic continuation to a holomorphic function on $U$. For normal complex spaces Riemann's second theorem on the removal of singularities also holds: if  $\textrm{codim}_x A\geq 2$ at every point $x\in A$ , then the analytic continuation in question is possible without the requirement that the function is bounded. A reduced complex space $X$ is normal if and only if for every open set $U\subset X$ the restriction mapping of holomorphic functions
$$\Gamma(U,\Ol_X)\to \Gamma(U\setminus X_\sg, \Ol_X)$$
is bijective. 
A reduced complex space $X$ is a Stein space if and only if its normalization $X^\nu$ has this property.

\section{Geometric measure theory}

Let $X$ be a topological space and $\mathcal{M}$ a $\sigma-$algebra. A measure is usually defined as a $\sigma-$additive function $\mu:\mathcal{M}\to[0,+\infty]$ such that $\mu(\emptyset)=0$. However, we will call \emph{measure} a function $m$ defined on all the subsets of $X$, with values in $[0,+\infty]$, such that $m(\emptyset)=0$, which is only required to be countably subadditive\footnote{these set-functions are usually called \emph{outer measures}.}.

Given $(\mathcal{M},\mu)$, for any subset $A$ of $X$ we define
$$m(A)=\inf\{\mu(E)\ : A\subset E\in\mathcal{M}\}\;.$$
The sets $A$ for which
$$m(E)=m(E\cap A)+m(E\setminus A)\qquad \forall\;E\subseteq X\;,$$
are called \emph{$\mu-$measurable.}

A measure $m$ is called \emph{Borel regular} if every open set is $\mu-$measurable and if, for each $A\subseteq X$ there exists a Borel set $B\subseteq X$ with $A\subseteq B$ and $m(A)=,(B)$. 

Moreover, suppose $X$ is a locally compact Hausdorff space, then a measure $m$ is called a \emph{Radon measure} if the following conditions hold:
\begin{enumerate}
\item every compact set has finite $m$ measure;
\item every open set is $m-$measurable and, if $V\subseteq X$ is an open set, then 
$$m(V)=\sup\{m(K)\ :\ K\textrm{ is compact and }K\subseteq V\}\;;$$
\item for every $A\subseteq X$, 
$$m(A)=\inf\{m(V)\ :\ V\textrm{ is open and }A\subseteq V\}\;.$$
\end{enumerate}

\medskip

Let $\mathcal{L}^n$ be the Lebesgue measure on $\R^n$. 

We give now a general construction, due to Caratheodory, which permits to obtain many ''geometric'' measures, like the Hausdorff measure or the spherical measure.

Let $\mathcal{F}$ be a collection of sets in $\R^n$ and let $\zeta:\mathcal{F}\to[0,+\infty]$ be a function, called the \emph{gauge} of the measure we are going to construct. We define the \emph{preliminary measure} $\phi_\delta$, with $0<\delta\leq+\infty$, as follows: if $A\subseteq\R^n$, then 
$$\phi_\delta(A)=\inf\left\{\sum_{S\in\mathcal{G}}\zeta(S)\ :\ \mathcal{G}\subset\mathcal{F}\cap\{\mathrm{diam}(S)\leq\delta\}\textrm{ and }A\subset\bigcup_{S\in\mathcal{G}}S\right\}\;.$$
By definition $\phi_{\delta_1}\geq\phi_{\delta_2}$ if $0<\delta_1\leq\delta_2$, thus we may set
$$\psi(A)=\lim_{\delta\to0}\phi_\delta(A)=\sup_{\delta>0}\phi_\delta(A)\;.$$
Clearly $\psi$ is a measure and it can be shown that every open set is $\psi-$measurable; indeed, one has
$$\phi_\delta(A\cup B)\geq\phi_\delta(A)+\phi_\delta(B)$$
with $0<\delta<\dist(A,B)$.

It is not difficult to show that, if $\mathcal{F}$ is the family of all Borel sets, then every subset of $\R^n$ is contained in a Borel set with the same $\phi_\delta$ measure i.e. $\psi$ is a regular Borel measure.

\medskip

For any non-negative real number $k$ let
$$
\Omega_k=\frac{[\Gamma(1/2)]^k}{\Gamma(k/2+1)}\;.
$$
where $\Gamma$ is the Gamma function. If $k$ is integer, $\Omega_k$ is the volume of the unit ball in the Euclidean $k-$space.  

Let $X$ be a metric space, $\mathcal{F}$ is the family of all subsets $S$ of $X$ and the gauge function 
$$
\zeta_1(S)=\Omega_k2^{-k}(\mathrm{diam}S)^k
$$
for $S\neq\emptyset$ and $\zeta_1(\emptyset)=0$ 

The resulting measure, obtained by the Caratheodory construction is nothing but the \emph{$k-$dimensional Hausdorff measure} on $X$.

We observe that, if $X=\R^n$, we obtain the same measure replacing $\mathcal{F}$ by the collection of open (closed or convex ) subssets.

It is immediate that $\H^0$ is the counting measure.

\begin{Prp}For $0\leq s<t<+\infty$ and $A\subseteq X$, we have that
\begin{enumerate}
\item $\H^s(A)<\infty$ implies that $\H^t(A)=0$;
\item $\H^t(A)>0$ implies that $\H^s(A)=\infty$.
\end{enumerate}\end{Prp}

The \emph{Hausdorff dimension} of a set $A\subseteq X$ is
\begin{eqnarray*}\dim_{\H}A&=&\sup\{s\ :\ \H^s(A)>0\}=\sup\{s\ :\ \H^s(A)=\infty\}\\
 &=&\inf\{t\ :\ \H^t(A)<\infty\}=\inf\{t\ :\ \H^t(A)=0\}\;.\end{eqnarray*}

\medskip

The Caratheodory method can be used to construct the \emph{$k-$dimensional spherical measure} $\mathcal{S}^k$: in this case  $\mathcal{F}$ is the set of all closed balls in $\R^n$ and $\zeta=\zeta_1$ as before. We have $\H^k\leq \mathcal{S}^k\leq 2^k\H^k$.

\medskip

The \emph{$k-$dimensional upper density} of a measure $m$ at a point $p$ is defined by 
$$\Theta^{*k}(m,p)=\limsup_{r\to0}\frac{m(B(p,r))}{\Omega_kr^k}\;,$$
where $B(x,d)$ is the ball of center $x$ and radius $d$. Similarly, 
$$
\Theta_*^k(m,p)=\liminf_{r\to0}\frac{m(B(p,r))}{\Omega_kr^k}\;
$$
defines the \emph{$k-$dimensional lower density}.
\subsection{Rectifiable sets}

Let $k$ be an integer with $1\leq k$. A set $A\subset X$ is said to be \emph{countably $k-$rectifiable} if 
$$S\subset S_0\cup\left(\bigcup_{j=1}^\infty f_j(A_j)\right)$$
where $\H^k(S_0)=0$, $A_j\subseteq\R^k$ and $f_j:A_j\to X$ are Lipschitz functions for $j=1,2, \ldots$.

Usually, we will also require $S$ to be $\H^k-$measurable and sometimes to have locally finite $\H^k$ measure. This definition differs from the one given \cite{federer1} because of the presence of the set $S_0$ (see \cite{ambrosio1}, \cite{krantz1}).

We note that, if $X=\R^n$, the maps $f_j$ can be thought to be defined all over $\R^k$, withouth affecting their Lipschitz constants in view of Kirszbraun's theorem (see \cite{federer1}).

We summarise some results about rectifiable sets in euclidean spaces, so we set $X=\R^n$.

\begin{Lmm} The set $S$ is countably $k-$rectifiable if and only if $S\subseteq\bigcup_j T_j$ where $\H^k(T_0)=0$ and each $T_j$, $j\leq 1$ is a $k-$dimensional, embedded $\Ci^1-$manifold of $\R^n$.\end{Lmm}

\begin{Prp} If $S$ is $\H^k-$measurable and countably $\H^k-$rectifiable, then $S=\bigcup_j S_j$ with $\H^k(S_0)=0$, $S_i\cap S_j=\emptyset$ if $i\neq j$ and for $j\geq1$ $S_j\subseteq T_j$ and $T_j$ is a $k-$dimensional, embedded $\Ci^1-$manifold of $\R^k$.\end{Prp}

\section{Classical theory of currents}\label{sec_class_cur}

We introduce the basic definitions, notations and facts about the classical theory of currents referring to the book by Federer \cite{federer1} for a systematic exposition of the theory. Some more accessible texts are \cite{krantz1, morgan1}.

\subsection{Mass, comass and flat norm for forms}
Let $Y$ be a differentiable manifold of dimension $n$. A {\it $k$-vector} at a point $x\in Y$ is an element of $\bigwedge^k T_xY$. A {\it $k$-vector field} on $Y$ is a smooth section of $\bigwedge^k TY$. A {\it $k$-covector} at a point at a point $x\in Y$ is an element of $\bigwedge^k T^\ast_xY$. A {\it $k$-covector} or (smooth) {\it $k$-form} on $Y$ is a smooth section of $\bigwedge^k T^\ast Y$.

A vector $v\in \bigwedge^k T_xY$ (a covector $\omega\in \bigwedge^k T^\ast_xY$) is said to be {\it simple} if it is the exterior product of $1$-vectors ($1$-covectors).

For any multiindex $I=\{i_1,\ldots,i_k\}$, $1\le i_s\le n$, $1\le s\le k$ we set $\vert I\vert=k$ and 
$$
\frac{\de}{\de x^I}=\frac{\de}{\de x^{i_1}}\wedge\cdots\wedge \frac{\de}{\de x^{i_k}}\;,
$$
$$
{\rm d}x^I={\rm d}x^{i_1}\wedge\cdots\wedge {\rm d}x^{i_k}\;.
$$
If 
 $$
\omega=\sum_{\vert I\vert=k}\omega_I{\rm d}x^I,\>\>\> v=\sum_{\vert I\vert=k}v^I\frac{\de}{\de x^I}\;,
$$
we define
\begin{equation}\label{COUPLE}
\langle \omega,v\rangle=\sum_{\vert I\vert=k}\omega_I v^I
\end{equation}
Assume now that $Y$ is a riemannian manifold with scalar product $g$ and let 
$$
\sum_{r,s=1}^n g_{rs}dx^r\otimes dx^s
$$
be the $1^{\rm st}$ fundamental form of $(Y,g)$. Denote by $dV_g$ the associated volume form.

Given a $k$-vector $v\in \bigwedge^k T_xY$, for any multindex $I=\{i_1,\ldots,i_k\}$ we define
$$
v_I=v_{i_1\cdots i_k}=g_{i_1s_1}\cdots g_{i_ks_k}v^{s_1\cdots s_k}\;;
$$
the quantities $v_I$ are the components of a $k$-covector. The real number 
\begin{equation}\label{LENGTH}
\vert v\vert=\left(\sum_{\vert I\vert=k}v_I v^I\right)^{1/2}\;.
\end{equation}
is, by definition, the {\it length} of $v$.
In a similar way, given a $k$-covector $\omega\in \bigwedge^k T^\ast_xY$ we define 
$$
\omega^I=\omega^{i_1\cdots i_k}=g^{i_1s_1}\cdots g^{i_ks_k}\omega_{s_1\cdots s_k}\;;
$$
where $\left(g^{ij}\right)$ the inverse of $\left(g_{ij}\right)$.
The {\it length} of $\omega$ is then
\begin{equation}\label{LENGTH1}
\vert \omega\vert=\left(\sum_{\vert I\vert=k}\omega_I \omega^I\right)^{1/2}\;.
\end{equation} 

The {\it mass} of a $k$-vector $v$ at $x$ is defined by
\begin{equation}\label{MASSA}
\Vert v\Vert=\inf\left\{\sum_\alpha \vert v_\alpha \vert:v_\alpha\> {\rm is\,\,simple\,\,\,and}\,\, v=\sum _\alpha v_\alpha\right\}\;.
\end{equation}
Similarly, the {\it comass} of a $k$-covector $\omega$ at $x$ is defined by
\begin{equation}\label{COMASSA}
\Vert \omega\Vert=\sup\left\{\sup \langle\omega, v\rangle:v\> {\rm is\,\,simple\,\,\,and}\,\, \vert v\vert\le1\right\}\;.
\end{equation}
\begin{Rem}\label{MACOMA}
The following holds\medskip
\begin{itemize}
\item[1)] $\Vert v\Vert=\sup\left\{ \langle\omega, v\rangle: \Vert\omega\Vert\le 1\right\}$;\medskip
\item[2)] the norms $\Vert\,\,\Vert$ are equivalent to the norms $\vert\,\,\vert$;\medskip
\item[3)] if $k=1$, $k=n-1$, then from (\ref{COMASSA}) it follows that
$$
\Vert \omega\Vert=\vert \omega\vert=\left(\sum_j\omega_j \omega^j\right)^{1/2}
$$
\item[4)] $\Vert\,\,\Vert=\vert\,\,\vert$ on the simple vectors and covectors.
\end{itemize}
\end{Rem}
Given a $k$-vector field on $Y$ the number 
\begin{equation}\label{MASSA1}
\Vert v\Vert=\sup\limits_{x\in Y}\Vert v(x)\Vert
\end{equation}
is called the {\it mass} of $v$. Similarly, 
\begin{equation}\label{MASSA2}
\Vert \varphi\Vert=\sup\limits_{x\in Y}\Vert \varphi(x)\Vert
\end{equation}
is called the {\it comass} of the $k$-form $\varphi$ on $Y$.

The {\it flat norm} of a $k$-form $\varphi$ is defined by
\begin{equation}\label{FLAT}
\Vert \varphi\Vert_{\rm flat}=\sup\limits_{x\in Y}\max\left\{\Vert\varphi(x)\Vert,\Vert {\rm d}\varphi(x)\Vert\right\}\;.
\end{equation}

\subsection{Currents}
Let us recall the main definitions about currents. Let $Y$ be a differentable manifold of real dimension $n$. We denote by $\mathcal D^k(Y)$ the space of the $k$-forms on $Y$ with compact support, endowed with the usual topology and we put $\mathcal D(Y)=\mathcal D^0(Y)$. A {\it $k$-current} on $Y$ is a continuous linear functional $T:{\mathcal D}^k(Y)\to\C$; $k$ is called the {\it dimension} of $T$, $n-k$ the {\it degree} of $T$. The space of all $k$-currents on $Y$ will be denoted by $\mathcal D(Y)'$. 

\medskip

The {\it differential} or {\it boundary} of $T$ is the $(k+1)$-current $dT$ defined by
$$
dT(\varphi)=T( d\varphi)
$$
for every $\varphi\in{\mathcal D}^{k-1}(Y)$.

\medskip

Let $\mathcal S$ the family of the closed subsets $C$ of $Y$ with the following property: if $x\in C\setminus Y$ there exists a neighbourhood of $x$ such that $T(\varphi)=0$ for every $k$-form $\varphi$ such that $\supp(\varphi)\Subset U$. If $T\neq 0$, then $\mathcal S\neq\emptyset$. The closed subset 
$$
\bigcap\limits_{C\in\mathcal S} C
$$
is called the {\it support} of $T$ and denoted $\supp T$.
 
 \medskip
 
 Let $f:Y\to Y'$ be a smooth map between differentiable manifolds. Let $T$ be a $k$-current on $T$ such that $f\vert_{{\supp}T}$ is a proper map. Given a $k$-form $\varphi\in{\mathcal D}^k(Y')$ consider a function $\chi\in C_c^{\infty}(Y)$ such that $\chi=1$ on a neighbourhood of $\supp T\cap f^{-1}(\supp \varphi)$; then the complex number $T(\chi f^\ast\varphi)$
 is independent of $\chi$ and $\varphi\mapsto T(\chi f^\ast\varphi)$ defines a $k$-current on $Y$ called the {\it direct image} of $T$ and denoted by $f_\ast T$. The operation $f_\ast$ commutes with the boundary i.e. $d\circ f_\ast=f_\ast\circ d$.

\medskip

For $\omega\in\mathcal{D}^m(Y)$, $m\leq k$, we define
$$T\llcorner\omega(\varphi)=T(\omega\wedge\varphi)$$
for any $\varphi\in\mathcal{D}^{k-m}(Y)$. As it is easily seen $T\llcorner\omega$ is a $(k-m)-$current.

\medskip

For $f\in \Ci^0(Y)$, $f\ge 0$, let
\begin{equation}\label{TOTVAR}
\Vert T\Vert(f)=\sup\left\{T(\varphi), \varphi\in\mathcal D^k(Y):\Vert\varphi(x)\Vert\le f(x)\right\};
\end{equation}
if $\Vert T\Vert(f)<+\infty$ for every $f\in \Ci_c^0(Y)$ the current $T$ is said to have {\it locally finite mass}. 

If $\Vert T\Vert(1)<+\infty$ we say that $T$ has a {\it finite mass}. Then $f\mapsto\Vert T\Vert(f)$ is a Radon measure with a density $ d\Vert T\Vert$ such that
$$
\Vert T\Vert(f)=\int_Yf d\Vert T\Vert
$$
for all continuous functions $f$ on $Y.$ By virtue of Radon-Nikodym theorem there exists $\Vert T\Vert$-a.e. a $k$-vector field $\sf T$ such that $\Vert T\Vert$-a.e. 
$\vert \sf T\vert=1$, $\sf T$ is simple and for every $k$-form $\varphi$ one has
$$
T(\varphi)=\int_Y\langle \varphi,{\sf T}\rangle\,d\Vert T\Vert.
$$
The {\it mass} of $T$ is 
$$
{\mass}(T)=\int_Yd\Vert T\Vert=\sup\left\{T(\varphi), \varphi\in\mathcal D^k(Y):\Vert\varphi(x)\Vert\le 1\right\}.
$$
If $T$ is locally of finite mass the {\it mass} of $T$ is defined for every compact subset $K$ of $Y$ by
$$
{\mass}_K(T)=\int_K d\Vert T\Vert=\sup\left\{T(\varphi), \varphi\in\mathcal D^k(Y):\Vert\varphi\Vert_K\le 1\right\}\;.
$$
The vector space of $k-$currents of locally finite mass is denoted by $\mathfrak{M}_{k,\loc}(Y)$, while $\mathfrak{M}_{k}(Y)$ denotes the space of $k-$currents of finite mass.

\paragraph{Example: currents of integration}
Let $V$ be a closed submanifold of a riemannian manifold $Y$. Then
$$
\varphi\mapsto \int_V \varphi
$$
is a $k$-current. If the $k$-volume of $V$ is finite then $T$ has a finite mass ${\mass}(T)=k$-{\rm volume of }$V$.

The same definition if $V$ is replaced by a locally finite integral polyedral $k$-chain $P$. Here polyedral $k$-chains are meant formal sums $\sum_{j\in\Z}n_j\sigma_j$ where $n_j\in\Z$, $j\in\Z$, and $\sigma_j$ are Lipschitz maps $\triangle_k\to Y$ where $\triangle_k$ is the oriented standard $k$-simplex. Currents of this type are said to be {\it currents of integration} and are denoted by $[P]$. 

\medskip

Let ${\sf C}^k(Y;\Z)$ be the group of the locally finite integral polyedral $k$-chains and by $\mathcal P_k(Y;\Z)$ the corresponding group of the integration currents. If $K$ is a compact subset of $Y$, ${\sf C}_K^k(Y;\Z)$ denotes the subgroup of chains of ${\sf C}^k(Y;\Z)$ with support in $K$ and  $\mathcal P_{K,k}(Y;\Z)$ the group of the corresponding currents. If we work with real coefficients the corresponding spaces of currents are denoted by $\mathcal P_k(Y;\R)$ and $\mathcal P_{k,K}(Y;\R)$ respectively. 

\bigskip

\paragraph{Normal currents}

A $k$-current $T$ which is locally of finite mass is said to be {\it locally normal} if $dT $ has a locally finite mass. For these currents we use the seminorm
\begin{equation}\label{SNORM}
{\mathbf{N}}_K(T)={\mass}_K(T)+{\mass}_K(dT)\;,  
\end{equation}
$K$ a compact in $Y$. The space of all locally normal currents is denoted by $\mathfrak{N}_{k,\loc}(Y)$. The space of the normal currents whose support is contained in a compact subset $K$ of $Y$ is denoted by $\mathfrak{N}_{k,K}(Y)$ and the one of normal currents ($T$ and $dT$ of finite mass) by $\mathfrak{N}_k(Y)$.

\paragraph{Rectifiable currents.}

$k$-currents which are limits in mass norm of integration currents are called {\it rectifiable}. We denote by $\mathfrak{R}_{k,K}(Y;\Z)$ the space of all $k$-current $T$ such that for every $\epsilon>0$ there exists $P\in \mathcal P_{k,K}(Y;\Z)$ with ${\mass}_K(T-[P])<\epsilon$.

A rectifiable $k$-current $T$ is an {\it integral current} if $dT$ is rectifiable. In particular, rectifiable current are normal. Spaces of rectifiable $k$-currents on a riemannian manifold $Y$ are denoted by $\mathfrak{I}_{k,K}(Y)$ and $\mathfrak{I}_k(Y)$, $K$ being a compact of $Y$.  

\subsection{Dual flat seminorm for currents.}

The {\it dual flat seminorm} of a $k$-current $T$ is defined by
 \begin{equation}\label{FLAT1}
 \Vert T\Vert_{K,\rm flat}=\sup\left\{\vert T(\varphi)\vert:\Vert \varphi\Vert_{K,\rm flat}\le 1\right\}
 \end{equation}
 If $\Vert T\Vert_{K,\rm flat}<+\infty$ then $\supp T\subset K$ is compact and 
 \begin{equation}\label{FLAT2}
 \Vert dT\Vert_{K,\rm flat}\le\Vert T\Vert_{K,\rm flat}\;.
 \end{equation}
It can be shown that
\begin{equation}\label{FLAT3}
\Vert T\Vert_{K,\rm flat}=\inf_S\left\{{\mass}(T- dS)+{\mass}(S):\supp S\Subset K\right\}
\end{equation}
(and by Hahn-Banach theorem this inf is attained). 

The $\Vert \Vert_{K,\rm flat}$-closure of $\mathfrak{N}_{k,K}(Y)$ is denoted by $\mathfrak{F}_{k,K}(Y)$. The union of all $\mathfrak{F}_{k,K}(Y)$ is denoted by $\mathfrak{F}_k(Y)$; its elements are called {\it flat $k$-chains}. A {\it locally flat $k-$chain} is a current $T$ such that $T\llcorner f\in \mathfrak{F}_k(Y)$ for every $f\in\Ci^\infty_c(Y)$.

\medskip

 Currents with support in $K$of the form
$$
T=R+ dS
$$
where $R$, $S$ are rectifiable, $T\in\mathfrak{R}_{k,K}(Y;\Z)$, $S\in \mathfrak{R}_{k+1,K}(Y;\Z)$ are called {\it integral flat $k$-chains} with support in $K$; the space of such currents is denoted by $\mathcal{F}_{k,K}(Y)$.  The {\it flat norm} of a integral flat $k$-chain $T\in\mathcal F_{k,K}(Y)$ is defined by (\ref{FLAT3}). 

\medskip

We recall some results on flat currents.

\begin{Prp} A flat current of finite mass is mass limit of normal currents.\end{Prp}
\noindent{\bf Proof: } Indeed, given $T\in\mathfrak{F}_m(Y)\cap\mathfrak{M}_m(Y)$ and $\epsilon>0$, we can find $Q_\epsilon$ normal such that $\flat{T-Q_\epsilon}<\epsilon$, so there exists a current $S_\epsilon$ with compact support such that
$$\epsilon>\flat{T-Q_\epsilon}=\mass(T-Q_\epsilon-dS_\epsilon)+\mass(S_\epsilon)\;.$$
As $Q_\epsilon$ is normal and $T$ is of finite mass, also $T-Q_\epsilon$ is of finite mass; this implies that $dS_\epsilon$ is of finite mass, so $S_\epsilon$ is normal.

Moreover, $Q_\epsilon+dS_\epsilon$ is normal and $\mass(T-Q_\epsilon-dS_\epsilon)<\epsilon$, which proves the thesis. $\Box$

\medskip

\begin{Prp}\label{prp_caratt_flat}Given $T\in\mathfrak{F}_{k,K}(Y)$, we can find a $k-$vectorfield $\xi$ and a $(k+1)-$vectorfield $\eta$ on $Y$, supported in $K$, with $L^1(K, dV_g)$ coefficients, such that 
$$T(\omega)=\int\limits_{\supp\omega}\langle\xi,\omega\rangle dV_g+\int\limits_{\supp\omega}\langle\eta, d\omega\rangle dV_g\;.$$
\end{Prp}
\noindent{\bf Proof: } If $K$ is contained in a coordinate chart, we can apply the result for $\R^n$, which can be found in \cite[4.1.18]{federer1}. For a generic compact, it is enough to write it as union of coordinate chart and employ a smooth partition of unity. $\Box$

\begin{Rem}On the other hand, it is obvious that every current of the form
$$T(\omega)=\int\limits_{\supp\omega}\langle\xi,\omega\rangle dV_g+\int\limits_{\supp\omega}\langle\eta, d\omega\rangle dV_g$$
is flat, as it is written as $A+dB$, with $A$, $B$ supported in $K$ and with finite mass. \end{Rem}

\chapter{Currents on singular complex spaces}
\epigraphhead[55]{\epigraph{There are, it has been said, two types of people in the world. There are those who, when presented with a glass that is exactly half full, say: this glass is half full. And then there are those who say: this glass is half empty. The world belongs, however, to those who can look at the glass and say: What's up with this glass? Excuse me? Excuse me? This is my glass? I don't think so. My glass was full! And it was a bigger glass!}{T. Pratchett - \emph{The Truth}}}

This chapter opens with a hurried presentation of the theory of local metric currents on metric spaces, which is explained with more ease and care in \cite{lang1}. We complete this sketchy picture with the detailed proof of the comparison theorem between metric and classical currents on a manifold and with a note on the homology theory given by normal metric currents.

The second section introduces the {\em Dolbeault decomposition} for metric currents i.e. decomposition in currents of type $(p,q)$ and, as a consequence, the definition of the operators $\de$ and $\debar$. The main properties of these operators are investigated, mainly in connection with the other operations on metric currents, i.e. contraction, pushforward, slicing.

The sheafs of metric currents we define in the third section are one of the main reasons why we work with local metric currents. 

An investigation on holomorphic currents, i.e. $\debar-$closed $(p,n)-$currents, is the subject of the last section. 

\section{Local metric currents}

We recall in this section the basic facts of the metric currents on a metric space, firstly developed by Ambrosio and Kirchheim in \cite{ambrosio1}. Our exposition will follow closely \cite{lang1}.

\medskip

Let $X,\ Y$ be metric spaces. We denote by $\Lip(X,Y)$ the set of all Lipschitz maps from $X$ to $Y$ and by $\Lip_L(X,Y)$ the set of $L-$Lipschitz maps i.e. maps $X\to Y$ with Lipschitz constant $L$; $\Lip_\loc(X,Y)$ the set of those maps $f:X\to Y$ such that for every $K\subset X$ compact, $f\vert_K\in\Lip(K,Y)$. If $Y=\C$ we simply we write $(X)$ instead of $(X,\C)$ in the previous symbols. Given $f\in \Lip(X,Y)$, we define
$$
\Lip(f)=\inf\{L\in[0,+\infty)\ :\ f\in\Lip_L(X,Y)\}\;.
$$
Recall that if $A\subset X$ and $f\in \Lip_L(A,\R^n)$, then there exists $\tilde{f}\in\Lip_{\sqrt{n}L}(X,\R^n)$ such that $f=\tilde{f}\vert_A$; moreover, every $f:X\to\C$ bounded and uniformly continuous can be approximated uniformly by Lipschitz functions.

\subsection{Definitions}

In what follows, $X$ will be a locally compact metric space. For a compact subset $K\subset X$ and $L\in [0,\infty)$ we put 
$$
\Lip_{K,L}=\{f\in\Lip_L(X)\ :\ \supp(f)\subset K\}
$$
The union $\D(X)$ of all such spaces is the algebra of all $f\in \Lip(X)$ whose support $\supp(f)$ is compact. We say that $f_j\to f$ in $\D(X)$ if and only if all $f_j$ belong to some $\Lip_{K,L}(X)$ and $f_j\to f$ pointwise on $X$. Similarly, we say that $\pi_j\to\pi$ in $\Lip_\loc(X)$ if and only if for every compact $K$ there exists a constant $L_K$ such that $\pi_j\vert_K\in\Lip_{K,L_K}(X)$ for all $j$ and $\pi_j\to\pi$ pointwise on $K$.

Let
$$
\D^0(X)=\D(X)\qquad \D^m(X)=\D(X)\times[\Lip_\loc(X)]^m
$$
endowed with the product topology.

\medskip

An \emph{$m-$dimensional local metric current} $T$ in $X$ is a function $T:\D^m(X)\to\C$ such that\renewcommand{\labelenumi}{(\arabic{enumi})}
\begin{enumerate}
\item (multilinearity)\label{T_multilinearity} $T$ is $(m+1)-$linear;
\item (continuity)\label{T_continuity} $T(f^j,\pi_1^j,\ldots,\pi_m^j)\to T(f,\pi_1,\ldots, \pi_m)$ if $(f^j,\pi^j)\to(f,\pi)$ in $\D^m(X)$;
\item (locality)\label{T_locality} if $m\geq1$, $T(f,\pi_1,\ldots,\pi_m)=0$ whenever some $\pi_i$ is constant in a neighborhood of $\supp(f)$.
\end{enumerate}
The vector space of all $m-$dimensional local metric currents will be denoted by $\D_m(X)$; we can endow it with a weak topology such that $T_n\to T$ if and only if
$$T_n(f,\pi_1,\ldots, \pi_m)\to T(f,\pi_1,\ldots, \pi_m)\qquad \forall\ (f,\pi_1,\ldots,\pi_m)\in\D^m(X)\;.$$

\medskip

An immediate consequence of the definition is the following strict locality property proved in \cite[Section 2]{lang1}.

\begin{Lmm}Let $m\geq1$, $T\in\D_m(X)$, then $T(f,\pi_1,\ldots,\pi_m)=0$ whenever some $\pi_i$ is constant on $\supp(f)$.\end{Lmm}

Given $(u,v)\in\Lip_\loc(X)\times[\Lip_\loc(X)]^k$ and $T\in \D_m(X)$, with $m\geq k\geq0$, we define the current $T\llcorner(u,v) \in\D_{m-k}(X)$ by
$$(T\llcorner(u,v))(f,g)=T(uf,v,g)=T(uf,v_1,\ldots, v_k,g_1,\ldots, g_{m-k})$$
for $(f,g)\in\D^{m-k}(X)$. We have that
$$(T\llcorner(1,v))\llcorner(1,w)=T\llcorner(1,v,w)\;.$$

\begin{Prp}Suppose $T\in\D_m(X)$, $m\geq 1$ and $(f,\pi_1,\ldots, \pi_m)\in\D^m(X)$, then
\begin{enumerate}
\item if $m\geq2$ and $\pi_i=\pi_j$ with $i\neq j$, we have
$$T(f,\pi_1,\ldots, \pi_m)=0\;.$$
\item for every $g,\ h\in\Lip_\loc(X)$,
$$T(f,gh,\pi_2,\ldots,\pi_m)=T(fg,h,\pi_2,\ldots,\pi_m)+T(fh,g,\pi_2,\ldots,\pi_m)\;.$$
\end{enumerate}
\end{Prp}
For the proof see \cite[Proposition 2.4]{lang1}.

We also recall the following analogue of the chain rule \cite[Theorem 2.5]{lang1}.

\begin{Teo}\label{teo_chain}Suppose $m,n\geq1$, $T\in\D_m(X)$, $f\in\D(X)$, and let $U\subset\R^n$ be an open set, $\pi=(\pi_1,\ldots, \pi_n)\in\Lip_\loc(X,U)$ and $g=(g_1,\ldots, g_m)\in[\Ci^{1,1}(U)]^m$. Then, if $n\geq m$,
$$
T(f,g\circ\pi)=\sum_{\lambda\in\Lambda(n,m)}T\left(f\det\left[(D_{\lambda(k)}g_i)\circ\pi\right]_{i,k=1}^m,\pi_{\lambda(1)},\ldots,\pi_{\lambda(m)}\right)\;.
$$
If $n<m$, then $T(f,g\circ\pi)=0$.\end{Teo}

\medskip

The \emph{support} $\supp(T)$ of a local metric current $T\in\D_m(X)$ is the intersection of all closed sets $C\subset X$ such that $T(f,\pi)=0$ whenever $(f,\pi)\in\D^m(X)$ and $\supp(f)\cap C=\emptyset$. We have the following properties (see \cite[Section 3]{lang1}).
\renewcommand{\labelenumi}{\roman{enumi}. }
\begin{Lmm}Let $T\in\D_m(X)$, then
\begin{enumerate}
\item $\supp(T)$ is the set of all $x\in X$ such that for every $\epsilon>0$ there exists $(f,\pi)\in\D^m(X)$ with $\supp(f)\subset B_\epsilon(x)$ and $T(f,\pi)\neq0$;
\item if $f\vert_{\supp(T)}=0$ then $T(f,\pi)=0$;
\item if $m\geq1$, $T(f,\pi)=0$ whenever some $\pi_i$ is constant on $\{f\neq0\}\cap\supp(T)$.\end{enumerate}
\end{Lmm}

Obviously, given $A\subset X$ closed and $T_A\in\D_m(A)$, we can define $T\in\D_m(X)$ by
$$T(f,\pi_1,\ldots, \pi_m)=T_A(f\vert_A,\pi_1\vert_A,\ldots,\pi_m\vert_A)\;.$$

\begin{Prp}\label{prp_supp_corr}Let $T\in\D_m(X)$ and $A\subset X$ be a locally compact set containing $\supp(T)$. Then there is a unique $T_A\in\D_m(A)$ such that
$$T_A(f,\pi)=T(\tilde{f},\tilde{\pi})$$
whenever $(f,\pi)\in\D^m(A), (\tilde{f},\tilde{\pi})\in\D^m(X)$ and $\tilde{f}\vert_A=f$, $\tilde{\pi}_i\vert_A=\pi_i$. Moreover, $\supp(T)=\supp(T_A)$.\end{Prp}

The \emph{boundary} of $T\in\D_m(X)$, $m\geq1$, is the local current $dT\in\D_{m-1}(X)$ defined by
$$dT(f,\pi_1,\ldots,\pi_{m-1})=T(\sigma,f,\pi_1,\ldots,\pi_{m-1})\qquad\forall\; (f,\pi_1,\ldots,\pi_{m-1})\in\D^{m-1}(X)$$
where $\sigma\in\D(X)$ is any functions such that $\sigma\equiv1$ on $\{f\neq0\}\cap\supp(T)$. One can check that this defines indeed a local metric current such that $\supp(dT)\subset\supp(T)$ and $(dT)_A=d(T_A)$ for every locally closed subset $A$. If $m\geq 2$, locality implies that $ddT=0$. The operator $d:\D_m(X)\to \D_{m-1}(X)$ is linear and weakly continuous.

\begin{Lmm}If $T\in\D_m(X)$ and $(u,v)\in\Lip_\loc(X)\times[\Lip_\loc(X)]^k$, with $m>k$, then
$$(dT)\llcorner(u,v)=T\llcorner(1,u,v)+(-1)^kd(T\llcorner(u,v))\;.$$
\end{Lmm}

Given $T\in\D_m(X)$ a locally compact set set $A\subset X$ containing $\supp(T)$, $Y$ and a locally compact metric space. For every proper map $F\in\Lip_\loc(A,Y)$ the \emph{pushforward of $T$ via $F$} is the current $F_\sharp T\in\D_m(Y)$ given by
$$F_\sharp T(f,\pi_1,\ldots,\pi_m)=T_A(f\circ F,\pi_1\circ F,\ldots,\pi_m\circ F)\;.$$
We have $\supp(F_\sharp T)=F(\supp(T))$ and, if $m\geq1$, $F_\sharp(dT)=d(F_\sharp T)$; moreover, if $Z$ is another locally compact metric space and $G\in\Lip_\loc(Y,Z)$ is proper, we have that $(G\circ F)_\sharp T=G_\sharp F_\sharp T$.

If $F\in\Lip_\loc(X,Y)$ is proper, the operator $F_\sharp:\D_m(X)\to\D_m(Y)$ is linear and weakly continuous.

\subsection{Mass}

Given $T\in\D_m(X)$ and an open set $V\subset X$, we define the \emph{mass} $M_V(T)$ of $T$ in $V$ as the least number $M\in[0,+\infty]$ such that
$$\sum_{\lambda\in\Lambda}T(f_\lambda,\pi^\lambda)\leq M$$
whenever $\Lambda$ is a finite set, $(f_\lambda,\pi^\lambda)\in\D(X)\times[\Lip_1(X)]^m$, $\supp(f_\lambda)\subset V$ and $\sum_{\lambda\in\Lambda}|f_\lambda|\leq 1$. $M_X(T)=M(T)$ will be called the \emph{total mass} of $T$. We set
$$M_{m,\loc}(X)=\{T\in\D_m(X)\ :\ M_V(T)<+\infty\ \forall\;V\Subset X\}$$
$$M_m(X)=\{T\in\D_m(X)\ :\ M(T)<+\infty\}$$
and we define an outer measure $\|T\|$ by
$$\|T\|(A)=\inf\{M_V(T)\ :\ V\subset X\textrm{ is open and }A\subset V\}$$
for every $A\subset X$. We note that the mass is lower semicontinuous with respect to weak convergence: if $T_n\to T$ weakly in $\D_m(X)$, then
$$M_V(T)\leq\liminf_{n\to\infty}M_V(T_n)$$
for every open set $V\subset X$. We have the following results (see \cite[Section 4]{lang1}).

\begin{Prp}For $m\geq0$, $(M_m(X), M)$ is a Banach space.\end{Prp}

\begin{Teo}Let $T\in\D_m(X)$, then
\begin{enumerate}
\item the function $\|T\|:2^X\to[0,+\infty]$ is a Borel regular outer measure;
\item we have $\supp(\|T\|)=\supp(T)$ and $\|T\|(X\setminus\supp(T))=0$;
\item for every open set $V\subset X$,
$$\|T\|(V)=\sup\{\|T\|(K)\ :\ K\subset X\textrm{ is compact and }K\subset V\}\;;$$
\item if $T\in M_{m,\loc}(X)$, then $\|T\|$ is a Radon measure and
$$|T(f,\pi)|\leq\prod_{i=1}^m\Lip(\pi_i\vert_{\supp(f)})\int_X|f|d\|T\|$$
for all $(f,\pi)\in\D^m(X)$.
\end{enumerate}
\end{Teo}

Let $\B^\infty_c(X)$ be the algebra of all bounded complex-valued Borel functions whose support is compact; from the previous Theorem it follows that every $T\in M_{m,\loc}(X)$, we can naturally extends to
$$T:\B^\infty_c(X)\times[\Lip_{\loc}(X)]^m\to\C$$
by a density argument.

\begin{Teo}The extended functional $T$ has the following properties:
\renewcommand{\labelenumi}{(\arabic{enumi}) }
\begin{enumerate}
\item (multilinearity) $T$ is $(m+1)-$linear on $\B^\infty_c(X)\times[\Lip_\loc(X)]^m$;
\item (continuity) $T(f^j,\pi^j)\to T(f,\pi)$ whenever $f,f^j\in\B^\infty_c(X)$, $\sup_j\|f^j\|_\infty<\infty$, $\bigcup_j\supp(f^j)\subset K\Subset X$, $f^j\to f$ pointwise on $X$ and $\pi^j\to\pi$ in $[\Lip_\loc(X)]^m$;
\item (locality) if $m\geq 1$, $T(f,\pi)=0$ whenever some $\pi_i$ is constant on the support of $f\in\B^\infty_c(X)$;
\item for all $(f,\pi)\in\B^\infty_c(X)\times[\Lip_\loc(X)]^m$
$$|T(f,\pi)|\leq\prod_{i=1}^m\Lip(\pi_i\vert_{\supp(f)})\int_X|f|d\|T\|\;.$$
\end{enumerate}\end{Teo}

This extension allows us to define $T\llcorner u$, more generally, for any locally bounded Borel function, in particular for $u=\chi_B$ a characteristic function of some Borel set $B\subset X$. we write $T\llcorner B$ instead of $T\llcorner\chi_B$. Clearly, given $(u,v)\in\B^\infty_\loc(X)\times[\Lip_\loc(X)]^k$ and $T\in M_{m,\loc}(X)$, the following holds for any open set $V\subset X$ 
$$M_V(T\llcorner(u,v))\leq\prod_{i=1}^k(v_i\vert_V)\int_V|u|d\|T\|$$
and $T\llcorner(u,v)\in M_{m-k,\loc}(X)$.

\begin{Lmm} Suppose that $T\in M_{m,\loc}(X)$, $Y$ is a locally compact metric space, $F\in\Lip_{\loc}(X,Y)$ is proper when restricted to $\supp(T)$. Then $F_\sharp T\in M_{m,\loc}(Y)$ and\renewcommand{\labelenumi}{\roman{enumi}. }
\begin{enumerate}
\item for all $(f,\pi)\in\B^\infty_c(Y)\times[\Lip_\loc(Y)]^m$ and $\sigma\in\B^\infty_c(X)$ such that $\sigma=1$ on $\{f\circ F\neq 0\}\cap\supp(T)$,
$$F_\sharp T(f,\pi)=T(\sigma(f\circ F), \pi\circ F)\;;$$
\item for every Borel set $B\subset Y$,
$$M((F_\sharp T)\llcorner B)\leq \Lip(F\vert_{F^{-1}(B)\cap\supp(T)})^m\|T\|(F^{-1}(B))\;.$$
\end{enumerate}
\end{Lmm}

\begin{Lmm}Let $T\in M_{m,\loc}(X)$ and $B\subset X$ either a $\sigma-$finite with respect to $\|T\|$ Borel set or an open set. Then $\|T\|(B)$ is the least number such that
$$\sum_{\lambda\in\Lambda}T(f_\lambda,\pi^\lambda)\leq \|T\|(B)$$
whenever $\Lambda$ is finite, $(f_\lambda,\pi^\lambda)\in\B^\infty_c(X)\times[\Lip_1(X)]^m$ and $\sum_{\lambda}|f^\lambda|\leq \chi_B$. Moreover $\|T\|\llcorner B=\|T\llcorner B\|$ and $\|T\|(B)=M(T\llcorner B)$.\end{Lmm}

A slightly different extension can be realized for a current $T\in M_m(X)$, namely we can extend it to
$$T:\B^\infty(X)\times[\Lip(X)]^m\to\C$$
with the following properties:\renewcommand{\labelenumi}{(\arabic{enumi}) }
\begin{enumerate}
\item (multilinearity) $T$ is $(m+1)-$linear;
\item (continuity) $T(f^j,\pi^j)\to T(f,\pi)$ whenever $(f,\pi),(f^j,\pi^j)\in\B^\infty(X)\times[\Lip_L(X)]^m$, $\sup_j\|f^j\|_\infty$ is finite, $(f^j,\pi^j)\to (f,\pi)$ pointwise on $X$;
\item (locality) if $m\geq 1$, $T(f,\pi)=0$ whenever some $\pi_i$ is constant on the support of $f$;
\item for all $(f,\pi)\in\B^\infty(X)\times[\Lip(X)]^m$
$$|T(f,\pi)|\leq\prod_{i=1}^m\Lip(\pi_i\vert_{\supp(f)})\int_X|f|d\|T\|\;.$$
\end{enumerate}
\renewcommand{\labelenumi}{\roman{enumi}. }
\medskip

For a locally compact metric space $X$, the Banach space $(M_m(X),M)$ is exactly the vector space of $m-$dimensional metric currents defined by Ambrosio and Kirchheim in \cite{ambrosio1}, endowed with the mass norm.

\subsection{Locally normal currents}

Given $T\in\D_m(X)$ and an open set $V\subset X$, we define
$$N_V(T)=M_V(T)+M_V(dT)$$
if $m\geq 1$ and $N_V(T)=M_V(T)$ if $m=0$, and let $N(T)=N_X(T)$. Let
$$N_{m,\loc}(X)=\{ T\in\D_m(X)\ :\ N_V(T)<+\infty\textrm{ for every open set }V\Subset X\}$$
and
$$N_m(X)=\{T\in\D_m(X)\ :\ N(T)<+\infty\}\;.$$
The elements of $N_{m,\loc}(X)$, $N_m(X)$ are called \emph{locally normal} currents, \emph{normal} currents, respectively. $N_m(X)$ is a Banach space with the norm $N$.

Given $T\in N_{m,\loc}(X)$ and $(u,v)\in\Lip_\loc(X)\times[\Lip_\loc(X)]^k$, where $m>k$, we have
$$d(T\llcorner(u,v))=(-1)^k((dT)\llcorner(u,v)-T\llcorner(1,u,v))\;,$$
so that
$$M_V(d(T\llcorner(u,v)))\leq \prod_{i=1}^k\Lip(v_i\vert_V)\left(\int_V|u|d\|dT\|+\Lip(u\vert_V)\|T\|(V)\right)$$
for every open set $V\subset X$. Therefore $T\llcorner(u,v)\in N_{m-k,\loc}(X)$. We note that push-forwards of locally normal currents are locally normal. We state now a usefull property of the locally normal currents (see \cite[Lemma 5.2]{lang1}, or \cite[Proposition 5.1]{ambrosio1} for the case of normal currents).

\begin{Prp}\label{prp_norm_cont0}Let $T\in N_{m,\loc}(X)$, $m\geq1$, then
\begin{enumerate}
\item for every $(f,g)\in\D^m(X)$ with $g_2,\ldots, g_m\in\Lip_1(X)$,
$$|T(f,g)|\leq\Lip(f)\int_{\supp(f)}|g_1|d\|T\|+\int_X|fg_1|d\|dT\|\;;$$
\item for all $(f,g), (\tilde{f},\tilde{g})\in\D(X)\times[\Lip_1(X)]^m$,
$$|T(f,g)-T(\tilde{f},\tilde{g})|\leq\int_X|f-\tilde{f}|d\|T\|+$$
$$+\sum_{i=1}^m\left(\Lip(f)\int_{\supp(f)}|g_i-\tilde{g}_i|d\|T\|+\int_X|f||g_i-\tilde{g}_i|d\|dT\|\right)\;.$$
\end{enumerate}\end{Prp}

Indeed, this proposition is true under slightly weaker assumptions.

\begin{Prp}\label{prp_norm_cont}
Let $T:\D^k(X)\to\C$ be a multilinear local functional, with locallly finite mass, for which the product rule holds; furthermore, let us suppose that $dT$ is of locally finite mass too. Then $T$ is  continuous.\end{Prp}
\noindent{\bf Proof: } Let $f\in \D(X)$, $\pi_1,\pi_1'\in\Lip_\loc(X)$ and $\pi=(\pi_2,\ldots,\pi_k)\in[\Lip_\loc(X)]^{k-1}$; then
$$\begin{array}{rcl}& &T(f,\pi_1,\pi)-T(f,\pi_1',\pi)\\&=&T(1,f\pi_1,\pi)-T(1,f\pi_1',\pi)-T(\pi_1,f,\pi)+T(\pi_1',f,\pi)\\
&=&dT(f\pi_1,\pi)-dT(f\pi_1',\pi)-T(\pi_1,f,\pi)+T(\pi_1',f,\pi).\end{array}$$
Using the locality property, $|T(f,\pi_1,\pi)-T(f,\pi_1',\pi)|$ can be esimated with
$$\left(\int_{\supp f}|f||\pi_1-\pi_1'|d\|dT\|+\Lip(f)\int_{\supp f}|\pi_1-\pi_1'|d\|T\|\right)\prod_{j=2}^k\Lip(\pi_j\vert_{\supp f})$$
Repeating this argument for every component, we have that
$$|T(f,\pi)-T(f',\pi')|\leq\int_U|f-f'|d\|T\|+$$
$$+\sum_{i=1}^k\left(\int_U|f||\pi_i-\pi_i'|d\|dT\|+\Lip(f)\int_{\supp f}|\pi_i-\pi_i'|d\|T\|\right)\prod_{j\neq i}\Lip(\pi_j\vert_{\supp f})$$
As $\supp f$ is compact, if the components of $\omega^i$ converge pointwise to the components of $\omega$ on $X$, they converge uniformly on $\supp f$, therefore, applying the previous estimate to the difference $T(\omega_i)-T(\omega)$ we obtain that $T(\omega_i)$ converges to $T(\omega)$, since the integrands in the right-hand side all converge uniformly to $0$ on $\supp f$. $\Box$

Finally, let us we give a characterization of top dimensional locally normal currents in Euclidean spaces. In order to formulate the result, we define the top dimensional current associated to a $L^1_\loc$ function in an Euclidean space: given $U\subset\R^n$ open and $u\in L^1_\loc(U)$, we define
$$[u](f,\pi_1,\ldots, \pi_n)=\int_{\supp f} uf \det\mathrm{Jac}(\pi)d\mathcal{L}^n\;.$$
By Rademacher's theorem, $\pi_1,\ldots, \pi_n$ are almost everywhere differentiable and their differential are bounded on $\supp f$, therefore the Jacobian is well defined and its product with $uf$ is integrable. The current $[u]$ is obviously multilinear, local and of locally finite mass. The continuity follows by a well known properties of $W^{1,\infty}_\loc$ functions. Therefore, $[u]\in M_{n,\loc}(\R^n)$.

\begin{Teo} Let $T\in N_{k,\loc}(U)$. Then there exists $u\in BV_\loc(U)$ such that $T=[u]$ and $\|dT\|=|Du|$.\end{Teo}

\subsection{Rectifiable currents and slicing}

Let $\H^k$ be the $k-$dimensional Hausdorff measure on $X$. A current $T\in M_{k,\loc}(X)$ is said to be \emph{locally rectifiable} if 
\begin{enumerate}
\item $\|T\|$ is concentrated on a countably $\H^k-$rectifiable set;
\item $\|T\|$ vanishes on $\H^k-$negligible Borel sets.
\end{enumerate}
The space of locally rectifiable currents is denoted by $\mathcal{R}_{k,\loc}(X)$, the elements of $\mathcal{R}_k(X)=\mathcal{R}_{k,\loc}(X)\cap M_k(X)$ are called \emph{rectifiable} currents.

We say that a locally rectifiable current $T$ is \emph{locally integer rectifiable} if for any $\phi\in\Lip(X,\R^k)$ and compact Borel  $B$ set in $X$ one has $\phi_\sharp(T\llcorner B)=[u]$ with $u\in L^1(\R^k, \Z)$. The space of such currents is denoted by $\mathcal{I}_{k,\loc}(X)$; the elements of $\mathcal{I}_k(X)=\mathcal{I}_{k,\loc}(X)\cap M_k(X)$ are called \emph{integer rectifiable} currents.

We define the spaces of \emph{locally integral} and \emph{integral} currents as follows.
$$I_{k,\loc}(X)=\{T\in\mathcal{I}_{k,\loc}(X)\ :\ dT\in\mathcal{I}_{k-1,\loc}(X)\}\;,$$
$$I_k(X)=I_{k,\loc}\cap N_k(X)\;:$$

\begin{Teo}Let $T\in M_{k,\loc}(X)$, $k\geq1$. Then $T\in\mathcal{R}_{k,\loc}(X)$ (resp. $T\in\mathcal{I}_{k,\loc}(X)$) if and only if there exist a sequence $\{K_i\}$ of compact sets in $\R^k$, a sequence $\{\theta_i\}$ of functions in $L^1(\R^k,\R)$ (resp. $L^1(\R^k,\Z)$) with $\supp \theta_i\subset K_i$ and a sequence $\{f_i\}$ of bi-Lipschitz maps $f_i:K_i\to X$ such that 
$$\|T\|(A)=\sum_i\|(f_i)_\sharp[\theta_i]\|(A)\qquad T(f,\pi)=\sum_i(f_i)_\sharp[\theta_i](f,\pi)$$
for every Borel set $A\subset X$ and for every $(f,\pi)\in\D^k(X)$.\end{Teo}

\medskip

Given $T\in M_{k,\loc}(X)$, $\pi\in\Lip_\loc(X,\R^m)$, we define the \emph{slice} $\langle T,\pi,x\rangle\in \D_{k-m}(X)$ by
$$\langle T,\pi,x\rangle(f,\eta)=\lim_{\epsilon\to0}T(f\rho_\epsilon\circ\pi, \pi, \eta)$$
where $\rho_\epsilon$ is any family of mollifiers, for every $x\in\R^m$ for which the limit exists.

\begin{Teo}\label{teo_ext_slice}If $T\in N_{k,\loc}(X)$, with $\supp T$ separable, $\pi\in Lip_\loc(X,\R^m)$, then
\begin{enumerate}
\item for $\mathcal{L}^m-$almost every $x\in\R^m$, the slice $\langle T,\pi,x\rangle$ exists and is locally normal and $d\langle T,\pi, x\rangle=(-1)^m\langle dT, \pi, x\rangle$;
\item for all $(f,g)\in\mathcal{B}^\infty_c\times[\Lip_\loc(X)]^{k-m}$, 
$$\int_{\R^m}\langle T,\pi,x\rangle(f,g)dx=T\llcorner(1,\pi)(f,g)\;;$$
\item for every $\|T\llcorner(1,\pi)\|-$measurable set $B\subset X$, 
$$\int_{\R^m}\|\langle T,\pi, x\|(B)dx=\|T\llcorner(1,\pi)\|(B)\;.$$
\end{enumerate}
\end{Teo}

\begin{Teo} Suppose $m,m'>1$, $k>m+m'$, $\pi\in\Lip(X,\R^m)$, $\pi'\in\Lip(X,\R^{m'})$, $T\in N_{k,\loc}(X)$ and $\supp T$ is separable. Then
$$\langle T,(\pi,\pi'), (x,x')\rangle=\left\langle \langle T, \pi, x\rangle,\pi', x'\right\rangle$$
for $\mathcal{L}^{m+m'}-$almost every $(x,x')\in\R^{m+m'}$.\end{Teo}

We give some results which allow us to recover informations on the rectifiability of a current from the rectifiability of the slices.

\begin{Teo} Let $T\in N_{k,\loc}(X)$, $k\geq1$ be such that $\supp T$ is separable and $\mathscr{P}$ a countable subset of $\Lip_1(X)$ which is dense in the uniform norm on compact sets. If, for each $\pi\in\mathscr{P}^k$, $\langle T,\pi, x\rangle\in \mathcal{R}_{0,\loc}(X)$  (resp. $\mathcal{I}_{0,\loc}(X)$) for $\mathcal{L}^k-$almost every $x\in\R^k$, then $T\in\mathcal{R}_{k,\loc}(X)$ (resp. $\mathcal{I}_{k,\loc}(X)$). \end{Teo}

\begin{Teo} Let $T\in N_{k,\loc}(X)$, $1\leq m\leq k$, be such that $\supp T$ is separable. Then 
\begin{enumerate}
\item if $T\in\mathcal{R}_{k,\loc}(X)$ (resp. $\mathcal{I}_{k,\loc}(X)$) and $\pi\in \Lip_\loc(X,\R^m)$, then $\langle T,\pi, x\rangle\in\mathcal{R}_{k-m,\loc}(X)$ (resp. $\mathcal{I}_{k-m,\loc}(X)$) for $\mathcal{L}^m-$almost every $x\in\R^k$;
\item conversely, if for each $\pi\in \Lip_\loc(X,\R^m)$, $\langle T,\pi, x\rangle\in\mathcal{R}_{k-m,\loc}(X)$ (resp. $\mathcal{I}_{k-m,\loc}(X)$) for $\mathcal{L}^m-$almost every $x\in\R^k$, then $T\in\mathcal{R}_{k,\loc}(X)$ (resp. $\mathcal{I}_{k,\loc}(X)$).
\end{enumerate}\end{Teo}

Finally, we give an important characterization of the space of integral currents.

\begin{Teo} $I_{k,\loc}(X)=\mathcal{I}_{k,\loc}(X)\cap N_{k,\loc}(X)$, for all $k\geq 0$. \end{Teo}

This result is often known as \emph{boundary rectifiability theorem}, as it can be stated as follows: if $T\in\mathcal{I}_{k\,loc}(X)$ and $dT$ has locally finite mass, then $dT\in\mathcal{I}_{k-1,\loc}(X)$. 

\subsection{Comparison theorem}\label{sec_comp}

With the notations introduced in section \ref{sec_class_cur} we state now the following
\begin{Teo} \label{teo_comp}Let $U\subset\C^N$ be an open set, $N\geq1$. For every $m\geq0$ there exists an injective linear map $C_m:\D_m(U)\to\Di_m(U)$ such that
$$C_m(T)(fdg_1\wedge\ldots\wedge dg_m)=T(f,g_1,\ldots, g_m)$$
for all $(f,g_1,\ldots, g_m)\in\Ci^\infty_c(U)\times[\Ci^\infty(U)]^m$. Moreover
\begin{enumerate}
\item for $m\geq1$, $d\circ C_m=C_{m-1}\circ d$;
\item for all $T\in\D_m(U)$, $\|T\|\leq\mass(C_m(T))\leq{N\choose m}\|T\|$;
\item the restriction of $C_m$ to $N_{m,\loc}(U)$ is an isomorphism onto $\mathfrak{N}_{m,\loc}(U)$;
\item the image of $C_m$ contains the space $\mathfrak{F}_{m,\loc}(U)$.\end{enumerate}
\end{Teo}
For the proof see \cite[Theorem 5.5]{lang1} and \cite[Theorem 11.1]{ambrosio1}.
Analogous result holds for manifolds:

\begin{Teo}\label{teo_comp_man}Let $U$ be an $N-$dimensional complex manifold, $N\geq 1$. For every $m\geq0$ there exists an injective linear map $C_m:\D_m(U)\to \Di_m(U)$ such that
$$C_m(T)(fdg_1\wedge\ldots\wedge dg_m)=T(f,g_1,\ldots, g_m)$$
for all $(f,g_1,\ldots, g_m)\in\Ci^\infty_c(U)\times[\Ci^\infty(U)]^m$. The following properties hold:
\begin{enumerate}
\item for $m\geq1$, $d\circ C_m=C_{m-1}\circ d$;
\item there exists a positive constant $c_1$ such that, for all $T\in\D_m(U)$
$$
c_1^{-2}\|T\|\leq\mass(C_m(T))\leq c_1^{2}{N\choose m}\|T\|;
$$
\item the restriction of $C_m$ to $N_{m,\loc}(U)$ is an isomorphism onto $\mathfrak{N}_{m,\loc}(U)$;
\item the image of $C_m$ contains the space $\mathfrak{F}_{m,\loc}(U)$.\end{enumerate}
\end{Teo}
\noindent{\bf Proof: }Fix a locally finite covering $\{U_j\}_{j\in\N}$ of relatively compact open sets with bi-Lipschitz coordinate charts $\phi_j:U_j\to\Omega_j\subseteq\C^N$.

Let $\{\rho_j\}_{j\in\N}$ be a smooth partition of unity subordinated to the covering $\{U_j\}_{j\in\N}$; for each $j\in\N$, the current $T_j=T\llcorner\rho_j$ is supported in $U_j$, therefore belongs to $\D_m(U_j)$.

The induced map $(\phi_j)_\sharp$ is an isomorphism between $\D_m(U_j)$ and $\D_m(\Omega_j)$; let $C_m^j$ be the linear injective map given by Theorem \ref{teo_comp} between $\D_m(\Omega_j)$ and $\Di_m(\Omega_j)$. Then $(\phi_j^{-1})_*$ is an isomorphism between $\Di_m(\Omega_j)$ and $\Di_m(U_j)$, which can be injected into $\Di_m(U)$.

Therefore, for every $j\in\N$, we have the map 
$$T\mapsto T_j\mapsto (\phi_j)_\sharp T_j\mapsto C_m^j((\phi_j)_\sharp T_j)\mapsto (\phi_j^{-1})_*C_m^j((\phi_j)_\sharp T_j)=R_m^j(T)\;.$$
All the intermediate steps are linear and injective, so the result is linear and injective. We set
$$C_m(T)=\sum_j R_m^j(T)\;,$$
which is well defined because the covering is locally finite; moreover, $C_m$ is linear and injective.

Given $(f,g_1,\ldots, g_m)\in\Ci^\infty_c(U)\times[\Ci^\infty(U)]^m$, we have that
$$T(f,g)=\sum_j T(\rho_j f, g)=\sum_j (\phi_j)_\sharp T_j ((\rho_j\circ\phi_j^{-1})\cdot(f\circ\phi_j^{-1}), g_1\circ\phi_j^{-1},\ldots, g_m\circ\phi_j^{-1})$$
where the sums are indeed finite, because $f$ has compact support. Now, by Theorem \ref{teo_comp} the last sum is equals
$$\sum_j C_m^j((\phi_j)_\sharp T_j)((\rho_j f)\circ\phi_j^{-1}dg_1\circ\phi_j^{-1}\wedge\ldots\wedge dg_m\circ\phi_j^{-1})=$$
$$\sum_j(\phi_j^{-1})_*C_m^j((\phi_j)_\sharp T_j)(\rho_j fdg_1\wedge\ldots \wedge dg_m)=\sum R_m^j(T)(fdg_1\wedge\ldots \wedge dg_m)=$$
$$C_m(fdg_1\wedge\ldots \wedge dg_m)\;.$$

\medskip

The conclusions of Theorem \ref{teo_comp} hold for $C_m^i$ and for the pushforward maps, so we just have to check that they still hold after contraction with a $0-$form and after a locally finite sum.

\smallskip

\noindent{\emph{i. }} We know that, for $m\geq1$, $d\circ C_m^j=C_{m-1}^j\circ d$. We also have
$$d(T_j)=d(T_\llcorner\rho_j)=dT\llcorner\rho_j+T\llcorner(\sigma_j,\rho_j)=(dT)_j+T\llcorner(\sigma_j,\rho_j)\;,$$
with $\sigma_j$ a compactly supported smooth function equal to $1$ on $\supp\rho_j$ and to $0$ outside $U_j$, so
$$R^j_{m-1}(dT)= (\phi_j^{-1})_*C_{m-1}^j((\phi_j)_\sharp (dT)_j)=(\phi_j^{-1})_*C_{m-1}^j((\phi_j)_\sharp (d(T_j)-T\llcorner(\sigma_j,\rho_j)))=$$
$$(\phi^{-1}_j)_*C_{m-1}^j(d(\phi_j)_\sharp T_j-((\phi_j)_\sharp T\llcorner(\sigma_j,\rho_j)))=(\phi_j^{-1})_*(dC_m^j((\phi_j)_\sharp T_j)-C_{m-1}^j((\phi_j)_\sharp T\llcorner (\sigma_j,\rho_j)))=$$
$$dR_m^j(T)-S_{m-1}^j(T)\;.$$
For a given classical form $fdg_1\wedge\ldots\wedge dg_{m-1}$, we have that $S_{m-1}^j(T)(fdg_1\wedge\ldots\wedge dg_{m-1})\neq0$ only for a finite number of $j\in\N$ (namely, those such that $fd\rho_j \neq0$), so
$$\sum_{j\ :\ fd\rho_j\neq0}\!\!\!\!S_{m-1}^i(T)(fdg_1\wedge\ldots\wedge dg_{m-1})=\sum(\phi_j^{-1})_* C_{m-1}^j((\phi_j)_\sharp T\llcorner (\sigma_j,\rho_j)))(fdg_1\wedge\ldots\wedge dg_{m-1})=$$
$$\sum T\llcorner (\sigma_j,\rho_j)(f,g_1,\ldots, g_{m-1})\;.$$
We can replace $\sigma_j$ with a $\sigma$, independent of $j$, defined by
$$\sigma=\max\{\sigma_j\ \vert\ fd\rho_j\neq0\}\;.$$
Thus the last sum is equal to
$$\sum T\llcorner (\sigma_j,\rho_j)(f,g_1,\ldots, g_{m-1})=T(\sigma f, \sum\rho_j, g_1,\ldots, g_{m-1})=0$$
 because $\sum\rho_j$ is constantly equal to $1$ on the support of $\sigma f$ (which coincides with $\supp f$).

Therefore $d\circ R_{m-1}^j=R_{m}^j\circ T$ and then obviously $d\circ C_{m-1}=C_m\circ d$.

\smallskip

\noindent{\emph{ii. }} Upon taking a refinement of our open covering, we can assume that there exists a positive constant $c_1$ such that $ \Lip(\phi_j), \Lip(\phi_j^{-1})\leq c_1$. Denoting $\mu$ the (metric) mass measure of $T$ and $m_j$ the (classical) mass measure of $R_m^j(T)$ we have
$$\|T\llcorner\rho_j\|(B)=(\mu\cdot\rho_j)(B)\;,$$
for every Borel set $B$, and
$$(\phi_j)_\sharp(\mu\cdot \rho_j)(B')\leq c_1(\mu\cdot\rho_j)(\phi_j^{-1}(B'))\;,$$
so
$$m_j(B)\leq c_1^2{N\choose m}(\mu\cdot\rho_j)(B)\;.$$
Summing on $j$ and denoting $\mu'$ the mass measure of $C_m(T)$, we obtain
$$\mu'(B)\leq c_1^2{N\choose m}\mu(B)\;.$$
To obtain the other estimate, we observe that
$$C_m(T)\llcorner \rho_j=R_m^j(T)$$
because
$$C_m(T)\llcorner\rho_j(fdg_1\wedge\ldots\wedge dg_m)=C_m(T)(f\rho_jdg_1\wedge\ldots\wedge dg_m)=T(f\rho_j, g_1,\ldots, g_m)=$$
$$R_m^j(T)(fdg_1\wedge\ldots\wedge dg_m)\;;$$
therefore
$$C_m^{-1}=\sum_j(\phi_j)_\sharp\circ (C_m^j)^{-1}\circ(\phi_j^{-1})_*(T\llcorner\rho_j)$$
with $T$ a classical current. So we obtain the estimate
$$\mu(B)\leq c_1^2\mu'(B)\;.$$

\smallskip

\noindent{\emph{iii.} and \emph{iv. }} The class of locally normal currents is stable under pushforward and contraction by a smooth function and the same is true for locally flat currents. Therefore these two points follow easily from the corresponding ones in \ref{teo_comp}. $\Box$

\subsection{Homology of normal currents}

Given metric space $X$, we can consider the chain complex
$$\ldots\longrightarrow N_k(X)\stackrel{d}{\longrightarrow} N_{k-1}(X)\longrightarrow\ldots\longrightarrow N_1(X)\stackrel{d}{\longrightarrow}N_0(X)\longrightarrow 0$$
where $N_k(X)$ is the space of normal metric currents with compact support, and the associated homology
$$H_k(X)=\frac{\mathrm{Ker}\{d:N_k(X)\to N_{k-1}(X)\}}{\mathrm{Im}\{d:N_{k+1}(X)\to N_k(X)\}}$$
Obviously, if $f:X\to Y$ is a Lipschitz map, we obtain the pushforward operator $f_\sharp:N_k(X)\to N_k(Y)$ for every $k$ and, since  $f_\sharp$ and $d$ commute, we have an induced operator
$$H(f):H_k(X)\to H_k(Y)$$
such that $H(\mathrm{Id})=\mathrm{Id}$ and $H(f\circ g)=H(f)\circ H(g)$. In other words, $H$ is a covariant functor from the category of metric spaces with Lipschitz functions to the category of abelian groups. In what follows we will write $f_*$ instead of $H(f)$.

Moreover, if $A$ is closed subset of $X$,we define $N_k(X,A)$ setting
$$
N_k(X)/N_k(A).
$$
Since $d:N_k(X)\to N_{k-1}(X)$ sends $N_k(A)$ in $N_{k-1}(A)$ we can consider the relative homology groups $H_k(X,A)$ and we have the long exact sequence of the pair, as well as for singular homology
$$\ldots H_k(A)\to H_k(X)\to H_k(X,A)\stackrel{d'}{\to}H_{k-1}(A)\to H_{k-1}(X)\to\ldots $$
where $d'$ is an homomorphism of degree $-1$.

\begin{Prp} \label{prp_mayer_vietoris}
Let $\{U,V\}$ be an open covering of $X$, let $i_U, i_V$ be the inclusions of $U\cap V$ in $U$ and $V$ respectively and let $j_U, j_V$ be the inclusions of $U$ and $V$ respectively in $X$. Then the short sequence of chain complexes
$$0\to N_*(\bar U\cap \bar V)\stackrel{(i_U)_*\oplus(i_V)_*}{\longrightarrow}N_*(\bar U)\oplus N_*(\bar V)\stackrel{(j_U)_*-(j_V)_*}{\longrightarrow}N_*(X)\to0$$
is exact.
\end{Prp}
\noindent{\bf Proof: }
Given $T\in N_k(U\cap V)$ with $(i_U)_\sharp(T)=0$, for every form $(f,\pi)\in\D^k(\bar U)$ we have
$$T(f\vert_{U\cap V},\pi\vert_{U\cap V})=0$$
so $T(g,\eta)=0$ for every $(g,\eta)\in\D^k(\bar{U}\cap\bar{V})$, that is $T=0$.

Moreover, if $(j_U)_\sharp(T)=(j_V)_\sharp(S)$, with $T\in N_k(\bar U)$ and $S\in N_k(\bar V)$, then $\supp ((j_U)_\sharp(T))=\supp((j_V)_\sharp(S))\subseteq \bar U\cap \bar V$; this means that $T=(i_U)_\sharp R$ and $S=(i_V)_\sharp R$ with $R\in N_k(\bar U\cap \bar V)$.

Finally, given $T\in N_k(X)$, we can consider a partition of unity subordinated to the covering $\{U, V\}$, $\{\phi_U,\phi_V\}$. The current $T\llcorner\phi_U$ has support contained in $U$, therefore there is $S_1\in N_k(\bar U)$ such that $T_\llcorner\phi_U=(j_U)_\sharp S_1$; similarly, there is $S_2\in N_k(\bar V)$ such that $-T\llcorner\phi_V=(j_V)_\sharp S_2$. So, we have that $T=(j_U)_\sharp S_1-(j_V)_\sharp S_2$ and the exactness of the sequence follows. $\Box$

\medskip

By employing the usual techniques of homological algebra and Proposition \ref{prp_mayer_vietoris}, we can now prove the Mayer-Vietoris sequence theorem for the homology of normal currents.

\begin{Prp}\label{prp_excision}
Given a closed subset $A$ of $X$ and an open set $U$ such that $\bar{U}$ is contained in the interior of $A$, we have that the inclusion map $(X\setminus U, A\setminus U)\to (X,A)$ induces an isomorphism in homology.
\end{Prp}
\noindent{\bf Proof: }
The result follows from the exactness of Mayer-Vietoris sequence, in the same way as in singular homology. $\Box$

\medskip

\begin{Prp}\label{prp_homot_inv}
Two Lipschitz-homotopic Lipschitz maps 
$$f\sim g:X\to Y$$
induce the same homomorphism in homology
\end{Prp}
\noindent{\bf Proof: }
Let $H:X\times[0,1]\to Y$ be the Lipschitz homotopy between $f$ and $g$ and let us define the operator
$$K:N_k(X)\to N_{k+1}(Y)$$
by the following formula
$$K(T)(f,\pi_1,\ldots, \pi_{k+1})=\sum_{i=1}^{k+1}(-1)^{i+1}\int_0^1T\left(f\circ H \cdot \frac{\de \pi_i\circ H}{\de t},\ldots, \hat{\pi_i},\ldots\right)$$
Arguing like in \cite[Proposition 10.2]{ambrosio1}, we see that if $T\in N_k(X)$, $K(T)$ is also in $N_{k+1}(Y)$ and the following holds
$$d(K(T))=-K((dT))+g_\sharp T-f_\sharp T.$$
Consequently, if $dT=0$, we see that $g_\sharp T-f_\sharp T$ is in the image of $d:N_{k+1}(Y)\to N_k(Y)$, that is $f_*=g_*$ as applications between $H_*(X)$ and $H_*(Y)$. $\Box$

\medskip

\begin{Prp}\label{prp_dim_ax}
If $X$ is a metric space with only one point, we have
$$H_*(X)=\left\{\begin{array}{ll}\mathbb{K}&\textrm{ if }*=0\\
0&\textrm{otherwise}\end{array}\right.$$
where $\mathbb{K}$ is either $\R$ or $\C$.
\end{Prp}
\noindent{\bf Proof: }
The thesis is obvious, as $M_0(X)=\{\alpha \delta_x\ \vert\ \alpha\in\mathbb{K}\}\cong \mathbb{K}$. $\Box$

\medskip

The previous results mean that the functor $H_*$ satisfies the axioms of Eilenberg and Steenrod for homology, therefore $H_*(X)$ is the usual singular homology with real (or complex) coefficients, whenever $X$ is a CW-complex.

\subsection{Examples}\label{sec_ex1}

We give some examples of metric currents.

\paragraph{$0$-currents and Borel measures}

Let $X$ be a metric space, $\mu$ a locally finite Borel measure and $\psi \in L^1(X,\mu)$. We define the functional
$$T(f)=\int_X f\psi d\mu$$
for every $f\in\D^0(X)$. This functional is obviously multilinear; the locality property is empty and if $f_j\to f$ pointwise with $f_j,f\in \Lip_{K,L}$, then there exists $j_0$ such that for every $j>j_0$, $|f(x_0)-f_j(x_0)|\leq \epsilon$, with some given $x_0$, so 
$$|f(x)|, |f_j(x)|\leq |f(x_0)|+\epsilon+L\mathrm{dist}(x,x_0)$$
whence $T(f_j)\to T(f)$. Finally, we have
$$T(f)\leq\int_{\supp f}|f||\psi|d\mu$$
i.e. the mass of $T$ is the measure $|\psi|d\mu$, which is finite.

The same is true if $\psi\in L^1_\loc(X,\mu)$, but then $T$ will be of \emph{locally} finite mass.

\paragraph{Currents of integration on a manifold}

Let $M$ be a real $n-$dimensional riemannian manifold, with volume form $dV$; as a particular case of the previous example, for every $\psi\in L^1(M, dV)$, we have the $0-$current
$$T(f)=\int_M f\psi dV\;.$$
On the other hand, since every compactly supported section of $\Lambda^n T^*M$, with $L^\infty$ coefficients,  is a multiple of $dV$ by a compactly supported $L^\infty$ function (so, in particular, it belongs to $L^1$) we can also define the $n-$current
$$
T(f,g_1,\ldots, g_n)=\int_M fdg_1\wedge\ldots\wedge dg_n.
$$
This functional is clearly multilinear and local, and, in view Rademacher's theorem, we have
$$
|T(f,g)|\leq\int_M|f||\det(\nabla g_1,\ldots, \nabla g_n)|dV\leq\int_M |f|d V\prod_{j=1}^n\|\nabla g_j\|_{\infty,\supp f}
$$
$$
\leq \int_{\supp f}|f|dV\prod_{j=1}^n\Lip(g_j\vert_{\supp f}).
$$
Moreover, by Stokes' theorem
$$
dT(f,g_1,\ldots, g_{n-1})=\int_{bM}fdg_1\wedge\ldots\wedge dg_{n-1}dW
$$
where $dW$ is the volume form induced on $bM$ by the riemannian structure on $M$; in particular, $dT$ is of locally finite mass; therefore, $T$ is locally normal, hence continuous by Proposition \ref{prp_norm_cont}. We denote such a current $T$ by $[M]$. Clearly, $[fM]=[M]\llcorner f$ is again a metric $0-$current, for $f\in \mathcal{B}^\infty_\loc(M)$.

\paragraph{Currents associated to $k-$forms}

Let $M$ as before and let $\alpha$ be a metric $k-$form $(\psi,\eta_1,\ldots,\eta_k)$. We define the $(n-k)-$current $T$ by the following formula
$$T(f,g_1,\ldots, g_{n-k})=[M]\llcorner\alpha(f,g_1,\ldots, g_{n-k})=\int_{M}f\psi d\eta_1\wedge\ldots\wedge d\eta_k\wedge d g_1\wedge\ldots\wedge dg_{n-k}\;.$$
This is obviously a metric functional and we have
$$|T(f,g)|\leq \prod\Lip(g_j\vert_{\supp f})\int_{M}|f||\psi|dV\prod \Lip(\eta_j\vert_{\supp f})\;,$$ 
which means that $T$ is of locally finite mass.

Moreover, an easy computation gives
$$dT=[bM]\llcorner\alpha+(-1)^{k+1}[M]\llcorner(d\alpha)\;,$$
proving that $T$ is locally normal if and only if $\nabla\psi\in L^\infty$ and $bM$ is locally of finite volume.

If $M=\R^n$, this kind of currents can be denoted by $[\alpha]$ and then $d[\alpha]=(-1)^{k+1}[d\alpha]$.

\paragraph{Currents associated to a submanifold}

Let $M$ be as above and consider a $k-$dimensional submanifold $N\subset M$. Then the current
$$T(f,g_1,\ldots, g_k)=\int_N fdg_1\wedge\ldots \wedge dg_k$$
is nothing but the image of $[N]\in \D_k(N)$ under the map $i_\sharp$ induced by $i:N\to M$.

Obviously, we can consider more general objects than submanifolds. For example, if $M$ is a complex manifold, all the complex subspaces of $M$ induce closed currents of locally finite mass by integration on their regular part (cfr \cite{lelong1}).

\paragraph{Currents associated to a vector field}

Let $U$ be an open set of $\R^n$ and let $\xi$ be a compactly supported $k-$vector field with $L^1$ coefficients (with respect to the Lebesgue measure $d\mathcal{L}$). The classical $k-$current
$$T(\omega)=\int_U\langle\xi, \omega\rangle d\mathcal{L}$$
is flat (see Proposition \ref{prp_caratt_flat}), therefore is in the image of $C_k$, so there exists a metric current in $\D_k(U)$ which coincides with $T$ on every smooth form. This is the metric current associated with $\xi$ (with respect to the Lebesgue measure).

\section{Dolbeault decomposition}
\label{sec_dolb_dec}
In what follows, $X$ will be a (reduced) complex analytic space; we can endow $X$ with a metric space structure in several situations, for example, when $X$ can be embedded in $\C^N$ or in $\CP^N$, when $X$ is a K\"ahler space or when $X$ is a Kobayashi-hyperbolic space. The aim of this section is to endow the space of metric currents on $X$ with a complex structure, that is to define a Dolbeault decomposition in $(p,q)-$currents.

\subsection{$(p,q)-$currents}

Given an open subset $U$ of the analytic space $X$ consider the $*-$algebra $\A(U)$ generated by $\Ol(U)$ i.e the subset of $\Lip_\loc(U)$ of complex valued functions of the form 
$$\sum_{i=1}^k f_i\overline{g}_i\qquad f_i,\ g_i\in\Ol(U)\;.
$$
 The elements of
$$
\A^m(U)=\D(U)\times[\A(U)]^m\;
$$
are called \emph{analytic forms} on $U$.
We remark that $\A(U)$ is dense in $\Lip_\loc(U)$ with respect to the topology defined before.

Given $(f,\pi)\in\A^m(U)$, we say that $(f,\pi)$ is of \emph{pure type $(p,q)$} if $p+q=m$ and there is a partition $I\cup J=\{1,\ldots, m\}$, with $|I|=p$ and $|J|=q$, such that $\pi_i$ is holomorphic and $\pi_j$ is antiholomorphic for every $i\in I$, $j\in J$. The vector space of such forms will be denoted by $\A^{p,q}(U)$.

\begin{Prp}\label{prp_forme_pq} Given $T\in\D_m(U)$, for every $(f,\pi)\in\A^m(U)$, there exist $(f_i,\pi^i)\in\A^{p_i,q_i}(U)$ for every $p_i+q_i=m$ such that
$$T(f,\pi)=\sum_{p_i+q_i=m}T(f_i,\pi^i)\;.$$
\end{Prp}
\noindent{\bf Proof: } Since an element $h\in\A(U)$ is of the form
$$h=\sum_{l=1}^kf_i\overline{g}_i$$
with $f_i,\ g_i\in\Ol(U)$. By property (\ref{T_multilinearity}) of local metric currents, we can reduce ourselves to the case when
$$
(f,\pi)=(f,f_1\overline{g}_1,\ldots, f_m\overline{g}_m)\;;
$$
then, applying Theorem \ref{teo_chain}, we obtain
$$T(f,\pi)=\sum_{p+q=m}\sum_{I,J}(-1)^{|IJ|}T(ff_{j_1}\cdots f_{j_q}\overline{g}_{i_1}\cdots \overline{g}_{i_p},f_{i_1},\ldots, f_{i_p},\overline{g}_{j_1},\ldots, \overline{g}_{j_q})$$
where $I\cup J=\{1,\ldots, m\}$ is a partition, with $|I|=p$, $|J|=q$, and $(-1)^{|IJ|}$ is the sign of the permutation  $(1,\ldots, m)\mapsto (I,J)$. $\Box$

\begin{Rem}\label{rem_forme_pq_gen} The same conclusions of Proposition \ref{prp_forme_pq} hold true  if $T:\A^m(U)\to\C$ is multilinear, alternating, respecting the product rule.\end{Rem}

The behaviour of a current on forms of pure type completely determines the current. Indeed

\begin{Prp}\label{prp_0_su_pq}If $T\in\D_m(U)$ and $T(f,\pi)=0$ whenever $(f,\pi)\in\A^{p,q}(U)$ with $p+q=m$, then $T=0$.\end{Prp}
\noindent{\bf Proof: } By Proposition \ref{prp_forme_pq}, we know that $T(f,\pi)=0$ for every $(f,\pi)\in\A^m(U)$. Moreover, given $(f,\pi)\in\D^m(U)$, we can find a sequence $\{(f^j,\pi^j)\}\subset\A^m(U)$ such that $(f^j,\pi^j)\to(f,\pi)$ in $\D^m(U)$. Therefore
$$T(f,\pi)=\lim_{j\to\infty}T(f^j,\pi^j)=0\;.$$
$\Box$

\medskip

We say that a current $T\in\D_m(U)$ is of \emph{bidimension} $(p,q)$ if $T(f,\pi)=0$ for every $(f,\pi)\in\A^{r,s}(U)$ with $(r,s)\neq(p,q)$; the space of $(p,q)-$currents will be denoted by $\D_{p,q}(U)$. Equivalently, mimicking the property (\ref{T_locality}) of local metric currents, we can say that $T\in\D_m(U)$ is of bidimension $(p,q)$ if
$$T(f,\pi_1,\ldots,\pi_m)=0$$
whenever there exists $I\subset\{1,\ldots, m\}$, with $|I|>p$, such that $\pi_i\vert_{\supp(f)}$ is holomorphic for every $i\in I$, or $J\subset\{1,\ldots, m\}$, with $|J|>q$, such that $\pi_i\vert_{\supp(f)}$ is antiholomorphic for every $j\in J$.

\subsection{Properties}\label{ssc_prop_pq}

It is easy to see that $\D_{p,q}(U)$ is closed for the weak topology on $\D_m(U)$ and, by Proposition \ref{prp_0_su_pq}, $\D_{p,q}(U)\cap \D_{r,s}(U)=\{0\}$ as soon as $(p,q)\neq(r,s)$.

\medskip

A \emph{Dolbeault decomposition} of a current $T\in \D_m(U)$ is
$$T=T_1+\ldots+ T_k$$
where $T_i\in\D_{p_i,q_i}(U)$, $p_i+q_i=m$.

\medskip

A Dolbeault decomposition for a current $T\in\D_m(U)$, if it exists, is unique: if $T=T_1+\ldots+T_k=S_1+\ldots+S_k$ with $T_i,\ S_i\in\D_{p_i,q_i}(U)$, then for every $(f,\pi)\in\A^{p_i,q_i}(U)$
$$T_i(f,\pi)-S_i(f,\pi)=T(f,\pi)-T(f,\pi)=0$$
so $T_i=S_i$, by Proposition \ref{prp_0_su_pq}.

\medskip

If $T$ admits a Dolbeault decomposition, we will denote by $T_{p,q}$ its unique component of bidimension $(p,q)$.

In general, given $T\in\D_m(U)$, we define $T_{p,q}:\A^m(U)\to\C$ as the alternating $(m+1)-$linear maps which respects the product rule and such that
$$T_{p,q}(f,\pi)=T(f,\pi)\qquad\forall\;(f,\pi)\in\A^{p,q}(U)$$
$$T_{p,q}(f,\pi)=0\qquad\forall\;(f,\pi)\in\A^{r,s}(U),\ (r,s)\neq(p,q)$$
By Remark \ref{rem_forme_pq_gen}, this defines uniquely $T_{p,q}$ on $\A^m(U)$.

$T$ admits a Dolbeault decomposition if $T_{p,q}$ extends to a current in $\D_{p,q}(U)$ for every bidimension such that $p+q=k$.

\medskip

We observe that, if $T=T_1+\ldots+T_k$ is the Dolbeault decomposition for $T\in\D_m(U)$, then $\supp(T)\supseteq\supp(T_i)$ for $1\leq i\leq k$.

Obviously, if $F\in\mathrm{Hol}(X,Y)$ is proper on $\supp(T)$, for $T\in\D_{p,q}(X)$, then $F_\sharp T\in\A_{p,q}(Y)$.

\medskip

Now, we turn our attention to the behaviour of metric currents on the regular part of an analytic space. As it is remarked in \cite{lang1}, after Theorem 2.5, if $U$ is an open set in some $\R^n$, $\D_m(U)=\{0\}$ if $m>n$. The same obviously holds also if $U$ is an open set of a $n-$dimensional real manifold. In particular, $\D_m(U)=0$ for every $U\Subset X_\rg$ as soon as $m>2n$, where $n=\dim_\C X_\rg$. We can actually say a little more.

\begin{Prp} If $n=\dim_C X_\rg$, then $\D_m(U)=\{0\}$ for every open set $U\subseteq X$ and every $m>2n$.\end{Prp}
\noindent{\bf Proof: } Let $U\subseteq X$ be an open set and  $T\in\D_m(U)$ with $m>2n$. As a consequence of the previous considerations (or directly of  Theorem \ref{teo_chain}), if $(f,\pi)\in \D^m(U)$ and $\supp(f)\Subset X_\rg$, then $T(f,\pi)=T\llcorner U_\rg(f,\pi)=0$.

Therefore $\supp(T)\subseteq X\setminus X_\rg=X_\sg$, so, by Proposition \ref{prp_supp_corr}, $T=i_\sharp T^1$, where $i:X_\sg\to X$ is the inclusion and $T^1\in \D_m(X_\sg\cap U)$; now, we note that $X_\sg$ is a complex analytic space of dimension strictly less that $n$, so we conclude, by induction, that $T=0$. $\Box$

\medskip

The same kind of result can be obtained for bidimension:

\begin{Prp}\label{prp_zero_high_dim}If $n=\dim_\C X_\rg$, $U\subseteq X$ is an open set and $T\in \D_{p,q}(X)$ with $p>n$ or $q>n$, then $T=0$.\end{Prp}
\noindent{\bf Proof: } It is enough to show that $T$ is zero when applied to analytic forms of pure type $(p,q)$. Take $(f,\pi)\in \A^{p,q}(U)$ and suppose that $\pi_1,\ldots,\pi_p$ are holomorphic and $\pi_{p+1},\ldots, \pi_{p+q}$ are antiholomorphic. Suppose also that $p>n$ (the case $q>n$ is analogous) and that $\supp(f)\Subset U_\rg$ is contained in complex chart $V$.

We can apply Theorem \ref{teo_chain} to compute $T(f,\pi)$, using complex coordinates on $V$. Since $p>n$, there will be some coordinate function which appears more than once, so $T(f,\pi)=0$. As long as $\supp(f)$ is compact in $X_\rg$, we can cover it by a finite number of coordinate charts and still obtain that $T(f,\pi)=0$. This shows that $T\llcorner X_\rg=0$, because of Proposition \ref{prp_0_su_pq}. So, $T=i_\sharp T^1$ with $i:X_\sg\to X$ the inclusion and $T^1\in\D_{p,q}(X_\sg\cap U)$ and, by induction on the dimension, the thesis follows. $\Box$.

\medskip

\begin{Prp}\label{prp_pushf_pq}If $F:X\to Y$ is a proper holomorphic map of complex spaces and $T\in\D_{p,q}(X)$, then $F\sharp T\in \D_{p,q}(Y)$.\end{Prp}

\subsection{Density argument}
Up to now, we have defined the $(p,q)-$components of a current as functionals $T_{p,q}:\A^m(U)\to\C$. We observe that $T_{p,q}$ is clearly multilinear and it is a simple matter to check that it is local. The problem now is to investigate when these maps are restrictions of actual metric currents. This is equivalent to the continuity of $T_{p,q}:\A^m(U)\to\C$. Indeed we have     
\medskip

\begin{Teo}\label{teo_est_G}Let $T:\A^m(U)\to\C$ be a multilinear functional, local and continuous on $\A^m(U)$; then there exists a unique metric current which, restricted to $\A^m(U)$, concides with $T$. We will denote this unique extension with $T$.\end{Teo}
\noindent{\bf Proof: }Without loss of generality we may assume that $U$ is biholomorphic to an analitic set $A$ in some open set $\Omega$ of $\C^N$.

We observe the following
\begin{Lmm} \label{lmm_cont_0}Let $T:\A^m(U)\to\C$ be a multilinear, local, continuous functional and let $(f,\pi^i)$, $(f,\eta^i)$ two sequences of metric forms in $\A^m(U)$ such that $(f,\pi^i)-(f,\eta^i)\longrightarrow 0$, then $T(f,\pi^i)-T(f,\eta^i)\longrightarrow 0$.\end{Lmm}
It's enough to consider the difference componentwise. 
 
 Next, define $T(f,\pi)$ for $(f,\pi)\in\D^m(U)$ taking a sequence $\{(f,\pi^i)\}_i\in \A^m(U)$ such that $(f,\pi^i)\to(f,\pi)$ in $\D^m(U)$ and $\Lip(\pi^i_j)\to\Lip(\pi_j)$ for $j=1,\ldots,m$. The sequence of complex numbers $\{T(f,\pi^i)\}_i$ is then bounded in $\C$ therefore we can choose two sequences $\mu_i$ and $\nu_i$ of natural numbers in such a way that $T(f,\pi^{\mu_i})$ and $T(f,\pi^{\nu_i})$ are convergent and 
$$(f,\pi^{\mu_i})-(f,\pi^{\nu_i})\longrightarrow0.$$
Then, by Lemma \ref{lmm_cont_0}, we have
$$T(f,\pi^{\mu_i})-T(f,\pi^{\nu_i})\longrightarrow 0$$
and we can put
$$T(f,\pi)=\lim_{i\to\infty}T(f,\pi^i)$$
If $(f,\eta^i)_i$ is another sequence of metric forms in $\A^m(U)$ which converges to $(f,\pi)$ (even without the condition on the Lipschitz constants), we can proceed in the same way to obtain that $T(f,\pi^i)-T(f,\eta^i)\longrightarrow 0$. This shows that the definition of $T(f,\pi)$ does not depend on the approximating sequence $(f,\pi)$.

The extension of $T$ is obviously multilinear and it enjoys the locality property (it's enough to approximate a form with a constant coefficient with forms in $\A^m(U)$ with the same constant coefficients, which is possible because $\A$ contains the constants). Thus, in order to end the proof we have to show that $T$ is continuous.

Suppose that the sequence $\{(f_i,\pi^i)\}_i\subset\D^m(U)$ converges to$(f,\pi)\in\D^m(U)$. For every fixed $i\in\N$, there is a sequence $\{(f_i,\pi^{i,\nu})\}_{\nu}$ in $\A^m(U)$ which converges to $(f_i,\pi^i)$, together with the Lipschitz constants of the components, on every compact. Fix  compact $K\subset U$ and suppose that $\Lip(\pi^i_k\vert_K),\Lip(\pi_k\vert_K)\leq C$. We can choose an integer $\nu(i)$ such that $\Lip(\pi^{i,\nu}_k\vert_K)\leq2C$, if $\nu>\nu(i)$. Therefore, for every monotonically increasing map $\mu:\N\to\N$, $\mu\geq\nu$, the sequence $
(f_i,\pi^{i,{\mu(i)}})$ converges to $(f,\pi)$ and, by the definition of $T(f,\pi)$, 
$$\lim_{i\to\infty}T(f_i,\pi^{i,{\mu(i)}})=T(f,\pi).$$
Now, by the definition of $T(f_i,\pi^i)$,  we can choose $\mu(i)$ such that  
$$|T(f_i,\pi^i)-T(f_i,\pi^{i,{\mu(i)}})|\leq 2^{-i}$$ 
and then
$$|T(f_i,\pi^i)-T(f,\pi)|\leq|T(f_i,\pi^i)-T(f_i,\pi^{i,{\mu(i)}})|+|T(f_i,\pi^{i,{\mu(i)}})-T(f,\pi)|$$
which tends to $0$ as $i$ approaches infinity. $\Box$

\begin{Rem} \label{rem_algebra_ext}Lemma \ref{lmm_cont_0} and Theorem \ref{teo_est_G} still hold if we replace $\A^m(U)$ with $\D(U)\times[\mathscr{G}(U)]^m$, where $\mathscr{G}(U)$ is any subalgebra of $\Lip_\loc$ containing the constants and such that every element  $u\in\Lip_\loc$ can be approximated by elements of $\mathscr{G}$ with Lipschitz constants on any compact $K$ bounded by $\Lip(u\vert_K)+\epsilon$, for any $\epsilon$.\end{Rem}

\medskip
\noindent{\bf Example } Let $T\in M_1(\C)$ be the metric current defined by
$$T(f,\pi)=\int_{S^1}fd\pi\;.$$
When $fd\pi$ is a smooth compactly supported $1-$form, this integral defines a classical flat current, therefore it can be extended to a metric current, which turns out to be of finite mass; $T$ can be written as
$$T(f,\pi)=\int_{S^1}\frac{i}{2}\pair{fd\pi, z\de_z-\bar{z}\de_{\bar{z}}}d\H^1$$
so
$$T(f,\pi)=\int_{S^1}\frac{i}{2}\pair{fd\pi, z\de_z}d\H^1-\int_{S^1}\frac{i}{2}\pair{fd\pi, \bar{z}\de_{\bar{z}}}d\H^1$$
$$=T_{1,0}(fd\pi)+T_{0,1}(fd\pi)\;.$$
Now, let us consider the two metric forms $(f,z)$ and $(f,z^2\bar{z})$, which coincide on $S^1=\supp T=\supp T_{1,0}=\supp T_{0,1}$; we have
$$T_{1,0}(f,z)=\frac{i}{2}\int_{S^1}zf(z)d\H^1$$
$$T_{1,0}(f,z^2\bar{z})=i\int_{S^1}zf(z)d\H^1\;.$$
Therefore, the conclusions of Proposition \ref{prp_supp_corr} don't hold for $T_{1,0}$, implying that it cannot be extended to a metric current.

\medskip

\noindent{\bf Example } Let $T\in M_{2,\loc}(\C)$ be the local metric current given by
$$T(f,\pi_1,\pi_2)=\int_{\C}\frac{1}{i\pi z}f\det(\nabla\pi)dz\wedge d\bar{z}\;.$$
We know that $T$ is a local metric current by \cite[Theorem 2.6]{lang1}; let $S=dT$ be its boundary. We compute the classical $(0,1)-$component of $S$, obtaining that 
$$S_{0,1}(f,\pi)=Cf(0)\frac{\de\pi}{\de z}(0)$$
for some constant $C$. It is easy to see that $S_{0,1}$ doesn't fullfill the conclusions of Proposition \ref{prp_supp_corr}.

The examples given in \ref{sec_ex1} can be defined on complex manifolds or on complex spaces; in particular, given a complex space $X$ of pure dimension $n$ (either as an abstract space of as a subspace of a bigger complex space $M$), the metric current
$$[X](f,g_1,\ldots, g_{2n})=\int_{X_\rg}f\det(\nabla g_1,\ldots, \nabla g_{2n})d\H^{2n}\;,$$
with $\H^{2n}$ the $2n-$dimensional Hausdorff measure, is closed (as $bX=\emptyset$), hence normal.

\begin{Prp}The current $[X]$ is of bidimension $(n,n)$.\end{Prp}
\noindent{\bf Proof: }Every $(p,q)-$component of $X$ is such that $p+q=2n$, therefore, unless $p=q=n$, there is one between $p$ and $q$ greater than $n$. By Proposition \ref{prp_zero_high_dim}, $T_{p,q}$ is then zero. So $[X]=[X]_{n,n}$. $\Box$

\medskip

From this result, it follows that, given a metric form $\alpha$ of pure type $(p,q)$, the current $[X]\llcorner\alpha$ is of bidimension $(n-p,n-q)$. 

\begin{Prp}If $T\in M_{k}(X)$ is represented by the integration on $X_\rg$ against a form with $L^1_\loc(X,\H^{2n})$ coefficients, then it admits a Dolbeault decomposition. \end{Prp}
\noindent{\bf Proof: } If $T$ is representable by such a form, then also every $T_{p,q}$ is, therefore it is enough to show that every functional which admits such a representation is a metric current. This can be done arguing as in Proposition \ref{prp_cont_inc}. $\Box$

\subsection{$\partial$ and $\debar$}

Let us suppose that $T$ is a $(p,q)-$current whose boundary admits a Dolbeault decomposition. Then we can define $\de T$ and $\debar T$ as follows.

Write $dT=S_1+\ldots+S_h$ with $S_i\in M_{p_i,q_i}(U)$ where $p_i+q_i=p+q-1=m$ since $dT\in M_{p+q-1}(U)$. If $(f,\pi)$ is a $m-$form of pure type $(p_i,q_i)$, then
$$S_{p_i,q_i}(f,\pi)=dT(f,\pi)=T(1,f,\pi)$$
if $p_i>p$ or $q_i>q$,and consequently $T(1,f,\pi)=0$, since $T$ is a $(p,q)-$current. Therefore, we can only have two cases: $p=p_i$ and $q-1=q_i$ or $p-1=p_i$ and $q=q_i$ i.e.
$$dT=S_{p,q-1}+S_{p-1,q}$$
and we put
$$\de T=S_{p-1,q}\qquad \debar T=S_{p,q-1}$$
Therefore, if a current $T$ admits a decomposition in $(p,q)$ components, we can define $\de T$ and $\debar T$ setting
$$
\de T=\sum_{i=1}^h\de T_i\qquad \debar T=\sum_{i=1}^h\debar T_i
$$
where $T=T_1+\ldots+T_h$ is the $(p,q)$ decomposition.

\medskip

\begin{Prp}\label{prp_pushf_debar}
If $H:X\to Y$ is a holomorphic map between complex spaces, $U\subseteq X$ an open set, then, for every current $T\in \D_m(U)$ for which $\de T$ and $\debar T$ are defined, the following hold:
$$H_\sharp \de T=\de H_\sharp T\qquad H_\sharp\debar T=\debar H_\sharp T$$
\end{Prp}
\noindent{\bf Proof: }
By Proposition \ref{prp_pushf_pq}, the pushforward of a $(p,q)-$current is a $(p,q)-$current. If $T\in M_{p,q}(U)$, 
$$d(H_\sharp T)=S_{p-1,q}+S_{p,q-1}=\de(H_\sharp T)+\debar(H_\sharp T)$$
by definition, but 
$$d(H_\sharp T)=H_\sharp(dT)=H_\sharp(\de T+\debar T)=H_\sharp(\de T)+H_\sharp(\debar T)$$
and so, by the uniqueness of decomposition, the thesis follows. $\Box$

\medskip

Moreover, it is easy to check that $C_m\circ\de=\de\circ C_m$ and $C_m\circ\debar=\debar\circ C_m$, where $C_m$ is the map given by Theorems \ref{teo_comp}, \ref{teo_comp_man}.

\begin{Prp}\label{prp_formule_de}
We have that $\de^2=\debar^2=0$ and $\de\debar=-\debar\de$
\end{Prp}
\noindent{\bf Proof: }
By the locality property, we have that $d^2T=0$, therefore
$$0=(\de+\debar)(\de+\debar)T=\de^2 T+(\de\debar T+\debar\de T)+\debar^2 T.$$
Since the right hand side is a decomposition in $(p,q)$ components, every term has to be zero:
$$\de^2T=0\qquad \debar^2 T=0\qquad\de\debar T+\debar\de T=0$$
and, as we didn't make any assumption on $T$, the thesis follows. $\Box$

\medskip

We give a formula for $\debar T$, for analytic subsets of $\C^n$.

\begin{Lmm}\label{lmm_eq_debar} Suppose $X$ is an analytic subset of some open set $\Omega\subseteq \C^n$. Given $T\in \D_m(X)$ such that $\debar T$ exists as a metric current and $(f,g_1,\ldots, g_{m-1})=(f,g)\in\D^{m-1}(U)$ with $\Ci^2$ coefficients, we have 
\begin{equation}\label{eq_debar}\debar T(f,g)=\sum_{j=1}^nT\left(\frac{\de f}{\de\bar{z}_j}, \bar{z}_j, g\right)+\sum_{k=1}^{m-1}(-1)^{k}\sum_{j=1}^n T\left(f, \bar{z}_j,\ldots, \frac{\de g_k}{\de\bar{z}_j},\ldots\right)\;.\end{equation}
\end{Lmm}
\noindent{\bf Proof: } The formula clearly holds for a classical current, integrating by parts. We already observed that $C_{m-1}\circ\debar=\debar\circ C_m$ and we know that there exists a metric current $S=\debar T$. The thesis follows. $\Box$

\medskip

We can define a multilinear, local functional of finite mass by (\ref{eq_debar}) and denote it by $\debar T$. If
$$W_m(X)=\{T\in \D_m(X)\ :\ \debar T\in\D_{m-1}(X)\}$$
then $\debar W_m(X)\subseteq W_{m-1}(X)$, by Proposition \ref{prp_formule_de}. The space $W_{p,q}(X)$ is defined in the same way but unfortunately, we don't have any decomposition theorem for $W_m$ in terms of $W_{p,q}$.

\medskip

\begin{Prp} \label{prp_int_parti}Suppose $U\subset X$ can be embedded as an analytic set into $\C^n$. Given $T\in W_m(U)$, $(u,v)\in\D^k(U)$ with $\Ci^\infty(U)$ coefficients, we have that
$$\debar (T\llcorner(u,v))=$$
$$(-1)^{k}\left((\debar T)\llcorner(u,v)- \sum_{j=1}^n T\llcorner(\de u/\de \bar{z}_j, \bar{z}_j, v)-\sum_{h=1}^k(-1)^h\sum_{j=1}^nT(u,\bar{z}_j,\ldots, \de v_h/\de\bar{z}_j,\ldots)\right)$$
where $z_1,\ldots, z_n$ are the coordinates of the embedding in $\C^n$.
\end{Prp}
\noindent{\bf Proof: }By Lemma \ref{lmm_eq_debar}
$$((\debar T)\llcorner(u,v))(f,g)=(\debar T)(uf,v,g)=\sum_{j=1}^nT\left(\frac{\de (uf)}{\de \bar{z}_j}, \bar{z}_j, v,g\right)+$$
$$+\sum_{h=1}^k(-1)^h\sum T\left(uf, \bar{z}_j,\ldots, \frac{\de v_h}{\de\bar{z}_j},\ldots, g\right)+\sum_{h=1}^{m-k}(-1)^{h+k}\sum_{j=1}^nT\left(uf, \bar{z}_j, v, \ldots, \frac{\de g_{h}}{\de\bar{z}_j},\ldots\right)$$
for every $(f,g)\in\D^{m-k}(U)$ with $\Ci^2$ coefficients.

We have
$$T\left(\frac{\de uf}{\de \bar{z}_j}, \bar{z}_j, v,g\right)=T\left(u\frac{\de f}{\de\bar{z}_j}, \bar{z}_j, v,g\right)+T\left(f\frac{\de u}{\de \bar{z}_j}, \bar{z}_j, v,g\right)$$
and we notice that
$$\sum_{j=1}^nT\left(u\frac{\de f}{\de\bar{z}_j}, \bar{z}_j, v,g\right)=(-1)^k\sum_{j=1}^nT\left(u\frac{\de f}{\de\bar{z}_j}, v,\bar{z}_j,g\right)=$$
$$=(-1)^k\sum(T\llcorner(u,v))\left(\frac{\de f}{\de\bar{z}_j},\bar{z}_j,g\right)\;.$$
So, again by Lemma \ref{lmm_eq_debar}, we obtain
$$\sum_{j=1}^nT\left(u\frac{\de f}{\de\bar{z}_j}, \bar{z}_j, v,g\right)+\sum_{h=1}^{m-k}(-1)^{h+k}\sum_{j=1}^nT\left(uf, \bar{z}_j, v, \ldots, \frac{\de g_{h}}{\de\bar{z}_j},\ldots\right)=$$
$$=(-1)^k\debar (T\llcorner(u,v))(f,g)\;.$$

Therefore we have
$$\debar (T\llcorner(u,v))=$$
$$(-1)^{k}\left((\debar T)\llcorner(u,v)- \sum_{j=1}^n T\llcorner(\de u/\de \bar{z}_j, \bar{z}_j, v)-\sum_{h=1}^k(-1)^h\sum_{j=1}^nT(u,\bar{z}_j,\ldots, \de v_h/\de\bar{z}_j,\ldots)\right)\;.$$
Now, the algebra of $\Ci^2$ functions satisfies the hypotheses of Remark \ref{rem_algebra_ext}, therefore the current $(\debar T)\llcorner(u,v)$ is uniquely determined by this formula, which therefore holds for every $(f,g)\in\D^{m-k}(U)$. $\Box$

\medskip

A similar formula holds for $\de T$.

%\begin{Rem}The boundary of an integer rectifiable current whose boundary has locally finite mass is again integer rectifiable and normal; therefore, the $\de$ and $\debar$ of an integer rectifiable $(p,q)-$current are again rectifiable.\end{Rem}

\subsection{Rectifiability and slicing}

Let $T$ be a metric current whose boundary admits a Dolbeault decomposition. Then $\supp \debar T, \supp \de T\subseteq \supp dT$ and, moreover, $\|\debar T\|_A, \|\de T\|_A\leq C\| dT\|_A$ for every $A\subset X$. In particular, if $dT$ is rectifiable, then $\debar T$ and $\de T$ are too.

\medskip

The slices of a current $T$ through a map $\pi:X\to\R^{n}$ are defined by
$$\langle T,\pi, x\rangle(f,\eta)=\lim_{\epsilon\to 0} T(\rho_{\epsilon, x}f, \pi, \eta)=\lim_{\epsilon\to 0}(-1)^{n(k-1-n)}\pi_\sharp(T\llcorner(f,\eta))(\rho_{\epsilon, x}, x_1,\ldots, x_n)$$
with $\rho_{\epsilon, x}$ any family of smooth approximations of  $\delta_x$. If $\pi:X\to\R^{2n}\cong\C^n$, we can write the slices as
$$\langle T,\pi, x\rangle=\lim_{\epsilon\to 0} \pi_\sharp(T\llcorner(f,\eta))(\rho_{\epsilon, x}, z_1,\bar{z}_1,\ldots, z_n,\bar{z}_n)\;.$$

\begin{Prp}The operators $\de$ and $\debar$ commute with the slicing through holomorphic maps.\end{Prp}
\noindent{\bf Proof: }Let $T\in N_{k,\loc}(X)$ be a locally normal current, such that $\de T$ and $\debar T$ are again locally normal. Let $\pi:X\to\C^n$ be a holomorphic map and suppose that $U\subset X$ is biholomorphic to an analytic subset of an open set of some $\C^N$. Let $(f,\eta)$ be a $(k-2n-1)-$metric form with $\Ci^2$ coefficients supported in $U$. We treat only the case of $\debar T$, the proof for $\de T$ being analogous.

Let $z_1,\ldots, z_n$ and $w_1,\ldots, w_N$ be holomorphic coordinates in $\C^n$ and $\C^N$, respectively. The slice $\langle \debar T, \pi, x\rangle$ exists for a.e. $x$, by Theorem \ref{teo_ext_slice} and we have
$$\langle \debar T,\pi, x\rangle(f,\eta)=\lim_{\epsilon\to 0}\pi_\sharp((\debar T)\llcorner(f,\eta))(\rho_{\epsilon,x}, z_1,\bar{z}_1,\ldots, z_n, \bar{z}_n)\;.$$
Now, we set  $\widetilde{\eta}^{jh}$ to be the $(k-2n-1)-$tuple differing from $\eta$ only in the $h-$th component, which is $(-1)^h\de \eta_h/\de\bar{w}_j$. By Proposition \ref{prp_int_parti}
$$(\debar T)\llcorner(f,\eta)=(-1)^{k-1}\debar(T\llcorner (f,\eta)) + \sum_{j=1}^NT\llcorner(\de f/\de \bar{w}_j, \bar{w}_j, \eta)+\sum_{h,j}T\llcorner(f, \bar{w}_j,\widetilde{\eta}^{jh})$$
and we note that $T\llcorner(f,\eta)$ is a $2n+1-$form, so 
$$\pi_\sharp\debar T\llcorner(f,\eta)=\debar\pi_\sharp T\llcorner(f,\eta)=0$$
by Proposition \ref{prp_pushf_debar}, as $\pi$ is holomorphic. It follows
$$\pi_\sharp((\debar T)\llcorner(f,\eta))(\rho_{\epsilon, x},z_1,\ldots, \bar{z}_n)=\sum_{j=1}^N\pi_\sharp(T\llcorner(\de f/\de\bar{w}_j, \bar{w}_j, \eta))(\rho_{\epsilon, x},z_1,\ldots, \bar{z}_n)+$$
$$+\sum_{j,h}\pi_\sharp(T\llcorner(f,\bar{w}_j,\widetilde{\eta}^{jh}))(\rho_{\epsilon,x},z_1,\ldots,\bar{ z}_n)$$
$$=\sum_{j=1}^NT\left(\frac{\de f}{\de \bar{w}_j}\cdot(\rho_{\epsilon,x}\circ\pi), \pi_1, \bar{\pi}_1,\ldots, \pi_n,\bar{\pi}_n, \eta\right)+\sum_{j,h}T\left(f\cdot(\rho_{\epsilon,x}\circ\pi), \pi_1,\ldots, \bar{\pi}_n, \widetilde{\eta}^{jh}\right)$$
$$=\sum_{j=1}^NT\llcorner(\rho_{\epsilon,x}\circ\pi, \pi_1,\ldots, \bar{\pi}_n)(\de f/\de\bar{w}_j,\bar{w_j},\eta)+\sum_{j,h}T\llcorner(\rho_{\epsilon,x}\circ\pi,\pi_1,\ldots,\bar{\pi}_n)(f,\bar{w}_j,\widetilde{\eta}^{jh})\;.$$
Now, again by Lemma \ref{lmm_eq_debar}, we have
$$\debar\langle T,\pi, x\rangle(f,\eta)=\sum_{j=1}^N\langle T,\pi, x\rangle\left(\frac{\de f}{\de\bar{w}_j}, \bar{w}_j, \eta\right)+\sum_{j,h}\langle T,\pi,x\rangle(f,\bar{w}_j,\widetilde{\eta}^{jh})$$
$$=\sum_{j=1}^N\lim_{\epsilon\to0} T\llcorner(\rho_{\epsilon,x}\circ\pi, \pi_1,\ldots, \bar{\pi}_n)\left(\frac{\de f}{\de\bar{w}_j}, \bar{w}_j,\eta\right)+\sum_{j,h}\lim_{\epsilon\to0} T\llcorner(\rho_{\epsilon,x}\circ\pi, \pi_1,\ldots, \bar{\pi}_n)\left(f, \bar{w}_j,\widetilde{\eta}^{jh}\right)$$
and by the previous computation this is equal to
$$\lim_{\epsilon\to 0}\pi_{\sharp}((\debar T)\llcorner(f,\eta))(\rho_{\epsilon, x}, \pi_1,\ldots,\bar{\pi}_n)=\langle \debar T,\pi, x\rangle(f,\eta)\;.$$
Therefore, for every $(f,\eta)$ with $f\in\Ci^2_c$
$$\langle \debar T,\pi, x\rangle(f,\eta)=\debar\langle T,\pi, x\rangle (f,\eta)\;.$$
This means that the functional $\debar\langle T,\pi, x\rangle$ can be extended to all the metric forms as a metric current, by defining it equal to $\langle \debar T,\pi, x\rangle$. $\Box$

\medskip

As already said, the previous proof works also for the $\de$:
$$\de\langle T,\pi, x\rangle=\langle \de T,\pi, x\rangle\;.$$
Moreover, if $\de\debar T$ is a metric current, then
$$\de\debar\langle T,\pi, x\rangle=\langle \de\debar T,\pi, x\rangle\;,$$
or, which is the same,
$$dd^c\langle T,\pi,x\rangle=\langle dd^c T,\pi, x\rangle\;.$$

\subsection{Positive currents}\label{ssc_pos_cur}

Let $T$ be a metric $(p,p)-$functional. We say that $T$ is \emph{positive} if, given $\pi_1,\ldots, \pi_p\in\Ol(X)$, 
$$T(f,\pi_1,\overline{\pi}_1,\ldots, \pi_p,\overline{\pi}_p)\geq0$$ for every compactly supported Lipschitz function $f\geq0$ on $X$.

\begin{Prp} Let $X$ be a complex analytic set in $\C^N$ and $T$ a local $(p,p)-$metric current on $X$. If $T$ is positive, then $T$ is of locally finite mass.\end{Prp}
\noindent{\bf Proof: } Let $\chi_K$ be the indicatrix function for the compact $K$ in $X$ and define
$$m_I(K)=T(\chi_K,z_{i_1},\overline{z}_{i_1},\ldots,z_{i_p},\overline{z}_{i_p})\;.$$
Then, for any real valued Lipschitz function $f$ with compact support,
$$|T(f,z_{i_1},\ldots, \overline{z}_{i_p})|\leq\|f\|_\infty T(\chi_{\supp f},z_{i_1},\ldots,\overline{z}_{i_p})=\|f\|_\infty m_I(\supp f)$$
and the same holds for complex valued functions, separating real and imaginary part.

Moreover, given multi-indeces $I$ and $J$ of length $p$, we have that
$$T(f,z_{i_1},\bar{z}_{j_1},\ldots, z_{i_p}, \bar{z}_{j_p})=\sum_{s\in S} \alpha_sT(f, h_1,\bar{h}_1,\ldots, h_n,\bar{h}_n)$$
where
$$h_m=\frac{1}{2}(z_{i_m}+i^{\beta_s}z_{j_m})$$
where $S$, $\alpha_s$, $\beta_s$ are an appropriate indexes set and appropriate constants, respectively.

This implies that
$$f\mapsto T(f,z_{i_1}, \ldots, \bar{z}_{j_p})$$
is a complex measure and with a little more effort we obtain that its total variation on the compact $K$, which we denote by $m_{IJ}(K)$, is bounded by
$$\sum_{I\cup J\supset M\supset I\cap J}m_M(K)\;.$$

Therefore, given $g_1,\ldots, g_{2p}\in \mathrm{Lip}_{\mathrm{loc}}(X)$, we have
$$|T(f,g_1,\ldots, g_{2p})|=\left|\sum_{IJ} T\left(f\frac{\de g}{\de(z_I, \bar{z}_J)}, z_{i_1},\bar{z}_{j_1}\ldots,z_{i_p}, \overline{z}_{j_p}\right)\right|\leq\|f\|\prod\mathrm{Lip}(g_j)\sum_{I,J} m_{I,J}(\supp f)\;.$$
This means that
$$M_K(T)\leq \sum_{I,J} m_{IJ}(\supp f)<+\infty\;.$$
$\Box$

\medskip

\begin{Rem}It is easy to deduce that, if each $m_I(K)$ is uniformly bounded as $K$ varies through the compact sets of $X$, then $T$ is of finite mass.\end{Rem}

\begin{Rem}The hypotheses of the previous Proposition can obviously weakened assuming only that $T$ is a metric functional on local metric forms, multilinear, alternating and satisfying the product rule with respect to differentials.\end{Rem}

\medskip

Classical examples of positive currents (which are also metric examples, of course) are the integration on a complex analytic set and the $i\de\debar$ of a plurisubharmonic function.

Most of the results about classical positive currents are of local nature, concerning usually extensions through small sets, vanishing or representation by integration on analytic sets.

Let $T$ be a metric current on a complex space $X$. For every $x\in X$ there exists a neighborhood of $x$ which is biholomorphic to an analytic subset $A$ of some open domain $\Omega\subset\C^n$, therefore, locally, $T$ can be viewed as a current on $\Omega$, supported in $A$ belonging to the image of the pushforward operator induced by the inclusion $A\hookrightarrow\Omega$.

Positivity in the metric sense traslates quite obviously into positivity in the classical sense, therefore any local result holding for classical positive currents, holds also for metric positive currents on analytic spaces.

\begin{Teo}[Skoda-El Mir]Let $E\subset X$ be a closed complete pluripolar set (i.e. for every $x_0\in X$ there exists a neighborhood $U$ and a locally integrable plurisubharmonic function $u$ such that $E\cap U=\{z\in U\ :\ u(z)=-\infty\}$) and suppose $T$ is a positive $(p,p)-$current on $X\setminus E$. Assume that $T$ has finite mass on a neighborhood of every point of $E$. Then the trivial extension of $T$ to $X$ is closed.\end{Teo}
\noindent{\bf Proof: } Let $x_0\in X$ be fixed and consider a neighborhood $U$ of $x_0$ such that
\begin{enumerate}
\item $U$ is biholomorphic to an analytic subset $A$ of some open domain $\Omega$ in $\C^n$;
\item $U\cap E=\{z\in U\ :\ u(z)=-\infty\}$ for some $u\in L^1(U)_\loc\cap Psh(U)$.
\end{enumerate}

Shrinking $U$ and $\Omega$ if necessary, suppose $A$ is defined in $\Omega$ by the equations $g_1=\ldots=g_k=0$, with $g_1,\ldots, g_k\in\Ol(\Omega)$. Let $E'$ be the image of $E\cap U$ through the biholomorphism with $A$ and denote by $T'$ and $u'$ the pushforwards of $T\llcorner U$ and $u$. $T'$ will be of locally finite mass around the points of $E'$ and $u'$ will be a plurisubharmonic function on $A$, i.e. the restriction to $A$ of a plurisubharmonic function $\widetilde{u}$ on $\Omega$ which is not identically $-\infty$ on $A$. By a classical result, $\widetilde{u}$ is locally integrable in $\Omega$, with respect to the standard Lebesgue measure on $\C^n$.

The set $E'$ is closed and complete pluripolar in $\Omega$, as
$$E'=\{z\in\Omega\ :\ \log(|g_1|^2+\ldots+|g_k|^2+e^{\widetilde{u}})=-\infty\}\;.$$

Therefore, by the classical Skoda-El Mir theorem, $T'$ extends trivially through $E'$ to a closed positive current on $\Omega$, supported in $A$. As it is locally normal, it is locally flat, therefore metric, by \ref{teo_comp}.

So, the trivial extension of $T$ is locally metric, closed and positive, hence it is metric, closed and positive on $X$. $\Box$

\medskip

With the same kind of argument, we obtain the following results.

\begin{Teo}[Second theorem of support] Let $M\subset X$ be a $\mathrm{CR}-$submanifold with $\mathrm{CR}-\dim M=p$ such that there is a submersion $p:M\to Y$ of class $\Ci^1$ whose fibers $F_t=p^{-1}(t)$ are connected and are the integral manifolds of the holomorphic tangent space $TM\cap JTM$. Then any closed $(p,p)-$current $T$ on $X$ of locally finite mass with $\supp T\subseteq M$ can be written as
$$T=\int_Y[F_t]d\mu(t)$$
with a unique complex measure $\mu$ on $Y$. $T$ is positive  if and only if $\mu$ is positive.\end{Teo}

\begin{Cor} Let $A$ be an analytic subset of $X$ with global irreducible components $A_j$ of pure dimension $p$. Then any closed current $T$ of bidimension $(p,p)$ with locally finite mass and $\supp T\subseteq A$ is of the form $T=\sum m_j[A_j]$.\end{Cor}

\begin{Teo}[Lelong-Poincaré equation]Let $f$ be a meromorphic function on $X$ which does not vanish identically on any component of $X$. Let $\sum m_jZ_j$ be the divisor of $f$. Then  the function $\log|f|$ is locally integrable on $X$ and satisfies the equation
$$\frac{i}{\pi}\de\debar\log|f|=\sum m_j[Z_j]$$
as a metric current.
\end{Teo}
\noindent{\bf Proof: } The current $(i/\pi)\de\debar\log|f|$ is positive, hence of locally finite mass, and closed, hence locally normal, therefore it is a local metric current.

Clearly, the thesis holds in any relatively compact open set of $X_\rg$, therefore 
$$T=\frac{i}{\pi}\de\debar\log|f|-\sum m_j[Z_j]$$
is supported in $X_\sg$, but then by Proposition \ref{prp_zero_high_dim}, $T=0$. $\Box$

\section{Sheaves of currents}

Let $X$ be a locally compact, paracompact metric space. We consider the presheaf
$$U\mapsto \D_m(U)\qquad U\subseteq X \ \mathrm{open.}$$
If $V\subseteq U$ are two open sets of $X$, we have the map
$$\rho^V_U:\D_m(U)\to \D_m(V)$$
defined by
$$\rho^V_U(T)(f,\pi)=T\llcorner\chi_V(f,\sigma\pi)$$
for every $(f,\pi)\in\D^m(V)$, where $\sigma$ is any locally Lipschitz function in $U$ equal to $1$ on $\supp(f)$ and $0$ outside $V$.

\begin{Prp}\label{prp_local_id}
Let $\{U_i\}_{i\in I}$ be a family of open subsets of $X$ and $U$ their union. 
If $T,\ S\in \D_m(U)$ and $\rho^{U_i}_U(T)=\rho^{U_i}_U(S)$ for every $i\in I$, then $T=S$.
\end{Prp}
\noindent{\bf Proof: } Let us consider a partition of unit $\{\phi_i\}$ subordinated to the open covering $\{U_i\}$ (we can assume, by possibly refining the covering, that each $U_i$ is relatively compact). Define $T_i=T\llcorner\phi_i$ and $S_i=S\llcorner\phi_i$. Then, for every $(f,\pi)\in\D_m(U)$, we have
$$T_i(f,\pi)=\rho^{U_i}_U(T_i)(f\vert_{U_i},\pi\vert_{U_i})$$
and
$$S_i(f,\pi)=\rho^{U_i}_U(S_i)(f\vert_{U_i},\pi\vert_{U_i}).$$
So 
$$T(f,\pi)=\sum_{i\in I}T_i(f,\pi)=\sum_{i\in I}\rho^{U_i}_U(T_i)(f\vert_{U_i},\pi\vert_{U_i})=\sum_{i\in I}\rho^{U_i}_U(T)(f\phi_i\vert_{U_i},\pi\vert_{U_i})=$$
$$=\sum_{i\in I}\rho^{U_i}_U(S)(f\phi_i\vert_{U_i},\pi\vert_{U_i})=\sum_{i\in I}\rho^{U_i}_U(S_i)(f\vert_{U_i},\pi\vert_{U_i})=\sum_{i\in I}S_i(f,\pi)=S(f,\pi)$$
which is the thesis. $\Box$

\medskip

\begin{Prp}\label{prp_gluing_sec}
Given a collection $\{U_i\}_{i\in I}$ of open sets of $X$ and a collection $\{T_i\}_{i\in I}$ of local metric currents $T_i\in \D_m(U_i)$ such that 
$$\rho^{U_i\cap U_j}_{U_i}(T_i)=\rho^{U_i\cap U_j}_{U_j}(T_j)$$
for every pair of indexes $i,j\in I$, there is $T\in \D_m(U)$ with $U=\bigcup U_i$ such that $T_i=\rho^{U_i}_U(T)$.
\end{Prp}
\noindent{\bf Proof: }
We can suppose that each $U_i$ is relatively compact and locally finite and choose a partition of unity $\{\phi_i\}_{i\in I}$ subordinated to the covering $\{U_i\}_{i\in I}$ of $U$.
For $(f,\pi)\in\D^m(U)$, we set
$$T(f,\pi)=\sum_{i\in I}T_i\llcorner\phi_i(f,\pi)$$
The sum is finite, for $f$ has compact support, and $T$ is obviously multilinear and local. Continuity with respect to the topology of $\D^m(X)$ follows since we can restrict ourselves to a compact subset of $U$: if $(f_i,\pi^i)\to(f,\pi)$, the set $\bigcup\supp(f_i)$ is compact.

Therefore, $T\in \D_m(U)$ and is quite simple to check that $\rho^{U_i}_U(T)=T_i$ for every $i$. $\Box$

\medskip

This result combined with Proposition \ref{prp_local_id}, shows that the presheaf
$$U\mapsto \D_m(U)$$
with the obvious restriction morphisms, is a canonical presheaf and therefore there exists a sheaf $\D_m$ whose sections on an open set $U$ are precisely $\D_m(U)$.

\subsection{Locally finite mass}

We can define a sub-presheaf of $\D_m$ by setting
$$U\mapsto M_{m,\loc}(U)\qquad U\subseteq X\textrm{ open}$$
and by considering the same restriction morphisms.

Obviously, Proposition \ref{prp_local_id} still holds. In order to obtain the gluing property, we have to check that the current $T$ constructed in the proof of Proposition \ref{prp_gluing_sec} is of locally finite mass, as soon as the initial currents $T_i$ are.

\medskip

Let us suppose that, for every finite collection $\{(f_\lambda,\pi^\lambda)\}_{\lambda\in\Lambda}\subset\D(X)\times[\Lip_1(X)]^m$ such that $\supp(f_\lambda)\subset U_i$ and $\sum_\lambda |f_\lambda|\leq 1$, we have
$$\sum T_i(f_\lambda,\pi^\lambda)\leq M_i<+\infty\;.$$
Given a finite collection $\{(f_\lambda, \pi^\lambda)\}_{\lambda\in\Lambda} \subset \D(X)\times[\Lip_1(X)]^m$ with $\supp(f_\lambda)\subset U$ and $\sum_\lambda |f_\lambda|\leq 1$, we have that, for every $i$, the collection $\{(f_\lambda\phi_i,\pi^\lambda)\}$ is as described above, so
$$\sum_\lambda T(f_\lambda,\pi^\lambda)=\sum_\lambda\sum_i T_i(\phi_if_\lambda,\pi^\lambda)\;.$$
Now, let us fix $V\Subset U$ and suppose that $\supp(f_\lambda)\subset V$; then only a finite number of $\phi_i$ are not vanishing on $V$. Let $J=\{i\in I\ :\ U_i\cap V\neq \emptyset\}$ and set $M_J=\sum_{j\in J} M_j<+\infty$; then
$$\sum_\lambda T(f_\lambda,\pi^\lambda)\leq M_J$$
where $M_J$ does not depend on the collection $\{(f_\lambda,\pi^\lambda)\}$, but only on $V$. Thus $T\in M_{m,\loc}(U)$.

Summing up, we have the following proposition.

\begin{Prp}The assignment
$$U\mapsto M_{m,\loc}(U)$$
defines a canonical presheaf, with an associated sheaf $\mathscr{M}_{m}$. The sections of $\mathscr{M}_{m}$ are the metric $m-$currents with locally finite mass.\end{Prp}

In the same way, we can obtain the same result for normal currents.

\begin{Prp}The assignment
$$U\mapsto N_{m,\loc}(U)$$
defines a canonical presheaf, thus associated to a sheaf $\mathscr{N}_{m}$, whose sections are the locally normal metric $m-$currents.\end{Prp}

The sheaves $\D_{p,q}$, $\mathscr{M}_{p,q}$, $\mathscr{N}_{p,q}$ are defined exactly in the same way.

\subsection{Properties}

\begin{Prp} $\D_m$, $\mathscr{M}_m$, $\mathscr{N}_m$, together with their $(p,q)$ analogues, are fine sheaves. \end{Prp}
\noindent{\bf Proof: } Let $\{U_i\}$ be an open covering of $X$, $\{r_i\}$ a collection of Lipschitz functions which form a partition of unity subordinated to $\{U_i\}$. Given $T\in\D_m(X)$ (or $\mathscr{M}_m$ or $\mathscr{N}_m$), we can consider the currents $T_i=T\llcorner r_i$;  then $T_i\in\D_m(U_i)$ (or $\mathscr{M}_m$ or $\mathscr{N}_m$) and $\sum T_i=T$. So $\D_m$ (or $\mathscr{M}_m$ or $\mathscr{N}_m$) is a fine sheaf. $\Box$

\medskip

\begin{Prp}$\D_m$, $\mathscr{M}_m$, $\mathscr{N}_m$, together with their $(p,q)$ analogues, are soft sheaves. \end{Prp}
\noindent{\bf Proof: } We can either apply the previous proposition and observe that fine sheaves on paracompact Hausdorff spaces are soft, or notice that, given a closed set $F\subset X$, the map $i:F\to X$ induced by the inclusion is proper and Lipschitz, therefore if $T\in \D_m(F)$, $i_\sharp T\in \D_m(X)$ (and similarly for $\mathscr{M}_m$ and $\mathscr{N}_m$.) $\Box$

\medskip

In particular, the sheaves of currents are acyclic. We note that the operator $d$ is defined between $\D_m$ and $\D_{m-1}$ or between $\mathscr{N}_m$ and $\mathscr{N}_{m-1}$, so we have the following resolutions of the constant sheaf on a complex space $X$, with $\dim_\C X_\rg=n$:
$$0\to\C\to\D_{2n}\xrightarrow{d}\D_{2n-1}\xrightarrow{d}\cdots\xrightarrow{d}\D_1\xrightarrow{d}\D_0\to0$$
$$0\to\C\to\mathscr{N}_{2n}\xrightarrow{d}\mathscr{N}_{2n-1}\xrightarrow{d}\cdots\xrightarrow{d}\mathscr{N}_1\xrightarrow{d}\mathscr{N}_0\to0$$

\subsection{Dolbeault complex}\label{ssc_dolb_cpx}

We define the sheaf $\mathscr{F}_{p,q}$ as the sheaf of $(p,q)$-currents whose $\debar$ is again a metric current.

\begin{Prp}$\mathscr{F}_{p,q}$ is a fine and soft sheaf and therefore acyclic.\end{Prp}
\noindent{\bf Proof:}The arguments of the previous section apply to $\mathscr{F}_{p,q}$ as well. $\Box$

\medskip

On a complex space $X$ of dimension $n$, we have the resolution
$$0\to K_p\to\mathscr{F}_{n-p,n}\xrightarrow{d}\mathscr{F}_{n-p,n-1}\xrightarrow{d}\cdots\xrightarrow{d}\mathscr{F}_{n-p,1}\xrightarrow{d}\mathscr{F}_{n-p,0}\to0$$
where $K_p=\ker\{\debar:\mathscr{F}_{n-p,n}\to\mathscr{F}_{n-p,n-1}\}$, is the analogue of the $(p,0)-$holomorphic forms on a complex manifold. It follows that
$$H^{q}(\Omega, K_p)=\frac{\ker\{\debar:\mathscr{F}_{n-p,n-q}(\Omega)\to\mathscr{F}_{n-p,n-q-1}(\Omega)\}}{\mathrm{img}\{\debar:\mathscr{F}_{n-p,n-q+1}(\Omega)\to\mathscr{F}_{n-p,n-q}(\Omega)\}}$$
for every $\Omega\subseteq X$.

\section{Holomorphic currents}

Given a complex space $X$ of pure dimension $n$, a current $T\in\D_{p,n}(X)$ of locally finite mass such that $\debar T=0$ is called a \emph{holomorphic $p-$current}.

\begin{Prp} Let $T\in \D_{p, n}(X)$, with locally finite mass, such that $\debar T$ is a metric current and $Y\subset X$ an analytic set with $\dim Y\leq n-2$ and $X_\sg\subseteq Y$. If $T\llcorner (X\setminus Y)$ is a holomorphic $p-$current on $X\setminus Y$, then $T$ is a holomorphic $p-$current on $X$.\end{Prp}

\noindent{\bf Proof: } $\supp \debar T$ must be contained in $Y$. Since $\debar T$ is a $n+p-1$ current and $\dim Y\leq n-2<n+p-1$ for every $p\geq0$, we must have $\debar T=0$. $\Box$

\medskip

\begin{Rem} If $\dim Y=n-1$ and $p>0$, we have the same conclusion, namely that if $T$ is holomorphic outside $Y$, then it is holomorphic on all of $X$. The only remaining case is $p=0$, $\dim Y=n-1$.\end{Rem}

This latter case is non trivial, because, for instance, the $(0,1)-$current 
$$T(f,g)=\int_{\C}\frac{f(z)}{z}\frac{\de g(z)}{\de\bar{z}}dz\wedge d\bar{z}$$
is locally finite mass and $\debar (T\llcorner\C^*)=0$, but $\debar T=2i\pi\delta_0$.

\subsection{Characterization by growth conditions}

Suppose $X$ is an analytic subspace of some complex hermitian manifold $M$; let $\star$ be the Hodge isomorphism on differential forms and $\xi\mapsto\omega_\xi$ be the isomorphism between the exterior powers of the tangent bundle and the exterior powers of the cotangent bundle. Given a holomorphic $p-$current $T$ on $X$, we know that there exists a ($dV_X-$summable) $(p,n)-$vector field $\xi$ on $X_\rg$ such that 
$$T(f,g)=\int_{X_\rg}\langle \xi, fdg_1\wedge\ldots\wedge dg_{n+p}\rangle dV_X\;,$$
with $(f,g)\in\D^{p,n}(X_\rg)$; such a vector field $\xi$ is holomorphic on $X_\rg$. We note that
$$\int_{X_\rg}\langle\xi,fdg_1\wedge\ldots\wedge dg_{n+p}\rangle dV_X=\int_{X_\rg}\omega_\xi\wedge\star(fdg_1\wedge\ldots\wedge dg_{n+p})=$$
$$(-1)^{(n+p)^2}\int_{X_\rg}\star(\omega_\xi)\wedge(f dg_1\wedge \ldots\wedge dg_{n+p})\;.$$
The form $\omega=(-1)^{(n+p)^2}\star\omega_\xi$ is the form associated to the current $T$ on $X_\rg$.

\begin{Teo}If $T$ is a holomorphic $p-$current on $X$, $\omega\in\Omega^{n-p}(X_\rg)$  the associated form, then for every compact $K\subset X_\sg$, and every $(f,g)\in\D^{p+n}(X)$, we have
\begin{eqnarray}\label{eq_int_omega1}\int_{K^\epsilon}\omega\wedge dg_1\wedge\ldots\wedge dg_{n+p}&\xrightarrow[\epsilon\to0]{}&0\\\label{eq_int_omega2}
\int_{bK^\epsilon}\omega\wedge( dg_1\wedge\ldots\wedge dg_{n+p-1})_{p,n-1}&\xrightarrow[\epsilon\to0]{}&0\;.\end{eqnarray}\end{Teo}

\noindent{\bf Proof: } $\omega$ is holomorphic hence $L^1_\loc$ on $X\setminus Y$ with respect to $dV_X$. So $T\llcorner(X\setminus Y)$ is locally flat on $X\setminus Y$, and consequently the representation by the form $\omega$ holds for any form with bounded coefficients and compact support.

Since $T$ is locally finite mass on $X$, given a compact set $K\subset Y$ and denoting by $K_\epsilon$ the $\epsilon-$neighbourhood of $K$ in $X$, we have 
$$\lim_{\epsilon\to 0} T\llcorner(1-\mathds{1}_{K_\epsilon})(f,g_1,\ldots, g_{n+p})=T(f(1-\mathds{1}_{K}),g_1,\ldots, g_{n+p})\;,$$
for every $(f,g)\in\D^{p+n}(X)$. Therefore, $\omega\wedge dg_1\wedge\ldots\wedge dg_{n+p}$ has to be locally integrable, which is equivalent to (\ref{eq_int_omega1}).

Moreover, given $(f,g)\in\D^{p+n-1}(X)$, with $f$, $g_i$ smooth, let $K\supset\supp f\cap X_\sg$. Then we have
$$0=\debar T(f,g)=\int_{X}\omega\wedge\debar(fdg_1\wedge\ldots\wedge dg_{n+p-1})=\lim_{\epsilon\to0}\int_{bK^\epsilon}f\omega\wedge(dg_1\wedge\ldots\wedge dg_{n+p-1})_{p,n-1}$$
which is equivalent to (\ref{eq_int_omega2}), by approximation. $\Box$

\medskip

In the particular case of $p=0$, $\dim X_\sg=n-1$, equation (\ref{eq_int_omega2}) becomes
$$\int_{bK^\epsilon}\omega\wedge dg_1\wedge\ldots\wedge dg_{n-1}\xrightarrow[\epsilon\to0]{}0\;.$$

\begin{Prp}Let $\omega\in\Omega^{n}(X_\rg)$ be a holomorphic $n-$form satisfying (\ref{eq_int_omega1}) and (\ref{eq_int_omega2}). Then $T_\omega\in\D_{0,n}(X_\rg)$, given by integration against $\omega$, extends to a holomorphic $T\in\D_{0,n}(X)$.\end{Prp}

\noindent{\bf Proof: }Condition (\ref{eq_int_omega1}) ensures that $\omega$ is in $L^1_\loc(X,dV_X)$, therefore $T$ is a well defined metric current on $X$ of locally finite mass. Moreover, as $T$ is of bidimension $(0,n)$, $\debar T= dT$, so
$$\debar T(f,g)=dT(f,g)=T(\sigma, f,g)=\lim_{\epsilon\to0}T(\sigma(1-\mathds{1}_{K^\epsilon}),f,g)\;,$$
because $T$ is locally finite mass. By Stokes theorem
$$T(\sigma(1-\mathds{1}_{K^\epsilon}),f,g)=\int_{bK^\epsilon}f\omega\wedge dg_1\wedge\ldots\wedge dg_{n-1}$$
which goes to $0$ as $\epsilon\to0$ by (\ref{eq_int_omega2}). Therefore $\debar T=0$. $\Box$

\medskip

It is worth noting that the only first condition guarantees that $T_\omega$ extends to a metric current on $X$.

\begin{Rem} We note that, if $\dim X_\sg\leq n-2$, then condition (\ref{eq_int_omega2}) is useless: $T$ being a metric current, also $\debar T=dT$ is and, since $T\llcorner X_\rg$ is holomorphic, $\supp\debar T\subseteq X_\sg$, thus, by Proposition \ref{prp_zero_high_dim}, $dT=\debar T=0$.\end{Rem}

\begin{Teo}Let $\omega\in\Omega^{n-p}(X_\rg)$ be a holomorphic $(n-p)-$form satisfying (\ref{eq_int_omega1}) and (\ref{eq_int_omega2}). Then the current $T_\omega\in\D^{p,n}(X_\rg)$, given by integration against $\omega$, extends to $T\in\D^{p,n}(X)$, holomorphic. \end{Teo}

\noindent{\bf Proof: } The extension $T$ is of locally finite mass, as observed in the previous Proposition. 

Now, let us assume that $(f,g)\in\D^{p+n-1}(X)$ is of class $\Ci^2$; then
$$\debar T(f,g)=\int_{X_\rg}\omega\wedge\debar(fdg_1\wedge\ldots\wedge dg_{n+p-1})=\lim_{\epsilon\to0}\int_{bK^\epsilon}\omega\wedge(fdg_1\wedge\ldots\wedge dg_{n+p-1})_{p,n-1}$$
where $K=\supp f\cap X_\sg$. By (\ref{eq_int_omega2}), we have $\debar T(f,g)=0$. For a generic metric form, we obtain the result by approximation. $\Box$

\medskip

\begin{Rem}We can substitute $K^\epsilon$ with any sequence of compact sets shrinking onto $K$.\end{Rem}

\medskip

The two conditions (\ref{eq_int_omega1}) and (\ref{eq_int_omega2}) relate the growth of the coefficients of $\omega$ to the decay of the volume element near $K$ and the growth of the coefficients of $d\omega$ to the decay of the volume of the $bK_\epsilon$'s. Obviously this description is quite rough because it doesn't take into account the differentials $dg_i$s.

\medskip

\paragraph{The smooth case} On a manifold, if $K$ is $(2n-2k)-$dimensional, the volume of $K_\epsilon$ goes to $0$ as $\epsilon^{2k}\H^{2n-2k}(K)$ and the differentials of the local coordinates are bounded away from $0$. In the same way, the volume of $bK_\epsilon$ goes down as $\epsilon\H^{2n-2k}(K)$. 

Holomorphic metric currents coincide with the classical holomorphic forms, in the smooth case. This is an easy consequence of the comparison theorem \ref{teo_comp_man}.

However, we note that the only meromorphic forms which satisfy (\ref{eq_int_omega1}) and (\ref{eq_int_omega2}) are the holomorphic ones. A meromorphic function can only have $1-$codimensional singularities and we can assume that such a singular set is locally described by $g=0$, $g\in\Ol$. Outside a set of (real) codimension at least $4$, we can take $g$ as a local coordinate and write the meromorphic function as $g^m\cdot h$, with $h\in \Ol$ independent of $g$.

The first condition and Fubini's theorem imply that $g^m$ has to be integrable in a neighbourhood of $0$, so $m\geq -1$; the second condition means that the supremum of $f$ on $bK_\epsilon$, which is $C\epsilon^m$, multiplied by the volume of $bK_\epsilon$, that is $\epsilon\H^{2n-2k}(K)$, has to go to $0$, so $C\H^{2n-2k}(K)\epsilon^{m-1}\to0$, which means $m\geq0$.

This holds outside a set of complex codimension at least $2$, which cannot contain the zero of a single holomorphic function.

In conclusion, in the smooth case, the two conditions imply that the extension of the current through an analytic subset is indeed a global holomorphic current.

\subsection{Examples}

We give now some examples of the previous characterization.

\paragraph{The cusp}

Let us consider the complex curve $X=\{z^3=w^2\}$ in $\C^2$ with parametrization $\phi:t\mapsto (t^2, t^3)$. This map is a biholomorphism between $X_\rg$ and $\C^*$, therefore every $\omega\in\Omega^{1-p}(X_\rg)$ corresponds to $\widetilde{\omega}\in\Omega^{1-p}(\C^*)$. For $p=0$ we have $\widetilde{\omega}=h(t)dt$. Conditions (\ref{eq_int_omega1}) and (\ref{eq_int_omega2}) for $\phi_*(h(t)dt)$ on balls centered at $(0,0)$, the only singular point of $X$, write
$$\int_{|t|<\epsilon'}|h(t)||dt\wedge\phi^*d\bar{z}|=\int_{|t|<\epsilon'}|h(t)||2\bar{t}dt\wedge d\bar{t}|\to0\;,$$
$$\int_{|t|<\epsilon'}|h(t)||dt\wedge\phi^*d\bar{w}|=\int_{|t|<\epsilon'}|h(t)||3\bar{t}^2dt\wedge d\bar{t}|\to0\;.$$
 By the first one the function $|h(t)||t|$ has to be in $L^1_\loc(\C,\mathcal{L})$ and this can happen only if $h(t)=g(t)t^{-2}$, with $g\in\Ol(\C)$. 
 The second one the implies 
$$\int_{|t|=\epsilon'}|h(t)||dt|\to0\;$$
which can happen only if $h(t)$ is indeed holomorphic on $\C$.

In conclusion, the holomorphic $(0,1)-$currents on $X$, are the pushforwards of holomorphic $(0,1)-$currents on $\C$ through the parametrization map $\phi$.

\bigskip

The exceptional divisor of the desingularization of a complex curve is made up of isolated points, therefore considering the projective curve instead of the affine one doesn't affect the computations above: the pullbacks of the two metrics are locally equivalent.

Moreover, also the holomorphic functions on a deleted neighborhood are the same. 

Both these considerations are false in higher dimension.

\medskip

Every complex curve admits a smooth normalization, by Ruckert's Parametrization Theorem; this observation will be developed in Chapter IV.

\paragraph{The affine quadratic cone}

Let $X=\{z^2=xy\}$ be a complex surface in $\C^3$, with $(0,0,0)$ as the only singular point. We consider on it the metric induced by this embedding in $\C^3$.

We have the parametrization $\phi:(s,t)\mapsto(st, t/s, t)$, with $(s,t)\in\C^*\times\C$, whose image is $X\setminus\{y=0\}\cup\{(0,0,0)\}$. The preimage of the singular point is $E=\{t=0\}$ and the pullback of a $(0,2)-$forms is of the form
$$\left(\psi_1(st,t/s, t)\frac{\bar{t}}{\bar{s}}+\psi_2(st,t/s, t)\bar{t}-\psi_3(st,t/s,t)\frac{\bar{t}}{\bar{s}^2}\right)d\bar{s}\wedge d\bar{t}\;.$$
The preimage of a ball around $(0,0,0)$ of radius $\epsilon$ is contained, in the $(s,t)$ plane, in the set
$$C(\epsilon)=\left\{(s,t)\in\C^2\ :\ |t|<\epsilon, \frac{|t|}{\epsilon}<|s|<\frac{\epsilon}{|t|}\right\}\;.$$
Therefore a holomorphic function on $C(\epsilon)\setminus E$ is a sum of monomials $s^nt^m$ and a holomorphic $(2,0)-$form is a sum of $s^nt^mds\wedge dt$.

Given $h(s,t)ds\wedge dt$, the first condition we have to check is
$$\int_{C(\epsilon)}|h(s,t)|t|(1+|s|+|s|^{-2})ds\wedge d\bar{s}\wedge dt\wedge d\bar{t}\to0$$
as $\epsilon\to0$.

If we assume $h(s,t)=s^nt^m$, then we have
$$\int_{a<|s|<b}|s|^n(1+|s|^{-2}+|s|^{-1})ds\wedge d\bar{s}=2\pi\int_{a}^b|s|^{n-3}(|s|^{2}+|s|^{1}+1)d|s|=$$
$$2\pi\left(\frac{b^{n+2}-a^{n+2}}{n+2}+\frac{b^{n}-a^{n}}{n}+\frac{b^{n+1}-a^{n+1}}{n+1}\right)\;.$$
Here we assume that none among $\{n+1, n, n-1\}$ is equal to $-1$. Otherwise a logarithm will appear. 
Now, $a=|t|/\epsilon$ and $b=\epsilon/|t|$, so, for example, 
$$b^n-a^n=\frac{\epsilon^{2n}-|t|^{2n}}{\epsilon^n|t|^n}\;.$$
Then
$$2\pi\int_{0}^\epsilon|t|^m|t|\frac{\epsilon^{2n}-|t|^{2n}}{\epsilon^n|t|^n}|t|d|t|=\int_0^\epsilon |t|^{2+m-n}\epsilon^{n}-|t|^{2+m+n}\epsilon^{-n}d|t|$$
which converges if and only if $2+m-n\geq0$ and $2+m+n\geq 0$. Integrating also the other terms we obtain the other conditions $1+m-n\geq0$ and $3+m+n\geq 0$, $m-n\geq0$ and $m+n+4\geq 0$. These six conditions are satisfied if and only if $m\geq \max\{n, -n-2\}$.
 
\medskip

To check (\ref{eq_int_omega2}), we need to compute the pullbacks through $\phi$ of the $(0,1)-$forms on $X$, i.e. restrictions of $(0,1)-$forms on $\C^3$. We notice that the boundary of $C(\epsilon)$ is
$$bC(\epsilon)=\{(s,t)\ :\ |t|\leq\epsilon, |s|=|t|/\epsilon\}\cup\{(s,t)\ :\ 0<|t|\leq\epsilon, |s|=\epsilon/|t|\}\;.$$
Now,
$$\phi^*d\bar{x}=\bar{s}d\bar{t}+\bar{t}d\bar{s}\;,\qquad \phi^*d\bar{y}=\frac{d\bar{t}}{\bar{s}}-\frac{\bar{t}}{\bar{s}^2}d\bar{s}\;,\qquad \phi^*d\bar{z}=d\bar{t}\;.$$
So we need that
$$\int_{bC(\epsilon)}|s|^n|t|^m|t|(1+|s|^{-2})|d\bar{s}\wedge dt\wedge d\bar{t}|=4\pi^2\int_0^\epsilon |t|^{m+1-n}\frac{\epsilon^{2n+2}-|t|^{2n+2}}{\epsilon^{n+1}}d|t|\xrightarrow[\epsilon\to0]{}0\;.$$
 This implies $m+1-n\geq0$, $m+3+n\geq0$; moreover, the result is going to $0$ with $\epsilon$ if and only if $m+3\geq0$.
 
 We also need that
 $$\int_{bC(\epsilon)}|s|^n|t|^m(|s|+1+|s|^{-1})|ds\wedge dt\wedge d\bar{t}|\to 0\;.$$
 This gives $m-n-1\geq0$, $m+n+1\geq0$ for the integrability and $m+1\geq0$ for the limit to be $0$.
 
\medskip

We conclude that $m+1\geq0$ and $m\geq\max\{n,-n-2\}$.

%\paragraph{The projective quadratic cone}
%
%Let $X=\{[x,y,z,t]\in\CP^3\ :\ xy=z^2\}$, with the Fubini-Study metric induced by the immersion in $\CP^3$. In this case, a deleted neighborhood $U_\rg$ of $[0,0,0,1]$ is biholomorphic to the product $\mathbb{D}^*\times Q$, with $Q$ a complex nonsigular conic; if we denote by $\pi:\mathbb{D}^*\times Q\to\mathbb{D}^*$ the projection, we can write every holomorphic function on $U_\rg$ as $h\circ\pi$, with $h\in\Ol(\mathbb{D}^*)$, because $Q$ is compact, so every holomorphic function on it is constant.
%
%Let $\phi:\C^2\to\CP^3$ be given by
%$$\phi(s,t)=[st, t/s, t, 1]\;.$$
%With respect to the Fubini-Study metric, in the chart $t\neq0$, the distance $\rho([u,1],[v,1])$ between two points is given by
%$$\cos\rho([u,1],[v,1])=\frac{|(u,v)+1|}{\sqrt{|u|^2+1}\sqrt{|v|^2+1}}$$
%so
%$$\cos(\rho([0,1],[\phi(s,t),1])=\frac{1}{\sqrt{|st|^2+|t/s|^2+|t|^2+1}}\;.$$
%Then the sets
%$$\left\{|st|^2+|t/s|^2+|t|^2+1< 1+3\epsilon^2\right\}$$
%are a collection of open sets, shrinking down to $[0,0,0,1]$. These sets are the same $C(\epsilon)$ described before. Their volume goes to zero as $\epsilon$, because they are comparable with balls of radius $\sqrt{\epsilon}$.
%
%\textcolor{red}{To be finished}

\subsection{Poincar\'e lemma}

Let us suppose that the space $X$, with $\dim X_\rg=n$, can be contracted locally at every point $x\in X$ by a Lipschitz contraction; that is, for every $x\in X$ there exist a neighbourhood $U_x$ and a Lipschitz  map $H:[0,1]\times U_x\to U_x$ such that $H(0,y)=y$, $H(1,y)=x$ for every $y\in U_x$.

Let $T\in M_{p,n}(X)$ be a normal holomorphic current, i.e. $\debar T=0$ and $dT=\de T\in M_{q-1,n}(X)$; $T$ defines a normal holomorphic current on $U_x$ for every $x$.

By \cite{ambrosio1}, the functional $S=T\vert_{U_x}\times [0,1]$ is well defined and is indeed a normal metric current on $U_x$; we note that
$$dS=-(dT)\times[0,1]+T\;.$$

If we suppose $dT=0$, we get $dS=T$; moreover, $dS=\de S+\debar S=T_{p,n}+ T_{p,n+1}$, but by \ref{prp_zero_high_dim}, no nontrivial component of bidimension $(p,n+1)$ can be present, which means that $T=dS=\de T$.

We can summarize these considerations in the following result.

\begin{Prp}Let $X$ be a complex space, which is locally Lipschitz contractible, then given a holomorphic $(p,n)-$current $T$ with $dT=0$, for evert $x\in X$ we can find a neighbourhood $U_x$ and a current $S\in M_{p+1, n}(U_x)$ such that $dS=T$ on $U_x$.\end{Prp}

\chapter{Analysis on singular complex spaces}
\epigraphhead[60]{\epigraph{They both savoured the strange warm glow of being much more ignorant than ordinary people, who were only ignorant of ordinary things.}{T. Pratchett - \emph{Equal Rites}}}

We give a definition of the Sobolev space $W^{1,2}$ on singular spaces, characterizing its elements as those functions which are locally $W^{1,2}$ on the regular part and whose gradient is square-integrable on the whole space. 

We compute the capacity of the singular set with respect to this Sobolev space and find it to be zero; as a consequence, we obtain an approximation result with functions supported away from the singularity.

The definitions of the local analogue of this space and of the corresponding spaces of vector fields appear in the third section; the corresponding approximation results for vector fields do not hold in general, so we introduce an additional hypothesis (which boils down essentially to the requirement that such a density holds).

Under this hypothesis, we replicate H\"ormander's $L^2$ techniques, paying attention also to the regularity results we obtain for the solutions. The existence on a given space of, let us say, holomorphic vector fields of  lower regularity than expected, can be employed to show that in such a space the density hypothesis does not hold. 

\section{Sobolev spaces}

Let $X$ be a complex analytic space embedded in $\C^N$, with $\dim_\C X_\rg=n$; we consider on $X$ the Hausdorff measure $\H^{2n}$ induced by the euclidean distance in $\C^N$.

We will say that a function $f:X\to\R$ belongs to $\Ci^\infty(X)$ if there exist a neighborhood $U$ of $X$ in $\C^N$ and a function $F:U\to\R$ such that $F\in\Ci^\infty(U)$ and $f=F\vert_X$; similarly, a function $f:X\to\R$ belongs to $\Ci^\infty_c(X)$ if $f\in\Ci^\infty(X)$ and $\supp f$ is compact in $X$.

A complex-valued (or vector-valued) function will be said to be of class $\Ci^\infty$ (resp. $\Ci^\infty_c$) on $X$ if its components are in $\Ci^\infty(X)$ (resp. $\Ci^\infty_c(X)$). If not otherwise stated, all the scalar functions will be complex-valued and by $TX_\rg$ we will denote the complexification of the real tangent bundle, i.e. $T_\R X_\rg\otimes_\R\C$.

Given $\Omega\subseteq X$, we will write $\Omega_\rg$ for $\Omega\cap X_\rg$ and $\Omega_\sg$ for $\Omega\cap X_\sg$.

\subsection{Definition and properties}

We consider the space $\Gamma(X_\rg, TX_\rg)$ of smooth vector fields on $X_\rg$ and we call \emph{regular} a vector field $Y$ such that $Y(\phi)\in L^\infty(X,\H^{2n})$ for every $\phi\in\mathcal{C}^\infty_c(X)$ and the following estimate holds:
$$\|Y(\phi)\|_\infty\leq M(\|\phi\|_{\infty,X}+\||\nabla\phi|\|_{\infty,X})$$
for some positive constant $M$, for every $\phi\in\mathcal{C}^\infty_c(X)$.
Let us denote by $R(X)$ the vector space of regular vector fields on $X_\rg$. For an open set $\Omega\subseteq X$ the space $R(\Omega)$ is defined in the same manner.

\begin{Rem}\label{rem_reg}If we denote by $Z'_1,\ldots, Z'_N$ the projections (with respect to the euclidean structure) of the coordinate fields in $\C^N$ on $TX_\rg$, we have that $Z_i', \overline{Z_i'}\in R(X)$ for $i=1,\ldots, N$. Therefore, the projection on $TX_\rg$ of every smooth vector field in $\C^N$ is in $R(X)$.\end{Rem}

Given an open set $\Omega\subset X$, let $W^{1,2}(\Omega)$ be the subspace of $L^2(\Omega,\H^{2n})$ defined by the following condition: $f\in W^{1,2}(\Omega)$ if and only if for every vector field $Y\in R(\Omega)$ there exists $g\in L^2(\Omega,\H^{2n})$ such that
\begin{equation}\label{eq_def_sob0}\int_{\Omega_\rg}Y(\phi) fd\H^{2n} = -\int_{\Omega_\rg}\phi gd\H^{2n}\qquad \forall\; \phi\in\Ci^\infty_c(\Omega)\end{equation}

%\begin{Rem} The definition of regular fields is motivated by the fact that we want the inclusion $\Ci^\infty_c(\Omega)\subseteq W^{1,2}(\Omega)$ for every open set $\Omega\subset X$.\end{Rem}

\begin{Rem} We note that the given definition implies that the map from $R(X)$ to $L^2(X,\H^{2n})$, $Y\mapsto g=g(Y)$, induced by $f$ is linear. So we can find an element $d f$ of $\Gamma(\Omega_\rg, T^*\Omega_\rg)$ such that 
$$\langle d f, Y\rangle=g(Y)\qquad\forall\; Y\in R(\Omega)$$
\end{Rem}

Therefore, we have a map $X_\rg\ni x\mapsto d_x f\in \Lambda^1(T\Omega_\rg)$ so that $\langle d f, Y\rangle$ is in $L^2(X_\rg, \H^{2n})$ for every $Y\in R(X)$. If we fix a Riemannian metric $g$ on $X_\rg$, we have an isomorphism between $T\Omega_\rg$ and $T^*\Omega_\rg$ and allows us to define the vector field $\nabla_g f$ 
$$g(\nabla_g f, Y)=\langle d f, Y\rangle$$
\label{pg_norma_eucl}We now consider the application $x\mapsto \|\nabla f_x\|_{g,x}$, where the norm is taken with respect to the chosen metric. In what follows, we will denote by $\nabla f$ the vector field obtained taking $g$ equal to the metric induced by the inclusion of $X$ in $\C^N$ and $|\nabla f|$ equal to the function $x\mapsto\|\nabla f_x\|_{g,x}$ (with $g$ again induced by the inclusion).

Now, we can write condition (\ref{eq_def_sob0}) more conveniently as
\begin{equation}\label{eq_def_sob}\int_{\Omega_\rg}Y(\phi) fd\H^{2n} = -\int_{\Omega_\rg}\phi \langle \nabla f, Y\rangle d\H^{2n}\end{equation}
and we say that $f$ is in $W^{1,2}(\Omega)$ if $f\in L^2(\Omega,\H^{2n})$ and, for every $Y\in R(X)$, $\langle \nabla f, Y\rangle \in L^2(\Omega,\H^{2n})$ and (\ref{eq_def_sob}) holds.

We will repeatedly use a known result about the measure of the intersection of a ball and a complex analytic set, whose proof can be found in \cite[Cor. 2]{lelong1}.

\begin{Prp}\label{prp_lelong}For every $p\in X$, there exist a ball $B=B(p,r)$ and a finite number $k(B)$ such that
$$\H^{2n}(X\cap B(q,\rho))\leq k(B) \rho^{2n}$$
for every $B(q,\rho)\subseteq B$.\end{Prp}

\medskip

From now on, we will say that $B$ is a ball in $\Omega$ if $B=B(q,\rho)\cap\Omega$; we set
$$k(\Omega)=\sup_{B\textrm{ ball in }\Omega}k(B)$$
In what follows we will always suppose that $k(\Omega)<\infty$, as it is the case when $\Omega$ is relatively compact, for instance.

We will also need the following property of Sobolev functions in Euclidean domains. 

\begin{Prp}\label{prp_max_Sob} Let $\Omega\subset X_\rg$ be an open set  such that $\overline{\Omega}\cap X_\sg=\emptyset$ and $f\in W^{1,2}(\Omega)$. Then the functions
$$f^+=\max\{f,0\}\ \textrm{ and }\ f^-=-\min\{f,0\}$$
belong to $W^{1,2}(\Omega)$ and
$$\nabla f^+=\left\{\begin{array}{cc}\nabla f&\textrm{a.e. on }\{f>0\}\\0&\textrm{a.e. on }\{f\leq0\}\end{array}\right.$$
$$\nabla f^-=\left\{\begin{array}{cc}0&\textrm{a.e. on }\{f\geq0\}\\-\nabla f&\textrm{a.e. on }\{f<0\}\end{array}\right.\;.$$
\end{Prp}
\noindent{\bf Proof: }If $\Omega$ is a coordinate chart the statement follows from the classic Sobolev theory in $\R^k$ \cite[Sec.~4.2-Thm.~4]{evans1}.

In the general case, we can covering $\Omega$ with coordinate charts $\{\Omega_j\}$ and construct a partition of unity $\{\chi_j\}$ subordinated to this open covering; as the closure of $\Omega$ doesn't intersect $X_\sg$, we can ask that, for every $p\in \overline{\Omega}$, there exist only a finite number of open sets $\Omega_j$ such that $p\in\overline{\Omega_j}$; moreover, we can assume $\|\nabla \chi_j\|\leq L$ for every $j$, for some positive real constant $L$.

Now, the function $f_j=f\chi_j$ belongs to $W^{1,2}(\Omega)$ and $W^{1,2}(\Omega_j)$ and, as $\Omega_j$ is biholomorphic to some open set of $\C^n$, we know that $f_j^+$ and $f_j^-$ both belong to $W^{1,2}(\Omega_j)$ and their gradients are given by the formulas above. Then, obviously $f_j^+$ and $f_j^-$ belong to $W^{1,2}(\Omega)$; moreover we have 
$$f^+=\sum_j f_j^+\qquad f^-=\sum_jf_j^-$$
as both sums are locally finite. 

We observe that, locally, 
$$\nabla f^+=\sum_j\nabla f_j^+\qquad f^-=\sum_jf^-_j\qquad\ \ \mathrm{a.e.}$$
as these sums are locally finite too. But then we have immediately that $|\nabla f^+|\leq|\nabla f|$ and $|\nabla f^-|\leq|\nabla f|$ almost everywhere, i.e. $f^+, f^-\in W^{1,2}(\Omega)$. The formulas for their gradients follow at once from the analogue results for $f^+_j$ and $f^-_j$. $\Box$

\medskip

As a corollary, we have that $h=\max\{f,g\}\in W^{1,2}(\Omega)$ for every $f,g\in W^{1,2}(\Omega)$ and
$$\nabla h=\left\{\begin{array}{cc}\nabla f&\textrm{a.e. on }\{f\geq g\}\\\nabla g &\textrm{a.e. on }\{f\leq g\}\end{array}\right.$$

\begin{Teo}\label{teo_caratt_W}A function $f\in L^2(\Omega,\H^n)$ belongs to $W^{1,2}(\Omega)$ if and only if $f\in W^{1,2}(U)$ for every $U\Subset \Omega_\rg$ and $|\nabla f|\in L^2(\Omega,\H^{2n})$.\end{Teo}
\noindent{\bf Proof: } One implication is obvious: if $f\in W^{1,2}(\Omega)$, then $f\in W^{1,2}(U)$ for every $U\Subset\Omega$ so, in particular, for every $U\Subset\Omega_\rg$; moreover, by Remark \ref{rem_reg}, $\langle Z'_i,\nabla f\rangle\in L^2(X,\H^{2n})$, i.e. $|\nabla f|\in L^2(X, \H^{2n})$.

Conversely, for every $\phi\in\Ci^\infty_c(\Omega)$ and $\epsilon>0$,  since $\H^{2n-1}(\Omega_\sg)=0$ and $\Omega_\sg$ is closed in $\Omega$, we can cover $\Omega_\sg\cap\supp\phi$ with finitely many balls $\{B_k(p_k,\rho_k)\}_{k=1}^m$, such that $\rho_k<\epsilon$ for every $k$ and $\sum_h \rho_h^{2n-1}<\epsilon$. Let $B'_k$ be the open ball with center $p_k$ and radius $2\rho_k$; for $\epsilon$ small, the balls $B_k'$ have compact closure in $\Omega$.

We remark that 
$$\H^{2n}\left(\bigcup_k B'_k\right)\leq \sum_k\H^{2n}(B'_k)\leq k(\Omega)s^{2n}\sum_k\rho_k^{2n}\leq k(\Omega)2^{2n}\epsilon\sum_k\rho_k^{2n-1}\leq k(\Omega)2^{2n}\epsilon^2\;.$$

We can construct a smooth cut-off function $\phi_k:\C^N\to[0,1]$ such that $\supp \phi_k\subset\Omega\setminus B_k$ and $\phi_k\equiv 1$ in $\Omega\setminus B'_k$ and $\|\nabla\phi_k\|_\infty\leq 2/\rho_k$.

Let now $f$ be a bounded function, belonging to $W^{1,2}(U)$ for every $U\Subset\Omega_\rg$ and whose gradient is in $L^2(\Omega_\rg,\H^{2n})$. Then, using Proposition \ref{prp_lelong}, for every $Y\in R(X)$ we have

$$\left|\int_{\Omega_\rg}Y(\phi) fd\H^{2n} - \int_{\Omega_\rg}Y\left(\phi\prod_k\phi_k\right) fd\H^{2n} \right|$$
$$\leq\int_{\bigcup_k B'_k}|Y(\phi)f|d\H^{2n}+\|\phi\|_\infty\sum_k\int_{B'_k\setminus B_k}|f\nabla\phi_k|d\H^{2n}$$
$$\leq\|Y(\phi)\|_\infty\|f\|_\infty\H^{2n}\left(\bigcup_k B'_k\right)+\|\phi\|_\infty\|f\|_\infty\sum_k\|\nabla\phi_k\|_\infty\H^{2n}(B'_k\setminus B_k)$$
$$\leq\|Y(\phi)\|_\infty\|f\|_\infty k(\Omega)2^{2n}\epsilon^2+2^{2n}k(\Omega)\|\phi\|_\infty\|f\|_\infty\sum_k\|\nabla\phi_k\|_\infty\rho^{2n}_k$$
$$\leq\|Y(\phi)\|_\infty\|f\|_\infty k(\Omega)2^{2n}\epsilon^2+2k(\Omega)2^{2n}\|\phi\|_\infty\|f\|_\infty\epsilon.$$

Similarly, we obtain

$$\left|\int_{\Omega_\rg}\phi \langle \nabla f, Y\rangle d\H^{2n}-\int_{\Omega_\rg}\phi\prod_k\phi_k \langle\nabla f, Y\rangle d\H^{2n}\right|\leq2\|\phi\|_\infty\int_{\Omega_\rg\cap\bigcup_k B'_k}\!\!\!\!\!\!\!\!\!\!\!\!\!\!\!\!\!\!|\langle\nabla f, Y\rangle|d\H^{2n}$$
$$\leq2\|\phi\|_\infty\|\langle\nabla f, Y\rangle\|_2\H^{2n}\left(\bigcup_k B'_k\right)^{1/2}\leq 2^{n+1}\|\phi\|_\infty\|\langle \nabla f, Y\rangle\|_2\sqrt{k(\Omega)}\epsilon$$

Now, we know that
$$\int_{\Omega_\rg}Y\left(\phi\prod_k\phi_k\right) fd\H^{2n}=-\int_{\Omega_\rg}\phi\prod_k\phi_k \langle \nabla f, Y\rangle d\H^{2n}$$
because 
$$\displaystyle\supp \left(\phi\prod_k\phi_k\right)\Subset \Omega\setminus\bigcup_kB'_k\subset\Omega_\rg\;.$$ 
We know that the sides of the last equality tend to the sides of (\ref{eq_def_sob}), so we can conclude that equality (\ref{eq_def_sob}) holds for $f$, so $f\in W^{1,2}(\Omega)$.

\medskip

Now, let $f$ be a function in $W^{1,2}(U)$ for every $U\Subset\Omega_\rg$, with $L^2$ gradient, and let 
$$f_M=\min\{f^+, M\}-\min\{f^-,M\}.$$ 
By Proposition \ref{prp_max_Sob}, we know that $f_M\in W^{1,2}(U)$ and $\nabla f_M=\nabla f$ where $|f|\leq M$ and $\nabla f_M=0$ elsewhere. 

Then, since $f_M\in W^{1,2}(U)$ for every $U\Subset\Omega_\rg$ and $|\nabla f_M|\in L^2(\Omega,\H^{2n})$, from the previous computations we infer $f_M\in W^{1,2}(\Omega)$. Given again $Y\in R(\Omega)$, for every $\phi\in \Ci^\infty_c(\Omega)$, we have 
$$\int_\Omega Y(\phi)f_M d\H^{2n}\longrightarrow \int_\Omega Y(\phi)f d\H^{2n}$$
as $M\to\infty$, because $f_M\to f$ in $L^2-$norm and $Y(\phi)\in L^\infty(\Omega,\H^{2n})$.  Moreover
$$\int_\Omega Y(\phi)f_M d\H^{2n}=-\int_\Omega \phi \langle \nabla f_M, Y\rangle d\H^{2n}$$
and $|\nabla f_M|\leq|\nabla f|$, so, given $N\in \N$, with $N>M$,
$$\int_{\Omega}|\nabla f_M-\nabla f_N|^2d\H^{2n}\leq 4\int_{M\leq |f|\leq N}|\nabla f|^2d\H^{2n}$$
and this goes to zero when $M, N$ go to infinity, so 
$$\int_{\Omega}|\langle \nabla f_M-\nabla f_N, Y\rangle|^2 d\H^{2n}\leq \int_{\Omega} |\nabla f_M-f_N|^2d \H^{2n}\int_{\Omega}|Y|^2d\H^{2n}$$
goes to $0$. Therefore $\{\langle \nabla f_M,Y\rangle\}_M$ is Cauchy and converges to $\langle \nabla f, Y\rangle$, whence
$$\int_\Omega \phi \langle \nabla f_M, Y\rangle d\H^{2n}\longrightarrow\int_\Omega \phi\langle \nabla f, Y\rangle d\H^{2n}$$
for every $\phi\in\Ci^\infty_c(X)$.

In conclusion, we have that
$$\int_\Omega Y(\phi)f d\H^{2n}=-\int_\Omega \langle \nabla f, Y\rangle\phi d\H^{2n}$$
and this shows that $f\in W^{1,2}(\Omega)$. $\Box$

\bigskip

\begin{Cor}We have the inclusion $\Ci^\infty_c(\Omega)\subseteq W^{1,2}(\Omega)$.\end{Cor}
\noindent{\bf Proof: } Obviously, if $f\in\Ci^\infty_c(\Omega)$, then $f\in\Ci^\infty(\overline{U})$ for every $U\Subset\Omega_\rg$, so $f\in W^{1,2}(U)$; moreover, with the notation of Remark \ref{rem_reg}, $|\nabla f|^2$ is bounded by
$$|Z_1'(f)|^2+\ldots+|Z_N'(f)|^2+|\overline{Z_1}'(f)|^2+\ldots+|\overline{Z_N}'(f)|^2$$
From that remark it follows also that $Z_i'(f)$ and $\overline{Z_i}'(f)$ are in $L^\infty$. So $|\nabla f|^2\in L^\infty(\Omega)$, therefore $|\nabla f|\in L^2(\Omega,\H^{2n})$. $\Box$

\bigskip

\begin{Prp}\label{prp_max_in_W}Let $f,g\in W^{1,2}(\Omega)$ and set $h=\max\{f,g\}$, then $h\in W^{1,2}(\Omega)$ and
$$\nabla h=\left\{\begin{array}{cc}\nabla f&\textrm{a.e. on }\{f\geq g\}\\\nabla g &\textrm{a.e. on }\{f\leq g\}\end{array}\right.\;.$$\end{Prp}
\noindent{\bf Proof: }As $f,g\in W^{1,2}(\Omega)$, by Theorem \ref{teo_caratt_W} we know that $f,g\in W^{1,2}(U)$ for every $U\Subset\Omega_\rg$ and $|\nabla f|,|\nabla g|\in L^2(\Omega,\H^{2n})$. Fix $U\Subset\Omega_\rg$, then $h\vert_U\in W^{1,2}(U)$, because of the corollary to Proposition \ref{prp_max_Sob}; moreover, if  $\Omega^1=\{f\geq g\}$ and $\Omega^2=\{g\geq f\}$, then we now that
$$\nabla h=\nabla f\qquad\H^{2n}-\textrm{a.e. in every }U\Subset\Omega^1_\rg=\Omega^1\cap\Omega_\rg$$
$$\nabla h=\nabla g\qquad\H^{2n}-\textrm{a.e. in every }U\Subset\Omega^2_\rg=\Omega^2\cap\Omega_\rg$$
therefore, by $\sigma-$additivity,
$$\nabla h=\nabla f\qquad\H^{2n}-\textrm{a.e. in }\Omega^1$$
$$\nabla h=\nabla g\qquad\H^{2n}-\textrm{a.e. in }\Omega^2$$
so that we have
$$\int_{\Omega}|\nabla h|^2d\H^{2n}\leq\int_{\Omega^1}|\nabla f|^2d\H^{2n}+\int_{\Omega^2}|\nabla g|^2d\H^{2n}\leq\|\nabla f\|_{2,\Omega}^2+\|\nabla g\|_{2,\Omega}^2<+\infty$$
i.e. $|\nabla h|\in L^2(\Omega,\H^{2n})$. In conclusion, by Theorem \ref{teo_caratt_W}, we have that $h\in W^{1,2}(\Omega)$. $\Box$

\subsection{Capacity of the singular set}
The capacity of a set $E\subseteq \Omega\subseteq X$ with respect to the space $W^{1,2}(\Omega)$ is defined by
$$
C_{1,2}^\Omega(E)=\inf \{\|\nabla f\|_{2}^2,\ f\in W^{1,2}(\Omega),\ E\subset\Int\{f\geq1\}\}
$$
where $\Int A$ means the topological interior of the set $A$, relative to $\Omega$. We remark that, by requiring $E\subset\Int\{f\geq1\}$, we actually ask that there exists a function in the equivalence class of $f$ in $W^{1,2}(\Omega)$ for which that inclusion holds.

\begin{Teo}\label{teo_cap_0}If $E\subseteq\Omega$ and $\H^{2n-2}(E)<+\infty$, then $C^\Omega_{1,2}(E)=0$.\end{Teo}
\noindent{\bf Proof: }Given an open set $U$, $E\subset U\subset \Omega$, we can find a countable family of balls
$B_k=B(p_k,\rho_k)$ whose union contains $E$, with
$B(p_k,2\rho_k)\subseteq U$ and
$$\sum_{k=1}^\infty c_{2n-2}\rho^{2n-2}_k\leq \H^{2n-2}(E)+1\;.$$
To this aim, for every $h\in\Z$ we define
$$E_h=E\cap\left\{p\in \Omega\ :\ \frac{1}{2^{h+1}}<\mathrm{dist}(p,b\Omega)\leq\frac{1}{2^h}\right\}\;.$$
Since the closure of $E_h$ is compact in $\Omega$, by the very
definition of spherical Hausdorff measure we can find a countable
covering of $E_h$ by balls $B(q_j,s_j)$ with $B(q_j,2s_j)\subset U$,
with $\sum_jc_{2n-2}s_j^{n-2}\leq\H^{n-2}(E_h)+2^{-|h|-1}$. Since
the $E_h$'s give a partition of $E$, we obtain the required family of balls $B(p_k,r_k)$ collecting all these balls.

Now, we want to prove that given $U$ as above there exists a
function $u\in W^{1,2}(\Omega)$ with the following properties:
\renewcommand{\labelenumi}{\roman{enumi}. }
\begin{enumerate}
\item[(i)] $0\leq u\leq 1$ and $u=0$ a.e on $\Omega\setminus U$;
\item[(ii)] $E\subset\Int\{u\geq 1\}$;
\item[(iii)] $\displaystyle{\int_{\Omega}|\nabla u|^2d\H^{2n}\leq
K(\Omega)(\H^{2n-2}(E)+1)}$,
\end{enumerate}
where $K(\Omega)$ depends only on $\Omega$.

To this aim, we choose $u_k:\C^N\to[0,1]$ piecewise smooth (and
radially linear) such that
$$u_k\vert_{B(p_k,\rho_k)}\equiv 1\;,\qquad u_k\vert_{X\setminus B(p_k,2\rho_k)}\equiv 0\;.$$
We observe that $\|\nabla u_k\|_\infty\leq \rho_k^{-1}$ on $\C^N$.
Moreover, we set
$$v_m=\max\{u_1,\ldots,u_m\}.$$
The function $v_m$ obviously belongs to $W^{1,2}(\Omega)$, as
$$|\nabla v_m|^2\leq \sum_{k=1}^m|\nabla u_k|^2\;$$
is non-negative and has the property (i). Moreover, applying Proposition~\ref{prp_lelong} once more, we
obtain
\begin{eqnarray}\label{upboundvm}
\int_{\Omega}|\nabla v_m|^2d\H^{2n}&\leq&
\sum_{k=1}^m\H^{2n}(\Omega_\rg\cap B(p_k,2\rho_k))\rho_k^{-2} \nonumber\\
&\leq&
k(\Omega)\sum_{k=1}^m 2^n\rho^{2n-2}_k\nonumber \\
&\leq&
K(\Omega)(\H^{2n-2}(E)+1),
\end{eqnarray}
 where $K(\Omega)=2^{2n}k(\Omega)/c_{2n-2}$ depends only on $\Omega$.

The function $u=\sup_mv_m=\lim_mv_m$ belongs to
$L^2(\Omega,\H^{2n})$ and satisfies (i) by monotone convergence. The
function $u$ is identically equal to $1$ on the union of the balls
$B(p_i,\rho_i)$, so that since this union contains $E$ also
condition (ii) is satisfied. Since $|\nabla v_m|$ is bounded in
$L^2(\Omega,\H^{2n})$, because of \eqref{upboundvm}, we obtain that
$u\in W^{1,2}(\Omega)$ and $\nabla v_m\to\nabla
u$ weakly in $L^2$. Then, the lower semicontinuity of the $L^2$ norm under
weak convergence ensures condition (iii).

\medskip

Now, we consider a non-increasing sequence of open neighborhoods of
$E$, $\{U_j\}_{j\in\N}$ such that $U_j\subset\Omega$ and
$\H^{2n}(U_j)\to 0$ as $j\to \infty$ and  construct a function
$w_1:\Omega\to [0,1]$, vanishing outside $U_1$, as described before. Let $V_1$ be an open set containing $E$ and such that
$V\subseteq\{w_1=1\}$ and construct a non-negative function $w_2:\Omega\to [0,1]$, vanishing
outside this open set$V_1\cap U_2$. By iterating this procedure, we obtain a sequence of
non-negative functions $\{w_j\}$, vanishing outside smaller and
smaller neighborhoods of $E$, equal to $1$ on some other
neighborhoods of $E$, uniformly bounded in the $W^{1,2}$ norm.

We remark that $\nabla w_j$ can be nonzero only in $U_j\cap
V_{j-1}\setminus V_j$, so the gradients are mutually orthogonal in
$L^2$.

We set
$$S_j=\sum_{k=1}^j\frac{1}{k}$$
$$g_j=\frac{1}{S_j}\sum_{k=1}^j\frac{w_k}{k}$$
The functions $g_j$ belong to $W^{1,2}(\Omega)$ and $E\subset\Int\{g_j\geq 1\}$; moreover
$$\|\nabla g_j\|^2\leq\frac{1}{S_j^2}\sum_{k=1}^j\frac{C^2}{k^2}\xrightarrow[j\to\infty]{}0$$
So, by the definition of $W^{1,2}(\Omega)-$capacity of $E$, we have that
$$C^\Omega_{1,2}(E)\leq\|\nabla g_j\|_2^2\xrightarrow[j\to\infty]{}0$$
i.e. $C^\Omega_{1,2}(E)=0$. $\Box$

\begin{Rem} Since the family $U_j$ above is non-increasing we also proved that for any
open set $U$, $E\subset U\subset\Omega$, there exists
a Sobolev function with arbitrarily small $W^{1,2}$ norm identically
equal to 1 on an open set containing $E$ and vanishing almost everywhere
on $\Omega\setminus U$.\end{Rem}

\bigskip

\begin{Teo}\label{teo_approx_W}$\Ci^{\infty}(\overline{\Omega})$ is dense in  $W^{1,2}(\Omega)$.\end{Teo}
\noindent{\bf Proof: } Fix $f\in W^{1,2}(\Omega)\cap L^\infty$ and apply Theorem\ref{teo_cap_0} with 
$E=\Omega_{\rm sing}$. We find a sequence of functions $g_j\in W^{1,2}(\Omega)$ with values in $[0,1]$ strongly convergent to $0$ in
$W^{1,2}(\Omega)$,  identically equal to $1$ in open
sets $V_j$ containing $E$ and with supports contained in open sets $U_j\subseteq\Omega$, with $\H^{2n}(U_j)\to 0$.
Now, the function $(1-g_j)f$ is supported in $\Omega\setminus V_j$,
so we can find $h_j\in \Ci^{\infty}(\overline{\Omega})$ satisfying $\|(1-g_j)f-h_j\|_{W^{1,2}(\Omega)}\leq 1/j$.
(ref needed)

Since $|g_j|\leq 1$Now, we have that
$$\|g_jf\|_{2}^2\leq\int_{U_j}|f|^2d\H^{2n}\rightarrow 0$$
Moreover
$$\|\nabla(g_jf)\|_2\leq\|g_j\nabla f\|_2+\|f\nabla g_j\|_2\leq
\int_{U_j}|\nabla f|^2d\H^{2n}
+\|f\|_\infty\|\nabla g_j\|_2\rightarrow 0.$$
So $h_j\to f$ in $W^{1,2}(\Omega)$.

In the case when $f$ is possibly unbounded, we can consider the functions $f_M$ defined \ref{teo_caratt_W}. It is immediate to see that
$|f_M|\leq |f|$ and that $|\nabla f_M|\leq|\nabla f|$ $\H^{2n}$-a.e. thus, by the dominated convergence theorem,
$f_M\to f$ in $W^{1,2}(\Omega)$. In order to ends the proof we approximate each $f_M$ with
$\Ci^{\infty}(\overline{\Omega})$-functions and then apply the classical diagonal method. $\Box$

\section{Generalizations}

We present some extensions and variations over the previous definitions and results; most of the proofs will only be sketched, as they resemble closely the corresponding ones from the previous pages.

\subsection{Restricting to a sub-bundle}\label{ssc_subbundle}

We can generalize the previous results, restricting further the space of test vector fields and consequently substituting the integrability condition on $|\nabla f|$ with a request on the appropriate element of the dual.

\medskip

Let $F$ be a sub-bundle of $TX_{\rg}$ on $X_\rg$ and let $R_F(X)=R(X)\cap \Gamma(X_\rg, F)$ be the space of regular sections of $F$. Then, for any $\Omega\subseteq X$, we can define the space $W^{1,2}_F(\Omega)$  as the space of functions $f$ belonging to $L^2(\Omega,\H^{2n})$ such that, for every $Y\in R_F(X)$ there exists $g\in L^2(\Omega,\H^{2n})$ satisfying (\ref{eq_def_sob0}) holds.

\medskip

The map which sends $Y$ to $g$ is still linear, so we can find an element $d^Ff\in \Gamma(\Omega_\rg, F^*)$ such that 
$$\langle d^Ff, Y\rangle=g(Y)$$
then, by duality, we can find a section $\nabla^F_g f$ of $F$ itself such that
$$g=\langle \nabla^F_gf, Y\rangle.$$
As before, we will drop the subscript $\nabla^F_gf$ whenever $g$ is the metric induced by the Euclidean metric in $\C^N$.
\medskip

In order to obtain an analogue of Theorem \ref{teo_caratt_W}, we need to adapt Proposition \ref{prp_max_Sob} to this new setting. Looking carefully at the proof given in \cite{evans1}, one can notice that, in fact, the argument exploited there gives a more precise result; namely, if $U$ is an open set in $\R^M$, $f, g$ belong to $L^2(U)$ and $Y$ is a smooth vector field on $U$ such that, for every $\phi\in\Ci^\infty_c(U)$, we have
$$\int_{U}fY(\phi)dx=-\int_Ug\phi dx$$
then also the positive and negative parts of $f$, denoted by $f^+$ and $f^-$, have the same property:
$$\int_{U}f^+Y(\phi)dx=-\int_{U\cap\{f\geq0\}}g\phi dx$$
$$\int_{U}f^-Y(\phi)dx=\int_{U\cap\{f\leq0\}}g\phi dx$$

So, given $U\Subset \Omega_\rg$, for every $f,h\in W^{1,2}_F(U)$, $\max\{f,h\}\in W^{1,2}_F(U)$. Moreover, in the proof of Theorem \ref{teo_caratt_W} we fixed the vector field at the beginning, therefore the same argument shows the following.

\begin{Teo}\label{teo_caratt_F}A function $f\in L^2(\Omega,\H^{2n})$ belongs to $W^{1,2}_F(\Omega)$ if and only if $f\in W^{1,2}_F(U)$ for every $U\Subset\Omega_\rg$ and $|\nabla^Ff|\in L^2(\Omega, \H^{2n})$.\end{Teo}

Moreover, as $W^{1,2}(\Omega)\subseteq W^{1,2}_F(\Omega)$, Theorem \ref{teo_cap_0} still holds for the capacity with respect to $W^{1,2}_F(\Omega)$:
$$
C^{\Omega,F}_{1,2}(E)=\inf \{\|\nabla^F f\|_{2}^2,\ f\in W^{1,2}_F(\Omega),\ E\subset\Int\{f\geq1\}\}
$$
thus, if $\H^{2n-2}(E)$ is finite, then $C^{\Omega, F}_{1,2}(E)=0$.

As another consequence of the inclusion $W^{1,2}(\Omega)\subseteq W^{1,2}_F(\Omega)$, we have that also Theorem \ref{teo_approx_W} holds: indeed, since the functions $g_{j}$ employed in the proof of the Theorem belong to $W^{1,2}_F(\Omega)$ as well, therefore, given $f\in W^{1,2}_F(\Omega)\cap L^\infty$, $(1-g_{j}f)$ belongs to $W^{1,2}_F(\Omega)$ and approximates $f$ as $j\to+\infty$. The approximation by bounded functions is granted by the analogue of Proposition \ref{prp_max_Sob}.

So $\Ci^\infty(\overline{\Omega})$  is dense in $W^{1,2}_F(\Omega)$ with respect to the norm $\|\cdot\|_2+\|\nabla^F\cdot\|_2$.

\medskip

Let us consider the decomposition of the complexified tangent $TX_\rg$ into the holomorphic and the anti-holomorphic factors:
$$TX_\rg=T'X_\rg\oplus T''X_\rg$$
Letting $F=T'X_\rg$, we obtain the space $W^{1,2}_{\de}(\Omega)$, $d^F$ is the $\de$ operator and $\nabla^F=\nabla^{(1,0)}$ is the $(1,0)-$component of the gradient. Similarly, letting $F=T''X_\rg$, we obtain the space $W^{1,2}_{\debar}(\Omega)$, $d^F$ is the $\debar$ operator and $\nabla^F=\nabla^{(0,1)}$ is the $(0,1)-$component of the gradient. Hence, Theorem \ref{teo_caratt_F} holds in these two particular cases, together with the density theorem for smooth functions on $\overline{\Omega}$.

\bigskip

\subsection{An alternative proof of Theorem \ref{teo_caratt_W}}\label{ssc_alt_proof}

We produce an alternative proof of the mentioned Theorem that doesn't make use of the approximation by bounded functions, employing the fact that $\H^{2n-2}(\Omega_\sg)<+\infty$. With the notations of Theorem \ref{teo_caratt_W} and the choice of $B_k$'s and $B'_k$'s made in the proof of Theorem \ref{teo_cap_0}, we notice that
$$\int_{\Omega_rg}\!\!\!\!\!Y(\phi)fd\H^{2n}=\int_{\Omega_\rg}\!\!\!\!\!\prod \phi_kY(\phi)fd\H^{2n}+\int_{\Omega_\rg}\!\!\!\!\!\left(1-\prod\phi_k\right)Y(\phi)fd\H^{2n}$$
$$=\int_{\Omega_\rg}\!\!\!\!\!Y\left(\phi\prod \phi_k\right)fd\H^{2n}-\int_{\Omega_\rg}\!\!\!\!\!\phi Y\left(\prod \phi_k\right)fd\H^{2n}+\int_{\Omega_\rg}\!\!\!\!\!\left(1-\prod\phi_k\right)Y(\phi)fd\H^{2n}\;.$$
Now, 
$$\left|\int_{\Omega_\rg}\!\!\!\!\left(1-\prod\phi_k\right)Y(\phi)fd\H^{2n}\right|\leq\int_{\bigcup B'_k}\!\!\!\!\!|Y(\phi)f|d\H^{2n}\leq\|Y(\phi)\|_\infty\left(\H^{2n}\left(\bigcup B'_k\right)\right)^{1/2}\|f\|_2$$
which goes to $0$ with $\epsilon$. Moreover
$$\int_{\Omega_\rg}\!\!\!\!\!Y\left(\phi\prod \phi_k\right)fd\H^{2n}=-\int_{\Omega_\rg}\!\!\!\!\!g\phi\prod\phi_kd\H^{2n}$$
and
$$\left|\int_{\Omega_\rg}\!\!\!\!\!g\phi\left(1-\prod\phi_k\right)d\H^{2n}\right|\leq\|\phi\|_{\infty}\|g\|_2\left(\H^{2n}\left(\bigcup B'_k\right)\right)^{1/2}$$
which, as before, goes to $0$ with $\epsilon$.

We need to show that
$$\int_{\Omega_\rg}\!\!\!\!\!\phi Y\left(\prod \phi_k\right)fd\H^{2n}\to 0$$
 as well.  We have
$$\left|\int_{\Omega_\rg}\!\!\!\!\!\phi Y\left(\prod \phi_k\right)fd\H^{2n}\right|\leq\|\phi\|_\infty\left(\sum\|\nabla\phi_k\|_2\right)\left(\int_{\bigcup B'_k}\!\!\!\!\!|f|^2d\H^{2n}\right)^{1/2}$$
which goes to $0$ as $\epsilon\to 0$ and $\bigcup B'_k$ converges to $\Omega_\sg$.

\subsection{$W^{1,p}$}

All the theory just developed can be transferred withouth many changes to the spaces $W^{1,p}$. Let $p>1$ be fixed.

Namely, $W^{1,p}(\Omega)$ is the subspace of $L^p(\Omega,\H^{2n})$ defined as follows: $f\in W^{1,p}(\Omega)$ if and only if for every vector field $Y\in R(\Omega)$ there exists $g\in L^p(\Omega,\H^{2n})$ such that
\begin{equation}\label{eq_def_sob_p0}\int_{\Omega_\rg}Y(\phi)fd\H^{2n}=-\int_{\Omega_\rg}\phi gd\H^{2n}\qquad \forall\; \phi\in\Ci^\infty_c(\Omega)\;.\end{equation}
We proceed to define $df\in \Gamma(\Omega_\rg, T^*\Omega_\rg)$ and, specifying a metric, $\nabla f\in \Gamma(\Omega_\rg, T\Omega_rg)$, so that we can rewrite equation (\ref{eq_def_sob_p0}) as
\begin{equation}\label{eq_def_sob_p}\int_{\Omega_\rg}Y(\phi)fd\H^{2n}=-\int_{\Omega_\rg}\phi\langle \nabla f, Y\rangle d\H^{2n}\;.\end{equation}

The analogue of Proposition \ref{prp_max_Sob} holds, because the result is true on $\R^k$ (see again \cite[Sec. 4.2-Thm. 4]{evans1}); so we have the following characterization of $W^{1,p}(\Omega)$.

\begin{Teo}\label{teo_caratt_Wp} A function $f\in L^p(\Omega,\H^{2n})$ belongs to $W^{1,p}(\Omega)$ if and only if $f\in W^{1,2}(U)$ for every $U\Subset \Omega_\rg$ and $|\nabla f|\in L^p(\Omega,\H^{2n})$.\end{Teo}
\noindent{\bf Proof: } As in the result for $p=2$, one implication is obvious. 

Using the same notations and ideas of the proof of Theorem \ref{teo_caratt_W}, we define $\{B_k\}$, $\{B'_k\}$ and $\{\phi_k\}$ in the same way; now, if $f$ is also assumed to be bounded, we have for every $Y\in R(\Omega)$

$$
\int_{\Omega}Y(\phi)fd\H^{2n}-\int_{\Omega}Y\left(\phi\prod\phi_k\right)fd\H^{2n}\xrightarrow[\epsilon\to0]{}0\qquad\forall\;\phi\in\Ci^\infty_c(\Omega)\;,
$$
with exactly the estimates given in the proof for $p=2$. The proof of the other estimate runs as follows:

$$\left|\int_{\Omega_\rg}\phi \langle \nabla f, Y\rangle d\H^{2n}-\int_{\Omega_\rg}\phi\prod_k\phi_k \langle\nabla f, Y\rangle d\H^{2n}\right|\leq2\|\phi\|_\infty\int_{\Omega_\rg\cap\bigcup_k B'_k}\!\!\!\!\!\!\!\!\!\!\!\!\!\!\!\!\!\!|\langle\nabla f, Y\rangle|d\H^{2n}$$
$$\leq2\|\phi\|_\infty\|\langle\nabla f, Y\rangle\|_p\H^{2n}\left(\bigcup_k B'_k\right)^{1/q}\leq 2^{n+1}\|\phi\|_\infty\|\langle \nabla f, Y\rangle\|_p\sqrt[q]{k(\Omega)\epsilon}$$

where $q$ is such that $p^{-1}+q^{-1}=1$. The rightmost quantity goes to $0$ as $\epsilon\to 0$, because $q>0$; now the proof goes on as in the case $p=2$, for a bounded $f$.

To obtain also the result for unbounded $f$, we apply the same approximation procedure, defining the functions $f_M$ and noting that, as $f_M\to f$ in $L^p$, 

$$\int_{\Omega}Y(\phi)f_Md\H^{2n}\to\int_{\Omega}Y(\phi)fd\H^{2n}\;.$$

Moreover, we know that $|\nabla f_M|\leq |\nabla f|$, so, by Minkowski inequality, we have
$$\int_{\Omega}|\nabla f_M-\nabla f|^pd\H^{2n}\leq\left\{\left(\int_{M\leq|f|\leq N}\!\!\!\!|\nabla f_M|^pd\H^{2n}\right)^{1/p}\!\!\!\!+\left(\int_{M\leq |f|\leq N}\!\!\!\!|\nabla f_N|^pd\H^{2n}\right)^{1/p}\right\}^p$$
$$\leq2^p\int_{M\leq|f|\leq N}\!\!\!\!|\nabla f|^pd\H^{2n}$$
which goes to $0$ as $M\leq N$ go to infinity. We can therefore conclude as in the case $p=2$.  $\Box$

\medskip

\begin{Rem}We can also avoid the use of the approximation by bounded functions, adapting the alternative proof of this result given in Section \ref{ssc_alt_proof}, but we need to strenghten the hypothesis: let $q$ be such that $p^{-1}+q^{-1}=1$, then, if $\H^{2n-q}(\Omega_\sg)<+\infty$, we have
$$\sum\|\nabla\phi_k\|^q_q\leq k(\Omega)\sum2^n\rho_k^{2n-q}\leq K(\Omega)(\H^{2n-q}+1)<+\infty$$
and we obtain the same result.
\end{Rem}

With the same proof, we get the inclusion of the smooth test functions in our space.

\begin{Cor}We have the inclusion $\Ci^{\infty}_c(\Omega)\subseteq W^{1,p}(\Omega)$.\end{Cor}

Obviously, also $W^{1,p}(\Omega)$ has a lattice structure.

\begin{Prp}Let $f,g\in W^{1,p}(\Omega)$.  Then $h=\max\{f,g\}\in W^{1,p}(\Omega)$ and
$$\nabla h=\left\{\begin{array}{cc}\nabla f&\textrm{a.e. on }\{f\geq g\}\\\nabla g &\textrm{a.e. on }\{f\leq g\}\end{array}\right.\;.$$\end{Prp}

\medskip

We define the capacity of a set $E\subseteq \Omega\subseteq X$ in the same way:

$$C_{1,p}^\Omega(E)=\inf \{\|\nabla f\|_{p}^p,\ f\in W^{1,p}(\Omega),\ E\subset\Int\{f\geq1\}\}$$

with the same cautions pointed out before. The corresponding result for sets of zero capacity is the following.

\begin{Teo} If $E\subseteq\Omega$ and $\H^{2n-p}(E)<+\infty$, then $C_{1,p}^\Omega(E)=0$.\end{Teo}
\noindent{\bf Proof: }With the same notation of Theorem \ref{teo_cap_0}, we have
$$\sum_{k=1}^\infty c_{2n-p}\rho_k^{2n-p}\leq \H^{2n-p}(E)+1\;.$$
The function $u$ we want to construct will have properties (i), (ii) and
$$\textrm{(iii-p) } \int_{\Omega}|\nabla u|^pd\H^{2n}\leq K(\Omega)(\H^{2n-p}(E)+1)\;.$$
The proof now goes on exactly as in the case $p=2$; indeed, as we supposed $p>1$, the series $\sum_k k^{-p}$ is convergent, so the result follows in the same way. $\Box$

\medskip

As a corollary, we again obtain the density of smooth functions

\begin{Teo}$\Ci^\infty(\overline{\Omega})$ is dense in $W^{1,p}(\Omega)$.\end{Teo}

As a proof of the fact that full generality is seldom really achievable, we remark that it is possible to develop a $W^{1,p}$ theory restricted to a sub-bundle, obtaining analogous results to those presented in Section \ref{ssc_subbundle}.

\section{Currents and vector fields}

\subsection{Classical currents in $\C^N$}

We use the notation of \ref{sec_class_cur}.

Let $\mathcal{D}'(\C^N)$ be the space of distributions on $\C^N$.  Given $T\in \mathfrak{M}_{p,q}(\C^N)$ and fixed multi-indeces $I$ and $J$ of length $p$ and $q$ respectively, we have that the functional
$$\Ci^\infty_c(\C^N)\ni\phi\mapsto T_{IJ}(\phi)=T(\phi dz_I\wedge d\overline{z}_J)$$
is a distribution; moreover, as $\mass(T)<+\infty$,
$$|T(\phi dz_I\wedge d\overline{z}_J)|\leq\|\phi\|_\infty \mass(T)$$
so $T_{IJ}$ is indeed a measure.

Therefore, we can write
$$T=\sum_{|I|=p,\ |J|=q} T_{IJ}dz_I\wedge d\overline{z}_J$$
meaning that, if $\omega$ is a $(p,q)-$form given by
$$\omega=\sum_{|I|=p,\ |J|=q}\omega_{IJ}dz_I\wedge d\overline{z}_J$$
then
$$T(\omega)=\sum_{|I|=p,\ |J|=q} T_{IJ}(\omega_{IJ})\;.$$

\medskip

Given a Borel set $A\subseteq\C^N$ and a measure $\mu$ we denote by $\EL^2_{p,q}(A,\mu)$  the space of currents $T\in \mathfrak{M}_{p,q}(\C^N)$ such that, for any multi-indeces $I$ and $J$, $\vert I\vert=p$, $\vert J\vert=q$, one has $T_{IJ}=f_{IJ}\mu\llcorner A$, with $f_{IJ}\in L^2(A,\mu)$. We endow $\EL^2_{p,q}(A,\mu)$ with the norm
$$\|T\|_2=\left(\sum_{|I|=p, |J|=q}\int_{A}|f_{IJ}|^2d\mu\right)^{1/2}$$
In what follows $\mu$ will be the Hausdorff measure $\H^{2n}$ on $X$ and $A$ an open set $\Omega$ in $X$, so we will write $\EL^2_{p,q}(\Omega)$ for $\EL^2_{p,q}(\Omega, \H^{2n})$.

\begin{Prp} \label{prp_cont_inc}The inclusion $i:(\EL^2_{p,q}(\Omega), \|\cdot\|_2)\to(\mathfrak{M}_{p,q}(\C^N), \mass(\cdot))$ is a bounded linear map.\end{Prp}
\noindent{\bf Proof:} The map $i$ is obviously linear. Moreover, given $T\in \EL^2_{p,q}(\Omega)$ and a smooth, compactly supported $(p,q)-$form $\omega$ on $\C^N$, we have that
$$|T(\omega)|=\left|\sum_{I,J}\int_{\Omega_\rg}\omega_{I,J}f_{I,J}d\H^{2n}\right|\leq \sum_{I,J}\|f_{I,J}\|_2\|\omega_{I,J}\|_2$$
$$\leq\H^{2n}(\Omega)^{1/2}\sum_{I,J}\|f_{I,J}\|_2\|\omega_{I,J}\|_\infty\leq\H^{2n}(\Omega)^{1/2}\|T\|_2\|\omega\|_\infty$$
so that
$$\mass(i(T))\leq\|T\|_2\H^{2n}(\Omega)^{1/2}$$
whence $\|i\|\leq\H^{2n}(\Omega)^{1/2}$. $\Box$

We set $H_{p,q}(\Omega)=i^{-1}(\mathfrak{F}_{p,q}(\Omega))$, so $(H_{p,q}(\Omega), \|\cdot\|_2)$ is a Hilbert space.

\subsection{Vector fields}

A $1-$vector field of class $W^{1,2}_{\loc}$ (resp. $L^2_{\loc}$) on $\Omega$ is a (non necessarily continuous) section $\xi:\Omega_\rg\to T\Omega_\rg$ such that for every $\phi\in\Ci^\infty_0(\Omega)$ we have $\xi(\phi)\in W^{1,2}(\Omega)$ (resp. $L^2(\Omega,\H^{2n})$).

A $k-$vector field of class $W^{1,2}_{\loc}$ (resp. $L^2_{\loc}$) on $\Omega$ is a (non necessarily continuous) section $\xi:\Omega_\rg\to \Lambda^k T\Omega_\rg$ such that
 for every $\phi\in\Ci^\infty_0(\Omega)$ the contraction of $\xi$ with $\phi$ is a $(k-1)-$vector field of class $W^{1,2}_{\loc}$ (resp. $L^2_{\loc}$) on $\Omega$.
 
\medskip
 
In terms of global coordinates of $\C^N$, the coefficients of $\xi$ are functions in $W^{1,2}(\Omega')$ for every $\Omega'\Subset\Omega$, because the coefficient of $\de^k/\de x_I$ are
$$\langle \xi, dx_{i_1}\wedge\ldots\wedge dx_{i_k}\rangle=\xi(x_{i_1},\ldots, x_{i_k})$$
and the coordinate functions can be extended to compactly supported functions outside $\Omega'$.

\medskip

The divergence of a $k-$vector field $\xi$ of class $W^{1,2}_{\loc}$ on $\Omega$ is a $(k-1)-$vector field  $\chi$ of class $L^2_{\loc}$ on $\Omega$ such that for every $\phi_1,\ldots, \phi_k\in \Ci^\infty_0(\Omega)$ we have
$$\int_{\Omega_\rg}\xi(\phi_1,\ldots,\phi_k)d\H^{2n}=(-1)^k\int_{\Omega_\rg}\phi_1\chi(\phi_2,\ldots, \phi_k)d\H^{2n}$$
and we write $\chi=\mathrm{div}\xi$. This vector field $\chi$ obviously exists locally in the regular part. Moreover, as the components of $\xi$ are in $W^{1,2}(\Omega')$ for every relatively compact open set $\Omega'\subset\Omega$, we can show that $\chi$ is of class $L^2_{\loc}$ on $\Omega$.

\begin{Prp} For every $k-$vector field $\xi$ of class $W^{1,2}_{\loc}$ on $\Omega$ and every smooth compactly supported $k-$form $\omega$ on $\C^N$ such that $\supp \omega\cap X\Subset\Omega$, we have that
$$\int_{\Omega_\rg}\langle \xi, d\omega\rangle d\H^{2n}=(-1)^k\int_{\Omega_\rg}\langle \mathrm{div}\xi, \omega\rangle d\H^{2n}$$
 \end{Prp}
 
 \noindent{\bf Proof: } Let $f_1,\ldots, f_k$ be smooth compactly supported functions on $\C^N$ such that $g_i=f_i\vert_\Omega$ and the $(k-1)-$form
 $$\omega=f_1df_2\wedge\ldots\wedge df_k;$$
 we have $d\omega=df_1\wedge df_2\wedge\ldots\wedge df_k$ and 
 $$\int_{\C^N}\langle \xi ,d\omega\rangle d\H^{2n}=\int_{\Omega_\rg}\!\!\langle \xi ,d\omega\rangle d\H^{2n}=\int_{\Omega_\rg}\!\!\xi (g_1,\ldots, g_k)d\H^{2n}$$
 $$=(-1)^k\int_{\Omega_{\rg}}\!\!g_1\mathrm{div}\xi (g_2,\ldots, g_k)d\H^{2n}=(-1)^k\int_{\Omega_\rg}\!\!\langle \mathrm{div}\xi , \omega\rangle d\H^{2n}$$
 $$=(-1)^k\int_{\C^N}\langle \mathrm{div}\xi ,\omega\rangle d\H^{2n}.$$
 By linearity, the thesis follows. $\Box$
 
 \medskip
 
 If $\xi$ is a $k-$vector field of class $W^{1,2}_{\loc}$ on $\Omega$, then $\mathrm{div}\xi$ is a $(k-1)-$vector field of class $L^2_{\loc}$ on $\Omega$. Since $\Omega_\rg$ is a complex manifold, we have a natural Dolbeault decomposition for vector fields, so, if $\xi$ is of type $(p,q)$, we can write
 $$\mathrm{div}\xi=\div\xi+\dbiv\xi$$
 where $\div\xi$ is the $(p-1,q)-$component and $\dbiv\xi$ is the $(p,q-1)-$component. Obviously, both these vector fields are still of class $L^2_{\loc}$ on $\Omega$. By linearity, $\div$ and $\dbiv$ extend to arbitrary vector fields.
 
 \medskip
 
 We will say that $\xi$ is of class $W^{1,2}$ (resp. $L^2$) on $\Omega$ if it is of class $W^{1,2}_{\loc}$ (resp. $L^2_{\loc}$) on $\Omega$ and $|\xi|,\ |\mathrm{div}\xi|\in L^2(\Omega,\H^{2n})$ (resp. $|\xi|\in L^2(\Omega,\H^{2n})$), where $|\cdot|$ is the norm described on page \pageref{pg_norma_eucl}. We will write $\|\xi\|_2$ for the $L^2-$norm of the function $|\xi|$; occasionally, we will use the notation $\|\xi\|_{2,\Omega}$ to stress that the integration is on the set $\Omega$.
 
\begin{Rem}\label{rem_caratt_W}Repeating the proof of Section \ref{ssc_alt_proof}, we can show that $\xi$ is of class $W^{1,2}_\loc$ (resp. $W^{1,2}$) if it is of class $W^{1,2}_\loc$ on $\Omega_\rg$ and $|\xi|$, $|\mathrm{div}\xi|\in L^2_\loc(\Omega,\H^{2n})$ (resp. $L^2(\Omega,\H^{2n})$).\end{Rem}
 
\medskip

We can identify $\xi$ with an element of $\EL_k(\Omega)$ or, if $\xi$ is of type $(p,q)$, with an element of $\EL_{p,q}(\Omega)$ and this identification preserves the norms. Indeed the $\EL_{p,q}(\Omega)$-norm 
$$\|\xi\|_{\EL}=\sqrt{\displaystyle\sum_I\int_{\Omega_\rg}\langle \xi, dx_{i_1}\wedge\ldots\wedge dx_{i_k}\rangle^2d\H^{2n}}$$
whereas, if we consider the $L^2-$norm of the function $x\mapsto |\xi_x|$, we obtain
$$\||\xi_x|\|_2=\sqrt{\displaystyle\int_{\Omega_\rg}|\xi_x|^2d\H^{2n}}=\sqrt{\displaystyle\int_{\Omega_\rg}\sum_I \langle \xi, dx_{i_1}\wedge\ldots\wedge dx_{i_k}\rangle^2dH^{2n}}=\|\xi\|_{\EL}$$
Thus the two norms coincide.
 
\begin{Prp}\label{prp_w_normal} If $\xi$ is a $(p,q)-$vector field of class $W^{1,2}$ on $\Omega$, then the current
$$T_\xi(\omega)=\int_{\Omega_\rg}\langle\omega,\xi\rangle d\H^{2n}\qquad \forall\ \omega\in\Di^{p+q}(\C^N)$$
belongs to $\mathfrak{N}_{p,q}(\Omega)$.\end{Prp}

\noindent{\bf Proof: } The functional $T_\xi$ is obviously linear. Moreover, for any smooth $(p,q)-$form $\omega$, we have that
$$|T_\xi(\omega)|\leq\|\omega\|_{2,\Omega}\|\xi\|_{2,\Omega}\leq\sqrt{\H^{2n}(\Omega)}\|\omega\|_\infty\|\xi\|_{2,\Omega}$$
so the current $T_\xi$ is of finite mass. And, using the previous proposition, we then obtain that
$$|T_\xi(d\omega)|=\left|\int_{\Omega_\rg}\langle d\omega,\xi\rangle d\H^{2n}\right|=\left|\int_{\Omega_\rg}\langle \omega,\mathrm{div}\xi\rangle d\H^{2n}\right|\leq\sqrt{\H^{2n}(\Omega)}\|\omega\|_\infty\|\mathrm{div}\xi\|_{2,\Omega}$$
which implies that also $dT_\xi$ is of finite mass. This proves that $T_\xi$ is a normal $(p,q)-$current.
$\Box$

\medskip

By the map $i$ defined in Proposition \ref{prp_cont_inc}, we have that the vector field $\xi$ corresponds to the current $T_\xi$.

\begin{Prp}\label{prp_equiv_corr_cv}$H_{p,q}(\Omega)$ is isomorphic to the space of $(p,q)-$vector fields of class $L^2$ on $\Omega$.\end{Prp}
\noindent{\bf Proof: }If $T\in H_{(p,q)}(\Omega)$, by definition the coefficients $T_I$ are in $L^2(\Omega,\H^{2n})$; moreover, since $T$ is a flat current, the associated vector field is tangent to its support, that is, it is tangent to $\Omega_\rg$. 

On the other hand, let $\chi$ be a $(p,q)-$vector field of class $L^2$ on $\Omega$. We can find $(p,q)-$vector fields $\chi_j$ of class $L^2$ on $\Omega$, such that $\chi_j\to\chi$ with respect to the $L^2-$norm and that $\supp\chi_j\cap\Omega_\sg=\emptyset$.

Now, $\chi_j$ being supported in a closed set of $\Omega_\rg$, we can approximate its coefficients with respect to local coordinates by using $W^{1,2}(\Omega)$ functions. So, we construct $(p,q)-$vector fields $\xi_{j,k}$ of class $W^{1,2}$ on $\Omega$ as follows.

Take an open locally finite coordinate covering $\{U_\alpha\}_\alpha$ of $\Omega_\rg$ and a partition of unity $\{\phi_\alpha\}_\alpha$ subordinated to $\{U_\alpha\}_\alpha$ and on every $U_\alpha$ consider the vector field $\chi_{j,\alpha}=\phi_\alpha\chi_j$. Then, hoosing a basis $X_{1,\alpha},\ldots, X_{2n,\alpha}$ for the trivialization of $TU_\alpha$, we have that
$$\chi_{j,\alpha}=\sum_i f^i_{j,\alpha}X_{i,\alpha}$$
with $f^i_{j,\alpha}\in L^2(U_\alpha, \H^{2n})$.
Let $k$ be a positive integer and $K\subset\Omega_\rg$ a compact set such that $\H^{2n}(\Omega_\rg\setminus K)\leq 1/k$, and take functions $g^i_{j,k,\alpha}$ in $W^{1,2}(U_\alpha)$ such that:
\begin{enumerate}
\item $\displaystyle{\sum_{U_\alpha\cap K\neq\emptyset}\sum_i\|f^i_{j,\alpha}-g^i_{j,k,\alpha}\|_2\leq 1/k}$
\item $\displaystyle{\sum_{U_\alpha\cap K=\emptyset}\sum_i\|\nabla g^i_{j,k,\alpha}\|_2\leq k}.$
\end{enumerate}
The locally finite sum 
$$\xi_{j,k}=\sum_\alpha\sum_ig^i_{j,k,\alpha}X_{i,\alpha}$$
converges to a vector field of class $W^{1,2}$ on $\Omega$. 

We obviously have that $\xi_{j,k}\to \chi_j$ with respect to the $L^2-$norm, so, letting $k,j$ go to infinity, we have an approximation of $\chi$ by vector fields $\xi_m$ of class $W^{1,2}$. This means that $T_{\xi_m}\to T_\chi$, because of the Proposition \ref{prp_cont_inc}, so $T_\chi$ is a limit of normal currents in the mass norm, that is, is flat. Hence, as $T_\chi\in\mathfrak{F}_{p,q}(\C^N)$, we have that $\chi$ corresponds to an element of $H_{p,q}(\Omega)$. $\Box$

%\begin{Rem}As $T_\xi$ is normal, $T_\xi\in \mathfrak{F}_{p,q}(\C^N)$,  so it is a metric current and, being supported in $\Omega$, belongs to $M_{p,q}(\Omega)$ and it is flat; therefore $T_\xi\in i(H_{p,q}(\Omega))$.\end{Rem}

\medskip

In view of Proposition \ref{prp_equiv_corr_cv}, we will identify the space of $(p,q)-$vector fields with $L^2$ coefficients and $H_{p,q}(\Omega)$.

In what follows we will be interested also in another class of vector fields, namely, the $k-$vector fields $\xi$ of class $L^2$ on $\Omega$ which have an $L^2$ divergence in a weak sense. We will denote by $D_{\mathrm{div}}$ the set of such vector fields; obviously, the vector fields of class $W^{1,2}$ defined above, belong to $D_{\mathrm{div}}$.

\section{Graph norm density}

The main goal of this section is to obtain some density results in the spirit of \ref{teo_approx_W} for vector fields.

\medskip

Given $\xi$ of class $W^{1,2}$ on $\Omega$, we consider the vector field $\xi_{j,m}=(1-g_{j,m})\xi$ (with the notation of Theorem \ref{teo_approx_W}). Then 
$$\|\xi-\xi_{j,m}\|_2=\|g_{j,m}\xi\|_2\leq \|\xi\|_{2,\supp g_{j,m}}\xrightarrow[j,m\to\infty]{}0$$
and, if $|\xi|\in L^\infty$,
$$\|\mathrm{div}(\xi-\xi_{j,m})\|_2=\|\mathrm{div}g_{j,m}\xi\|_2\leq\||\xi|\nabla g_{j,m}\|_2+ \|\mathrm{div}\xi\|_{2,\supp g_{j,m}}\xrightarrow[j,m\to\infty]{}0$$
So, as before, it is enough to show that the vector fields with $L^\infty$ modulus are dense among those of class $W^{1,2}$ on $\Omega$. The same computations and conclusions hold for the set $D_{\mathrm{div}}$.

\medskip

We denote by $D_{\debar}$ the domain of the operator $\dbiv:H_{p,q}(\Omega)\to H_{p,q-1}(\Omega)$.

\begin{Prp}\label{prp_dense_dom}The set $D_{\debar}$ is dense in $H_{p,q}(\Omega)$.\end{Prp}

\noindent{\bf Proof: }If $\xi$ is a $(p,q)-$vector field of class $W^{1,2}(\Omega)$, then $\dbiv \xi\in H_{(p,q)}(\Omega)$ and, because of the density of  $W^{1,2}$ in $L^2$, we have the thesis.$\Box$

\medskip

We can define the adjoint $\dbad$ of $\dbiv$ as a densely defined operator:
$$\dbad:H_{0,0}(\Omega)\to H_{0,1}(\Omega)$$
by the formula
$$(\dbiv \xi, \eta)=(\xi,\dbad \eta)$$
for $\xi\in D_{\debar}$. If $\eta$ is of class $W^{1,2}(\Omega)$, then, given a coordinate chart $U\Subset\Omega_\rg$ with local coordinates $z_1,\ldots, z_n$, we have that
$$\dbad \eta=-\sum_{j=1}^n\frac{\de \eta}{\de z_j}\de_j$$
on $U$. This holds on every coordinate chart of the regular part and these writings patch together; let $\{U_k\}$ be a collection of such charts (such that every point of $\Omega_\rg$ belongs to at most $M$ of them) and let $z_1^{(k)},\ldots, z_n^{(k)}$ be coordinates on $U_k$, then
$$\sum_k\int_{U_k}\sum_{j=1}^n\left|\frac{\de \eta}{\de z_j^{(k)}}\right|^2d\H^{2n}\leq\sum_{k}\int_{U_k}|\nabla\eta|^2d\H^{2n}\leq M\|\nabla\eta\|_2^2$$
Indeed, this expression for $\dbad \eta$ makes sense as soon as we have some global upper bound for the right-hand side.

\subsection{Weighted $L^2$ spaces}

Given $\Omega\subset X$ (not necessarily bounded) and $\phi\in\Ci^1(\Omega)$, we define the norm $\|\cdot\|_{\phi}$ by
$$\|g\|_{\phi}=\left(\int_{\Omega_\rg}g^2e^{-\phi}d\H^{2n}\right)^{1/2}$$
and we consider the spaces $L^2(\Omega,\phi)$ whose elements are the functions $f\in L^1_{\mathrm{loc}}(\Omega,\H^{2n})$ such that $\|f\|_{\phi}<+\infty$.

\begin{Prp}\label{prp_pesi_cpt} If $\Omega'\Subset \Omega$ and $\phi\in\Ci^1(\Omega)$, then $L^2(\Omega')=L^2(\Omega',\phi)$.\end{Prp}
\noindent{\bf Proof: }  Because of the compactness  there exist two positive constants $c_1, c_2$ such that the following estimate
$$c_1\leq e^{-\phi}\leq c_2$$
holds on $\Omega'$ and
 $$c_1\|\cdot\|_2\leq\|\cdot\|_\phi\leq c_2\|\cdot\|_2.$$
 
That means that the $L^2-$norm is finite if and only if the weighted $L^2-$norm is finite. Therefore the two Hilbert spaces coincide and the norms on them are equivalent. $\Box$

\medskip

We define
$$W^{1,2}_{\mathrm{loc}}(\Omega)=\{f\in L^1_{\mathrm{loc}}(\Omega,\H^{2n})\ \vert\ f\in W^{1,2}(V)\ \forall\ V\Subset\Omega\}$$
and, given $\phi,\psi\in\Ci^1(\Omega)$, we set
$$W^{1,2}(\Omega,\phi,\psi)=\{f\in W^{1,2}_{\mathrm{loc}}(\Omega)\ \vert\ f\in L^2(\Omega, e^{-\phi}\H^{2n}),\ |\nabla f|\in L^2(\Omega, e^{-\psi}\H^{2n})\}$$
Given $\varphi\in \Ci^\infty_0(\Omega)$, its support $\Omega'$ is compact, so Proposition \ref{prp_pesi_cpt} applies, therefore, for every $f\in W^{1,2}(\Omega,\phi,\psi)$ and every $Y\in R(X)$, we have that
$$\int_{\Omega_\rg}\!\!\!\!\!\!fY(\varphi)d\H^{2n}=\int_{\Omega'_\rg}\!\!\!\!\!\!fY(\varphi)d\H^{2n}=-\int_{\Omega'_\rg}\!\!\!\!\!\!\varphi\langle \nabla f, Y\rangle d\H^{2n}=-\int_{\Omega_\rg}\!\!\!\!\!\!\varphi\langle \nabla f, Y\rangle d\H^{2n}$$
In the same way, we define $H_{p,q}(\Omega,\phi)$ to be the set of $(p,q)-$vector fields of class $L^1_{\mathrm{loc}}$ on $\Omega$ whose modulus belongs to $L^2(\Omega,\phi)$. The operator
$$\dbiv:H_{p,q}(\Omega,\phi)\to H_{p,q-1}(\Omega,\psi)$$
is defined by saying that $\xi=\dbiv \eta$ if this relation holds on every compact set of $\Omega$ and $|\xi|\in L^2(\Omega,\psi)$. We will denote by $D_{\debar}$
the domain of the operator $\dbiv$. The same definitions and considerations can be carried on for the operator $\mathrm{div}$ and its domain will be denoted by $D_{\mathrm{div}}$.

Obviously, $W^{1,2}(\Omega,\phi,\psi)$ is contained both in $D_{\mathrm{div}}$ and in $D_{\debar}$, so both $\mathrm{div}$ and $\dbiv$ have dense domain in $H_{p,q}(\Omega,\phi)$.

\medskip

Once again, using the techniques of the Theorem \ref{teo_approx_W}, we can approximate every vector field $\xi\in H_{p,q}(\Omega,\phi)$ with $|\xi|\in L^\infty(\Omega)$ by elements supported away from $\Omega_\sg$. 

\medskip

Endowed with the scalar product 
$$\left(\xi,\eta\right)_{\phi}=\int_{\Omega}\langle \xi_x,\eta_x\rangle_x e^{-\phi(x)}d\H^{2n}(x)$$
where $\langle \cdot,\cdot\rangle_x$ is the euclidean scalar product on $T_x\C^N$, the space $H_{p,q}(\Omega,\phi)$ becomes a Hilbert space.

We can therefore define the adjoint  operator $\debar^*:H_{p,q-1}(\Omega,\psi)\to H_{p,q}(\Omega,\phi)$. 
We note that, if $\xi\in H_{0,q}(\Omega,\phi)$ is of class $W^{1,2}(\Omega,\phi,\psi)$, given an open coordinate set $U\Subset\Omega_\rg$ with coordinates $z_1,\ldots, z_n$, we have that, on $U$, 
$$\dbad\xi=-e^{\phi}\sum_{j_1<\ldots<j_q}\sum_{j\neq j_k}\frac{\partial e^{-\psi}\xi_{j_1,\ldots,j_q}}{\partial z_j}(-1)^{m(j)}\partial_{\overline{z}_j}\wedge\partial_{\overline{z}_{j_1}}\wedge\ldots\wedge\partial_{\overline{z}_{j_q}}$$
In fact, gluing these local expressions we obtain a vector field $\dbad\xi$ on $\Omega_\rg$ (and therefore a.e. on $\Omega$) which locally satisfies
$$(\xi, \dbiv\eta)_\psi=(\dbad\xi, \eta)_\phi$$
and has bounded norm in $L^2(\Omega,\phi)$. Therefore it represents an element of $H_{0,q+1}(\Omega,\phi)$ and we can use boundedness to show that the equation above holds globally on $\Omega$.

\medskip

Clearly, if $\xi$ has $W^{1,2}(\Omega,\phi,\psi)-$coefficients from Proposition $\ref{prp_w_normal}$ easily follows that the associated current $T_\xi$ is locally normal, . 

\medskip

Now, given three weight functions $\phi_1,\ \phi_2,\ \phi_3$, we can consider the situation
$$\xymatrix{H_{p,q+1}(\Omega,\phi_1)\ar[r]^\dbiv&H_{p,q}(\Omega,\phi_2)\ar@/^1pc/[l]^{\dbad}\ar[r]^\dbiv&H_{p,q-1}(\Omega,\phi_3)\ar@/^1pc/[l]^{\dbad}}\;.$$
In $H_{p,q}(\Omega)$, we have the set $D_{\debar}\cap D_{\debar^*}$, intersection of the domains of $\dbiv$ and $\dbad$, which contains the vector fields with smooth coefficients. We define the norm
$$D_{\debar}\cap D_{\debar^*}\ni\eta\mapsto\|\eta\|_{\phi_2}+\|\dbiv\eta\|_{\phi_1}+\|\dbad\eta\|_{\phi_3}$$
called the \emph{graph norm}. 

\subsection{Density hypothesis}\label{ssec_dens_hyp}

We say that $\Omega$ is \emph{\boh} if 
\begin{enumerate}
\item every $k-$vector field of class $W^{1,2}$ on $\Omega$ can be approximated in the $W^{1,2}-$norm by $k-$vector fields with bounded euclidean norm;
\item every $k-$vector field in $D_{\debar}\cap D_{\debar^*}$ can be approximated in the graph norm by $k-$vector fields with bounded euclidean norm.
\end{enumerate}

\begin{Prp}\label{prp_reg=boh} If $\Omega\Subset X_{\rg}$, then $\Omega$ is \boh.\end{Prp}
\noindent{\bf Proof: } Since $\Omega$ is relatively compact inside $X_{\rg}$, which is a manifold, we can find a finite number of coordinate charts $U_1,\ldots, U_k$ of $X_{\rg}$ such that $\Omega\Subset U_1\cup\ldots\cup U_k$ and construct a partition of unity $\psi_1,\ldots, \psi_k$ subordinated to that covering.

Now, given a $k-$vector field $\xi$, we can consider $\xi_j=\psi_j\xi$, which is supported in $U_j$ and then obtain, by convolution, a smooth approximation of $\xi_j$, say $\{\xi_{j,m}\}_m$.

As there is only a finite number of open sets $U_j$, we can ask that, for every $j$ and $m$, the following holds:
$$\|\xi_j-\xi_{j,m}\|_2+\|\nabla(\xi_j-\xi_{j,m})\|_2\leq\frac{1}{m}$$
so that, letting
$$\xi_m=\xi_{1,m}+\ldots+\xi_{k,m}$$
we have
$$\|\xi-\xi_m\|_2\leq\frac{k}{m}\xrightarrow[m\to\infty]{}0$$
The same procedure works for the graph norm. $\Box$

\begin{Prp}\label{prp_away-boh}If the $k-$vector fields of class $W^{1,2}$ compactly supported in $\Omega_\rg$ are dense in the $k-$vector fields of class $W^{1,2}$ on $\Omega$, then $\Omega$ is {\boh} and the same statement holds for the graph norm.\end{Prp}
\noindent{\bf Proof: }We carry on the proof for the $W^{1,2}-$norm, but it works in the same way for the graph norm. Let $\xi$ be a $k-$vector field of class $W^{1,2}$ on $\Omega$ and let $\xi_j$ be a sequence of $k-$vector fields of class $W^{1,2}$, compactly supported in $\Omega_\rg$, approximating $\xi$ in the $W^{1,2}-$norm. Then, by Proposition \ref{prp_reg=boh}, we can find a sequence $\xi_{j,m}$ of smooth (hence, bounded) $k-$vector fields approximating $\xi_j$. Finally, by a diagonal procedure, we find a sequence $\xi_{j,m(j)}$ of bounded $k-$vector fields approximating $\xi$ in the $W^{1,2}-$norm. $\Box$

\begin{Prp}\label{boh}If $\Omega$ is \boh, then every vector field of class $W^{1,2}$ (respectively, in $D_{\debar}\cap D_{\debar^*}$) can be approximated by vector fields supported away from $\Omega_\sg$.\end{Prp}
\noindent{\bf Proof: } Let $\xi$ be a bounded $k-$vector field of class $W^{1,2}$ on $\Omega$ (respectively, in $D_{\debar}\cap\ D_{\debar^*}$) and let $g_{j,m}$ be the functions defined in the proof of Theorem \ref{teo_approx_W}. Then $\|g_{j,m}\xi\|_2\to 0$ since $\H^{2n}(\supp g_{j,m})\to 0$ and
$$\|\nabla(g_{j,m}\xi)\|_2\leq\|g_{j,m}\nabla \xi\|_2+\|\xi\llcorner\nabla g_{j,m}\|_2$$
where the first term goes to zero as the previous one and the second one is bounded from above by  $\|\xi\|_\infty\|\nabla g_{j,m}\|_2$ which also goes to zero. The same computations hold for $\dbiv$ and $\dbad$.

So, the vector fields $(1-g_{j,m})\xi$ tend to $\xi$ and are supported away from $\Omega_\sg$; now, given a generic vector field $\xi$ of class $W^{1,2}$ on $\Omega$, we can approximate it with bounded vector fields and again by the diagonal procedure, we have the thesis. $\Box$

\begin{Prp}\label{prp_density}Given a collection of smooth functions $\{\rho_n\}_{n\in\N}$, subordinated to a compact exhaustion of $\Omega$,  such that
$$e^{-\phi_3}|\debar \rho_n|^2\leq e^{-\phi_2}$$
$$e^{-\phi_2}|\debar \rho_n|^2\leq e^{-\phi_1}$$
the elements compactly supported in $\Omega_\rg$ are dense in $D_{\debar}\cap D_{\debar^*}$ with respect to the graph norm.\end{Prp}
\noindent{\bf Proof: }By virtue of Proposition \ref{boh} we may assume that $\xi$ is supported away from $\Omega_\sg$. This means that also $\rho_n\xi$ is supported away from $\Omega_\sg$, therefore we can apply the computations in \cite{hormander1}. $\Box$

\medskip

Collecting together Propositions \ref{boh} and \ref{prp_density}, we have the converse of Proposition \ref{prp_away-boh}: if $\Omega$ is \boh, then the vector fields compactly supported in $\Omega_\rg$ are dense in the suitable sense.

Following \cite[Section 5.2]{hormander1}, we choose an hermitian metric on $\Omega_\rg$ with $f\in\Ci^\infty(\Omega_\rg)$, such that $|\debar\rho_n|\leq 1$ for every $n$. The norms in what follows will always be with respect to this new metric and its associated volume form $dV=fd\H^{2n}$ with $f\in\Ci^\infty(\Omega_\rg)$.

\medskip

\begin{Teo}Let $\Omega\subset X$ be a \boh\  Stein open set. Then the equation $\dbiv \xi=\eta$ has a solution $\xi$ of class $L^2(\Omega,\loc)$ for any $\eta$ of class $L^2(\Omega,\loc)$.\end{Teo}
\noindent{\bf Proof: } By \cite{narasimhan1}, as $\Omega$ is Stein, we can find a strictly plurisubharmonic exhausting function, which is real analytic.

By \cite[Theorem 5.2.3]{hormander1}, applied to the open set $\Omega_\rg$, there exists a continuous function $C$ on $\Omega_\rg$ such that
$$\int (\lambda- C)|\eta|^2e^{-\phi}dV\leq4(\|\dbad\eta\|_\phi^2+\|\dbiv\eta\|_\phi^2)$$
for every $(p,q)-$vector field $\eta$ with coefficients $\Ci^\infty_c$ on $\Omega_\rg$ and for every $\phi\in\Ci^2(\Omega_\rg)$, where $\lambda$ is the lowest eigenvalue of the Levi form of $\phi$.

If $\phi$ is a strictly plurisubharmonic exhausting function for $\Omega$, we can modify it by composing with a convex, monotonically increasing function so that $(\lambda-C)>4$; then we have
$$\|\eta\|_\phi^2\leq \|\dbad\eta\|_\phi^2+\|\dbiv \eta\|_\phi^2\;.$$

By Proposition \ref{prp_density}, this estimate holds for any $\eta$ in $D_{\debar}\cap D_{\debar^*}$, therefore, by a standard result in functional analysis (e.g. see \cite[Lemma 4.1.1]{hormander1}), we have that the equation $\dbiv \xi=\eta$ has a solution of class $L^2(\Omega,\loc)$ for every $\dbiv-$closed $\eta$ of class $L^2(\Omega,\loc)$. $\Box$

\medskip

Let $\theta$ be the adjoint operator of $\dbiv$ between $L^2$ and $L^2$; $\theta$ and $\dbad$ differ only for multiplication operators: $\dbad=M'\theta M=M_1\theta+M_0$. Therefore we have density results also for the norm $\|\cdot\|+\|\dbiv\cdot\|+\|\theta\cdot\|$.

\begin{Prp}If $\eta$, $\dbiv\eta$, $\theta\eta$ are compactly supported and of class $L^2(\Omega)$, then $\eta$ is of class $W^{1,2}(\Omega)$.
\end{Prp}
\noindent{\bf Proof: } We use a variation of the proof of Lemma 4.2.3 in \cite{hormander1}. If $\eta$ is smooth and compactly supported in $\Omega_\rg$, we can suppose without loss of generality that its support is contained in a coordinate chart and we have the following estimate
$$\sum_{I,J}\sum_{j=1}^n\int\left|\frac{\de \eta_{I,J}}{\de\bar{z}_j}\right|^2d\H^{2n}\leq 2\|\theta\eta\|^2+\|\dbiv\eta\|^2\;,$$
with respect to some local coordinates (cfr \cite[Formula 4.2.7]{hormander1}). By density, we extend this estimate to every $\eta$ with compact support  in $\Omega$, belonging to $D_{\debar}\cap D_{\theta}$. Therefore we conclude that $\de \eta_{I,J}/\de\bar{z}_j\in L^2(\Omega)$ for every $I,J$ and every $j$.

Modifying the proof of Lemma 4.2.4 in \cite{hormander1} with the same argument we just used, we conclude that also $\de \eta_{I,J}/\de z_j\in L^2(\Omega)$ for every $I,J$ and every $j$. The thesis follows. $\Box$

\medskip

Employing the previous Lemma, we can repeat \emph{verbatim} the proof of \cite[Theorem 4.2.5]{hormander1}, obtaining the following result.

\begin{Teo}\label{teo_reg_w}Let $\Omega$ be a Stein \boh\ open set and let $0\leq s\leq\infty$. Then the equation $\debar\xi=\eta$ has a solution $\xi$ of class $W^{s+1,2}(\Omega,\loc)$ for every $\xi$ of class $W^{s,2}(\Omega,\loc)$ such that $\dbiv\eta=0$. Every $L^2(\Omega,\loc)$ solution of the equation has this property if $\xi$ is of bidimension $(p,n)$.\end{Teo}

We note that on a singular space the Sobolev embeddings do not hold, so we cannot conclude that there is a smooth solution when the datum is smooth; the best we can say is that the solution will be smooth on $\Omega_\rg$ and with every derivative in $L^2_\loc$ (namely, locally square integrable around singular points).

\medskip

\paragraph{Example} Let us consider the (reducible) singular space $X=\{zw=0\}$ in $\C^2$. The $(0,1)-$vector-fields
$$\xi_1=\sum_{j=1}^\infty \frac{t^n}{2^{n+1}}\frac{\de}{\de\bar{t}}\qquad \xi_2=\sum_{j=1}^\infty\frac{t^n}{3^{n+1}}\frac{\de}{\de\bar{t}}$$
are holomorphic on the unit punctured disc $\DD^*$, that is, solutions to the equations $\dbiv\xi=0$ on $\DD^*$. Therefore, their pushforwards through the maps 
$$f_1:t\to(t,0)\qquad f_2:t\to (0,t)$$
defined from $\DD^*$ to $X_\rg$, are holomorphic vector-fields on an open set of $X_\rg$. Their sum $\xi$, extended to $0$ at the singular point, is a bounded (hence $L^2$) solution of the equation $\debar\xi=0$ on an open set of $X$, but the vector-field $\xi$ isn't of class $\Ci^\infty$ around the singular point. Indeed, it is not even of class $\Ci^0$: let $\omega=\rho\cdot(d\bar{z}+d\bar{w})$, with $\rho\in\Ci^\infty_c(\C^2)$ such that $\rho\equiv1$ in a neighborhood of the origin, then, near $(0,0)$
$$\langle \xi, \omega\rangle=\left\{\begin{array}{lcl}(2-z)^{-1}&\textrm{ if }&w=0,\ z\neq0\\
(3-w)^{-1}&\textrm{ if }&z=0,\ w\neq0\end{array}\right.$$
which doesn't extend continuously to $(0,0)$.

\medskip

\paragraph{Example} Let $X=\{z^2=w^3\}$ be the cusp in $\C^2$; we consider the parametrization
$$f:t\mapsto (t^3,t^2)$$
and denote by $\xi=f_*\eta$ the pushforward through $f_*$ of the $(0,1)-$vector-field
$$\eta=\frac{1}{t}\frac{\de}{\de\bar{t}}\;.$$
For any given $\phi\in\Ci^\infty_c(\C^2)$, we have that
$$f^*(\xi(\phi))=f^*f_*\eta(f^*\phi)=\eta(f^*\phi)=\eta(\phi\circ f)$$
and
$$\int_{X}|\xi(\phi)|^2d\H^2=-\frac{i}{2}\int_{\C}|\eta(\phi(t^3,t^2))|^2(4|t|^2+9|t|^4)dt\wedge d\bar{t}=$$
$$-\frac{i}{2}\int_{\C}\frac{1}{|t|^2}\frac{\de\phi(t^3, t^2)}{\de\bar{t}}|t|^2(4+9|t|^2)dt\wedge d\bar{t}<+\infty$$
which means that $\xi$ is of class $L^2_\loc$. Moreover it is clear that $\dbiv \xi$ is zero outside $X_\sg$, which means that if an open set of $X$ containing the origin is \boh, $\dbiv \xi=0$ on $X$.

But it is a simple matter of computation to show that  $\dbiv\eta=\delta_0$ on $\C$; therefore
$$\dbiv f_*\eta=f_*\dbiv\eta=\delta_{(0,0)}\;.$$
Therefore, no open set of $X$ containing the origin is \boh.

\chapter{Some applications and examples}
\epigraphhead[55]{\epigraph{If you must hold yourself up to your
children as an object lesson (which is not at all necessary),
hold yourself up as a warning and not as an example.}{G. B. Shaw - \emph{A Treatise on Parents and Children}}}

In this chapter, we collect some examples of complex spaces, such as locally reducible spaces, complex curves, hypersurfaces in $\C^n$, where we can apply the theory of metric currents to the Cauchy-Riemann equation.

The results themselves are not at all striking, but they are useful to understand some characteristics of metric currents. In particular, the dependence on the geometric structure of the space, the ``growth conditions'' which are required by the mass and continuity properties, the  ``flatness'' of metric currents (i.e. the fact that a metric current supported in some closed set is the pushforward through inclusion of a current on that set).

Moreover, in the case of hypersurfaces, we only show that if we are able to solve the (classical) $\debar-$equation in $L^p$ on the complement of the hypersurface (for $p$ small enough), then we can solve the $\debar-$equation in the metric currents on the hypersurface. We also present a possible approach to the equation in $L^p$, showing its equivalence to an estimate on compactly supported forms.

\section{Locally reducible spaces}\label{sec_normal_cr}

Some of the material exposed here appears also in \cite{mongodi1}.

\subsection{Structure theorem}

Let $X=L_1\cup\ldots\cup L_m$ be the union of linear subspaces $L_i$ of $\C^n$, with dimension $k_i$, such that $L_i\not\subseteq L_j$ whenever $i\neq j$. Obviously, $X$ is an analytic subset of $\C^n$.

Let $X_1$ be the singular set of $X$ and  let us suppose that we have indexed the subspaces in such a way to have $\dim L_i\geq\dim L_{i+1}$ for every $i$; we also suppose that $X$ isn't contained in any proper subspace of $\C^n$. 

Now, we consider a set $\mathcal{B}=\{L_1\}\cup\{L_{i_1},\ldots, L_{i_k}\}$ such that
$$\bigoplus_{L\in\mathcal{B}}L=\C^n\qquad L\not\subseteq\bigoplus_{\stackrel{L'\in\mathcal{B}}{L\neq L'}}L'\quad \forall\ L\in\mathcal{B}$$
We have the projections 
$$\pi_1:X\to L_1\qquad\mathrm{and}\qquad \pi_{i_h}:X\to L_{i_h}$$
and the inclusions $j_1$ and $j_{i_h}$. 

Let
$$S=\bigcup_{L\in\mathcal{B}}L\qquad S'=\bigcup_{L\not\in\mathcal{B}}L\;.$$

Given $T\in \D_m(X)$, we consider the currents
$$T_1=(j_1\circ \pi_1)_\sharp T$$
$$T_{i_h}=(j_{i_h}\circ\pi_{i_h})_\sharp T$$
and the difference
$$R=T-T_1-T_{i_1}-\ldots-T_{i_k}$$

Then for $(f,\xi)\in\D^m(X)$ such that $\supp(f)\subseteq (L_1\cup L_{i_1}\cup\ldots\cup L_{i_k})\setminus X_1$ we can find $f_L\in \Lip_b(L)$ for $L\in\mathcal{B}$ such that 
$$f=\sum_{L\in\mathcal{B}}f_L\;.$$

Therefore we have
$$T(f,\xi)=\sum_{L\in\mathcal{B}} T(f_L,\xi)=\sum_{L\in\mathcal{B}} T\llcorner\chi_{L\setminus X_1}(f_L,\xi)$$ $$=T_1\llcorner\chi_{L_1\setminus X_1}(f_{L_1},\xi)+\sum_{h=1}^k T_{i_h}\llcorner\chi_{L_{i_h}\setminus X_1}(f_{L_{i_h}},\xi)=T_1(f_{L_1},\xi)+\sum_{h=1}^k T_{i_h}(f_{L_{i_h}},\xi)$$
$$=T_1(f,\xi)+\sum_{h=1}^k T_{i_h}(f,xi)$$
and
$$\supp(R)\subseteq \bigcup_{L\not\in\mathcal{B}}L\cup X_1=X_R$$
Now, as $\dim L_1\geq \dim L_i$ for every $i$, the maximum of the dimensions of the irreducible components of $X_R$ is less or equal to $\dim L_1$ so we can repeat our argument on $X_R$ to obtain a decomposition of $R$. Eventually, the remainder will have support contained in $X_1$, whose irreducible components have dimension strictly less than $\dim L_1$.

Thus, we obtain a decomposition of $T=T_1+\ldots+ T_N$ with the following properties: $T_i$ is the pushforward through an holomorphic map $h_i$ of a current $\widetilde{T}_i$ on some $\C^{k_i}$.

\medskip

Let us suppose that $T$ is a $\debar-$closed current of bidimension $(p,q)$ on $X$ and let us write $T=\sum (h_j)_\sharp Z_j$; since all the maps involved are holomorphic, $Z_j$ is of bidimension $(p,q)$ and $\debar Z_j=0$ for every $j$. 
Then, for every $j$, we can solve the equation $\debar V_j=Z_j$ on $\C^{k_j}$ with a metric current $V_j$ (e.g., by convolution with the Cauchy kernel) and then the current $U=\sum(h_j)_\sharp V_j$ satisfies
$$\debar U=\sum_{j=1}^N\debar (h_j)_\sharp V_j=\sum_{j=1}^n(h_j)_\sharp \debar V_j=\sum_{j=1}^N(h_j)_\sharp Z_j=T\;.$$

\bigskip
This allows us to prove the following results

\begin{Prp}\label{prp_loc_struct}Let $X$ be a complex space which is locally biholomorphic to a union of linear subspaces of $\C^N$, e.g. if $X$ can locally be realized as a normal crossings divisor, $\Omega\subseteq X$ an open set and $T\in \D_m(\Omega)$. For every $x\in\Omega$, we can find an open set $\omega\ni x$, $\omega\Subset\Omega$, holomorphic maps $h_i:V_i\to \Omega$, for $i=1,\ldots, k$, with $V_i\subset\C^{n_i}$ and metric currents $T_i\in\D_m(V_i)$ such that the current
$$T-\sum_{i=1}^k(h_i)_\sharp T_i$$
has support disjoint from $\overline{\omega}$.\end{Prp}
\noindent{\bf Proof: } It is enough to choose an open set $U$ containing $x$, such that there exists a biholomorphism of $U$ with a union of linear subspaces of $\C^N$, then we can apply the decomposition shown above for the current $T\llcorner\sigma$, with $\sigma$ Lipschitz and supported in $U$. We now set $\omega=\mathrm{Int}(\supp\sigma)$. $\Box$

\medskip

\begin{Teo}\label{teo_loc_deb} Let $X$ be as in Proposition \ref{prp_loc_struct}. Then given $\Omega\subseteq X$ open, $T\in \D_m(\Omega)$ such that $\debar T=0$ and $x\in\Omega$, we can find an open set $\omega$, $\omega\ni x$, $\omega\Subset\Omega$ and a current $S\in\D_{m+1}(U)$ such that $T-\debar S$ has support disjoint from $\overline{\omega}$.\end{Teo}
\noindent{\bf Proof: }By Proposition \ref{prp_loc_struct} (and with the same notation), we find $\omega$ such that $T-\sum (h_i)_\sharp T_i$ has support disjoint from its closure; moreover, we know that $\debar T$ has support disjoint from $\omega$ if and only if $\debar (h_i)_\sharp T_i$ does. So, we can solve $\debar S_i=T_i$ in $V_i$, by convolution, and we set
$$S=\sum_{i=1}^k (h_i)_\sharp S_i\;.$$
By Proposition \ref{prp_pushf_debar}, we have that
$$\debar S(f,\pi)=\sum_{i=1}^k\debar((h_i)_\sharp S_i))(f,\pi)=\sum_{i=1}^k(h_i)_\sharp(\debar S_i)(f,\pi)=\sum_{i=1}^k(h_i)_\sharp(T_i)(f,\pi)=T(f,\pi)$$
for every $(f,\pi)\in\D^m(U)$ with $\supp f\subset\omega$. So $T-\debar S$ has support disjoint from $\overline{\omega}$. $\Box$

\subsection{Holomorphic currents}

We now investigate the solutions of $\debar T=0$ when $T$ is a $(p,n)-$current on a complex space $X$, with $\dim_\C X_\rg=n$, which is locally biholomorphic to a union of linear subspaces of some $\C^N$. In the smooth case, these currents correspond to $(n-p,0)-$forms with holomorphic coefficients.

\begin{Prp}\label{prp_loc_hol_curr}Let $X$ be as in Proposition \ref{prp_loc_struct}, $\Omega\subseteq X$ an open set. Then, given a $\debar-$closed current $T\in \D_{p,n}(\Omega)$ and a point $x\in \Omega$, we can find an open set $\omega$, $\omega\ni x$, $\omega\Subset\Omega$, 
a finite number of holomorphic maps $h_i:V_i\to\Omega$, with $V_i\subseteq\C^{n}$ open sets, and holomorphic $(n-p)-$forms $f_i\in\Omega^p(V_i)$ such that 
$$T-\sum_{i=1}^k(h_i)_\sharp[f_i]$$
has support disjoint from $\overline{\omega}$.\end{Prp}
\noindent{\bf Proof: } By Proposition \ref{prp_loc_struct} (and with the same notation), we find $\omega$ such that $T-\sum (h_i)_\sharp T_i$ has support disjoint from its closure; moreover, we know that $\debar T$ has support disjoint from $\omega$ if and only if each $(h_i)_\sharp T_i$ does. As $\dim_\C X_\rg=n$, every $V_i$ is an open set in $\C^n$; we have that $\debar T_i$ is zero on $h_i^{-1}(\omega)$ so, as $T_i$ is of bidimension $(p,n)$, we can find a $(n-p,0)-$form $f_i$ with holomorphic coefficients such that $T_i=[f_i]$ in $h_i^{-1}(\omega)$.

By pushforward, we have the thesis. $\Box$

\medskip

We remark that the currents described in the previous Proposition are locally flat in any local affine embedding of $X$, so the sheaf $K_p$ defined in \ref{ssc_dolb_cpx} coincides with the kernel of $\debar:\D_{p,n}\to\D_{p,n-1}$.

\begin{Prp}Let $X$ be as in Proposition \ref{prp_loc_struct} and $\pi:X^\nu\to X$ is the normalization of $X$. Then $K_p=\pi_*\Omega^{n-p}_{X^\nu}$ where $\Omega^{n-p}_{X^\nu}$ is the sheaf of holomorphic $(n-p)-$forms.\end{Prp}
\noindent{\bf Proof: } Given $\Omega\subset X$ which is biholomorphic to a union of open neighborhoods of $0$ in $\C^n$, the preimage $\pi^{-1}(\Omega)\subset X^\nu$ is the disjoint union of these neighborhoods and $X^\nu$ is smooth. By Proposition \ref{prp_loc_hol_curr}, we know that $\debar-$closed $(p,n)-$currents on $\Omega$ are the holomorphic currents on $\pi^{-1}(\Omega)$. So the thesis follows. $\Box$

\medskip

We remark that $\Omega\subseteq X$ is Stein if and only if $\pi^{-1}(\Omega)$ is Stein, so
$$\frac{\ker\{\debar:\mathscr{F}_{n-p,n-q}(\Omega)\to\mathscr{F}_{n-p,n-q-1}(\Omega)\}}{\mathrm{img}\{\debar:\mathscr{F}_{n-p,n-q+1}(\Omega)\to\mathscr{F}_{n-p,n-q}(\Omega)\}}=H^q(\Omega, K_p)=H^q(\pi^{-1}(\Omega), \Omega^{p}_{X^\nu})=0\;.$$

\medskip

We are now ready to give a global version of Proposition \ref{teo_loc_deb}. 
\begin{Teo} Let $X$ be a Stein space with completely reducible singularities, $T\in \D_{n,n-1}(X)$ a $\debar-$closed metric current; then there exists a metric current $S$ of bidimension $(n,n)$ such that $\debar S=T$.\end{Teo}
\noindent{\bf Proof: } We consider an open covering $\{V_i\}_{i\in\N}$ of $X$ such that we can solve the Cauchy-Riemann equation on every $V_i$; we obtain a collection $\{S_i\}_{i\in\N}$ of metric currents with $S_i\in \D_{n,n}(\overline{V_i})$ and $\debar S_i=T$ in $\D_{n,n-1}(V_i)$.

On the sets $V_{ij}=V_i\cap V_j$, we have that $R_{ij}=S_i-S_j$ is a $\debar-$closed $(n,n)$ metric current in $M_{n,n}(\overline{V_{ij}})$; lifting the covering $\{V_i\}$ to a covering $\{\Omega_j\}$ of the normalization $Y$ of $X$, we can find holomorphic functions $f_{ij}\in\Ol(\overline{\Omega_{ij}})$ such that $R_{\nu(i)\nu(j)}=\pi_\sharp[\Omega_{ij}]_\llcorner f_{ij}$.

We now recall that the normalization of a Stein space is Stein, so $H^1(Y,\Ol)=0$, therefore there exist functions $f_i\in\Ol(\Omega_i)$ such that $f_{ij}=f_i-f_j$. Defining $R_i=\pi_\sharp[\Omega_i]\llcorner f_i\in \D_{n,n}(\overline{V_i})$, on $\overline{V_{ij}}$ we have 
$$R_{ij}=R_i-R_j$$
so
$$S_i-R_i=S_j-R_j\;.$$
Thus we can define a metric current $S$ such that $\debar S=T$.$\Box$

\subsection{Density properties}

Let $X$ be a complex space. Suppose that there exist a complex manifold $Y$ and a finite morphism $\pi:Y\to X$ which is a biholomorphism on the regular part of $X$ and such that for every $y\in Y$ $D\pi_y$ is invertible. For instance, this is the case if $X$ has only normal crossings, the normalization $X^\nu$ is smooth and the canonical map $\pi:X^\nu\to X$ is locally invertible. Then

\begin{Prp}\label{prp_norm_boh}Every open set $\Omega\subseteq X$ is \boh.\end{Prp}
\noindent{\bf Proof: } In what follows, $\H^{2n}$ can be substituted by any measure of the form $e^{-\phi}\H^{2n}$ with $\phi\in\Ci^1(\Omega)$.

Let $\Omega'=\pi^{-1}(\Omega)$ and let $\{U_i\}_{i\in\N}$ be a collection of open sets of $Y$ such that
$$\Omega'\subseteq\bigcup_iU_i\;,$$ 
$\pi\vert_{U_i}:U_i\to\pi(U_i)\textrm{ is invertible}$ and every $U_i$ is biholomorphic to a bounded open set $V_i\subset\C^n$.

We denote by $\phi_i$ the inverse of $\pi\vert_{U_i}$ on $\pi(U_i)$. Let $E=\pi^{-1}(X_\sg)$; since $\pi$ is locally invertible, $\H^{2n}(E)=\H^{2n-1}(E)=0$ and $\H^{2n-2}(E\cap K)<+\infty$ for every $K\subset Y$ compact. We denote by $g_i$ the biholomorphism between $U_i$ and $V_i$.

Given a vector field $\xi$ on $\Omega$ of class $L^2$, we can define almost everywhere the pushforward $(\pi^{-1})_*\xi=\xi'$, by setting 
$$\xi'\vert_{U_i\setminus E}=(\phi_i)_*(\xi\vert_{\pi(U_i)\setminus X_\sg})$$
Let $\mu=\pi^*\H^{2n}$ be the pullback of the Hausdorff measure. Then on every set $U_i$ there exists a bounded non-vanishing smooth function $h_i$ such that $\mu=h_i \cdot g_i^*\mathcal{L}$, where $\mathcal{L}$ is the standard Lebesgue measure on $\C^n$.

Now, we construct a partition of unity $\{\chi_j\}$ subordinated to the covering $\{U_j\}$ and we consider the vector fields $\xi'_j=\chi_j\xi'$. As $\xi'_j$ is supported in $U_j$, we can regularize it by a convolution, obtaining a famili $\xi'_{j,\epsilon}$.

If $\{U_j\}$ is a finite family, we set 
$$\xi'_\epsilon=\sum_j\xi'_{j,\epsilon}\;;$$
if not, we set
$$\xi'_{1/n}=\sum_{j=1}^n\xi'_{j,1/n}\;.$$
In both cases, $\{\xi'_{1/n}\}$ is a sequence of bounded (smooth) vector fields which converge to $\xi'$ in the $L^2$ norm with respect to the measure $\mu$. Obviously, defining $\xi_n=\pi_*\xi'_{1/n}$, we obtain a sequence of bounded vector fields (because $\pi$ is locally invertible) which converges to $\xi$ in $L^2(\Omega,\H^{2n})$.

\medskip

Let us suppose that $\xi\in D_{\mathrm{div}}$, then on every $U_i$ the divergence of $\xi'_j$ (in principle, considered in a distributional sense) has finite $L^2$ norm with respect to the measure $\mu$; therefore, by the standard properties of convolution, $\|\mathrm{div}\xi'_j-\mathrm{div}\xi'_{j,1/n}\|_2\to 0$ as $n\to\infty$, so also $\|\mathrm{div}\xi-\mathrm{div}\xi_n\|_2\to0$.

In the same way, let us take $\xi\in D_{\debar}\cap D_{\debar^*}$; then we have that on every $U_i$, the vector fields $\dbiv \xi'_j$ and $\dbad\xi'_j$ have finite $L^2$ norms with respect to the appropriate pullback measures. Again, by the standard properties of convolution, we have
$$\|\dbiv\xi'_{j,1/n}-\dbiv\xi'_j\|+\|\dbad\xi'_{j,1/n}-\dbad\xi'_j\|\to0\;,\quad n\to\infty\;$$
whence $\|\dbad\xi_n-\dbad\xi\|_{\phi_3}+\|\dbiv\xi_n-\dbiv\xi\|_{\phi_1}\to0$ as $n\to\infty$. $\Box$

\bigskip

Therefore, by the conclusions of section \ref{ssec_dens_hyp}, we have that every vector field on $\Omega$ in $D_{\debar}\cap D_{\debar^*}$ can be approximated by vector fields compactly supported in the regular part $\Omega_\rg$. 

From now on, we assume that $\Omega$ is pseudoconvex.

For a vector field, which is compactly supported in $\Omega_\rg$, we can employ the techniques of \cite{hormander1} to show that
\begin{equation}\label{eq_apriori}\|\xi\|_{\phi_2}\leq C(\|\dbiv\xi\|_{\phi_1}+\|\dbad\xi\|_{\phi_3})\;.\end{equation}
By density, this estimate holds for every element of $D_{\debar}\cap D_{\debar^*}$, so, by \cite{hormander1}, we have that for every $\xi\in\ker\dbiv$ we can find $\eta\in H_{p,q+1}(\Omega,\phi_3)$ such that $\dbiv\eta=\xi$.

\begin{Teo}Given a $(p,q)-$vector field $\xi$ on $\Omega$ of class $L^1_{\loc}$ such that $\dbiv\xi=0$ in a distributional sense, then we can find $\eta$, a $(p,q+1)-$vector field on $\Omega$ of class $L^1_{\loc}$ such that $\dbiv \eta=\xi$.\end{Teo}
\noindent{\bf Proof: } We can find $\phi_1,\ \phi_2,\ \phi_3$ such that $\xi\in H_{p,q}(\Omega,\phi_2)$ and such that equation (\ref{eq_apriori}) holds. Therefore, we have a solution to the Cauchy-Riemann equation for $\xi$ in $H_{p,q+1}(\Omega,\phi_3)$. This solution $\eta$ is obviously at least of class $L^2_{\loc}$. $\Box$

\section{Complex curves}\label{co_cu}

Given a complex curve $X\subset\C^m$, let $\pi:\widetilde{X}\to\C^m$ be the normalization of $X$. For any $p\in X_\sg$ and any $q\in \pi^{-1}(p)$ we can find a local holomorphic parametrization containing $q$ and a holomorphic change of coordinates in $\C^m$ so that, locally, $\pi$ is given by holomorphic functions $(P_1(t),\ldots, P_m(t))$, $P_j:\mathbb{D}\to\C^m$, such that
$$\begin{array}{ll}P_i(t)=0&\textrm{ for }1\leq i\leq k-1\\
P_k(t)=t^{N_k}& \\
P_i(t)=\sum_{j\geq N_i}a_{ij}t^j&k+1\leq i\leq m\end{array}$$
where $1\leq N_k<N_{k+1}<\cdots< N_m$. We define $N(q)=N_k$ and 
$$N(p)=\sum_{q\in\pi^{-1}(p)} N(q)\;.$$
If $\omega$ is the standard hermitian metric on $\C^m$, locally $\pi^*\omega$ is given by
$$\frac{i}{2}h(t)dt\wedge d\bar{t}$$
with
$$h(t)=N(q)^2|t|^{2N(q)-2}+O(|t|^{2N(q)})\;.$$
We define $E=\pi^{-1}(X_\sg)$.

\subsection{Holomorphic currents on complex curves}\label{ssc_hol_cc}

We study in greater detail the structure of holomorphic currents on complex curves. As it was already observed the existence of a parametrization allows us to explicit the conditions (\ref{eq_int_omega1}) and (\ref{eq_int_omega2}) and determine in terms of multiplicities of the singular points which poles are permitted for the coefficients of a holomorphic current.

\medskip

\begin{Prp}\label{prp_laurent}Let $\omega$ be a $(1,1)-$form compactly supported in $\mathbb{D}^*$. We define
$$T(\omega)=\int_{\mathbb{D}^*}t^{-j}\omega$$
Then $(\pi_q)_\sharp T$ extends to a metric current on $\pi_q(\mathbb{D})$ if and only if $j\leq2N(q)-1$. Moreover $(\pi_q)_\sharp T$ extends to a holomorphic current on $\pi_q(\mathbb{D})$ if and only if $j\leq N(q)-1$.\end{Prp}
\noindent{\bf Proof: } Given a smooth compactly supported $(1,1)-$form $\xi$ on $\C^m$, we have 
$$T(\pi_q^*\xi)=\int_{\mathbb{D}^*}\frac{i}{2t^j}f(t)h(t)dt\wedge d\bar{t}$$
where $f(t)$ a $\Ci^\infty_c$ function on $\mathbb{D}$. We know that $t^{-j}h(t)\in L^1(\mathbb{D})$ if and only if $|t|^{-j}h(t)|t|$ is bounded, that is, if and only if $-j+2N(q)-2+1\geq0\Leftrightarrow j\leq 2N(q)-1$. So $(\pi_q)_*T$ fulfills (\ref{eq_int_omega1}) if and only if $j\leq 2N(q)-1$.

On the other hand, given $\eta$ a $(1,0)-$form, smooth and compactly supported in $\C^m$, we have that $\eta=\sum a_idz_i$ with $a_i\in \Ci^\infty_c(\C^m)$. By the assumptions on $\pi_q$ made at the beginning of this section, we have that $\pi_q^*dz_i=0$ for $1\leq i\leq k-1$, $\pi_q^*dz_k=N(q)t^{N(q)-1}dt$ and $\pi_q^*dz_i=O(t^{N(q)})dt$ for $k+1\leq i\leq m$. Therefore
$$\frac{1}{t^j}\pi_q^*\eta=N(q)a_k\circ\pi_qt^{N(q)-1-j}dt +O(t^{N(q)-j})dt$$
so, if $j\leq N(q)-1$, $\alpha=N(q)-1-j\geq0$ and 
$$\lim_{\epsilon\to0}\left|\int_{|t|=\epsilon}N(q)a_k\circ\pi_q t^\alpha dt\right|\leq\lim_{\epsilon\to0}2\pi\epsilon\max_{|t|\leq\epsilon}|N(q)a_k\circ\pi_q t^\alpha|=0;$$
similarly
$$\lim_{\epsilon\to0}\int_{|t|=\epsilon}O(t^{N(q)-j})dt=0\;.$$
If $N(q)=j$, then 
$$\lim_{\epsilon\to0}\int_{|t|=\epsilon}O(t^{N(q)-j})dt=0$$
but
$$\lim_{\epsilon\to0}\int_{|t|=\epsilon}\frac{1}{t}N(q)a_k\circ\pi_q dt=i\pi N(q)a_k\circ\pi_q(0)\;.$$
The same  happens for $j\geq N(q)+1$, with one of the terms included in $O(t^{N(q)-j})$.
In conclusion, (\ref{eq_int_omega2}) is satisfied if and only if $j\leq N(q)-1$. $\Box$

\medskip

As we already observed, every $\debar-$closed current in $\Di_{r,1}(X)$ can be restricted to a holomorphic current on the regular part of $X$; in particular, if $r=1$, then $T$ coincides with a holomorphic function on the regular part and, by the Proposition \ref{prp_laurent}, this holomorphic function has to satisfy some growth conditions on the singular points.

\begin{Teo} $T\in\Di_{1,1}(X)$ is $\debar-$closed if and only if $T=\pi_\sharp[\sigma]$ with $\sigma\in H^0(\widetilde{X}, \Ol(|E|-E))$.\end{Teo}
\noindent{\bf Proof: } $\sigma\in H^0(\widetilde{X}, \Ol(|E|-E))$ if and only if, around every point $q\in E$, $\sigma$ can be written as $t^{-N(q)+1}g(t)$ where $t$ is a local coordinate and $g(t)$ is holomorphic.

By Proposition \ref{prp_laurent}, the integration against $\sigma$ induces a holomorphic $(1,1)-$current on $X$ by pushforward if and only if, writing the Laurent series of $\sigma$ in a neighbourhood of $q\in E$, the negative powers of the local coordinate have exponents less (in absolute value) or equal to $N(q)-1$.

By the previous considerations, every $\debar-$closed $(1,1)-$current on $X$ can be represented, on the regular part, as integration against some holomorphic function $h_1\in\Ol(X_\rg)$; the pullback $\pi^*h_1$ defines, by integration, a $(1,1)-$current on $\widetilde{X}\setminus E$ whose pushforward can be extended. Therefore $\pi^*h_1$ is a section of $\Ol_{\widetilde{X}}(|E|-E)$. $\Box$

\medskip

The same calculations can be carried on for $(0,1)-$currents.

\begin{Prp}\label{prp_laurent1}Let $\omega$ be a $(0,1)-$form compactly supported in $\mathbb{D}^*$ and
$$T(\omega)=\int_{\mathbb{D}^*}t^{-j}dt\wedge\omega$$
Then $(\pi_q)_\sharp T$ extends to a metric current on $\pi_q(\mathbb{D})$ if and only if $j\leq N(q)$. Moreover $(\pi_q)_\sharp T$ extends to a holomorphic current on $\pi_q(\mathbb{D})$ if and only if $j\leq 0$.\end{Prp}
\noindent{\bf Proof: } Again, we impose the conditions given by (\ref{eq_int_omega1}) and (\ref{eq_int_omega2}).

Let $\eta$ be a $(0,1)-$form in $\C^m$, with $\eta=\sum\eta_id\bar{z}_i$, $\eta_i\in\Ci^\infty_c(\C^m)$; then
$$\pi^*_q\eta=N(q)\eta_k\circ\pi_q\bar{t}^{N(q)-1}d\bar{t}+O(\bar{t}^{N(q)}d\bar{t}$$
so $t^{-j}\pi^*_q\eta\in L^1(\DD^*)$ if and only if $N(q)-1-j\geq -1$, that is $j\leq N(q)$.

Let $f\in\Ci^\infty_c(\C^m)$; we observe that, for $j\leq1$,
$$\int_{|z|=\epsilon}O(|t|^{-j})dt\xrightarrow[\epsilon\to0]{}0\;.$$
Moreover, if $j\leq0$, we have
$$\left|\int_{|z|=\epsilon}t^{-j}f\circ\pi_q dt\right|\leq2\pi\epsilon\max_{|t|\leq\epsilon}|t^{-j}f\circ\pi_q|\xrightarrow[\epsilon\to0]{}0\;;$$
if $j=1$
$$\lim_{\epsilon\to0}\int_{|z|=\epsilon}t^{-1}f\circ\pi_q dt=i\pi f(\pi_q(0))\;$$
and if $j\geq1$, the integral of $O(|t|^{-j})$ does not converge. $\Box$

\medskip

We obtain an analogue of the previous theorem.

\begin{Teo} $T\in\Di_{0,1}(X)$ is $\debar-$closed if and only if $T=\pi_\sharp[\alpha]$ with $\alpha\in H^0(\widetilde{X}, \Omega^1)$.\end{Teo}

\begin{Rem} In Propositions \ref{prp_laurent} and \ref{prp_laurent1}, there is a value of $j$ for which the current $T$ is metric, not holomorphic, but $\debar T$  is finite. In both cases, $\debar T$ is a Dirac $\delta$ in $\pi_q(q)$; so, if $T$ is of bidegree $(0,1)$, $\debar T$ is again a metric current, while in the other case $\debar T$ isn't metric, because its support is a discrete set, but its dimension is greater than $0$.\end{Rem}

\subsection{Currents with $L^p$ coefficients}

In this section we are going to solve the Cauchy-Riemann equation for metric currents with $L^p$ coefficients for some $p$. To do so, we will focus on solving it on a neighborhood of a singular point, producing a local solution around every preimage of it through $\pi$.

We use the same notations as in \ref{co_cu}.
\begin{Prp}\label{prp_lp_desing}Given a $k-$vector field $\xi$ on $\widetilde{X}$, locally integrable on $\widetilde{X}\setminus E$, let $\eta=\pi_*\xi$. Then $\eta\in L^p_\loc(X,\H^2)$ if and only if for every $q\in\tilde{X}$ there exists a local coordinate $t$ such that $\xi\in L^p_\loc(\widetilde{X}, |t|^{(N(q)-1)(2+kp)}dV)$. \end{Prp}
\noindent{\bf Proof: } An easy computation shows that
$$|\eta|\circ\pi=|t|^{k(N(q)-1)}|\xi||1+O(t)|\;$$
therefore
$$\int_{\pi(\mathbb{D})}\!\!\!\!|\eta|^pd\H^2=C\int_{\mathbb{D}}(|\eta|\circ\pi)^p|t|^{2N(q)-2)}dV=C\int_{\mathbb{D}}|\xi|^p|1+O(t)|^p|t|^{kp(N(q)-1)+2N(q)-2}dV$$
and the thesis follows. $\Box$

By Proposition \ref{prp_lp_desing}, we can reduce the problem to solving the Cauchy-Riemann equation with weights. As before, we take $p\in X_\sg$ and $q\in\pi^{-1}(p)$.

\medskip

If $N(q)=1$, we obtain the usual $L^p$ spaces of the unit disc with respect to the Lebesgue measure and we know that the equation $\partial u/\partial \bar{t}=f$ can be solved in $L^p$ by means of the convolution with the Cauchy kernel (ref?).

\medskip

If $N(q)>1$, we recall the following theorem, proved by Fornaess and Sibony in \cite{fornsib}.

\begin{Teo}\label{teo_forn_sib}Let $\Omega\Subset\C$ be an open set, let $1<p\leq 2$, let $\phi:\Omega\to\R\cup\{-\infty\}$ be subharmonic and assume $f:\Omega\to\C$ is a measurable function with $\int_{\Omega}|f|^pe^{-\phi}<\infty$. Then there exists a measurable function $u:\Omega\to\C$ such that 
$$\left(\int_{\Omega}|u|^pe^{-\phi}\right)^{1/p}\leq\frac{5}{2(p-1)}(\mathrm{diam}\Omega)\left(\int_{\Omega}|f|^pe^{-\phi}\right)^{1/p}$$
and such that $\partial u/\partial \bar{t}=f$ in the sense of distributions.\end{Teo}

Set $\Omega=\mathbb{D}^*$ and $\phi(t)=-\alpha\log|t|$. Then  
$$\int_{\Omega}|f|^pe^{-\phi}=\int_{\mathbb{D}}|f|^p|t|^{\alpha}=C\int_{U}|\pi_*f|^pd\omega<+\infty$$
where $U$ is the image of the unit disc in $X$ through the chosen local parametrization. Then Theorem \ref{teo_forn_sib}applies, hence we can find $u\in L^p(\mathbb{D}^*, |t|^{\alpha})=L^p(\mathbb{D},|t|^{\alpha})$ such that $\partial u/\partial \bar{t}=f$ in the sense of distributions on $\mathbb{D}^*$.

\begin{Prp}\label{boh3} Let $x\in X_\sg$, $V$ be an open neighborhood of $x$ and $T\in \D_{0,0}(V)$ with $L^p_{\loc}$ coefficients, $p\in(1,2]$. Then we can find a smaller open set $U$ still containing $x$ and a current $S\in\D_{0,1}(U)$, with $L^p_{\loc}$ coefficients, such that $\debar S=T$ on $U$.\end{Prp}

\noindent{\bf Proof: } We consider, for every $q\in\pi^{-1}(x)$ a local parametrization $\pi_q:\mathbb{D}\to V$ and we set
$$U=\bigcup_{q\in\pi^{-1}(x)}\pi_q(\mathbb{D})\;.$$

\medskip

We notice that
$$f_q=(\pi^{-1}_q\vert_{\pi_q(\mathbb{D}^*)})_\sharp T_{\pi_q(\mathbb{D}^*)}$$
is a $(0,0)-$current on $\mathbb{D}^*$ with coefficients in $L^p(\mathbb{D},|t|^{2N(q)-2})$. 

If $N(q)=1$, $\pi^{-1}$ extends with non-vanishing Jacobian on $x$, so we can find $u_q\in L^p(\mathbb{D})$ such that $\partial u_q(t)/\partial\bar{t}=f_q$; therefore, if we set
$$U_q=u_q\frac{\partial}{\partial t}$$
we have $\debar U_q=f_q$ on $\mathbb{D}$ and $(\pi_q)_\sharp U_q\in\D_{0,1}(\pi_q(\mathbb{D}))$ has $L^p_\loc$ coefficients.

If $N(q)>1$, by Theorem \ref{teo_forn_sib}, we can find $u_q\in L^p(\mathbb{D},|t|^{2N(q)-2})\subseteq L^p(\mathbb{D}, |t|^{(N(q)-1)(2+p)})$ such that $\partial u_q(t)/\partial \bar{t}=f_q$ on $\mathbb{D}^*$ in the sense of distributions; therefore, if we set
$$U_q=u_q\frac{\partial}{\partial t}$$
we have $\debar U_q=f_q$ on $\mathbb{D}^*$ and $(\pi_q)_\sharp U_q\in\D_{0,1}(\pi_q(\mathbb{D}))$ has $L^p_\loc$ coefficients, by Proposition \ref{prp_lp_desing}.

Define $S_1=\sum_q (\pi_q)_\sharp U_q$ and observe that $S_1$ has $L^p_\loc$ coefficients on $U$; moreover $\supp(\debar S_1 - T)\subseteq\{x\}$. It follows that $\debar S_1-T=\alpha_x\delta_x$. Now, if $N(q)=1$ for every $q\in\pi^{-1}(x)$, then $\alpha_p=0$, because we are in the Scenario of section \ref{sec_normal_cr}; if there exists $q$ such that $N(q)>1$, then we consider the current 
$$R=(\pi_q)_\sharp\frac{\alpha_x}{t\pi}\frac{\partial}{\partial t}\;.$$ 
Obviously, $1/t\in L^p(\mathbb{D}, |t|^\beta)$, if $1<p\leq2$ and $\beta\geq2$, so $R$ has $L^p_\loc$ coefficients on $U$. Moreover
$$\debar R=(\pi_q)_\sharp\alpha_p\debar\left(\frac{1}{t\pi}\frac{\partial}{\partial t} \right)=\alpha_x(\pi_q)_\sharp\delta_q=\alpha_x\delta_x\;.$$
If we set $S=S_1+R$, we have $\debar S=T$ on $U$ and $S$ has $L^p_\loc$ coefficients on $U$. $\Box$

\begin{Prp} Let $x\in X_\sg$, $V$ be an open neighborhood of $x$ and $T\in \D_{0,0}(V)$ with $L^2_{\loc}$. Then we can find a smaller open set $U$ still containing $x$ and a current $S\in\D_{1,1}(U)$, with $L^2_{\loc}$ coefficients, such that $\debar S=T$ on $U$.\end{Prp}
\noindent{\bf Proof: } We use the notations of Proposition \ref{boh3}.

We notice that
$$f_q\frac{\partial}{\partial t}=(\pi^{-1}_q\vert_{\pi_q(\mathbb{D}^*)})_\sharp T_{\pi_q(\mathbb{D}^*)}$$
is a $(1,0)-$current on $\mathbb{D}^*$ with coefficients in $L^2(\mathbb{D},|t|^{4(N(q)-1)})$, by Proposition \ref{prp_lp_desing}. 

As in the previous case, if $N(q)=1$ we can solve $\partial u_q/\partial \bar{t}=f_q$ on $\mathbb{D}$ arguing as in Proposition \ref{boh3}. If $N(q)>1$, we can solve $\partial u_q/\partial \bar{t}=f_q$ on $\mathbb{D}^*$. In both cases the solution belongs to the space $L^2(\mathbb{D},|t|^{4(N(q)-1)})\subseteq L^p(\mathbb{D}, |t|^{6(N(q)-1)})$.

Set 
$$U_q=u_q\frac{\partial}{\partial t}\wedge\frac{\partial}{\partial\bar{t}}$$
so that $\debar U_q=f_q\partial/\partial t$ on $\DD^*$. If we consider the current
$$S=\sum_q(\pi_q)_\sharp U_q$$
then we have $\supp(\debar S-T)\subseteq\{p\}$ $\debar S=T$ on $U$ in view of Remark \ref{rem_caratt_W}. $\Box$

\subsection{Density for curves}

For algebraic curves, the problem of the density of bounded vector fields has been solved essentially in \cite{bruning1} by Br{\"u}ning, Peyerimhoff and Schr{\"o}der.  We describe briefly their approach.

\medskip

Let us consider the operator $\debar:\Di^{p,q}(X_\rg)\to\Di^{p,q+1}(X_\rg)$. The (affine or projective) metric on $X_\rg$ is not complete, so, in principle, we have different closed extensions of this operator; we consider only extensions which have $\Di^{p,q+1}(X_\rg)$ in the domain of their adjoint and we denote by $\debar_{\min}$ and $\debar_{\max}$ the minimal and maximal closed extensions with such property.

Namely
$$D_{\debar_{\min}}=\left\{s\in L^2_{p,q}(X,d\H^{2n})\ \vert\ \exists \{s_n\}\in\Di^{p,q}(X_\rg)\ :\begin{array}{lcr} s_n\to s&\mathrm{in}&L^2_{p,q}(X,d\H^{2n})\\ & & \\
\debar s_n\to \debar s& \mathrm{in}&L^2_{p,q}(X,d\H^{2n})\end{array}\right\}$$

\smallskip

$$D_{\debar_{\max}}=\left\{s\in L^2_{p,q}(X, d\H^{2n})\ \vert\ \debar s\in L^2_{p,q}(X, d\H^{2n})\right\}\;,$$
where $\debar s$ is understood in the sense of distributions.

\medskip

With the notation of \ref{ssc_hol_cc}, we set $h(t)=N(q)^2|t|^{2N(q)-2}+O(|t|^{2N(q)})$ and consider $q\in \pi^{-1}(X_\sg)$, with a neighborhood that will be identified with $\mathbb{D}$. Denote by $U$ (resp. $U^*$) the image of $\mathbb{D}$ (resp. $\mathbb{D}^*$) through $\pi_q$  and define the operators
$$\Phi_0:\Di^{0,0}(\mathbb{D}^*)\ni f\mapsto (h^{-1/2}f)\circ\pi_q\in\Di^{0,0}(U^*)$$
and
$$\Phi_1:\Di^{0,1}(U^*)\ni fd\bar{z}\mapsto f\circ\pi_q^{-1}\in\Di^{0,0}(\mathbb{D}^*)\;;$$

$\Phi_0$ and $\Phi_1$ are unitary operators so we can define the operator 
$$D_1:\Phi_1\debar\Phi_0:\Di^{0,0}(\mathbb{D}^*)\to\Di^{0,0}(\mathbb{D}^*)\;;$$ 
Observe that
$$D_1(f)=\frac{\de}{\de\bar{z}}(h^{-1/2}f)\;.$$
Passing to polar coordinates, we define another unitary operator
$$\Phi_2:\Omega^{0,0}((0,\delta)\times \mathbb{S}^1)\to\Omega^{0,0}(\mathbb{D}^*)$$
with a particular choice of the number $\delta=\delta(N(q))$. Finally we set
$$D_2=2e^{-i\phi}\Phi_2^{-1}D_1\Phi_2\;.$$
Then
$$D_2f(x,\phi)=\left(s\frac{\de f}{\de x}+ib\frac{\de f}{\de\phi}+cf\right)(x,\phi)\;;$$
moreover, if $r(p)=\sharp\pi^{-1}(p)$, we set $R=\sum_{p\in X_\sg}r(q)$ and we consider $D_2$ as an operator on
$$L^2([0,\delta], L^2(\mathbb{S}^1_R))=L^2([0,\delta], L^2(\overbrace{\mathbb{S}^1\times\ldots\times\mathbb{S}^1}^{R\ \mathrm{ times}}))\;.$$
Then $D_2$ can be written as
$$D_2=B_1(x)\partial_x+x^{-1}(\tilde{S}_0+\tilde{S}_1(x))$$
where $B_1(x)$ the multiplication by $a(x,\cdot)$ in every $L^2(\mathbb{S}^1)$ and
$$\widetilde{S}_0=\bigoplus_{\substack{p\in\pi^{-1}(X_\sg)\\ q\in\pi^{-1}(p)}}N(q)^{-1}i\partial_\phi-\frac{1}{2}\;.$$
$$\widetilde{S}_1(x)=x^{1/N}b_1(x^{1/N}, \phi)i\partial_\phi+x^{1/N}c_1(x^{1/N},\phi)\qquad\textrm{ on each }L^2(\mathbb{S}^1)\;.$$
With some more technical details to check, we can apply the results from \cite{bruning2} and obtain the following theorem, whose proof can be found in \cite[Theorem 3.1]{bruning1}.

\begin{Teo} All closed extension of $\debar$ between the minimal and the maximal are Fredholm operators and these extensions correspond bijectively to the linear subspaces of the finite dimensional vector space
$$W=D_{\debar_{\max}}/D_{\debar_{\min}}\;.$$
Moreover, if $\debar_V$ is the extension corresponding to $V\subset W$, we have
$$\mathrm{ind}\debar_V=\mathrm{ind}\debar_{\min}+\dim V\;.$$\end{Teo}

We denote by $\chi(\widetilde{X})$ the arithmetic genus of $\widetilde{X}$, i.e.
$$\chi(\widetilde{X})=\sum_{q\geq0}(-1)^q\dim H^{0,q}(\widetilde{X})\;.$$
Using again some results from \cite{bruning2}, we obtain the following equalities.
\begin{eqnarray*}
\mathrm{ind}\debar_{\min}&=&\chi(\widetilde{X})\\
\mathrm{ind}\debar_{\max}&=&\chi(\widetilde{X})+\sum_{\substack{p\in X_\sg\\q\in\pi^{-1}(p)}}(N(q)-1)\end{eqnarray*}
which are contained in Theorems 4.1, 4.2 in \cite{bruning1}.

Finally, keeping the same notations we have the following

\begin{Teo}Given a complex curve $X\subset\CP^n$, the density hypotheses hold if and only if for every $p\in X_\sg$ and $q\in\pi^{-1}(p)$ one has $N(q)=1$.\end{Teo}
\noindent{\bf Proof: } Obviously, if the density holds, $\debar_{\min}=\debar_{\max}$, therefore we need
$$\mathrm{ind}\debar_{\min}=\mathrm{ind}\debar_{\max}$$
that is
$$\sum_{\substack{p\in X_\sg\\q\in\pi^{-1}(p)}}(N(q)-1)=0$$
which can happen if and only if $N(q)=1$. 

On the other hand, if $N(q)=1$ for every $q\in\pi^{-1}(p)$, for every $p\in X_\sg$, then $X$ is a normal-crossing and the density for such spaces was shown in Section \ref{sec_normal_cr}. $\Box$

\medskip

\begin{Rem}An equivalent formula for $\mathrm{ind}\debar_{\max}$ was obtained by Pardon in \cite{pardon1}, by sheaf theoretic methods.\end{Rem}

Clearly, the main point in the proof of the previous theorem consists in showing that density cannot hold if there is a point with $N(q)>1$. We point out that, using Theorem \ref{teo_reg_w}, an alternative proof of this fact could be given by providing examples of holomorphic vector-fields which are not of class $W^{s,2}$ for every $s$. We did it explicitly for the cusp in the second example following the mentioned theorem.

\section{Hypersurfaces}

We give two definitions regarding the solvability of Cauchy-Riemann equation in $L^r$ on a Stein domain in $\C^n$. 
\medskip
Let $q\geq1$. We say that $\Omega$ is \emph{$(C',r)-$regular} if for every $\debar-$closed $(p,q)-$form $\eta$ with coefficients in $L^r(\Omega)$, there exists a $(p,q-1)-$form $\omega$ with coefficients in $L^r(\Omega)$ such that $\debar \omega=\eta$ and $\|\omega\|_{L^r}\leq C'\|\eta\|_{L^r}$.

\medskip

We say that $\Omega$ is \emph{compactly $(C,r')-$regular} if for every $\debar-$closed $(n-p,n-q)-$form $u$ with coefficients in $L^{r'}(\Omega)$ which vanishes almost everywhere outside a compact set $K\subset\Omega$ there exists a $(n-p,n-q-1)-$form $g$ with coefficients in $L^{r'}(\Omega)$ which vanishes almost everywhere outside another compact set $K'\subset\Omega$ and such that $\debar g=u$ and $\|g\|_{L^{r'}}\leq C\|u\|_{L^{r'}}$.

\bigskip

\noindent{\bf Examples } 1) By H\"ormander's results, every bounded open set in $\C^n$ is $(C,2)-$regular, with $C$ depending only on $\mathrm{diam}\Omega$.

2) Every bounded, smoothly bounded strictly pseudoconvex domain in $\C^n$ is $(C,r)-$regular for every $1\leq r\leq+\infty$, with $C$ depending on the Levi form of the boundary and on the fourth derivatives of a defining function (see \cite{Kerzman70, Kerzman71}).

\bigskip

We recall the classical Bochner-Martinelli-Koppelmann formula, referring to \cite{Thi11} for further details. Let
$$k_{BM}(z)=c_n\sum_{j=1}^n(-1)^j\frac{z_jdz_1\wedge\ldots\wedge dz_n\wedge d\bar{z}_1\wedge\ldots \widehat{d{\bar{z}}_j}\ldots\wedge d\bar{z}_n}{|z|^{2n}}$$
$$c_n=(-1)^{n(n-1)/2}\frac{(n-1)!}{(2\pi i)^n}$$
and consider the map $\pi:\C^n\times\C^n\to\C^n$ given by $\pi(z,\zeta)=(z-\zeta)$. Set $K_{BM}=\pi^*k_{BM}$ and denote by $K_{BM}^{p,q}$ the component of $K_{BM}$ of bidegree $(p,q)$ in $z$ and $(n-p,n-q-1)$ in $\zeta$.

\begin{Teo}\label{teo_KBM}Let $D\Subset\C^n$ be a domain with piecewise $\Ci^1$ boundary. Let $0\leq q\leq n$. Then every $(p,q)-$form $v$ of class $\Ci^1$ on $\overline{D}$ is represented on $D$ by
$$v(z)=\int_{bD} K^{p,q}_{BM}(z,\zeta)\wedge v(\zeta)+\debar_z\int_{D}K^{p,q-1}_{BM}(z,\zeta)\wedge v(\zeta)+\int_{D}K_{BM}^{p,q}(z,\zeta)\wedge \debar v(\zeta)\;.$$
In particular, if $q\geq1$, $v$ has compact support and $\debar v=0$, we have that
$$u(z)=\int_{D}K^{p,q-1}_{BM}(z,\zeta)\wedge v(\zeta)$$
is a $\Ci^{1}$ solution of the equation $\debar u= v$ on $D$.\end{Teo}

It is not difficult to extend this formula to the case when $v$ is a classical $\debar-$closed $(p,q)-$current with compact support in $D$: for $\epsilon$ small enough, the regularization $v\ast\rho_\epsilon$ is a smooth, compactly supported, $\debar-$closed $(p,q)-$form in $D$. We obtain a solution
$$u_\epsilon(z)=\int_{D}K^{p,q-1}_{BM}(z,\zeta)\wedge v_\epsilon(\zeta)\;.$$
Given $\phi$ a smooth $(n-p,n-q+1)-$form with compact support in $D$, we have
$$\int_{D}u_\epsilon(z)\wedge \phi(z)=\int_{D}K_{p,q-1}(z,\zeta)\wedge v_\epsilon(\zeta)\wedge \phi(z)=(-1)^{p+q+1}\int_{D}\omega(\zeta)\wedge v_\epsilon(\zeta)$$
where $\omega(\zeta)=\int_{D}K_{p,q-1}(z,\zeta)\wedge\phi(z)$. By the properties of the convolution, $v_\epsilon\to v$ weakly as distributions, therefore
$$\lim_{\epsilon\to0}\int_{D}u_\epsilon(z)\wedge\phi(z)=v(\omega)\;.$$
So, the current $u=\lim_{\epsilon\to0}u_\epsilon$ exists and is such that
$$u(\phi)=(-1)^{p+q+1}v(\omega)$$
with $\debar\omega=\phi$. Moreover, let $\psi$ be a smooth $(n-p,n-q)-$form with compact support such that $\debar\psi=\phi$, then $v(\psi-\omega)=0$, as $\debar v=0$ as a current. Therefore $u(\phi)=(-1)^{p+q+1}v(\psi)$, that is $\debar u= v$ in $D$. We have proved the following result.

\begin{Prp}\label{prp_KBM}Let $v$ be a $(p,q)-$current with compact support in $D$, with $\debar v=0$. Then the current
$$u=\lim_{\epsilon\to0}\int_{D}K^{p,q-1}_{BM}(z,\zeta)\wedge (v\ast\rho_\epsilon(\zeta))$$
is well-defined in $D$ (where the limit is understood in the weak sense) and such that
$$\debar u=v$$
as currents in $D$.\end{Prp}

We will denote the limit $u$ by $K_{BM}\sharp v$.

\medskip

\begin{Rem}\label{rem_reg_KBM}The integral
$$\int_{D}K^{p,q-1}_{BM}(z,\zeta)\wedge v_\epsilon(\zeta)$$
can be understood as a convolution between the coefficients of $v_\epsilon$ and some coefficients of $k_{BM}$; the latters are, around $0$, $O(|z|^{1-2n})$, therefore in $L^r_{\loc}$ for $r<1+(2n-1)^{-1}$. Therefore, if $v$ is of locally finite mass, its coefficients are locally finite Radon measures and the result of the convolution will be a form $u$ with $L^r_{\loc}$ coefficients .\end{Rem}

Let now $X$ be a divisor in $\C^n$, $x\in X$ a point and $U$ a neighborhood of $x$ in $\C^n$. Then $\Omega=U\setminus X$ is a Stein domain in $\C^n$. We have the following result of local solvability for the Cauchy-Riemann equation on $X$.

\begin{Teo} Let $T\in \Di_{0,q}(X)$ with $\supp T$ compact contained in $U$, $M(T)<+\infty$ and $\debar T=0$. If there is $\epsilon$ for which $\Omega$ is $(C,r)-$regular for $1<r<1+\epsilon$, then there exists $S\in \Di_{0,q+1}(U\cap X)$ such that $\debar S=T$ in $\Di_{0,q}(X\cap U)$.\end{Teo}
\noindent{\bf Proof: } Let $j:X\to\C^n$ be the inclusion. We set
$$\widetilde{T}=j_\sharp T$$
and define
$$\widetilde{S}=K_{BM}\sharp \widetilde{T}\;.$$
$\widetilde{S}$ is a (classical) $(0,q+1)-$current on $\C^n$ such that $\debar \widetilde{S}=\widetilde{T}$; we have $M(S)\leq CM(T)$ for some constant $C$ and $d\widetilde{S}=\debar\widetilde{S}=T$, so $\widetilde{S}$ is normal, hence metric.

By Remark \ref{rem_reg_KBM}, $\widetilde{S}$ has coefficients belonging to $L^r_{\loc}(\C^n)$ for some $r>1$. Therefore, $\widetilde{S}$ restricts to a $L^{r}$ form on $\Omega$, where $\debar\widetilde{S}=0$; by hypothesis, $\Omega$ is $(C,r)-$regular for $r$ close enough to $1$, so there exists $\widetilde{R}$ with coefficients in $L^r(\Omega)$ such that
$$\debar \widetilde{R}=\widetilde{S}\qquad\mathrm{ on }\ \Omega\;.$$
Moreover, $\widetilde{R}$ is a $(0,q+2)-$current, therefore $d\widetilde{R}=\debar\widetilde{R}=\widetilde{S}$. So, $\widetilde{R}$ is normal, hence metric on $\Omega$; we can extend it to a normal current on $U$ by inclusion.

We define
$$S=\widetilde{S}- d\widetilde{R}=\widetilde{S}- \debar\widetilde{R}\;.$$
We have $\debar S=\debar{\widetilde{S}}-\debar^2\widetilde{R}=T$ on $U$ and $\supp S\subseteq X\cap U$. As $S$ is normal, this implies that $S\in \Di_{0,q+1}(X\cap U)$. $\Box$

\medskip

The condition of $(C,r)-$regularity for $U\setminus X$ is easily fulfilled for $r\geq 2$. Indeed, for $r=2$, we already observed that every bounded open set is $(C,2)-$regular.

In general, for $r\geq 2$, let $v$ be a $\debar-$closed $(p,q)-$form with $L^r$ coefficients on $U\setminus X$, then for every $(n-p,n-q-1)$ smooth form $\phi$, with compact support in $U\setminus X$, we have that
$$\int_{U\setminus X} v\wedge \debar\phi=0\;.$$
Now, given a $(n-p,n-q-1)$ smooth form $\psi$, with compact support in $U$, we can find smooth functions $\rho_k$ which are $0$ near $X$ and $1$ far from it, approximating the characteristic function of $U$ in $L^{r'}$   (with $r'=r/(r-1)$) and such that
$$\|\nabla \rho_k\|_{r',U}\to_{k\to\infty}0\;.$$
These functions exist because $\H^{2n-r'}(X\cap U)<+\infty$ (see \cite{evans1}). Therefore we have
$$0=\int_{U\setminus X}v\wedge \debar(\rho_k\psi)=\int_{U\setminus X}v\wedge\debar\rho_k\wedge \psi+(-1)^{p+q+1}\int_{U\setminus X}v\wedge\rho_k\debar\psi$$
and
$$\lim_{k\to\infty}\left|\int_{U\setminus X}v\wedge\debar\rho_k\wedge \psi\right|\leq\lim_{k\to\infty}\|v\|_{r, U}\|\nabla\rho_k\|_{r', U}\|\phi\|_\infty=0,$$
so
$$\int_{U}v\wedge\debar\psi=\lim_{k\to\infty}\int_{U\setminus X}v\wedge\rho_k\debar\psi=\lim_{k\to\infty}\int_{U\setminus X}v\wedge\debar\rho_k\wedge \psi+(-1)^{p+q+1}\int_{U\setminus X}v\wedge\rho_k\debar\psi=0\;.$$
This means that the unique $L^r$ extension of $v$ to $U$ is $\debar-$closed in $U$, so we can solve $\debar u=v$ in $L^r$ in $U$ with an estimate $\|u\|_r\leq C\|v\|_r$. The solution can obviously be restricted to $U\setminus X$ and, as $X$ is of measure $0$,  the norms don't change.

\medskip

Therefore, the main problem is to determine if $U\setminus X$ is $(C,r)-$regular for small values of $r$ near $1$. We proceed to establish some equivalent formulation of the problem

\subsection{The Cauchy-Riemann equation in $L^r$}

We start with a lemma of Functional Analysis, which dates back to Fichera \cite{fichera1}. Let $E_1$, $E_2$ be Banach spaces on the field $\mathbb{K}$ which can denote both the real or the complex field, let moreover $E_1^\ast$, $E_2^\ast$ be their topological duals. Consider a vector space $V$ on $\mathbb{K}$ and linear maps $F_1:V\to E_1$, $F_2:V\to E_2$. 

\begin{Teo}\label{teo_banach}
Given a linear functional ${ S }\in E_1^\ast$ such that
\begin{equation}\label{COMP1}({ S }\circ F_1)({\ker}\,F_2)=0\;,\end{equation} 
the equation
\begin{equation}\label{EQ1}
{ S }\circ {F_1}={ T }\circ {F_2}
\end{equation}
has a solution ${ T }\in E^\ast_2$ 
if and only if there exists a constant $C>0$ such that for every  $v\in V$ we have
\begin{equation}\label{STIMA1}
\inf\limits_{w\in{\ker}\, F_2}\Vert F_1(v)+F_1(w)\Vert_{E_1}\le C\Vert F_2(v)\Vert_{E_2}.
\end{equation}
In such a situation, setting $W=\left\{{ T }_0\in E^\ast_2:{ T }_0\circ F_2=0\right\}$, we have 
\begin{equation}\label{STIMA2}
\inf\limits_{{ T }_0\in W}\Vert { T }+{ T }_0\Vert_{E^\ast_2}\le C\Vert { S }\Vert_{E^\ast_1}.
\end{equation}
\end{Teo}

\noindent{\bf Proof: }
1) Let us suppose, at first, that ${\ker}\,F_2\subseteq{\ker}\,F_1$. We set $V_0=V/{\ker}\,F_2$; the maps $F_1$ e $F_2$ descend to the quotient as $\widetilde F_1:V_0\to E_1$, $\widetilde F_2:V_0\to E_2$, $\widetilde F_2$ is injective and the estimate (\ref{STIMA1}) gives
\begin{equation}\label{STIMA3}
\Vert \widetilde F_1(v_0)\Vert_{E_1}\le C\Vert \widetilde F_2(v_0)\Vert_{E_2}.
\end{equation}
for every $v_0\in V_0$. Therefore, we can assume that $F_2$ is injective. 

We begin by proving that estimate (\ref{STIMA1}) implies the existence of a solution to equation (\ref{EQ1}). By the previous paragraph, we can assume that
\begin{equation}\label{STIMA4}
\Vert F_1(v)\Vert_{E_1}\le C\Vert F_2(v)\Vert_{E_2}
\end{equation}
holds for every $v\in V$.

As $F_2(v)\in F_2(V)$ uniquely determines $v$, we can define on $F_2(V)$ a linear functional ${ T }$ by the formula 
$$
{ T }\left(F_2(v)\right):={ S }\left(F_1(v)\right)\;.
$$
 Such a functional is continuous by the estimate (\ref{STIMA4}), as we have
$$
\vert { T }\left(F_2(v)\right)\vert=\vert { S }\left(F_1(v)\right)\vert\le\Vert { S }\Vert_{E^\ast_1}\Vert F_1(v)\Vert_{E_1}\le C\Vert { S }\Vert_{E^\ast_1}\Vert F_2(v)\Vert_{E_2}
$$
which leads to
$$
\Vert{ T }\Vert_{F_2(V)^\ast}\le C\Vert { S }\Vert_{E_1^\ast}.
$$
By Hahn-Banach's theorem ${ T }$ extends as a continuous linear functional on $E_2$ without increasing its norm; we will denote such extension again by ${ T }$. ${ T }$ is a solution of equation (\ref{EQ1}) and verifies (\ref{STIMA3}).

To show necessity of estimate (\ref{STIMA4}), we can suppose $F_1$ and $F_2$ surjective, possibly passing to the closures $\overline{F_1(V)}$, $\overline{F_2(V)}$. Equation (\ref{EQ1}) has now a unique solution ${ T }\in E^\ast_2$ for every ${ S }\in E^\ast_1$. We denote by $L$ the linear operator $E^\ast_1\to E^\ast_2$ defined by $L({ S })={ T }$. We want to show that its graph $\Gamma(L)$ is closed and from this the continuity of $L$ will follow. Let $({ S }_n,L({ S }_n))\subset\Gamma(L)$ be a sequence converging to $({ S },{ T })\in E_1^\ast\times E_2^\ast$. For every $v\in V$ we have
$$
\vert { S }_n(F_1(v))-{ S }(F_1(v)) \vert\le\Vert { S }_n-{ S }\Vert_{E^\ast_1}\Vert F_1(v)\Vert_{E_1}
$$
$$
\vert L({ S }_n)\left(F_2(v)\right)-{ T }\left(F_2(v)\right)\vert\le\Vert L({ S }_n)-{ T }\Vert_{E^\ast_2}\Vert F_2(v)\Vert_{E_2}
$$
so
$$
{ S }\left(F_1(v)\right)=\lim\limits_{n\to+\infty}{ S }_n\left(F_1(v)\right)=\lim\limits_{n\to+\infty}L({ S }_n)(F_2(v))={ T }\left(F_2(v)\right),
$$
i.e. $L({ S })={ T }$. This shows that $\Gamma(L)$ is closed, thus implying that $L$ is continuous.

Estimate (\ref{STIMA2}) follows from the continuity of $L$. Estimate (\ref{STIMA4}) is a consequence of Banach-Steinhaus's theorem. Indeed, set 
$$
\mathcal V=\left\{v\in V:F_2(v)\neq 0\right\}
$$
and, for every $v\in\mathcal V$, let us consider the functional $\delta_v:E^\ast_1\to\mathbb{K}$ given by 
$$
\delta_v({ S })={ S }\left(\frac{F_1(v)}{\Vert F_2(v)\Vert_{E_2}}\right).
$$
For every  $v\in\mathcal V$ and every ${ S }\in E^\ast_1$ we have
$$
\vert\delta_v({ S })\vert=\left\vert{{ S }\left(\frac{F_1(v)}{\Vert F_2(v)\Vert_{E_2}}\right)}\right\vert\le\left\vert { T }\left(\frac{F_2(v)}{\Vert F_2(v)\Vert_{E_2}}\right)\right\vert\le\Vert { T }\Vert\le C\Vert { S }\Vert
$$
so $\{\delta_v\}_v\in\mathcal V$ is a pointwise uniformly bounded family of functionals on $E^\ast_1$. Banach-Steinhaus's theorem implies that $\{\Vert\delta_v\Vert_{E^{\ast\ast}_1}\}_v\in\mathcal V$ is a bounded set, i.e. there exists a constant $C'\in\mathbb{K}$ such that
$$
\Vert\delta_v\Vert_{E^{\ast\ast}_1}=\frac{\Vert F_1(v)\Vert_{E_1}}{\Vert F_2(v)\Vert_{E_2}}\le C'
$$
for every $v\in\mathcal V$. This gives (\ref{STIMA4}).

\medskip

2) The general case can be reduced to the previous one by considering the Banach space $\widetilde E_1=E_1/\overline {F_1({\ker}\,F_2)}$ with the quotient norm. Its dual consists of linear continuous functionals on $E_1$ which vanish on $F_1({\ker}\,F_2)$ and consequently on $\overline {F_1({\ker}\,F_2)}$. If ${ S }$ verifies the hypothesis, it determines a functional $\widetilde{ S }\in {\widetilde E}^\ast_1$. Now, we can consider the linear map $\widetilde F_1:V\to \widetilde E_1$ given by $v\mapsto [F_1(v)]$, where $[F_1(v)]$ is the equivalence class of $F_1(v)$ in $\widetilde E_1$. It's easy to check that ${\ker}\, F_2\subseteq{\ker}\,\widetilde F_1$, so, by what we have proved before,
$$
\widetilde { S }\circ {\widetilde F_1}={ T }\circ {F_2}
$$
has a solution if and only if the following estimate holds
$$
\Vert\widetilde F_1\Vert_{\widetilde E_1}=\inf\limits_{w\in{\ker}\, F_2}\Vert F_1(v)+F_1(w)\Vert_{E_1}\le C\Vert F_2(v)\Vert_{E_2}.
$$
This concludes the proof of the theorem
$\Box$

\bigskip

We move on to the study of the Cauchy-Riemann equation in $L^r$, with or without the compact support assumption.

Let $\Omega$ be an open Stein domain in a complex manifold of complex dimension $n>1$ and $r>1$ a real number; set $r'=r/(r-1)$. In this section, $q$ will always denote a strictly positive integer.

\begin{Teo}\label{teo_cns}Let $C'$ be a positive real number. $\Omega$ is $(C',r)-$regular if and only if
\begin{equation}
\label{eq_cns}
\inf_{\debar \phi_0=0}\|\phi+\phi_0\|_{L^{r'}}\leq C\|\debar\phi\|_{L^{r'}}\qquad \forall \phi\in \D^{n-p,n-q}(\Omega)
\end{equation}
with $\phi_0\in \D^{n-p,n-q}(\Omega)$ as well and $C$ another positive real number.\end{Teo}

\noindent{\bf Proof: } We set
$$E_1=L^{r'}_{n-p,n-q}(\Omega)\qquad E_2=L^{r'}_{n-p,n-q+1}(\Omega)\;,$$
$$V=\D^{n-p,n-q}(\Omega)$$
and
$$F_1(\phi)=\phi\qquad F_2(\phi)=(-1)^{p+q+1}\debar \phi$$
for $\phi\in V$.

It is well known that
$$E_1^\ast=L^r_{p,q}(\Omega)\qquad E_2^\ast=L^r_{p,q-1}(\Omega)$$
with the duality pairing given by
$$(\alpha,\beta)=\int_{\Omega}\alpha\wedge\beta\;.$$

We note that if $\phi\in\ker F_2$, then $\debar\phi=0$ and, as $H^{n-p,n-q}_c(\Omega)=0$ for $q>0$, we have $\psi\in\D^{n-p,n-q+1}(\Omega)$ such that $\debar\psi=\phi$, so condition (\ref{COMP1}) reads
$$0=\int_{\Omega}\alpha\wedge \phi=\int_{\Omega}\alpha\wedge\debar\psi=0\qquad \forall \psi\in\D^{n-p,n-q+1}(\Omega)\;,$$
i.e. $\debar\alpha=0$ in the sense of distributions.

\bigskip

We now prove that equation (\ref{eq_cns}) is a sufficient condition for $(C',r)-$regularity. Let $\eta\in E_1^\ast$ such that $\debar\eta=0$ and consider the equation (\ref{EQ1}); if $\omega\in E_2^\ast$ is a solution, we have
$$\int_{\Omega}\eta\wedge \phi=(-1)^{p+q+1}\int_{\Omega}\omega\wedge\debar \phi\qquad\forall \phi\in\D^{p,q}(\Omega)\;.$$
This means that 
\begin{equation}\label{eq_vera}
\debar(\omega)=\eta
\end{equation}
holds in the sense of distributions. Moreover, from (\ref{STIMA2}), we get
\begin{equation}
\label{eq_stimavera}
\inf_{\debar\omega_0=0}\|\omega+\omega_0\|_{L^r}\leq C\|\eta\|_{L^r}\;.
\end{equation}

Therefore, for a given $\epsilon>0$, there exists $\omega_0$ such that $\debar\omega_0=0$ and $\|\omega+\omega_0\|_{L^r}\leq (C+\epsilon)\|\eta\|_{L^r}$, so $\omega_1=\omega+\omega_0$ is the required solution. By Theorem \ref{teo_banach}, if
$$\inf_{\debar\phi_0=0} \|\phi+\phi_0\|_{L^{r'}}\leq C\|\debar\phi\|_{L^{r'}}\qquad \forall \phi\in V$$
then (\ref{eq_vera}) is always solvable and the estimate (\ref{eq_stimavera}) holds.

\medskip

To show necessity of equation (\ref{eq_cns})), we note that, if (\ref{eq_vera}) has always a solution which satistfies $\|\omega_1\|_{L^r}\leq C'\|\eta\|_{L^r}$, then 
$$\inf_{\debar\omega_0=0}\|\omega_1+\omega_0\|_{L^r}\leq C'\|\eta\|_{L^r}$$
and therefore
$$\inf_{\debar\phi_0=0} \|\phi+\phi_0\|_{L^{r'}}\leq C'\|\debar\phi\|_{L^{r'}}\qquad \forall \phi\in V\;.$$
This concludes the proof. 
$\Box$

\medskip

\begin{Teo}If $\Omega$ is compactly $(C,r')-$regular, then it is $(C',r)-$regular.
\end{Teo}
\noindent{\bf Proof: } By Theorem \ref{teo_cns}, it is enough to show that (\ref{eq_cns}) holds. Therefore, let us consider $\phi\in\D^{n-p,n-q}(\Omega)$ and set $u=\debar\phi$. By hypothesis, there exists a $(n-p,n-q)-$form $g$ with $L^{r'}$ coefficients vanishing outside some compact set such that $\debar g=u=\debar\phi$ and $\|g\|_{L^{r'}}\leq C\|u\|_{L^{r'}}=C\|\debar\phi\|_{L^{r'}}$.

Set $\psi_{\epsilon}=(g-\phi)\star\rho_\epsilon$ and note that $\psi_\epsilon\in \D^{n-p,n-q}(\Omega)$ and $\debar\psi_\epsilon=0$. Moreover
$$\|\phi+\psi_{\epsilon}\|_{L^{r'}}\leq\|g\star\rho_\epsilon\|_{L^{r'}}+\|\phi-\phi\star\rho_\epsilon\|_{L^{r'}}\;,$$
so, for any given $\delta>0$, we can take $\epsilon$ small enough such that
$$\|\phi+\psi_{\epsilon}\|_{L^{r'}}\leq \|g\|_{L^{r'}}+\delta\leq C\|\debar\phi\|_{L^{r'}}+\delta\;.$$

Therefore
$$\inf_{\debar \phi_0=0}\|\phi+\phi_0\|_{L^{r'}}\leq \inf_{\epsilon}\|\phi+\psi_{\epsilon}\|_{L^{r'}}\leq\lim_{\delta\to0}C\|\debar\phi\|_{L^{r'}}+\delta=C\|\debar\phi\|_{L^{r'}}\;,$$
which is what we needed to show. $\Box$\

\medskip

\begin{Prp}\label{prp_exh}Let $\{U_j\}_{j\in\N}$ be an exhausting sequence for $\Omega$ made of relatively compact Stein domains and suppose that $U_j$ is $(C'_j,r)-$regular. Then for every $\debar-$closed $(n-p,n-q)-$form $u$ with coefficients in $L^{r'}(\Omega)$ which vanishes almost everywhere outside a compact set $K\subset\Omega$ there exists a $(n-p,n-q-1)-$form $g$ with coefficients in $L^{r'}(\Omega)$ which vanishes almost everywhere outside another compact set $K'\subset\Omega$.
\end{Prp}
\noindent{\bf Proof: } Let $u$ be as said and consider $u_{\epsilon}=u\star\rho_\epsilon$. For $\epsilon$ small enough, there is $j$ such that $\supp u\Subset U_j$ and $u_\epsilon\in\D^{n-p,n-q}(U_j)$; moreover $\debar u_\epsilon=(\debar u)\star\rho_\epsilon=0$. As $H^{n-p,n-q}_c(U_j)=0$, there exists $\phi_\epsilon\in\D^{n-p,n-q-1}(U_j)$ such that $u_\epsilon=\debar\phi_\epsilon$. By Theorem \ref{teo_cns}, we have that
$$
\inf_{\substack{\supp \phi_0\subset U_j\\\debar \phi_0=0}}\|\phi_\epsilon+\phi_0\|_{L^{r'}}\leq C_j\|\debar\phi_\epsilon\|_{L^{r'}}=C_j\|u_\epsilon\|_{L^{r'}}\leq C_j\|u\|_{L^{r'}}\;.
$$
Let $g_\epsilon=\phi_\epsilon+\phi_0$ be such that $\|g_{\epsilon}\|_{L^{r'}}\leq C_j\|u\|_{L^{r'}}+\epsilon$; as the forms $g_\epsilon$ have smooth coefficients, they can be viewed as forms with coefficients in $L^{r'}(\Omega)$, with the same estimates on the norms. The set $\{g_\epsilon\}$ is weakly-$\star$ compact in $L^{r'}_{n-p,n-q-1}(\Omega)$; let $g$ be its weak-$\star$ limit.

By definition
$$
\int_{\Omega}g_\epsilon\wedge \psi=\int_{\Omega}g\wedge \psi\qquad \forall \psi\in\D^{p,q+1}(\Omega)\;,
$$
Moreover, if $\psi=\debar \eta$ with $\eta\in\D^{p,q}(\Omega)$, we have
$$
\int_{\Omega}u_\epsilon\wedge\eta=\int_{U_j}u_\epsilon\wedge \eta=\pm\int_{U_j}g_\epsilon\wedge\debar\eta=\pm\int_{\Omega}g_\epsilon\wedge\debar\eta\;.
$$
By convolution, $u_\epsilon\to u$ in $L^{r'}-$norm, so the first integral tends to
$$\int_{\Omega}u\wedge \eta\;.$$
By the weak-$\star$ convergence the last integral tends to
$$\int_{\Omega}g\wedge\debar \eta\;.$$
These two being equal (the sign depending on the degrees of the forms), we have that
$$\debar g= u$$
in the sense of distributions.

Now, if $\psi\in\D^{p,q+1}(\Omega)$ is supported outside $U_j$, 
$$\int_{\Omega}g_\epsilon\wedge \psi=0$$
and so
$$\int_{\Omega}g\wedge \psi=0$$
thus implying that $g$ vanishes almost everywhere outside $U_j$.
 $\Box$
 
 \medskip
 
 \begin{Rem}By the semicontinuity of the norm under weak-$\star$ convergence, we obtain that
 $$\|g\|_{L^{r'}}\leq C_j\|u\|_{L^{r'}}\;.$$
 Therefore, if $\sup C_j=C<+\infty$, we can conclude that $\Omega$ is compactly $(C,r')-$regular.\end{Rem}
 
 \begin{Rem}Every Stein domain on a complex manifold has the exhaustion required in Proposition \ref{prp_exh}. Indeed, we can exhaust any Stein domain with smoothly bounded strongly pseudoconvex domains and by the work of Kerzmann (\cite{Kerzman70, Kerzman71}), Demailly and Laurent-Thi\`ebault (\cite{demailly1}) we have the $(C'_j,r)-$regularity of $U_j$.\end{Rem}

\medskip

For instance, if $\Omega$ is a weakly pseudoconvex open set in $\C^n$ and $r=r'=2$, H\"ormander's results imply that $\Omega$ is $(C,2)-$regular, where $C$ depends only on $n$ and the diameter of $\Omega$. Therefore, we can exhaust $\Omega$ with open sets $U_j$ and the constants $C_j$ will converge to the constant $C$. So, an open set in $\C^n$ is $(C,2)-$regular if and only if is compactly $(C,2)-$regular.

\medskip

Another example is given by bounded, smoothly bounded strongly pseudoconvex domains in a complex manifold: such domains are $(C,r)-$regular with a constant $C$ depending on the derivatives of a defining function up to the fourth order. As it is always possible to approximate such a domain with smaller ones so that the boundaries converge as $\Ci^4$ manifolds, the constants $C_j$ converge as well.

\medskip

The problem of determining if $\Omega=U\setminus X$ is $(C,r)-$regular led to some interesting representation formulas for the $\debar$ equation in the polydisc (see \cite{amar1}).

\chapter{Complex geometry in infinite dimension}
\epigraphhead[60]{\epigraph{Let's think the unthinkable, let's do the undoable. Let us prepare to grapple with the ineffable itself, and see if we may not eff it after all.}{D. Adams - \emph{Dirk Gently's Holistic Detective Agency}}}

This last chapter deals with some infinite dimensional problems in complex analysis and geometry. We begin by studying the Cauchy-Riemann equation for metric currents in Banach spaces, introducing a variation of the metric currents, which we call quasi-local currents.  An adaptation of the cone construction allows us to show that Cauchy-Riemann equation is always solvable in terms of quasi-local currents.

Then we move on to geometric problems, such as the study of positive currents, the characterization of currents of integration on an analytic set and the boundary problem for holomorphic chains.

The formulations, letting alone the solutions, of these problems depend heavily on the definition of analytic set we want to work on. We give the definition of finite-dimensional analytic set following Aurich, Ramis, Ruget and others.

However, this definition won't always fit our needs, as the characterization of holomorphic chains shows.

\section{Cauchy-Riemann equation in Banach spaces}

The results collected here appear also in \cite{mongodi1}.

In the section, we focus our attention on metric currents on infinite dimensional Banach spaces. We want to stress on the fact that the local version of metric currents developed by U. Lang does not make sense in this setting, due to the fact that the compact sets  have empty interior. In the next section we will propose a replacement for these local objects by quasi-local currents.

Here we examine the behavior of metric currents in relation with their projections on finite dimensional subspaces. In order to recover informations on the whole space from its finite dimensional subspaces, we introduce the following category of Banach spaces (see also \cite{noverraz1}).

A (complex) Banach space $E$  is said to have the \emph{projective approximation property (PAP for short)} if there exist a constant $a$ and an increasing collection $\{E_t\}_{t\in T}$ of finite dimensional subspaces of $E$ such that
\begin{enumerate}
\item[1)] $\{E_t\}_{t\in I}$ is a directed set for the inclusion;
\item [2)]$\displaystyle{E=\overline{\bigcup_{t\in I} E_t}}$;
\item[3)] for every $t\in I$ there exists a projection $p_t:E\to E_t$ with $\|p_t\|\leq a$.
\end{enumerate}

Every Banach space with a Schauder basis has the PAP. Two important cases of PAP Banach spaces with no Schauder basis are $\Ci(K)$, the space of continuous functions on a compact space K with the sup norm and $L^p(X,\mu)$, with $1\leq p\leq+\infty$, where $X$ is a locally compact space and $\mu$ being a positive Radon measure. In this section, we will work with Banach spaces having the PAP; we will endow the set $I$ of indeces with the partial ordering coming from the inclusion relation between subspaces.

\medskip

\begin{Prp}\label{prp_limlip}Let $f\in \Lip(E)$ and define $f_t=f\circ p_t$. Then $f_t\to f$ pointwise and $\Lip(f_t)\leq a\Lip(f)$, for every $t\in I$.
\end{Prp}
\noindent{\bf Proof: } By property 2) in the definition of PAP, for every $x\in E$ there exists a sequence $\{x_j\}\subset E$, with $x_j\in E_{t_j}$, such that $x_j\to x$. By property 1), for a given $j$, we have that if $t>t_j$, then $f_t(x_j)=f(x_j)$. Moreover, 
$$\Lip(f_t)\leq \Lip(f)\cdot\|p_t\|\leq a\Lip(f)\;;$$
so, given $t, t'\in I$, let $j$ be such that $t_j\leq t$ and $t_j\leq t'$, then $f_t(x_j)=f_{t'}(x_j)=f(x_j)$ and
$$|f_t(x)-f_{t'}(x)|=|f_t(x)-f_t(x_j)+f_{t'}(x_j)-f{t'}(x_j)|\leq 2a\Lip(f)\|x-x_j\|$$
which goes to $0$ as $j\to\infty$. $\Box$

\begin{Prp}\label{prp_conv_proj}Let $T\in M_k(E)$ and define $T_t=(\pi_t)_\sharp(T)\in M_k(E_t)$ for every $t\in I$ such that $\dim_{\C}E_t\geq k$. By means of the inclusion $i_t:E_t\to E$, we can consider $T_t$ as an element of $M_k(E)$ and then, $T_t\to T$ weakly.\end{Prp}
\noindent{\bf Proof: } Let $\mu_t$ be the mass of $T_t$ and $\mu$ the mass of $T$; then $\mu_t=(p_t)_\sharp \mu$. By \cite[Lemma 2.9]{ambrosio1}, the support of $\mu$ is a $\sigma-$compact set, therefore for every $\epsilon>0$ there exists a compact $K_\epsilon$ such that $\mu(E\setminus K_\epsilon)\leq \epsilon$. As $p_t\to \mathrm{Id}_E$ uniformly on every compact set (because of PAP), we have that $\mu_t\to \mu$ on $K_\epsilon$, which implies that
$$\int_E f\circ p_t d\mu=\int_{E}fd\mu_t\to \int_E fd\mu$$
for every $f\in \Lip_b(E)$. This result and Proposition \ref{prp_limlip} now imply
$$T_t(f,\pi)=T(f\circ p_t,\pi\circ p_t)\to T(f,\pi)$$
which is our thesis. $\Box$

\medskip

Let $\{E_t, p_t\}_{t\in I}$ be the countable collection of subspaces and projections given by PAP. We call it a \emph{projective approximating sequence} (PAS for short) if $p_t\circ p_s=p_{\min\{s,t\}}$.

We note that every separable Hilbert space or, more generally, every Banach space with a Schauder basis contains a PAS.

\begin{Teo}\label{teo_ext}Let us suppose that $\{E_t, p_t\}$ is a PAS in $E$. If we are given a collection of metric currents $\{T_t\}_{t\in I}$ such that 
\begin{enumerate}
\item $T_t\in N_k(E_t)$,
\item $(p_{t}\vert_{E_{t'}})_\sharp T_{t'}= T_{t}$ for every $t, t'\in I$ with $t'>t$,
\item $\|T_t\|\leq (p_t)_*\mu$ and $\|d T_t\|\leq (p_t)_*\nu$ for every $t\in I$ and some $\mu,\ \nu$ finite Radon measures on $E$. 
\end{enumerate}
then there exists $T\in N_k(E)$ such that $(p_t)_\sharp T=T_t$ for every $t\in I$. \end{Teo}

\noindent{\bf Proof: } We consider, in $\E^k(E)$, the subspaces $(p_t)^*\E^k(E_t)$: their union $\mathcal{P}^k$ is dense, with respect to pointwise convergence, with bounded Lipschitz constants. We define a functional $T:\mathcal{P}^k\to\C$, by setting $T(f,\pi)=T_t(f,\pi)$, with $t$ such that $(f,\pi)\in(p_t)^*\E^k(E_t)$. By hypothesis 2), this definition is well posed; the functional so defined is obviously multilinear and local on $\mathcal{P}^k$; moreover, by hypothesis 3), we have that there exists a finite Radon measure $\mu$ on $E$ such that
$$T(f,\pi)\leq \prod_{j=1}^k \Lip(\pi_j) \int_E|f|d\mu\qquad \forall\ (f,\pi)\in\mathcal{P}^k\;.$$
In the same way, we infer that there exists a finite Radon measure $\nu$ such that
$$dT(f,\pi)\leq\prod_{j=1}^{k-1} \Lip(\pi_j)\int_E|f|d\nu\qquad\forall\ (f,\pi)\in\mathcal{P}^{k-1}\;;$$
therefore, the current $T$ is also continuous on $\mathcal{P}^k$, being normal on this set of metric forms (see \cite[Proposition~5.1]{ambrosio1}). 

Extending $T$ by density, we obtain a multilinear, local, continuous functional on $\E^k(E)$, whose mass is bounded by $\mu$ and whose boundary's mass is bounded by $\nu$; thus, the extension is a normal current, which we denote again with $T$, and it is not hard to check that $T$ satisfies $(p_t)_\sharp T=T_t$ $\forall\ t\in T$. $\Box$

\medskip

We can substitute the request of the existence of a PAS and the compatibility condition (hypothesis \emph{(2)}) with an assumption on the existence of a global object. A \emph{metric functional} is a function $T:\E^k(E)\to\C$ which is subadditive and positively $1-$homogeneous with respect to every variable. For metric functionals, we can define mass, boundary and pushforward (see Section 2 of \cite{ambrosio1}).

\begin{Prp}Let $E$ be a Banach space with PAP. Suppose that $T:\E^k(E)\to\C$ is a metric functional with $T$ and ${\rm d}T$ of finite mass, such that $(p_t)_\sharp T\in N_k(E_t)$ for every $t\in I$. Then $T\in N_k(E)$.\end{Prp}
\noindent{\bf Proof: } We have 
$$\|(p_t)_\sharp T\|\leq a(p_t)_*\|T\|$$
$$\|d(p_t)_\sharp T\|\leq a(p_t)_*\|dT\|\;,$$ 
so, by  the previous Theorem there exists $\widetilde{T}\in N_k(E)$ such that $(p_t)_\sharp \widetilde{T}=T_t$. This means that $\widetilde{T}$ and $T$ coincide on the metric forms in $\mathcal{P}^k$; by density, we conclude that $\widetilde{T}=T$. $\Box$

\subsection{Bidimension}

The definition of bidimension given in Section \ref{sec_dolb_dec} is meaningful also for a complex Banach space. For a careful analysis of the notion of holomorphy in this context we refer the interested reader to the first chapters in \cite{noverraz1}. Here we only notice that Lipschitz holomorphic functions are not necessarily dense in the space of Lipschitz functions and we cannot work with local concepts as in Section \ref{sec_dolb_dec}, because the spaces of local currents do not make sense on a Banach space.

However, inspired by the links we found between the finite dimensional projections of a current and the current itself, we would like to give a different characterization of $(p,q)-$currents. 

We say that $T\in M_{k}(E)$ is \emph{finitely} of bidimension $(p,q)$ if every finite dimensional projection of it is a $(p,q)-$current.

\begin{Prp} $T\in M_k(E)$ is a $(p,q)-$current if and only if it is finitely so.\end{Prp}
\noindent{\bf Proof: } A projection $p:E\to V$ is a continuous complex linear operator, thus holomorphic, so $T\in M_{p,q}(E)$ implies $p_\sharp T\in M_{p,q}(V)$, so one implication is proved.

On the other hand, if $h\in \Ol(E)$, then $h\vert_{E_t}\in \Ol(E_t)$ for every $t\in T$; so, if $(f,\pi)\in \E^k(E)$ contains $p+1$ holomorphic differentials, then so does $(f\vert_{E_t}, \pi\vert_{E_t})\in \E^k(E_t)$. Therefore,
$$T(f,\pi)=\lim_{t\in T} T(f\circ p_t, \pi\circ p_t)=\lim_{t\in T} T_{t}(f,\pi)\;.$$
As $T_t$ is a finite dimensional projection, it is of bidimension $(p,q)$, so the right hand side is zero. The same argument applies when $(f,\pi)$ contains $q+1$ antiholomorphic differentials, giving us the desired conclusion. $\Box$

As an application of Theorem \ref{teo_ext}, we have the following result about the existence of a Dolbeault decomposition for $T\in M_k(E)$.

\begin{Prp}\label{prp_dec}Let us suppose that $\{E_t, p_t\}$ is a PAS in $E$. Let $T\in N_k(E)$; if $T_t$ has a Dolbeault decomposition in normal $(p,q)-$currents in $E_t$ for all $t\in T$, with a finite Radon measure $\nu$ (independent of $t$) whose pushforward dominates the boundaries' masses, then also $T$ admits a Dolbeault decomposition.\end{Prp}
\noindent{\bf Proof: } Let us fix a pair $(p,q)$ such that $p+q=k$ and let $S_t$ be the $(p,q)-$component of $T_t$; by hypothesis, $S_t\in N_{p,q}(E_t)$ and $\|d S_t\|\leq (p_t)_*\nu$, independently of $t$, and it is not hard to show that $\|S_t\|\leq C'\|T_t\|\leq C''\|T\|$, with $C',\ C''$ independent of $t$ (in particular, independent of $\dim E_t$).

Last thing to check is the compatibility condition (condition ii)) in Theorem \ref{teo_ext}), but this follows easily from the invariance of the bidimension under pushforward by holomorphic maps. Applying Theorem \ref{teo_ext}, we have the thesis. $\Box$

\begin{Rem} In general it is not easy to verify the hypotheses of Proposition \ref{prp_dec} for a current $T\in N_k(E)$; however, this result is an example of a general phenomenon: in a Banach space with the projective approximation property, it is often enough to check a certain property for finite dimensional subspaces in order to obtain that it holds for the whole space. For instance, any equality between currents holds in $E$ if and only if it holds finitely, namely whenever the currents are pushed forward through a finite rank projection.\end{Rem}

Employing the idea given in this Remark, we can show the following.

\begin{Cor}If $T\in M_k(E)$ admits a Dolbeault decomposition, then it is unique.\end{Cor}

\subsection{Quasi-local metric currents}

To partially overcome the problems related to the lack of a local theory of currents, we introduce a new definition, which somehow locates midways between the local and the global one.

Let $\E_{\q}(E)$ be defined as follows:
$$\E_{\q}(E)=\bigcup_{R>0}\{f\in\Lip(E)\ :\ \supp f\subset B(0,R)\}$$
where $B(0,R)$ is the ball of centre $0$ and radius $R$ in $E$. We say that a sequence $\{f_j\}\subset \E_{\q}(E)$ converges to $f\in \E_{\q}(E)$ if $f_j\to f$ pointwise, $\Lip(f_j), \Lip(f)\leq C$ and $\supp f_j, \supp f\subset B(0,R)$ for some $R$.

We also define $\Lip_{\q}(E)$ to be
$$\{ f\in \Ci^0(E)\ :\ f\vert_{B(0,R)}\in \Lip(B(0,R))\ \forall\ R>0\}\;.$$
A sequence $\{\pi_j\}$ in this space converges if it converges pointwise and the Lipschitz constants on any fixed ball are uniformly bounded in $j$ (but not necessarily with respect to the radius of the ball).

Finally, we define the spaces of quasi-local metric forms as
$$\E^k_\q(E)=\E_\q(E)\times[\Lip_\q(E)]^k\;.$$

\medskip

A \emph{quasi-local} $k-$dimensional metric current is a functional $T:\E^k_\q(E)\to\C$ satisfying the following
\begin{enumerate}
\item $T$ is multilinear
\item if $(f^i,\pi^i)\to (f,\pi)$ in $\E^k_\q(E)$ then $T(f^i,\pi^i)\to T(f,\pi)$
\item $T(f,\pi)=0$ whenever there is an index $j$ such that $\pi_j$ is constant on a neighborhood of $\supp f$
\item for every $R>0$ there is a finite Radon measure $\mu_R$ such that
$$|T(f,\pi)|\leq\prod_{j=1}^k\Lip(\pi_j\vert_{B(0,r)})\int\limits_{B(0,R)}|f|d\mu_R$$
for every $(f,\pi)\in\E^k_\q(E)$ with $\supp f\subset B(0,R)$.\end{enumerate}
We denote the space of such currents by $M_{k,\q}(E)$.

\medskip

The last condition can be rephrased as: there exits a Radon measure $\mu$ on $E$, such that $\mu(B(0,R))<+\infty$ for every $R>0$ and such that
$$|T(f,\pi)|\leq\prod_{j=1}^k\Lip(\pi_j\vert_{\supp f})\int\limits_{\supp f}|f|d\mu$$
for every $(f,\pi)\in\E^k_\q(E)$. If $\mu$ happens to be a \emph{finite} Radon measure, then $T$ is indeed a $k-$dimensional metric current in the sense of Ambrosio and Kirchheim.

\medskip

The definitions of boundary, pushforward, contraction are the same of the usual metric currents. The pushforward can be defined for \emph{quasi-local} proper Lipschitz map, that is any map which is Lipschitz on $B(0,R)$ for every $R>0$ and such that the preimage of any bounded set is bounded.

\begin{Rem} The projections $p_t$ are by no means quasi-local, so we cannot repeat verbatim the arguments of the previous section.\end{Rem}

The space $N_{k,\q}(E)$ is defined as the set of quasi-local currents whose boundary is again a quasi-local current, that is, has quasi-locally finite mass.

\medskip

By the mass condition, we can extend any $T\in M_{k,\q}(E)$ to a functional on $k+1-$tuples $(f,\pi)$ where $\pi\in [\Lip_{\textrm{q-loc}}(E)]^k$ and $f\in\mathcal{B}_b^{\infty}(E)$, i.e. the algebra of bounded Borel functions with bounded support in $E$. The basic properties of metric currents hold true also for this quasi-local variant. Namely, we have the following.

\begin{Prp}Given $T\in M_{k,\q}(E)$, we denote again by $T$ its extension to $\mathcal{B}^\infty_b(E)\times[\Lip_{\textrm{q-loc}}(E)]^k$; then
\begin{enumerate}
\item $T$ is multilinear in $(f,\pi)$ and
$$T(fd\pi_1\wedge\ldots\wedge d\pi_k)+T(\pi_1df\wedge\ldots\wedge d\pi_k)=T(\sigma d(f\pi_1)\wedge\ldots\wedge d\pi_k)$$
whenever $f,\pi\in \E_\q(E)$ and $\sigma\in\mathcal{B}^\infty_b(E)$ is equal to $1$ on the support of $f\pi_1$ and
$$T(fd\psi_1(\pi)\wedge\ldots\wedge d\psi_k(\pi))=T(f\det\nabla\psi(\pi)d\pi_1\wedge\ldots\wedge d\pi_k)$$
whenever $\psi=(\psi_1,\ldots,\psi_k)\in\mathcal{C}^1(\R^k,\R^k)$;
\item $$\lim_{i\to\infty}T(f^i,\pi^i_1,\ldots,\pi^i_k)=T(f,\pi)$$
whenever $f^i-f\to0$ in $L^1(E,\|T\|)$ and $\pi^i_j\to\pi_j$;
\item $T(f,\pi)=0$ if $\{f\neq 0\}=\bigcup B_i$ with $B_i\in\mathcal{B}(E)$ and $\pi_i$ constant on $B_i$.
\end{enumerate}
\end{Prp}

The definition of $(p,q)-$current given in Section \ref{sec_dolb_dec} can be applied also to quasi-local currents. We have the same results presented in \ref{ssc_prop_pq}. Indeed, given $(f,\pi)\in \E^k_\q(E)$, the differentials $\pi_1,\ldots, \pi_k$ can be approximated by analytic functions of a finite number of variables, so the proofs go on almost identically.

\subsection{Quasi-local solution to $\debar U=T$}

Given a function $a:E\to\C$, we set
$$a_t(x)=a(tx)\qquad \forall\ t\in\C\;.$$

Let $T\in N_{k}(E)$ be a $(0,k)-$current, with $\supp T$ bounded and $0\not\in\supp T$; we define the following $(k+1)-$dimensional metric functional
$$C_{\debar}(T)(f,\pi_1,\ldots,\pi_{k+1})=\frac{1}{2\pi i}\sum_{j=1}^{k+1}(-1)^{j+1}\int_{\C}T\left(f_t\frac{\de \pi_{jt}}{\de \bar{t}}d\widehat{\pi}_{jt}\right)\frac{dt\wedge d\bar{t}}{t-1}\;,$$
where $d\widehat{\pi}_j$ is the wedge product of all the differentials different from $\pi_j$.

\begin{Lmm}$C_{\debar}(T)$ is a multilinear, local metric functional, with quasi-locally finite mass.\end{Lmm}
\noindent{\bf Proof: } Multilinearity and locality are obvious. Let $B(0,r)$ be a ball containing the support of $T$ and let $B(0,d)$ be a ball disjoint from the support of $T$.

We have that, for $fd\pi\in\E^{k+1}_\q(E)$, with $\Lip(\pi_i)=1$ and $\supp f\subset B(0,R)$, the following holds
$$|C_{\debar}(T)(fd\pi)|\leq \frac{r(k+1)}{2\pi}\int_{|t|<R/d}\frac{|t|^k}{|t-1|}\int_{E}|f_t|d\|T\|dt\wedge d\bar{t}\;.$$
Moreover, given a bounded borel set $A$, 
$$\|C_{\debar}(T)\|(A)\leq \frac{r(k+1)}{2\pi}\int_{|t|<R/d}\frac{|t|^k}{|t-1|}\|T\|(A/t)dt\wedge d\bar{t}\;.$$
We want to estimate this quantity when $A=B(0,R)$, with $R>r$. We split the integral in $t$ in two parts: one small ball around the origin and the rest of the ball of radius $R/d$. In the small ball around the origin, we have
$$\int_{|t|<\epsilon}\frac{|t|^k}{|t-1|}\|T\|(B(0,R/t))dt\wedge d\bar{t}\leq c_1\epsilon^{2+k}M(T)\;.$$
On the rest of the outer ball, we have
$$\int_{\epsilon<|t|<R/d}\frac{|t|^k}{|t-1}\|T\|(B(0,R/t))dt\wedge d\bar{t}\leq (R/d)^{k+2}M(T)\;.$$
So, letting $\epsilon\to0$, we get
$$M_{B(0,R)}(C_{\debar}(T))\leq \frac{r(k+1)}{2\pi}(R/d)^{k+2}M(T)\;.$$
Therefore, the mass of $C_{\debar}(T)$ is quasi-locally finite. $\Box$

\begin{Prp}\label{prp_d_cone}We have
$$d C_{\debar}(T)=C_{\debar}(d T) + T$$
as quasi-local metric functionals.
\end{Prp}
\noindent{\bf Proof: } Let $fd\pi$ be a quasi-local metric $k-$form such that $f$, $\pi_1,\ldots, \pi_{k}$ have Lipschitz derivatives. Define the function
$$\phi(t)=T\left( f_td\pi_t\right)$$
and observe that, as $T$ is a quasi-local current, thus continuous, then
$$\frac{\de\phi}{\de\bar{t}}=T\left(\frac{\de f_t}{\de\bar{t}}d\pi_t\right)+\sum_{j=1}^k(-1)^{j+1}T\left(fd\frac{\de\pi_{jt}}{\de\bar{t}}\wedge d\widehat{\pi}_{jt}\right)\;.$$
By the definition of boundary this expression can be rewritten as
\begin{equation}\label{eq_bdry}T\left(\frac{\de f_t}{\de\bar{t}}d\pi_t\right)+\sum_{j=1}^k(-1)^{j}\left[T\left(\frac{\de\pi_{jt}}{\de\bar{t}}df_t\wedge d\widehat{\pi}_{jt}\right)-dT\left(f_t\frac{\de\pi_{jt}}{\de\bar{t}}d\widehat{\pi}_{jt}\right)\right]\;,\end{equation}
exactly as in the proof of Proposition 10.2 in \cite{ambrosio1}.

Given a generic form $fd\pi\in\E^{k}_\q(E)$, the conclusion still holds: it is enough to approximate $f$ and $\pi_j$ by
$$f^\epsilon(x)=\int_{\C}f(sx)\rho_\epsilon(s)ds\wedge d\bar{s}\;,\qquad \pi_j^\epsilon(x)=\int_{\C}\pi_j(sx)\rho_\epsilon(s)ds\wedge d\bar{s}\;,$$
where $\rho_\epsilon$ are convolution kernels, compactly supported, $w^*-$converging to $\delta_1$.

By Fubini's theorem we have
$$\lim_{\epsilon\to0}\frac{\de f^{\epsilon}_t}{\de\bar{t}}(x)=\frac{\de f_t}{\de\bar{t}}(x)\;,\quad\lim_{\epsilon\to0}\frac{\de \pi_{jt}^\epsilon}{\de\bar{t}}(x)=\frac{\de\pi_{jt}}{\de\bar{t}}(x)$$
for $\|T\|+\|dT\|-$a.e. $x$, for $\mathcal{L}^2-$a.e. $t$; therefore the derivative with respect to $\bar{t}$ of $t\mapsto T(f^\epsilon_t d\pi^\epsilon_t)$ converges for a.e. $t\in\C$.

We notice that the supports of convolutions, for $\epsilon$ small enough are not significantly distant from the supports of the original functions. 

Since $d(C_{\debar}(T))(fd\pi)+C_{\debar}(dT)(fd\pi)$ is equal to the integral of the expression in (\ref{eq_bdry}), multiplied by $(t-1)^{-1}$, we have 
$$d(C_{\debar}(T))(fd\pi)+C_{\debar}(dT)(fd\pi)=\int_{\C}\frac{\de\phi(t)}{\de\bar{t}}\frac{dt\wedge d\bar{t}}{t-1}=\phi(1)=T(fd\pi)\;,$$
because, as $\supp f$ is bounded and $\supp T$ has positive distance from $0$, the function $\phi(t)$ is compactly supported. $\Box$

\begin{Cor}$C_{\debar}(T)$ belongs to $N_{k+1,\q}(E)$ and it is of bidimension $(0,k+1)$.\end{Cor}
\noindent{\bf Proof: } Employing the previous Proposition, we repeat the proof of Proposition 10.2 in \cite{ambrosio1}, obtaining the continuity of $C_{\debar}(T)$ as a quasi-local metric current. Moreover, an easy calculation shows that, given $fd\pi\in E^{k}_\q(E)$ and $h\in\Ol(E)$, we have
$$C_{\debar}(T)(fdh\wedge d\pi)=0$$
so $C_{\debar}(T)$ is of bidimension $(0,k+1)$.$\Box$

\begin{Cor}We have $\debar C_{\debar}(T)=C_{\debar}(\debar T)+T$.\end{Cor}
\noindent{\bf Proof: } For a $(0,q)-$current $S$, $dS=\debar S$.$\Box$

In view of the previous results we have the following quasi-local solution for the Cauchy-Riemann equation.

\begin{Teo}\label{teo_0q} Let $T\in N_{k}(E)$ be a $(0,k)-$current, with $\debar T=0$, such that $\supp T$ is bounded with a positive distance from $0$. Then there exists a quasi-local metric $(0,k+1)-$current $U$ such that $\debar U=T$.\end{Teo}
\noindent{\bf Proof: } If $dT=\debar T=0$, then, letting $U=C_{\debar}(T)$, we have
$$\debar U=C_{\debar}(\debar T)+T=T$$
as quasi-local currents. $\Box$

\begin{Rem}We can control the mass of the solution $U$ on a ball by the mass of $T$ on that ball, with a constant depending only on the support of $T$, the dimension of $T$ and the radius of the ball.\end{Rem}

This result can at once be extended to $(p,q)-$currents, if we add the request that $\debar T=\de T=0$.

\begin{Cor} Let $T\in N_{k}(E)$ be a $(p,q)-$current, with $\debar T=\de T=0$, such that $\supp T$ is bounded with a positive distance from $0$; then there exists a quasi-local metric $(p,q+1)-$current $U$ such that $\debar U=T$.\end{Cor}
\noindent{\bf Proof: } We perform the same cone construction and we note that, by Proposition \ref{prp_d_cone}, $dC_{\debar}(T)=T$. Now, if $fd\pi\in E^{k}_\q(E)$ has $p+1$ holomorphic differentials, say $\pi_1,\ldots, \pi_{p+1}$, and $h\in \Ol(E)$, then 
$$T\left(f_t\frac{\de h_t}{\de\bar{t}}d\pi\right)=T(0d\pi)=0$$
$$T\left(f_t\frac{\de \pi_{jt}}{\de\bar{t}}d\widehat{\pi}_{jt}\right)=T(0d\widehat{\pi}_{jt})=0\qquad\textrm{ for } 1\leq j\leq p+1$$
and $d\widehat{\pi}_{jt}$, for $j\geq p+2$, contains $p+1$ holomorphic differentials, so $T(\sigma d\widehat{\pi}_{jt})=0$ for any $\sigma$ with bounded support. This means that $C_{\debar}(T)$ is a $(p,q+1)-$current; since $dC_{\debar}(T)=T$ is of bidimension $(p,q)$, we conclude that $dC_{\debar}(T)=\debar C_{\debar}(T)$.$\Box$

With some more effort, we can obtain the general result for $(p,q)-$currents.

\begin{Teo}Let $T\in N_k(E)$ be a $(p,q)-$current, with $\debar T=0$, such that $\supp T$ is bounded with a positive distance from $0$; then there exists a quasi local metric $(p,q+1)-$current $U$ such that $\debar U=T$. \end{Teo}
\noindent{\bf Proof: } We remark that, as $T$ is normal and $\debar T=0$, then $dT$ admits a Dolbeault decomposition, where $(dT)_{p-1,q}=\de T= dT$. Let $h_1,\ldots, h_p\in\Ol(E)$ be holomorphic functions and set $H=(h_1,\ldots, h_p)$; then $$S_H=T\llcorner(1,h_1,\ldots, h_p)$$
is a $(0,q)-$current such that 
$$\debar S_H=dS_H=(dT)\llcorner(1,h_1,\ldots, h_p)=(dT)_{p-1,q}\llcorner(1,h_1,\ldots, h_p)$$
and the last term is $0$ by the definition of $(p-1,q)-$current. Therefore $\debar S_H=0$.

Now, by virtue of Theorem \ref{teo_0q}, there exists a $(0,q+1)-$current $V_H$ such that $\debar V_H=S_H$;  for each $(f,H,\pi)\in\E^{k+1}_{\q}(E)$ with $H$ a $p-$tuple of holomorphic functions, we define the metric functional
$$U(f,H,\pi)=V_H(f,\pi)$$
and we set
$$U(f,\eta)=0$$
whenever $\eta$ contains at most $p-1$ holomorphic functions and at least $q+2$ antiholomorhpic functions.

It is easy to check that $U$ is then defined for all (quasi-local) $k+1-$form with either holomorphic or antiholomorphic differentials; this allows us to extend $U$ as a multilinear, local, alternating functional on the (quasi-local) $k+1-$forms with analytic coefficients.

Whenever $\supp f\subset B(0,R)$, we have
$$|U(f,H,\pi)|=|V_H(f,\pi)|\leq\prod_{j=1}^{q+1}\Lip(\pi_j)C(R)\int_{B(0,R)}|f|d\|S_H\|\leq$$
$$\leq \prod_{j=1}^p\Lip(h_j)\prod_{j=1}^q\Lip(\pi_j)C(R)\int_{B(0,R)}|f|d\|T\|\;.$$
Thus $\|U\|_{B(0,R)}\leq C(R)\|T\|_{B(0,R)}$, which means that the mass of $U$ is quasi-locally finite, wherever $U$ is defined.

Moreover, let $(f,\eta)$ be a $k-$form with each of $\eta_1,\ldots,\eta_k$ holomorphic or antiholomorphic and
$$f=\chi_{B(0,R)}\cdot g$$
for some $r<R$ and $g$ holomorphic or antiholomorphic. Then
$$dU(f,\eta)=U(\chi_{B(0,R)}, g,\eta)\;.$$
If $(g,\eta)$ contains more that $p$ holomorphic functions or more than $q+1$ antiholomorphic functions, then $dU(f,\eta)=0$.

If $\eta=(H,\pi)$, with $H=(h_1,\ldots, h_p)$ holomorphic, then
$$|dU(f,\eta)|=|U(\chi_{B(0,R)}, g,H,\pi)|=|V_H(\chi_{B(0,R)}, g,\pi)|=|(dV_H)(f,\pi)|$$
$$=|S_H(f,\pi)|=|T(f,H,\pi)|=|T(f,\eta)|$$
so, in this case $\|dU\|\leq\|T\|$.

If $\eta$ contains only $p-1$ holomorphic functions, then $g$ has to be holomorphic (otherwise $U(\chi_{B(0,R)},g,\eta)=0$); so we set $\eta=(\pi',H')$ and
$$|U(\chi_{B(0,R)},g,\pi',H')|=|V_{gH'}(\chi_{B(0,R)}, \pi')|$$
$$=|(dV_{gH'})(\pi'_1\chi_{B(0,R)}, \pi'_2,\ldots,\pi'_{q+1})|$$
$$=|S_{gH'}(\pi_1'\chi_{B(0,R)}, \pi_2',\ldots,\pi_{q+1}')|=|T(\pi'_1\chi_{B(0,R)}, g, H', \widehat{\pi'}_1)|$$
$$=|dT(f\pi_1', H', \widehat{\pi'}_1)-T(f, \pi_1',H', \widehat{\pi'}_1)|=|dT(f\pi_1', H', \widehat{\pi'}_1)|$$
where we have employed the definition of boundary and the definition of $(p,q)-$current. We can suppose, without loss of generality, that $\pi'_1(0)=0$, so that we get
$$\|dU\|\leq R\|dT\|\;.$$

\medskip

We infer $\|dU\|\leq \max\{\|T\|, R\|dT\|\}$, on the space of forms $(f,\eta)$ where every component is either holomorphic or antiholomorphic on $B(0,R)\supseteq \supp f$.

This allows us to extend $U$ as a quasi-local current on forms $(f,\eta)\in\E^{k+1}_{\q}(E)$ with $f=g\chi_{B(0,R)}$ for some $R$ and $g$ holomorphic or antiholomorphic; by linearity in the first component, we can allow $g$ to be analytic, then by density (and quasi-local finiteness of mass) we can extend $U$ to $\E^k_{\q}(E)$.

By the previous computations, $U$ is then a quasi-locally normal current, of bidimension $(p,q+1)$, with $\debar U=T$; moreover, the mass of $U$ is controlled, on a fixed ball, by the mass of $T$ and we have $\|U\|_{B(0,R)}+\|dU\|_{B(0,R)}\leq A(R)\|T\|+B(R)\|dT\|$.
$\Box$

\begin{Rem} The hypothesis that $0$ has positive distance from the support of $T$ can be avoided, by constructing the cone from a point different from the origin, as long as the support of $T$ is bounded.\end{Rem}

Boundedness of $\supp T$ seems much harder to get rid of and, to date, we do not even know if it is possible. In the same way, it is not apparent that one can improve the estimates on mass in order to obtain a metric current (not a quasi-local one) from the cone construction. In the one-dimensional case, this cone construction consists in the convolution with the Cauchy kernel, which does not in general give a compactly supported solution.

\medskip

Two natural questions arise:
\begin{enumerate}
\item are there conditions on $T$ ensuring that the solution to $\debar U=T$  obtained with the cone construction  has bounded support? or finite mass?
\item can the alleged solution with bounded support be obtained as a minimizer for the mass or the quasi-local mass among all quasi-local solutions to the Cauchy-Riemann equation?
\end{enumerate}

\section{Analytic sets in Hilbert spaces}

There are many possible definitions for a finite dimensional analytic set in an infinite-dimensional space. Here we adopt the following (see \cite{ruget1}). Let $H$ be a complex Hilbert space (or more generally a Banach space): a closed set $A\subset H$ will be called a {\em finite-dimensional analytic set} in $H$ if, locally in $H$, $A$ is an analytic subspace of a complex submanifold of $H$ of finite dimension .

\begin{Rem} This definition is equivalent to the one given by Douady in \cite{douady1} (see \cite{ruget1} for the details).\end{Rem}

Another possible (and equivalent) definition is the following (see \cite{aurich1}): $A\subset H$ is a finite-dimensional analytic set if for each $v\in H$ there exist a neighborhood $U$, another complex Hilbert space $H'$ and a holomorphic map $F:U\to H'$ whose differential has finite dimensional kernel such that $A\cap U=F^{-1}(0)$. This follows easily from the Implicit Function Theorem; this turns out to be equivalent to asking for a map whose differential has finite dimensional kernel and cokernel (i.e. a Fredholm map) such that $A\cap U=F^{-1}(0)$. 

\medskip

We recall a result from \cite{ruget1}.

\begin{Teo}\label{teo_hilb_img} Let $X$ be an analytic subset of finite dimension, $W$ a hilbertian manifold and $f:X\to W$ a proper holomorphic map. Then, $f(X)$ is a finite-dimensional analytic subset of $W$.\end{Teo}

\noindent{\bf Example} The properness assumption cannot be dropped as shown by the following example. Let us consider the space $\ell^2$ of square-summable sequences of complex numbers and consider the holomorphic map $g:\mathbb{D}^2\to\ell^2$ given by
$$g(z,w)=\{zw^n\}_{n\in\N}\;.$$
The preimage of $\{0\}$ is $\{z=0\}$, which is not compact. Let $X=g(\DD^2)$ and assume, by contradiction, that there exists a holomorphic function $\Phi:\ell^2\to\C$ vanishing on $X$; then $\Phi\circ g=0$ and consequently 
$$0=\Phi\circ g(z,w)=\sum_{n,m\geq0}\frac{\de^{n+m} \Phi\circ g}{\de z^n\de w^m}(0,0)\frac{z^nw^m}{n!m!}\;.$$
One has
$$0=\frac{\de \Phi\circ g}{\de z}(0,0)=\frac{\de\Phi}{\de e_0}(0)$$
$$0=\frac{\de^2\Phi\circ g}{\de w\de z}(0,0)=\frac{\de}{\de w}\left(\sum_j\frac{\de \Phi}{\de e_j}(g)\frac{\de g_j}{\de z}\right)(0,0)=$$
$$\left(\sum_{j,k}\frac{\de^2 \Phi}{\de e_k\de e_j}(g)\frac{\de g_k}{\de w}\frac{\de g_j}{\de z} + \sum_j\frac{\de\Phi}{\de e_j}\frac{\de^2 g_{j}}{\de w\de z}\right)(0,0)=\frac{\de\Phi}{\de e_1}(0,0)\;.$$
since
$$\frac{\de g_k}{\de w}(0,0)=0\qquad \forall k$$
and
$$\frac{\de^2g_k}{\de z\de w}(0,0)\neq 0 \Leftrightarrow g_k(z,w)=zw \Leftrightarrow k=1\;.$$
Proceeding in this way, we can show that $\de\Phi/\de e_j=0$ in $0\in \ell^2$, for all $j\in\N$. Therefore all the derivatives of $\Phi$ vanish at the origin, which means that no regular hypersurface of $\ell^2$ can contain a neighborhood of $0\in X$.

\bigskip

The major advantage of the given definition is that we can recover the local properties of finite-dimensional analytic sets in an infinite dimensional space from the usual local results on analytic sets in a complex manifold. In particular, let $X$ be a finite-dimensional analytic set in $H$, then
\begin{enumerate}
\item $X$ admits a local decomposition in finitely many irreducible components;
\item such a decomposition is given by the closure of the connected components of the regular part of $X$;
\item for every $x\in X$ there exist a finite-dimensional subspace $L$ of $H$ and an orthogonal projection $\pi:H\to L$ which realizes a neighborhood of $x\in X$ as a finite covering on $L$;
\item $X$ is locally connected by analytic discs;
\item if $X$ is irreducible, every nonconstant holomorphic function is open;
\item if $X$ is irreducible, the maximum principle holds;
\item if $X$ is compact and $H$ is holomorphically separable, then $X$ is finite.
\end{enumerate}

\bigskip

The behaviour of the analytic sets in a Banach space can vary wildely, depending on the properties of the space. We give here some examples.

\medskip

\noindent{\bf Example } Let $c_0$ be the vector space of sequences of complex numbers vanishing at infinity, i.e. $\{a_n\}\subset\C$ such that $\lim_{n\to\infty}a_n=0$; $c_0$ is a Banach space with the supremum norm. We consider the holomorphic map $f:\DD\to c_0$ given by
$$f(z)=\{z^n\}_{n\in\N}\;;$$ $f$ is a regular and injective holomorphic map; its image is contained in the unit ball of $c_0$ and if $\{z_k\}$ is a sequence converging to $b\DD$,  $\{f(z_k)\}_k$ does not converge, therefore $f$ is proper. Thus, $f(\DD)$ is a complex manifold of dimension $1$ in $E$, which is bounded.

\medskip

\noindent{\bf Example } Generalizing the previous example, we consider the Banach space of $p-$summable sequences of complex numbers $\ell^p$ and the holomorphic map $F:\DD^k\to\ell^p$ given by
$$F(z_1,\ldots, z_k)=\{z^I\}_{I}$$
where $I$ varies through all the multi-indeces of length $k$. We have that
$$|z^I|\leq (\max\{|z_1|,\ldots, |z_k|\})^{|I|}$$
and the number of multiindexes $I$ with $|I|=i_1+\ldots+i_k=m$ is
$${m+k-1\choose k-1}$$
which is less than $(2m)^k$ if $m$ is large enough. Therefore
$$\sum_I |z^I|^p\leq\sum_{m} (2m)^k(\max\{|z_1|,\ldots, |z_k|\})^m$$
which converges for $\max\{|z_1|,\ldots, |z_k|\}<1$.

Again, the map $F$ is regular, injective and proper with unbounded image $F(\DD^k)$. We observe that $F(\DD^k)$ provides an example of finite dimensional manifold not contained in any finite-dimensional linear subspace.

\subsection{Positive currents}

The definition of \emph{positive} current given in \ref{ssc_pos_cur} applies to currents in infinite dimensional spaces as well. We say that a current is \emph{finitely positive} if every finite dimensional projection of it is positive.

\begin{Prp}Let $\Omega\subset H$ be an open set. $T\in M_{2p}(\Omega)$ is positive if and only if it is finitely so.\end{Prp}
\noindent{\bf Proof: } Obviously, if $T$ is positive then every complex linear pushforward of it is positive. On the other hand, if $p_m:H\to\C^m$ is the projection on the first $m$ coordinates, then, by Proposition \ref{prp_limlip}, $f\circ p_m\to f$ pointwise and $\Lip(f\circ p_m)\leq\Lip(f)$.

Therefore, 
$$T(f,\pi_1,\ldots,\overline{\pi}_p)=\lim_{m\to\infty}T(f\circ p_m, \pi_1\circ p_m,\ldots, \overline{\pi}_p\circ p_m)=$$
$$\lim_{m\to\infty}(p_m)_\sharp T(f\vert_{p_m(H)},\pi_1\vert_{p_m(H)},\ldots, \overline{\pi}_p\vert_{p_m(H)})\geq0\;,$$
which is the thesis. $\Box$

\medskip

On an infinite dimensional complex space, we can define plurisubharmonic functions as follows (see also \cite{mujica1} for a discussion of the properties of such functions).

Let $\Omega\subset H$ be an open set and let $u:\Omega\to\R\cup\{-\infty\}$ be an upper semicontinuous function (not identically equal to $-\infty$) such that
$$u(a)\leq\int_{0}^{2\pi}u(a+e^{i\theta}b)d\theta$$
for every $a\in\Omega$ and every $b\in H$ such that $a+\lambda b\in\Omega$ for every $\lambda\in\C$ with $|\lambda|\leq1$.

Define $d^c=i(\de-\debar)$.

\begin{Prp}Let $T\in M_{2p}(\Omega)$ be a positive closed current with bounded support and $u:\Omega\to\R\cup\{-\infty\}$ a bounded plurisubharmonic function. Then $dd^c(T\llcorner u)$ is a closed, positive metric current with bounded support and the following estimate holds:
$$M(dd^c(T\llcorner u))\leq \|u\|_\infty M(T)\;.$$
\end{Prp}
\noindent{\bf Proof: } We note that the result is true for any finite-dimensional projection of $T$. Namely, if $p_m$ is as above, 
$$(p_m)_\sharp dd^c(T\llcorner u)=dd^c((p_m)_\sharp T\llcorner (u\circ p_m))=dd^c(T_m\llcorner u_m)$$
with $T_m$ a positive closed current with compact support in $\C^m$ and $u_m$ a bounded plurisubharmonic function. Then we know that $$dd^c(T_m\llcorner u_m)=T_m\llcorner (dd^c u_m)$$
in the sense of distributions and for every compact $K$  we have
$$M_K(dd^c(T_m\llcorner u_m))\leq C \|u_m\|_{\infty, K}M_K(T_m)\;.$$
By \cite[Theorem 1.4]{ambrosio3}, there exists a subsequence $T_{m_j}$ which converges in mass to $T$, therefore $M_B(T_m)\to M_B(T)$ for every bounded set $B$ in $\Omega$.

This means that, for $j$ big enough,
$$M_K(dd^c(T_{m_j}\llcorner u_{m_j}))\leq C\|u\|_{\infty, K}M_K(T)$$
so, again by \cite[Theorem 1.4]{ambrosio3}, we can find a subsequence converging to some $S\in M_{2p-2}(\Omega)$.

Such an $S$ is such that $(\pi_m)_\sharp S=T_m$ for infinitely many $m$, therefore $S$ coincides, as a metric functional, with $dd^c(T\llcorner u)$. Moreover, $S$ is positive and closed and
$$M_K(S)\leq \|u\|_{\infty,K}M_K(T)\;.$$
We note also that $\supp S\subseteq \supp T$. $\Box$

\medskip

In view of the previous Proposition, we will denote by $dd^cu\wedge T$ the current $dd^c(T\llcorner u)$.

Proceeding by induction, we can give a meaning to the writing
$$dd^cu_1\wedge\ldots\wedge dd^c u_p\wedge T$$
for a current $T$ in the hypotheses of the previous Proposition. Such a definition allows us to write an analogue of the Monge-Amp\`ere operator in Hilbert spaces.

\subsection{Currents of integration on analytic sets}

Let $V$ be a finite-dimensional analytic set in some open domain $U\subset H$; since, by definition, $V$ is locally contained in some finite-dimensional submanifold, we now that it is locally of finite volume. Therefore, we can define the current $[V]$ of integration on the regular part of $V$; this will be, in general, a quasi-local metric current. However, if $B$ is a bounded, closed subset of $U$, the current $[V\llcorner B]$ will be a rectifiable metric current. The following results gives an estimate for the mass of such a current, analogous to Wirtinger formula in the finite dimensional case.

\begin{Prp}\label{prp_vol_in_H}Let $H$ be a Hilbert space, with scalar product $\langle\cdot,\cdot\rangle$, $V$  an analytic set in an open set $U\subset H$, with $\dim_\C V_\rg=p$. Let $\Omega$ be a ball in $U$ and let $[V]$ be the current of integration associated to $V\cap \Omega$ in $\Omega$. Then
$$M([V])\leq\lim_{n\to\infty}\sum_{1\leq i_1<\ldots<i_p\leq n}[V]\left(\frac{i^p}{2^pp!}, z_{i_1},\bar{z}_{i_1},\ldots, z_{i_p},\bar{z}_{i_p}\right)<+\infty\;,$$
where $\{z_j\}_{j\in \N}$ are coordinate functions with respect to some orthonormal basis.
\end{Prp}
\noindent{\bf Proof: }
Given an orthonormal basis $\{e_n\}_{n\in\N}$, let $E_m=\mathrm{span}\{e_1,\ldots, e_m\}$ and $\pi_m:H\to E_m$ the orthogonal projection. We denote $[V]_m=(\pi_m)_\sharp [V]$ and observe that $[V]_m\to [V]$ weakly, so that by the semicontinuity of the mass we get
$$M([V])\leq \liminf_{m\to\infty}M([V]_m)\;,$$
where the masses are relative to $\Omega$ or to $\Omega_m=\Omega\cap E_m$ respectively.

On the other hand, the projections $\pi_m$ have norm $\|\pi_m\|\leq 1$, so, if $\mu$ is the mass measure of $[V]$, we have
$$|[V]_m(f,\eta)|\leq \prod \Lip(\eta_j)\int_{U_m}|f|d(\pi_m)_\sharp \mu\;.$$
This implies that the mass measure $\mu_m$ of $E_m$ is dominated by $(\pi_m)_\sharp \mu$, therefore
$$\mu_m\leq (\pi_m)_\sharp \mu\xrightarrow[m\to\infty]{}\mu$$
which means that 
$$M([V])=\mu(\Omega)\leq\liminf_{m\to\infty}\mu_m(\Omega)\leq \lim_{m\to\infty}(\pi_m)_\sharp \mu(\Omega)=\mu(\Omega)=M([V])\;.$$
Now, $E_m$ with the induced scalar product is the usual complex hermitian space $\C^m$ and the pushforward of an analytic chain is again an analytic chain. Therefore
$$M([V]_m)\leq[V]_m(\omega_m^p/p!)=\H^{2p}(V_m)\leq C_p M([V]_m)$$
where
$$\omega_m=\frac{i}{2}\sum_{j=1}^mdz_j\wedge d\bar{z}_j=\frac{i}{2}\sum_{j=1}^me^*_j\wedge \bar{e}^*_j$$
and $C_p$ is a constant depending only on $p$ (see \cite{ambrosio1}, after Remark 8.4).
Noticing that $[V]_m(\omega_m^p/p!)=[V](\omega_m^p/p!)$, we obtain
$$M([V])\leq \lim_{n\to\infty}\sum_{1\leq i_1<\ldots<i_p\leq n}[V]\left(\frac{i^p}{2^pp!}, z_{i_1},\bar{z}_{i_1},\ldots, z_{i_p},\bar{z}_{i_p}\right)\leq C_pM([V])<+\infty\;,$$
which is the thesis. $\Box$

\medskip

\begin{Rem}The current $[V]$ is obviously positive, of bidimension $(k,k)$ for some $k$ and its boundary is supported outside $\Omega$.\end{Rem}

Given an orthonormal basis $\{e_j\}_{j\in\N}$ and a multiindex $I=(i_1,\ldots, i_k)$, we denote by $\pi_I$ the orthogonal projection from $H$ onto $\mathrm{Span}\{e_{i_1},\ldots, e_{i_k}\}$.

\begin{Teo}Let $\Omega\subset H$ be a ball, $S$ be a rectifiable current in $\Omega$. Suppose that
\begin{enumerate}
\item $\supp dS\cap\Omega=\emptyset$;
\item $S$ is a $(k,k)$ positive current.
\end{enumerate}
Then $S$ can be represented as a sum with integer coefficients of integrations on the regular parts of analytic sets. \end{Teo}

The absence of the specification \emph{finite-dimensional} is not a misprint (see the Remark following the proof).

\medskip

\noindent{\bf Proof: } As $S$ is a metric current, we can define its pushforward through any Lipschitz map.  We note that $(\pi_I)_\sharp S=m_I[V_I]$ with $V_I=\pi_I(\Omega)$. By \emph{(ii)} we know that $m\geq0$ and by the fact that $dS\llcorner\Omega=0$, we deduce that $S$ is integral in $\Omega$, therefore $m\in\N$.

By Theorem \ref{teo_ext_slice}, for almost every $y\in V_I$ we can define $\langle S,\pi_I,y\rangle$; moreover, we can find a $\H^{2k}-$rectifiable subset $B$ of $\supp S$ and an integer multiplicity function $\theta(x)$ such that $S=[B]\llcorner\theta$; then
$$\langle S,\pi_I,y\rangle=\sum_{x\in\pi_I^{-1}(y)\cap B}\theta(x)[x]$$
and
$$\sum_{x\in\pi_I^{-1}(y)\cap B}\theta(x)=\langle S,\pi_I, y\rangle(1)=m\;.$$

Let us call $G_I\subset V_I$ the set of $y$ such that the slice exists; then for $j\not\in I$ and $z\in G$, we set
$$P_j(z,W)=\prod_{x\in\pi_I^{-1}(z)\cap B}(W-w_j(x))^{\theta(x)}$$
where $w_j(x)$ is the $j-$th coordinate of $x$ in the fixed orthonormal basis.

We note that
$$\sum_{x\in\pi_I^{-1}(z)\cap B}\theta(x)w_j(x)^s=\langle S,\pi_I, z\rangle(w_j^s)$$
is a holomorphic function of $z$, because $\debar S=0$, therefore by a classical argument $P_j(z,W)$ is a polynomial with coefficients in $\Ol(V_I)$ for every $j\not\in I$.

After removing a $\H^{2k}-$negligible set from $G_I$, we have that $P_j(z,w_j)=0$ for every $j\not\in I$ and every $x=(z,w)\in \pi_I^{-1}(G_I)\cap B$.

Let us define
$$X_I=\{P_j(z,w_j)=0\ ,\ j\not\in I\}\qquad X=\bigcup_I X_I\;.$$

We can look at $X_I$ as the zero locus of the map 
$$P_I:H\to \mathrm{Span}_\C\{e_j\ :\ j\not\in I\}=H_1$$
given by
$$P_I(z,w)=\sum_{j\not\in I}e_jP_j(z,w_j)\;.$$

In order to show that $P_I$ is well defined, we observe that, since $P_j(z,W)$ is a polynomial in $W$, for a fixed $z$ we have 
$$|P_j(z,W)|\leq \min\{d(W, w_j(x))^m,\ x\in\pi_I^{-1}(z)\cap B\}$$
with $m\geq 1$, if $W$ is close enough to some $w_j(x)$. Therefore, for $p=(z,w)$ in a neighborhood of $B$, we can write
$$\sum_{j\not\in I}|P_j(z,w_j)|^2\leq\sum_{j}|w_j-w_j(x)|^{2m}\leq\|p-x\|^{2m}<+\infty$$
where $x$ is the nearest point in $B$ to $p=(z,w)$. The map $P_I$ is obviously holomorphic, as its entries are polynomials in $w_j$ with coefficients holomorphic in the first $k$ coordinates.

This shows that $X_I$ is locally given as the zero locus of a holomorphic map between Hilbert space, therefore it is an analytic set. Moreover, let $x\in X_I$ be a smooth point and suppose that $\dim_\C T_x X_I>k$. By construction, $\pi_I\vert_{T_xX_I}:T_xX_I\to\C^k$ has maximum rank, therefore we can find $m$ such that, setting $J=I\cup\{m\}$, the projection $\pi_J$ restricted to $T_xX_I$ is surjective onto $\C^{k+1}$. This means that $P_m(z,W)\equiv 0$, but this is impossible. So, $\dim_\C T_xX_I=k$, i.e. the regular part of $X_I$ is a smooth $k-$dimensional complex manifold.

Now, let
$$B_I=\{x\in B\ :\ J_{2k}\pi_I(x)\neq 0\}$$
i.e. the set of points $x$ of $B$ such that $D\pi_I$ has rank $2k$ on the approximate tangent to $B$ at $x$; define also
$$C_I=(B\cap V) \setminus \pi_I^{-1}(G_I)\;.$$
We have
$$\int_{B_I\cap C_I} J_{2k}\pi_I(x)d\H^{2k}(x)=\int_{V_I\setminus G_I}\left(\int_{\pi_I^{-1}(y)\cap B} gdH^0\right)d\H^{2k}(y)=0$$
where $g$ is the characteristic function of $B_I\cap C_I$. Since $J_{2k}\pi_I>0$ on $B_I\cap C_I$, this means that $H^{2k}(B_I\cap C_I)=0$.

Obviously, $B=\bigcup B_I$ and
$$B\cap\left(\bigcup_I \pi_I^{-1}(G_I)\right)\subset X\;,$$
but
$$B\setminus \bigcup_I \pi_I^{-1}(G_I)=\bigcap_{I}(B\setminus \pi_I^{-1}(G_I))\subset \bigcap_{I}((B\setminus B_I)\cup(B_I\cap C_I))\subseteq\bigcup_I(B_I\cap C_I)=D\;.$$
Since $\H^{2k}(B_I\cap C_I)=0$, we also have $\H^{2k}(D)=0$ and $\H^{2k}(B\setminus X)=0$. Moreover, as $X$ is closed in $\Omega$, $\supp S\subset X$; therefore $S=[B\cap X]$.

\medskip

If we denote by $X_\rg$ the union of the regular parts of $X_I$, we have that $S\llcorner X_\rg$ is a $(k,k)-$current, positive and closed,  with support on a $k-$dimensional smooth complex manifold. Therefore, $S\llcorner X_\rg$ can be written as a series with integer coefficients of the currents of integration on the connected components of $X_\rg$.

There exists $r>0$ such that $\pi_I(V)$ contains a ball of radius $r$ for every $I$; therefore, the $\H^{2k}-$measure of the regular part of $X_I$ is uniformly bounded from below independetly of $I$. On the other hand $S\llcorner X_\rg$ is of finite mass; therefore it has to be a finite sum.

Finally, let us consider the rectifiable set $R=B\setminus X_\rg$. If we project it on the first $m$ coordinates, for $m\geq k+1$, we obtain that its image is the singular set of a $k-$dimensional analytic space, therefore $\H^{2k}-$negligible; by \cite[Theorem 8.2]{ambrosio2}
$$\int_{R}\mathbf{J}_{2k}({d^R\pi_m}_x)d\H^{2k}(x)=\int_{\C^m}\sharp\{x\in R\cap \pi_m^{-1}(y)\}d\H^{2k}$$
with $\pi_{m}:H\to\mathrm{Span}\{e_1,\ldots, e_m\}$ the orthogonal projection. Let us denote by $\eta(x)$ the approximate tangent to $R$ in $x$; then the $2k-$jacobian of $\pi_m$ on $\mathrm{Tan}^{(2k)}(R,x)$ is given by the projection of $\eta(x)$ on $\mathrm{Span}\{e_1,\ldots, e_m\}$.

We define $A_m=\{ x\in R\ :\ \mathbf{J}_{2k}({d^R\pi_m}_x)>0\}$ and we note that $A_{k+1}\cup A_{k+2}\cup\ldots=R$, up to $\H^{2k}-$neglibigle sets. But
$$\int\limits_{A_m}\mathbf{J}_{2k}({d^R\pi_m}_x)d\H^{2k}(x)=\int\limits_{\C^m}\sharp\{x\in A_m\cap \pi_m^{-1}(y)\}d\H^{2k}=\int\limits_{\pi_m(A_m)}\sharp\{x\in R\cap \pi_m^{-1}(y)\}d\H^{2k}=0$$
because $\pi_m(A_m)\subseteq \pi_m(R)$, which is $\H^{2k}-$negligible. Therefore $\H^{2k}(R)=0$, so $S\llcorner X_\rg=S$ and this concludes the proof. $\Box$

\medskip

\begin{Rem}We cannot show that $X$ is a finite-dimensional analytic space in the sense precised in the beginning of this section; indeed, the example discussed before, the map $f(z,w)=(zw^n)_{n}$, gives an analytic space which carries a current of integration which satisfies the hypotheses of the previous theorem but cannot be written as the integration on a finite-dimensional analytic space.
\end{Rem}

\section{Boundaries of holomorphic chains}

\begin{Prp}\label{prp_emb}Let $M$ be a compact $\Ci^1$ submanifold of a complex reflexive Banach space $E$ with complex structure $J$, such that $\dim_\R M=2p-1$ and $\dim_\C T_zM\cap JT_zM=p-1$ for every $z\in M$. Then there exists a complex linear map $F:E\to\C^n$ for some $n>0$ which restrict to an embedding of $M$ into $\C^n$.\end{Prp}
\noindent{\bf Proof: }Given $z\in M$, let $l_{1,z},\ \ldots, l_{p,z}$ linearly independent elements of $E^*$ such that
$$\ker l_{1,z}\cap\ldots\cap\ker l_{p-1,z}\cap T_{z}M\cap JT_zM=\{0\}$$
and
$$\ker l_{1,z}\cap\ldots\cap\ker l_{p-1,z}\cap \ker\mathsf{Re}\ l_{p,z}\cap T_{z}M=\{0\}\;;$$
both these conditions are open. By compactness, we can find finitely many $l_1,\ldots, l_N$ such that
for every point $z\in M$ there exists indexes $j_1<\ldots<j_p$ such that
$$\ker l_{j_1}\cap\ldots\cap\ker l_{j_{p-1}}\cap T_{z}M\cap JT_zM=\{0\}$$
and (viewing $E$ as a real vector space)
$$\ker l_{j_1}\cap\ldots\cap\ker l_{j_{p-1}}\cap \ker\mathsf{Re}\ l_{j_p}\cap T_{z}M=\{0\}\;.$$
This means that if we define $L:E\to\C^N$ by $L=(l_1,\ldots, l_N)$, we have that the differential $dL$ is always of real rank $2p-1$ on $M$ and it is complex linear on $T_zM\cap JT_zM$.

\medskip

Let $U_1,\ldots, U_h$ be the open sets and $l_1,\ldots, l_p, l_{p+1},\ldots, l_{2p},\ldots, l_{hp}$ be the maps constructed as above and $\{V_j\}_{j=1}^K$ a collection of open sets in $M$ such that for each $V_j$ there exists a $U_{\nu(j)}$ such that $V_j\Subset U_{\nu(j)}$ and $\bigcup V_j= M$.

For a fixed $j$, the set $L^{-1}(L(V_j))$ is a union of $\mu_j$ connected components which are relatively compact in some open sets $U_{k_1},\ldots, U_{k_{\mu_j}}$; therefore, there exist $\mu_j$ linear maps $f_j^1,\ldots, f_j^{\mu_j}$ such that for each connected component there is one map which separates it from the others, that is, a map which has different values on it and on the union of the others.

Now, consider the map $F=(l_1,\ldots, l_{hp},f_1^1,\ldots, f_K^{\mu_K})$. By the first part of the construction, $F$ has an injective differential on $M$; by the second part, it is globally injective on $M$. Therefore $F$ is a holomorphic embedding of $M$ into $\C^{n}$, where $n=hp+\mu_1+\ldots+\mu_K$, realized with a complex linear map. $\Box$

\bigskip

Let $M$ be a compact $\Ci^1$ submanifold of a reflexive complex Banach space $E$ with complex structure $J$, with $\dim_\R M=2p-1$. $M$ induces a metric current $[M]$ of dimension $2p-1$.

\begin{Prp} \label{prp_caratt_CR}The following are equivalent:
\begin{enumerate}
\item $\dim_\C (T_z M\cap (JT_z M))=p-1$ $\forall\ z\in M$;
\item $[M](\alpha)=0$ for every metric $(r,s)-$form $\alpha$ on $E$ with $r+s=2p-1$ and $|r-s|>1$;
\item $M$ is locally the graph of a CR-function: for every $z\in M$ there exists $U$ neighbourhood of $z$ in $E$ such that $M\cap U$ is the graph of a function $f:\widetilde{M}\to E'$, $\widetilde{M}$ a CR-submanifold of $\C^p$, $E'$ a closed subspace of $E$ such that $E=E'\oplus \C^p$ and $f$ a CR-function on $\widetilde{M}$.
\end{enumerate}
\end{Prp}
\noindent{\bf Proof:} \emph{1) $\Longrightarrow$ 2) } 
Let $\alpha=(f,g_1,\ldots, g_r, h_1,\ldots, h_s)$ and let $i:M\to E$ an embedding whose differential is complex linear when restricted to (the preimage of) $T_zM\cap JT_zM$; then $[M]=i_\sharp T$, with $T\in M_{2p-1}(M)$. By the comparison theorem for manifolds, $T$ is induced by a classical current $T'$ on $M$; but then, 
$$[M](\alpha)=T(f\circ i, g_1\circ i,\ldots, h_s\circ i)=T'(f\circ i d(g_1\circ i)\wedge \ldots\wedge d(h_s\circ i))\;.$$
The functions $g_j\circ i$ have complex linear differentials on $i_*^{-1}(TM\cap JTM)$, therefore if there are more than $p$ of them, their wedge product will vanish; the same holds for the differentials of the functions $\overline{h}_j\circ i$. So $[M](\alpha)=0$ if $|r-s|>1$.

\medskip

\emph{2) $\Longrightarrow$ 1) } Let $\rho:E\to\C^N$ be a finite-dimensional embedding for $M$, which is holomorphic on $E$. If 
$$\dim_\C T_{\rho(z)}\rho(M)\cap J_{\C^N}T_{\rho(z)}\rho(M)=\dim_\C T_zM\cap JT_zM<p-1$$
then there exists a metric $(r,s)-$form $\beta$ on $\C^N$ with $r+s=2p-1$ and $|r-s|>1$ such that
$$\int_{\rho(M)}\beta\neq0$$
so
$$\int_M\rho^*\beta\neq0$$
and $\rho^*\beta$ is a $(r,s)-$form on $E$ with $|r-s|>1$.

\medskip

\emph{3) $\Longrightarrow$ 1) } Let $f:\widetilde{M}\to E'$ be the given CR-function; we define $G:\widetilde{M}\to \widetilde{M}\times E'$ by $G(p)=(p,f(p))$.

Let $F:\C^p\oplus E'\to E$ be the given isomorphism; then $(F\circ G)_*T_p\widetilde{M}=T_{F(p,f(p))} M$ and, since $T_{p}\widetilde{M}$ contains a complex subspace of complex dimension $p-1$, so does the tangent space of $M$.

\medskip

\emph{1) $\Longrightarrow$ 3) } Let us fix $z\in M$ and let $H_z$ be the complex subspace of $T_zE$ of (complex) dimension $p$ containing $T_zM$. By reflexivity, $E=E^*$,  so we have a splitting of $E=H_z\oplus E'$ for some closed subspace $E'$. By construction, $\pi:E\to H_z$ is a local embedding when restricted to a neighbourhood $U$ of $z$ in $M$, because it has a maximum rank differential at $z$.

Let $\widetilde{M}$ be the image of $U$ trough $\pi$; we have the function $f:\widetilde{M}\to E'$ defined by $(p,f(p))\in U\cap \pi^{-1}(p)$. By construction,  $f_*\vert_{T_p\widetilde{M}\cap JT_p\widetilde{M})}$ is $\C-$linear, so $f$ is a CR-function and $U$ is its graph. $\Box$

\bigskip

Let $S$ be a $(2p-1)-$current with compact support in a complex manifold $X$; we say that $S$ is \emph{maximally complex} if $M_{r,s}=0$ for $|r-s|>1$.

\begin{Rem} $S_{r,s}$ in general won't be a metric current (see subsection 1.3 for an example). Nevertheless, the above definition makes sense for any current $S$, as we only require that the functional $M_{r,s}$ be zero for values $(r,s)$ such that $|r-s|>1$.\end{Rem}

\begin{Prp}\label{prp_slice_maxcmp}Let $M$ be a $(2p-1)-$current with compact support in $X$, $F:X\to Y$ a Lipschitz holomorphic map. Suppose that $M$ is maximally complex, then the same is true for $F_\sharp M$ and, if $p>\dim_\C Y$, for $\langle M, F, \zeta\rangle$, given that $M$ is flat and slices exist.\end{Prp}
\noindent{\bf Proof: }
Obviously, we have $(F_\sharp M)_{r,s}=F_\sharp(M_{r,s})$ (this is an equality between metric functionals only, not metric currents).

Moreover, if $\dim_\C Y<p$ and if $\langle M, F, \zeta\rangle$ exists for some $\zeta\in Y$, let $\{\rho_{\epsilon, \zeta}\}$ be a family of smooth approximations of $\delta_\zeta$. Then locally (with $\supp f$ contained in a manifold chart for $Y$)
$$\langle M, F, \zeta\rangle(f,\eta)=\lim_{\epsilon\to 0}M(f\rho_{\epsilon,\zeta}, F, \overline{F}, \eta)\;.$$
So, if $M_{r,s}=0$ for $|r-s|>1$ then also $\langle M, F, \zeta\rangle_{r-m, s-m}=0$ for $1<|r-s|=|(r-m)-(s-m)|$, with $m=\dim_\C Y$. $\Box$

\bigskip

A \emph{MC-cycle} in a complex Banach space $E$ is a maximally complex $(2p-1)-$dimensional closed metric current, with compact support.

The definition is meaningless for $p=1$; the notion of \emph{moment condition} which substitutes the maximal complexity for $1-$dimensional currents cannot be given that easily in a Banach space and it turns out to be not automatically satisfied by a maximally complex current of higher dimension. The philosophical reason is the greater distance, in Banach spaces, between local and global aspects.

The following Theorem follows easily from Proposition \ref{prp_slice_maxcmp} and from the slicing properties of rectifiable currents.

\begin{Teo} Let $M$ be a rectifiable MC-cycle of dimension $(2p-1)$ in a Banach space $E$ and consider an holomorphic Lipschitz map $F: E\to \C^m$. Then
\begin{enumerate}
\item $F_\sharp M$ is a rectifiable MC-cylce of dimension $(2p-1)$ in $\C^m$;
\item if $m<p-1$, $\langle M,F,\zeta\rangle$ is a rectifiable MC-cycle of dimension $2(p-m)-1$ in $E$.
\end{enumerate}
\end{Teo}

\begin{Rem}By Theorem \ref{teo_ext_slice}, the slice $\langle M, F, \zeta\rangle$ exists rectifiable for almost every $\zeta \in \C^m$.\end{Rem}

\begin{Teo} Let $M$ be a MC-cycle of dimension $(2p-1)$ in $E$. Then, for every linear projection $\pi:E\to \C^p$ and every $\phi\in\Ol(\overline{\supp M})$, we have
$$\debar[\pi_\sharp(\phi M)]^{0,1}=0\;.$$
Moreover, there is a unique integrable compactly supported function $c_\phi$ in $\C^p$ such that
$$\debar c_\phi=[\pi_\sharp(\phi M)]^{0,1}$$
and such a function can be obtained by convolution with the Cauchy kernel or the Bochner-Martinelli kernel.\end{Teo}
\noindent{\bf Proof: } We know that $dM=0$; since $M$ is maximally complex, we have
$$M=M_{p,p-1}+M_{p-1,p}$$
so
$$0=dM=(dM)_{p-2,p}+(dM)_{p-1,p-1}+(dM)_{p,p-2}\;.$$
In particular, this means that $(dM)_{p,p-2}=0$. Therefore
$$\debar[\pi_\sharp(\phi M)]^{0,1}=\debar[\pi_\sharp(\phi M)]_{p,p-1}=[d\pi_\sharp(\phi M)]_{p,p-2}=
[\pi_\sharp(d\phi M)]_{p,p-2}$$
$$=[\pi_\sharp(\phi dM)]_{p,p-2}+[\pi_\sharp(M\llcorner(1,\phi))]_{p,p-2}$$
but $M\llcorner(1,\phi)$ has non vanishing $(r,s)-$components only for $(r,s)=(p-1,p-1)$ or $(r,s)=(p-2,p-1)$, so $[\pi_\sharp(M\llcorner(1,\phi))]_{p,p-2}=0$ and then
$$\debar[\pi_\sharp(\phi M)]^{0,1}=[\pi_\sharp(\phi dM)]_{p,p-2}=\pi_\sharp(\phi (dM)_{p,p-2})=0\;.$$
We note that $\pi_\sharp(\phi M)$ is a metric current in $\C^p$, therefore it is also a classical one, consequently its component of bidegree $(0,1)$ is a classical current as well and by the previous computation is $\debar-$closed. By a standard convolution-contraction with either the Cauchy kernel or the Bochner-Martinelli kernel, we can fin a compactly supported integrable function $c_\phi$ as requested. $\Box$

\medskip

\begin{Teo}Let $M$ be a compact, oriented $(2p-1)-$manifold (without boundary) of class $\Ci^2$ embedded in $H$, and suppose that there exists an orthogonal decomposition $H=\C^p\oplus H'$ such that the projection $\pi:H\to\C^p$, when restricted to $M$, is an immersion with transverse self-intersections. Then if $M$ is an $MC-$cycle there exists a unique holomorphic $p-$chain $T$ in $H\setminus M$ with $\supp T\Subset H$ and finite mass, such that $dT=[M]$ in $H$.\end{Teo}
\noindent{\bf Proof: } Let $\mathfrak{m}=\pi(M)\subset\C^p$; for every $\lambda\in H'\setminus\{0\}$, we define $\pi^\lambda(z)=(\pi(z), \langle z,\lambda\rangle)\in\C^{p+1}$.

By the previous results, $M^{\lambda}$ satisfies the same hypotheses in $\C^{p+1}$, therefore by \cite[Theorem 6.1]{harvey1} we can solve the problem for $M^{\lambda}=\pi^{\lambda}(M)$, finding a holomorphic $p-$chain $T^{\lambda}$ in $\C^{p+1}\setminus M$ with the required properties. Following the proof of Theorem 6.1 in \cite{harvey1}, we write
$$\C^p\setminus\mathfrak{m}=U_0\cup U_1\cup\ldots\cup U_k$$
where the $U_j$ are connected components and $U_0$ is unbounded; $T^\lambda$ is locally on each $U_j$ union of graphs of holomorphic functions
$$F_j^{\lambda, h}:U_j\to\C\qquad h=1,\ldots, n_{\lambda,j}\;.$$

Given another $\lambda'\in H'\setminus\{0\}$, we can consider the $p-$chain $T^{\lambda'}$, which will be given by holomorphic functions
$$F_j^{\lambda', h}:U_j\to\C\qquad h=1,\ldots, n_{\lambda',j}\;;$$
however, we can also consider, in $\C^{p+2}$, the manifold $M^{\lambda,\lambda'}$ and the associated solution $T^{\lambda,\lambda'}$; denoting by $p$ and $p'$ the restrictions of $\pi^\lambda$ and $\pi^{\lambda'}$ to $\C^{p+2}$, we have
$$p_* T^{\lambda,\lambda'}=T^{\lambda}\qquad p'_*T^{\lambda,\lambda'}=T^{\lambda'}\;.$$

Since the differentials of $p,\ p'$ are of rank $2p-1$ on $M^{\lambda,\lambda'}$ and because $p$ and $p'$ are holomorphic, their differentials are at least of rank $2p$ on $M^{\lambda,\lambda'}$; this means that they are of rank $2p$ in a neighborhood of $M^{\lambda,\lambda'}$ in $M^{\lambda,\lambda'}\cup\supp T^{\lambda,\lambda'}$ (which is locally a $\Ci^2$ manifold with boundary by Lemma 6.8 in \cite{harvey1}, ), therefore $n_{\lambda,j}=n_{\lambda,\lambda',j}=n_{\lambda',j}$ for every $j$ and every $\lambda,\lambda'\in H'\setminus\{0\}$.

Let $\{\lambda_i\}_{i\in I}$ be an orthonormal basis for $H'$ and consider the holomorphic functions
$$F_j^{\lambda_i, h}:U_j\to \C\qquad j=1,\ldots, k,\ \ h=1,\ldots, n_j,\ \ i\in I$$
and define
$$F_j^h=\sum_{i\in I}\lambda_i F_j^{\lambda_i, h}\;.$$

\noindent{\emph{The function $F_j^h$ is well defined. }} For any finite subset of indices $J\subset I$, we can consider the projection $$p_J:H\to\C^p\oplus\mathrm{Span}\{\lambda_i\}_{i\in J}$$
and the pushforward $[M]_J=(p_J)_\sharp [M]$; the functions $\{F^{\lambda_i, h}_j\}_{i\in J}$ give a solution for the finite-dimensional problem with datum $[M]_J$, therefore $S_{J,j,h}=\sum_{i\in J}F^{\lambda_i, h}_j\lambda_i$ is a holomorphic function with values in a finite-dimensional vector space, such that
$$|S_{J,j,h}(z)|\leq R$$
where $R$ is such that $\supp [M]_J\subset \C^p\times B(0,R)$, $B(z,r)$ being the ball with center $z$ and radius $r$ in $\mathrm{Span}\{\lambda_i\}_{i\in J}$.

Now, let us take $I=\N$ and fix $\epsilon>0$. By compactness, we can find $I'\subset I$ finite and set
$$V_\epsilon=\C^p\oplus\mathrm{Span}\{\lambda_i\}_{i\in I'}$$
so that $d(M,V_\epsilon)<\epsilon$; let $H'_\epsilon$ be the topological complement of $V_\epsilon$ in $H$, then the projection of $M$ on $H'_\epsilon$ lies in a ball of radius $\epsilon$ around $0$. Now, for any finite subset $J\subset I$ such that $\min J>\max I'$, we have that 
$$|S_{J,j,h}(z)|\leq \epsilon\;,$$
showing that the sequence of maps from $U_j$ to $H'$
$$\left\{\sum_{i=0}^m F^{\lambda_i,h}_j(z)\lambda_i\right\}_{m\in I}$$
is a Cauchy sequence with respect to the supremum norm on $U_j$. Therefore the limit $F^h_j(z)$ is well defined and continuous on the closure of $U_j$, because every element of the sequence is.

%We know, by construction, that $M^{\lambda_i}$ is contained in the closure of
%$$(U_1\cup\ldots\cup U_k)\times D(0,R_i)$$
%for some $R_i$. By Proposition \ref{prp_caratt_CR}, $M$ is locally a graph of a $CR$ function of class $\Ci^2$, so for every $p\in M$ there exists $U$ a neighborhood of $p$ in $H$, $\widetilde{M}\subset\C^p$ a CR manifold of class $\Ci^2$ and $f:\widetilde{M}\to H'$ such that $M\cap U$ is the graph of $f$.
%
%Since $f$ is at least of class $\Ci^2$, it can be extended to a tubular neighborhood of $\widetilde{M}$ and then to all $\C^p$ to a $\Ci^2$ function which can written as a sum $f=\sum\lambda_j f_j$. Then for every $z_0\in \widetilde{M}$ we have
%$$\|f_j\|_{\infty, \widetilde{M}}\leq |f_j(z_0)|d(\widetilde{M})\|f'_j\|_{\infty, \widetilde{M}}\leq |f_j(z_0)|d(\widetilde{M})\|f'\|_{\infty, \widetilde{M}}\;,$$
%where $d(\widetilde{M})$ is a constant depending on the geometry of $\widetilde{M}$. Hence
%$$\sum_{j\in I}\|f_j\|^2_{\infty, \widetilde{M}}\leq d(\widetilde{M})^2\|f'\|_{\infty, \widetilde{M}}^2\sum|f_j(z_0)|^2<+\infty\;.$$
%
%Therefore, by compactness, $M$ is contained in 
%$$(D(0,R_1)\times\ldots\times D(0,R_i)\times \ldots)\times\C^p$$
%with $\sum R_i^2<+\infty$. By projection, we have $M^{\lambda_i }\subset D(0,R_i)\times \C^p$, so $\supp T^{\lambda_j}\subset D(0,R_i)\times \C^p$ too. This follows from the construction for $T^{\lambda_i}$ given in \cite{harvey1}.
%
%Then we have
%$$\|F_j^{\lambda_i, h}\|_{\infty, U_j}\leq R_i$$
%which implies
%$$\sum_{i\in I}\|F_j^{\lambda_i, h}\|^2_{\infty, U_j}\leq \sum_{i\in I}R_i^2<+\infty\;.$$

\medskip

\noindent{\emph{The function $F_j^h$ is holomorphic. } } Indeed, for any $\lambda\in H'$, we write 
$$\lambda=\sum_{i\in I}\alpha_i \lambda_i$$
and
$$\langle F_j^h(z), \lambda\rangle=\sum_{i\in I}\alpha_i F_j^{\lambda_i, h}(z)\;.$$
We now observe that
$$\left|\sum_{i\in I}\alpha_i F_j^{\lambda_i, h}(z)\right|\leq\sqrt{\sum_{i\in I}|F_j^{\lambda_i,h}|^2}\sqrt{\sum_{i\in I}|\alpha_i|^2}\leq\|\lambda\|_{H'}\sqrt{\sum_{i\in I}|F_j^{\lambda_i, h}(z)|^2_{\infty, U_j}}<\|\lambda\|_{H'}\||F^j_h|\|_{\infty, U_j}$$
which is finite, and this implies that the sequence of holomorphic functions
$$\left\{\sum_{i=0}^m\alpha_iF_j^{\lambda_i,h}\right\}_{m\in I}$$
converges uniformly on $U_j$. The limit is then holomorphic, so $F_j^h$ is holomorphic.

\medskip

\noindent{\emph{The function $F_j^h$ extends $\Ci^1$ to the boundary. } } By \cite{harvey1}, there exist sets $A\subset\mathfrak{m}$ and $A_i\subset M^{\lambda_i}$ with $\pi(A_i)=A$, which are $\H^{2p-1}-$negligible and such that outside them we have $\Ci^1$ regularity for $\supp T^{\lambda_i}\cup M^{\lambda_i}$ and for the functions $F_j^{\lambda_i,h}$. Let us consider $p\in \mathfrak{m}\cap \overline{U_j}\setminus A$; for each $i\in I$, one of the following two cases can occur:
\begin{enumerate}
\item $F_j^{\lambda_i, h}(p)\not\in M^{\lambda_i}$,
\item $F_j^{\lambda_i,h}(p)\in M^{\lambda_i}$.
\end{enumerate}
In the former, $F_j^{\lambda_i,h}$ extends holomorphically through $p$, whereas in the latter we can find a relatively compact neighborhood $V$ of $p$ in $\mathfrak{m}$ such that $F_{j}^{\lambda_i,h}$ coincides on $V$ with some $CR$ function $f:V\to M^{\lambda_i}$. In both cases, $F_j^{\lambda_i,h}$ is of class $\Ci^1$ near $p$. Let $U$ be an open set with $\Ci^1$ boundary in $U_j$ such that $bU_j\cap bU= V$.

The restrictions of the derivatives of $F^{h}_j$ to $bU_j$ are continuous, when we derive in a direction tangent to $TbU$; however, by the Cauchy-Riemann equations, we can control the normal derivative with the tangential ones, therefore also the normal derivative of $F^{h}_j$ is a continuous function when restricted to $bU$.

We note that from this follows that the image of $bU$ through one of these maps is a compact set in $H'$ and we can replicate the previous argument, obtaining that the sequence
$$\left\{\sum_{i=0}^m \frac{\de}{\de z_s}F^{\lambda_i,h}_j(z)\lambda_i\right\}_{m\in I}$$
is a Cauchy sequence with respect to the supremum norm on $U$.

Therefore, the limit is continous on the closure of $U$, thus implying that
$$\left\|\left|\frac{\de}{\de z_s}F^h_j\right|\right\|_{\infty, U}<+\infty\;.$$
Moreover, on $bU\cap bU_j=V$, $F^h_j$ coincides with $f$ and we can cover $\H^{2p-1}-$almost all of $bU_j$ with open sets where $F^h_j$ coincides with some CR-functions realizing $M$ as a graph. Therefore, as $M$ is a compact $\Ci^1$ manifold,
$$\left\|\left|\frac{\de}{\de z_s}F^h_j\right|\right\|_{\infty, bU_j}<+\infty$$
hence
$$\left\|\left|\frac{\de}{\de z_s}F^h_j\right|\right\|_{\infty, U_j}<+\infty\;.$$
%By a well known result, on a point of $bU$, we can estimate the normal derivative of the real and imaginary part of $F_{j}^{\lambda_i,h}$ with the tangetial derivatives; therefore
%$$\left\|(F_j^{\lambda_i, h})'\right\|_{\infty, U}\leq C\left\|(F_j^{\lambda_i,h}\vert_{bU})'\right\|_{\infty,bU}\leq \left|(F_j^{\lambda_i,h}\vert_{bU})'(z_0)\right|d(bU)\|\nabla^2(F_j^{\lambda_i,h}\vert_{bU_j})\|_{\infty, bU}\;.$$
%Summation on $i\in I$ gives
%$$\sum_{i\in I}\left\|(F_j^{\lambda_i, h})'\right\|_{\infty, U}^2<+\infty\;.$$

\medskip

\noindent{\emph{The current of integration on the graph of $F_j^h$ has finite mass. } } By the previous paragraph, there exists a constant $C_{h,j}$ such that
$$|\nabla F^h_j(z)|^2=\sum_{i\in I}|\nabla F^{\lambda_i, h}_j(z)|^2\leq C_{h,j}\qquad\textrm{ for every }z\in U_j\;.$$
It is easy to show that there exists a polinomial $g_p(X)$ such that
$$\sum_i |a_i|\leq S<+\infty\Longrightarrow \sum_{|J|=p}\prod_{i\in J}|a_i|\leq g_p(S)<+\infty\;.$$
Therefore
$$\sum_{|J|=p}\prod_{i\in J}|\nabla F^{\lambda_i, h}_{j}(z)|\leq g_p(C_{h,j})<+\infty\qquad\textrm{ for every }z\in U_j\;.$$
We consider the $(p,p)-$form
$$\eta^h_j(z)=\sum_{|J|=p}\bigwedge_{i\in J}dF^{\lambda_i, h}_j(z)\wedge d\overline{F}^{\lambda_i, h}_j(z)$$
which is well-defined by the previous estimates and note that
$$\|\eta^h_j\|_{\infty, U_j}\leq g_p(C_{h,j})\;.$$
Let $\{w_i\}_{i\in I}$ be coordinates for the basis $\{\lambda_i\}_{i\in I}$, i.e. $w_i(v)=\langle v,\lambda_i\rangle$ for $v\in H'$, and denote by $T_{h,j}$ the (alleged) current of integration on the graph of $F^h_j$. Then
$$T_{h,j}(1, w_{i_1}, \overline{w}_{i_1},\ldots, w_{i_p},\overline{w}_{i,p})=\int\limits_{U_j}dF^{\lambda_{i_1}, h}_j(z)\wedge d\overline{F}^{\lambda_{i_1}, h}_j(z)\wedge\ldots\wedge dF^{\lambda_{i_p}, h}_j(z)\wedge d\overline{F}^{\lambda_{i_p}, h}_j(z)\;.$$
Therefore, by Proposition \ref{prp_vol_in_H}, we have
$$M(T_{h,j})\leq \mathcal{L}^{2p}(U_j)\sum_{p'=0}^pg_{p'}(C_{h,j})<+\infty\;.$$
We have to sum all the values from $0$ to $p'$ because we apply the formula of Proposition \ref{prp_vol_in_H} in $H$ and not in $H'$, so we have to consider also the $p-$tuples of coordinates coming in part from $\C^p$ and in part from $H'$.

\medskip

As the $F_j^h$ are a finite number of functions, we can consider the metric functional of integration on their graphs and denote it by $T$. $T$ is a holomorphic $p-$chain in $H\setminus M$, it has finite mass and its support is contained in a product of discs, therefore it is relatively compact in $H$. Moreover, for $\H^{2p-1}-$almost every point in $M$ there is a neighborhood where $\supp T\cup M$ is a $\Ci^1$ manifold.

This implies that $T$ is a metric rectifiable $(p,p)-$current in $H$. We note that for any finite-dimensional projection $p:H\to\C^m$, we have that $d(p_\sharp T)=p_\sharp[M]$; it is an easy application of Theorem \ref{teo_ext} to show that this implies $dT=[M]$.  Finally, it is not difficult to see that the map $x\mapsto (x,F_j^h(x))$ is proper into $H\setminus M$, which is an hilbertian manifold, hence by Theorem \ref{teo_hilb_img} its image is a finite dimensional complex space in $H\setminus M$.  
$\Box$

\medskip

\setlength{\epigraphwidth}{4.7cm}
\chapter*{Un-conclusions\markboth{UN-CONCLUSIONS}{}}
\epigraphhead[60]{\epigraph{Everything will be ok in the end. If it's not ok, it's not the end.}{Anonymous}}
\addcontentsline{toc}{chapter}{Un-conclusions}

What we tried to present in the previous pages cannot honestly be called a theory of metric currents on complex spaces: such a broad generality is far from being achieved. Instead, we spotted some possible applications of metric currents to the complex context, highlighting the characteristics of metric current who made these particular examples sound sensible to us.

In order to complete this picture, many gaps need to be plugged.

\medskip

On the side of the Cauchy-Riemann equation on singular spaces, we weren't able to tackle the problem on a general singularity. The following questions remain unanswered:
\begin{enumerate}
\item is there a way to obtain a structure theorem for metric currents on a general singular space, maybe in connection with some growth conditions on the singular part? could this help to solve the $\debar$ equation?
\item can we generalize the $L^p$ methods used for curves on a higher dimensional space? can we use the condition given in Theorem \ref{teo_cns} to study Cauchy-Riemann equation in $L^p$ on singular spaces, maybe with the aid of the representation formulas given in \cite{ander2}?
\item is it possible to finalize the method of solution we propose for hypersurfaces?
\end{enumerate}

\medskip

Another interesting question is related to the study of Sobolev spaces on singular complex spaces: can we employ the remarks made at the end of Chapter 3 to characterize the spaces where the density hypotheses hold? 

\medskip

Obviously, the last chapter raises the largest number of questions, starting form the very definition of finite-dimensional analytic set. In the characterization of holomorphic chains, we end up considering a wider class than the one given by the usual definition; which properties does this bigger class of spaces enjoy? Can we, at least in Hilbert spaces, give rise to a geometrically meaningful theory with this enlarged class?

Another question is related to the problem of boundaries of holomorphic chains; we solved the problem with a stronger hypothesis than the usual one, regarding the existence of a projection with only transversal crossings. Can we achieve an analogue of the classical result by Harvey and Lawson, where no mention of this requirement is made? 

Finally, it would be interesting to study whether the Monge-Amp\`ere equation can be solved with respect to a current and which geometrical applications could spring from this.

\medskip

There are many more doubts and questions which can be given voice, but these seems to us, if not the most interesting, the most related to what we presented in this thesis. 

As usual, the uncovered points and unanswered questions seem to be way more entertaining than what we already did.

\bibliography{bibsing}{}
\bibliographystyle{siam}

\end{document}